\documentclass[12pt]{article}

\usepackage{setspace}
\usepackage[utf8]{inputenc}

\usepackage{wrapfig}

\usepackage{booktabs}

\usepackage[ruled,linesnumbered]{algorithm2e}
\SetKwFor{While}{while}{:}{end while}
\SetKwFor{ForEach}{for each}{:}{end}

\usepackage{booktabs,caption}
\usepackage[flushleft]{threeparttable}

\usepackage[nointegrals]{wasysym}

\usepackage{adjustbox}

\usepackage{geometry}
\usepackage[hang,flushmargin]{footmisc}
\usepackage{refcount}

\usepackage[english]{babel}

\usepackage{changepage}
\usepackage{epsfig}
\usepackage{amsmath}
\usepackage{multicol}
\usepackage{tikz}

\usepackage{rotating}
\usetikzlibrary{positioning}

\usetikzlibrary{decorations.markings}
\usetikzlibrary{lindenmayersystems}
\usetikzlibrary{patterns}
\usetikzlibrary[shadings]
\usetikzlibrary{fadings}
\tikzfading
[
  name=fade out,
  inner color=transparent!0,
  outer color=transparent!100
]

\usepackage{mathtools}

\usepackage{amssymb}

\usepackage{array}

\usepackage{cancel}

\usepackage{xcolor} % Required for color customization
\usepackage{graphicx} % Required for inserting images

\usepackage[numbers]{natbib} % numeric citations like [3]
\usepackage[colorlinks=true,
            linkcolor=black,    % internal links (e.g., ToC, refs)
            citecolor=blue,     % citations
            urlcolor=blue]      % external URLs
           {hyperref}

\usepackage{url}

\usepackage{imakeidx}
\makeindex

\definecolor{lejla}{RGB}{120,0,133}

\usepackage{tocloft}
\cftsetindents{section}{0em}{2.5em}
\cftsetindents{subsection}{2.5em}{3em}
\cftsetindents{subsubsection}{5.5em}{3.5em}

\makeatletter
\let\oldl@section\l@section
\renewcommand{\l@section}[2]{\bfseries\oldl@section{#1}{#2}}
\makeatother

\renewcommand{\contentsname}{Table of Contents}

\newcommand{\diagonalStripes}[5]{%
  \begin{scope}
    \clip (#1,#2) rectangle ++(#3,#4);
    \begin{scope}[rotate=45, scale=#5]
      \foreach \i in {-20,...,20} {
        \pgfmathtruncatemacro{\c}{mod(\i,2)}
        \ifnum\c=0
          \fill[blue!50!white, opacity=0.5] (\i,-40) rectangle ++(1,80);
        \else
          \fill[green!50!white, opacity=0.5] (\i,-40) rectangle ++(1,80);
        \fi
      }
    \end{scope}
  \end{scope}
}

\title{Reconfiguration of Hamiltonian paths on grid graphs}
\author{}
\date{}

\begin{document}
\pagenumbering{roman}
\setcounter{page}{1}
\begin{titlepage}
    \centering
    \vspace*{0.5in}
    {\huge\bfseries Reconfiguration of Hamiltonian Cycles and Paths in Rectangular Grid Graphs \par}
    \vspace{1in}
    {\large Albi Kazazi\par}
    \vspace{1in}
    {\large A dissertation submitted to the Faculty of Graduate Studies\\
    in partial fulfillment of the requirements for the degree of\\
    Doctor of Philosophy\par}
    \vspace{0.75in}
    {\large Graduate Program in Mathematics and Statistics\\
    York University, Toronto, Ontario\par}
    \vspace{0.75in}
    {\large September 2025\par}
    \vspace{0.75in}
    \copyright\ \large Albi Kazazi, 2025
    \thispagestyle{empty}
\end{titlepage}
\cleardoublepage

% ABSTRACT
\setcounter{page}{2}  
\begin{center}
    {\Large \bfseries Abstract} \\[1em]
\end{center}
\addcontentsline{toc}{section}{Abstract}  % ADD THIS LINE
\noindent An \textit{\(m \times n\) grid graph} is the induced subgraph of the square lattice whose vertex set consists of all integer grid points 
\(\{(i,j) : 0 \leq i < m,\ 0 \leq j < n\}\). Let $H$ and $K$ be Hamiltonian cycles in an $m \times n$ grid graph $G$. We study the problem of reconfiguring $H$ into $K$ using a sequence of local transformations called \textit{moves}. A \textit{box} of $G$ is a unit square face. A box with vertices $a, b, c, d$ is \textit{switchable} in $H$ if exactly two of its edges belong to $H$, and these edges are parallel. Given such a box with edges $ab$ and $cd$ in $H$, a \textit{switch move} removes $ab$ and $cd$, and adds $bc$ and $ad$. A \textit{double-switch move} consists of performing two consecutive switch moves. If, after a double-switch move, we obtain a Hamiltonian cycle, we say that the double-switch move is \textit{valid}.
\null 
\noindent We prove that any Hamiltonian cycle $H$ can be transformed into any other Hamiltonian cycle $K$ via a sequence of valid double-switch moves, such that every intermediate graph remains a Hamiltonian cycle.
\null 
\noindent This result extends to Hamiltonian paths. In that case, we also use single-switch moves and a third operation, the \textit{backbite move}, which enables the relocation of the path endpoints.

% ACLNOWLEDGEMENTS 

% Version 1: Traditional

\newpage

% Version 2: Complaints and Grievances

%\newpage
%\section*{Acknowledgements}
%\addcontentsline{toc}{section}{Acknowledgements}  % Changed chapter to section
%\noindent \underline{\textbf{COMPLAINTS AND GRIEVANCES}}

%\null 

%\noindent The big one goes to my supervisor Neal Madras. His insightful and patient feedback has been invaluable. Through his mentorship, my time in the program has been one of growth, mathematical and beyond. I’m lucky and grateful to have had him as my advisor. If it isn't clear already, he has clearly set the bar carelessly high - how is anyone else supposed to measure up?? very annoying.

%\null 

%\noindent The next one is for Ada Chan. As my master's survey paper supervisor, she was very helpful, answering all my questions, and it was at that time that I realized that math research is a actually a lot of fun and so I should go for a PhD. As you can see, she was instrumental to the Neal-high-bar situation I have to contend with now. just as culpable!

%\null 

%\noindent aaand, I wouldn't be doing a master's at all if it wasn't for Paul Szeptycki and Youness Lamzouri. Attending their classes as an undergrad made me think that math is easy. Basically they were guilty of the same misdeed as Neal above - setting the bar carelessly high.

%\null 

%\noindent Lastly a quick one for my family and all their support. I will look back at this time and miss it, and it's your fault.

%\null 

%\noindent There are other people to blame, I'm sure, but I have restraint. $\square$

\newpage
\section*{Acknowledgements}
\addcontentsline{toc}{section}{Acknowledgements}  % Changed chapter to section

I would like to thank my supervisor Neal Madras for his insightful and patient feedback. Through his mentorship, my time in the program has been one of growth, mathematical and beyond. I'm lucky and grateful to have had him as my advisor.

\null 

Thanks to Ada Chan who, as my master's supervisor, was very helpful, answering all my questions. It was that experience that made me decide to pursue a PhD.

\null

Thanks to Paul Szeptycki and Youness Lamzouri, who taught early courses that made math feel approachable and worth the time.

\null

Thanks to my family for all their support.

\newpage
\begingroup
\setstretch{1.5}
\renewcommand{\contentsname}{\huge Table of Contents}
\renewcommand{\cftsecfont}{\normalsize}
\renewcommand{\cftsecpagefont}{\normalsize}
\addcontentsline{toc}{section}{Table of Contents}  % Move here
\tableofcontents
\endgroup
\newpage
\pagenumbering{arabic}
\setcounter{page}{1}
\phantomsection
\addcontentsline{toc}{section}{Introduction}
\section*{Introduction}
\null

Reconfiguration problems consist of a set of configurations, or a state space, and a set of reconfiguration rules that can be used to change one configuration into another. We call an application of a rule a move or a step. Given this data, we may construct a reconfiguration graph, where configurations are vertices, and two vertices are adjacent if there is a move that can change the configurations they represent into one another. Some popular problems include the knight's tour on a chessboard, the 15-puzzle with sliding blocks, and the Rubik's cube.

\null

\noindent Many well-known problems have been studied from a reconfiguration perspective, in particular with respect to the connectivity of the reconfiguration graph. One example is the reconfiguration of the $k$-colorings of a graph, where a $k$-coloring is an assignment of one of $k$ colors to each vertex so that adjacent vertices receive different colors, and the allowed move is to recolor one vertex at a time~\cite{celaya2016reconfiguring,bonsma2009finding, bonamy2014reconfiguration, feghali2020reconfiguration, johnson2016finding}. Another example involves dominating sets of size at most $k$, where a dominating set is a set of vertices such that every vertex is either in the set or adjacent to a vertex in the set. Here, the allowed move is to add or remove a single vertex~\cite{haas2014k, suzuki2016reconfiguration}. A related example involves independent sets. An independent set is a set of vertices in a graph such that no two vertices in the set are adjacent. We place a token on each vertex of an independent set of size at most $k$ and we allow one token to be slid at a time along an edge of the graph, provided that the resulting set is also an independent set of size $k$~\cite{demaine2015linear, bonamy2018recoloring}. For surveys on reconfiguration problems, see~\cite{van2013complexity,nishimura2018introduction}.

\null

\noindent\textbf{Hamiltonian cycles: history and difficulty.} The study of Hamiltonian paths and cycles in graphs dates back to Euler's 1736 resolution of the Königsberg bridge problem, which asked whether one could walk through the city crossing each of its seven bridges exactly once. This work established the foundation for graph theory. A century later, Hamilton introduced his icosian game in 1857, asking whether one could traverse all vertices of a dodecahedron exactly once and return to the starting point. Euler's work on the bridges problem led to a complete and easily verifiable characterization of graphs admitting Eulerian trails—graphs where every edge is traversed exactly once. Hamilton's problem, by contrast, remains notoriously difficult. Determining whether a graph contains a Hamiltonian cycle is NP-complete for general graphs, and few useful characterizations exist. Dirac's theorem provides some  partial insight: it states that every graph with $n \geq 3$ vertices where each vertex has degree at least $n/2$ must contain a Hamiltonian cycle. Even for special graph classes, Hamiltonicity questions remain open; Barnette's conjecture \cite{barnetteconjecture} -- that every 3-connected cubic planar bipartite graph is Hamiltonian -- has stood unresolved since 1969. Hamiltonian cycle problems arise naturally in optimization through the traveling salesman problem, in recreational mathematics through the knight's tour on a chessboard, and as canonical examples in computational complexity theory, making their study both practically and theoretically significant.

\null 

\noindent\textbf{Reconfiguration of Hamiltonian cycles.} The reconfiguration problem for Hamiltonian cycles asks the following: Given any two Hamiltonian cycles $H$ and $K$ in a graph and an allowed type of move, is there a sequence  $H_0, H_1, ..., H_r$ of Hamiltonian cycles, where $H=H_0$, $K=H_r$, such that for each $j \in \{1, \dots r\}$ $H_{j}$ is obtained by a single application of the allowed move on $H_{j-1}$? Often, the allowed move consists in adding and removing one or two pairs of edges from a given Hamiltonian cycle, in such a way that we obtain a new Hamiltonian cycle. For general graphs, this is a hard problem \cite{takaoka2018complexity}. 

In this thesis, we restrict our attention to the Hamiltonian cycle and path reconfiguration problem on a class of graphs called grid graphs, which are subgraphs of the square lattice. An \textit{\(m \times n\) grid graph $G$} is the induced subgraph of the square lattice whose vertex set consists of all integer grid points 
\(\{(i,j) : 0 \leq i < m,\ 0 \leq j < n\}\) with edges between vertices at distance 1. We call the unit-square faces of a grid graph \textit{boxes}. We show that for an $m \times n$ grid graph, the reconfiguration graph of Hamiltonian cycles and the reconfiguration graph of Hamiltonian paths is connected under a small set of suitable moves.

\begingroup
\setlength{\intextsep}{0pt}
\setlength{\columnsep}{10pt}
\begin{wrapfigure}[]{r}{0cm}
\setlength{\intextsep}{0pt}
\setlength{\columnsep}{20pt}
\begin{adjustbox}{trim=0cm 0cm 0cm 0cm}
\begin{tikzpicture}[scale=1.25]

\begin{scope}[xshift=0cm]
   \draw[gray,very thin, step=0.5cm, opacity=0.4] (0,0) grid (3.5,2.5); 

\draw[gray,very thin, step=0.5cm, opacity=0.4] (0,0) grid (2.5,2); 

\draw[blue, line width=0.5mm] (0,2.5)--++(0,-2.5)--++(0.5,0)--++(0,0.5)--++(0.5,0)--++(0,-0.5)--++(0.5,0)--++(0,1)--++(0.5,0)--++(0,-1)--++(0.5,0)--++(0,1)--++(0.5,0)--++(0,-1)--++(0.5,0)--++(0,1.5)--++(-2,0)--++(0,0.5)--++(2,0)--++(0,0.5)--++(-2.5,0)--++(0,-1.5)--++(-0.5,0)--++(0,1.5)--++(-0.5,0); 

%\node[] at (0.75, 0.25) [scale=0.8]{\small{X}};

\node[below] at (1.75, 0) [scale=1]{\begin{tabular}{c} Fig. I.1. A 1-complex Hamiltonian \\ cycle on an $8 \times 6$ grid graph. \end{tabular}};;
\end{scope}

\end{tikzpicture}
\end{adjustbox}
\end{wrapfigure}

\noindent The question of whether an $m \times n$ grid graph has a Hamiltonian path or cycle was first studied by Itai et al. in \cite{itai1982hamilton}. They showed that for an $m \times n$ grid graph to have a Hamiltonian cycle, it is necessary and sufficient that at least one of $m$ and $n$ is even. A \textit{polyomino graph} is a superclass of $m \times n$ grid graphs, where every edge belongs to a box. We give a constructive definition at the start of Chapter 1. Umans and Lenhart \cite{umans1997hamiltonian} gave a polynomial-time algorithm to find a Hamiltonian
cycle in polyomino graphs, if one exists. Chen et al. gave an efficient algorithm for constructing Hamiltonian paths in rectangular grid graphs \cite{chen2002efficient}.

Nishat and Whitesides \cite{nishat2017bend} introduced the ``flip'' and ``transpose'' moves described below, and a complexity measure called ``bend complexity" for Hamiltonian cycles in rectangular grid graphs. Roughly, a 1-complex Hamiltonian cycle is one in which every vertex of $G$ is connected to the boundary via a straight line.  They prove that using these two moves, it is possible to reconfigure any pair of 1-complex Hamiltonian cycles of $G$ into one another.

\noindent We dispense with the need for bend complexity constraints, proving that any Hamiltonian cycle in a rectangular grid graph can be reconfigured into any other, using a more general move, which we call a double-switch move.

\endgroup 

\null

\noindent Let $H$ be a Hamiltonian cycle of an $m\times n$ grid graph $G$. A box of $G$ with vertices $a, b, c, d$ is considered \textit{switchable} in $H$ if it has exactly two edges in $H$, and these edges are parallel. Let $abcd$ be a switchable box with edges $ab$ and $cd$ in $H$. We define a \textit{switch move} on the box $abcd$ in $H$ as follows: remove edges $ab$ and $cd$ and add edges $bc$ and $ad$. If $X$ is a switchable box in $G$ in $H$, we denote a switch move as $\textrm{Sw}(X)$.

\begingroup

\setlength{\intextsep}{0pt}
\setlength{\columnsep}{20pt}
\begin{center}
\begin{adjustbox}{trim=0cm 0cm 0cm 0cm}

\end{adjustbox}
\end{center}

\noindent A \textit{double-switch move} is a pair of switch operations where we first switch $X$ and then find another switchable box $Y$ and switch it, and denote it by $X \mapsto Y$. If, after a double-switch move, we obtain a new Hamiltonian cycle, we call the move a \textit{valid move}. See Figure I.2

\endgroup 

\begingroup
\setlength{\intextsep}{0pt}
\setlength{\columnsep}{20pt}
\begin{wrapfigure}[]{r}{0cm}
\setlength{\intextsep}{0pt}
\setlength{\columnsep}{20pt}
\begin{adjustbox}{trim=0cm 0cm 0cm 0cm}
%
\end{adjustbox}
\end{wrapfigure}

\null 

\noindent Let $X=abcd$ and $Y=dcef$ be boxes sharing the edge $cd$ of $G$. Assume that the edges $ab, fd, dc$ and $ce$ belong to $H$, and that the edges $fe,ad$ and $bc$ do not. A \index{flip|textbf}\textit{flip} move consists in removing the edges $fd, ce$ and $ab$, and adding the edges $ad, bc$ and $fe$. Effectively, this is the same as the move $X \mapsto Y$, obtained by first switching $X$ and then switching $Y$. See Figure I.3.

\endgroup

\begingroup
\setlength{\intextsep}{0pt}
\setlength{\columnsep}{20pt}
\begin{wrapfigure}[]{l}{0cm}
\setlength{\intextsep}{0pt}
\setlength{\columnsep}{20pt}
\begin{adjustbox}{trim=0cm 0cm 0cm 0cm}

\end{adjustbox}
\end{wrapfigure}

\noindent Consider the four boxes $X=abcd$, $Y=cbef$, $Z=cfgh$ and $W=dchi$ that are incident on the vertex $c$. Note that $X$ and $Y$ share the edge $cb$, $Y$ and $Z$ share $cf$, $Z$ and $W$ share $ch$, and $W$ and $X$ share $cd$. Assume that the edges $ab, be, ef, fc, cd$ and $hg$ belong to $H$ and that the edges $ad, bc$ and $fg$ do not. A \index{transpose|textbf}\textit{transpose} move consists in switching $X$ and then switching $Z$. See Figure I.4.

\noindent Nishat in \cite{nishat2020reconfiguration} showed that flip and transpose moves are always valid. The more general double-switch moves are sufficient for constructing algorithms that reconfigure arbitrary Hamiltonian cycles in grid graphs. This comes at the added cost of verifying the validity of each move. We provide such reconfiguration algorithms and prove the existence of all required moves.

\null

%We call e box of $G$ with exactly three edges in $H$ a \textit{leaf}. A flip move is a double switch move such that the switched boxes share an edge of $G$. Note that of them must necessarily be a leaf and the other a switchable box. A transpose move is a double switch move such that the switched boxes share exactly one vertex of $G$ and such that there is a leaf that is adjacent to both of them. The following theorem is the main result of this paper.

\noindent \textbf{Theorem I.1.} Let $H$ and $K$ be any two Hamiltonian cycles in an $m \times n$ grid graph $G$ with $n \geq m$. Then there exists a sequence of at most $n^2m$ valid double-switch moves that reconfigures $H$ into $K$.

\null 

\noindent See \cite{video} for an illustration.

\null

\noindent We also give an algorithm by which we are able to reconfigure any two Hamiltonian paths in an $m \times n$ grid graph. This requires a move that can relocate the end-vertices of a path, but switch moves do not relocate end-vertices. To address this, we use a move called a \textit{backbite} move, first introduced by Mansfield in 1982  ~\cite{mansfield1982monte}.

\begingroup
\setlength{\intextsep}{0pt}
\setlength{\columnsep}{20pt}
\begin{wrapfigure}[]{r}{0cm}
\begin{adjustbox}{trim=0cm 0cm 0cm 0cm}

\end{adjustbox}
\end{wrapfigure}

Let $H$ be a Hamiltonian path $v_1, \ldots, v_r$ of an $m\times n$ grid graph $G$, and let $v_s$ be adjacent to $v_1$, $s\neq 2$. If we add the edge $\{v_1,v_s\}$, we obtain a cycle $v_1, \ldots, v_s, v_1$, and a path $v_s, \ldots, v_r$. Now, if we remove the edge \{$v_{s-1},v_s\}$, we obtain a new Hamiltonian path $H'=(H \setminus \{v_{s-1},v_s\}) \cup \{v_1,v_s\}$. This operation is called a backbite move. See Figure I.5.

\endgroup

\null

\noindent \textbf{Theorem I.2.} Let $H$ and $K$ be any two Hamiltonian paths in an $m \times n$ grid graph $G$ with $n \geq m$. Then there exists a sequence of at most $n^2m+O(n^2)$ valid switch, double-switch, and backbite moves that reconfigures $H$ into $K$.

\null 

%\noindent See \cite{video1} \textcolor{red}{Make this video} for an illustration.

\setlength{\intextsep}{00pt}
\setlength{\columnsep}{20pt}
\begin{center}
\begin{adjustbox}{trim=0cm 0cm 0cm 0cm}

\end{adjustbox}
\end{center}

\noindent \textbf{Motivation.}  A self-avoiding walk is a walk in a lattice where every vertex is unique. A Hamiltonian path in a grid graph is an example of a self-avoiding walk. Madras and Slade in \cite{madras2013self} present a comprehensive and rigorous study of self-avoiding walks. One application of Theorems I.1 and I.2 is in chemical physics, drawing from the theory of self-avoiding walks. Researchers in \cite{oberdorf2006secondary}, \cite{jacobsen2008unbiased}, \cite{deutsch1997long}, and \cite{mansfield1982monte} use Monte Carlo methods to study statistical properties of polymer chains, which they abstracted as cycles and paths in the three-dimensional square lattice. They use self-avoiding walks to model how a flexible polymer chain is arranged in a liquid solution. A polymer chain's concentration is the fraction of vertices of the lattice that are occupied by the vertices (monomers) of the polymer. The authors consider maximally concentrated polymers (high-density polymers), where all the space is occupied by the polymer. These can be naturally represented as Hamiltonian paths or cycles. They view the set of Hamiltonian cycles of a rectangular grid graph as the state space of a Markov chain, with the double-switch move being the transition mechanism. Given a Hamiltonian cycle (a state in the state space), we choose two switchable boxes at random and perform a double-switch move. If the move is valid, then the new state is the resulting Hamiltonian cycle. Otherwise, we remain at the initial state and choose another pair of switchable boxes. The idea is that after a sufficiently large number of transitions, we obtain a sequence of many different states, which represents a reasonable random sample of the entire state space.

The goal of these methods is to obtain a uniform random sample of all Hamiltonian cycles of the grid. It follows from the theory of Markov chains that for the Monte Carlo simulation to generate such a distribution of Hamiltonian cycles, the chain must be \index{reversible (Markov chain)|textbf}\textit{reversible} and \index{irreducible (Markov chain)|textbf}\textit{irreducible}. Reversibility is satisfied when every move has an inverse, which is true here since the inverse of a double-switch move is obtained by switching the same pair of boxes again. Irreducibility is satisfied when we can reconfigure between any two Hamiltonian cycles using valid moves. Irreducibility ensures that the limit 
$$
\lim_{n\to\infty} P^n(i,j) = \pi(j)
$$
exists, and that it satisfies
$$
(1)\quad \sum_i \pi(i) = 1, \qquad 
(2)\quad \pi(j) = \sum_i \pi(i) P(i,j),
$$
where $\pi(j)$ is the long-term probability of Hamiltonian cycle $j$ occurring. Two other conditions, \index{aperiodic (Markov chain)|textbf}\textit{aperiodicity} (the chain can ``stay put" with positive probability, which happens here when moves are rejected) and finiteness of the state space (obviously true for a finite grid graph) are necessary and clearly satisfied but warrant mention for completeness. Furthermore, $\pi$ is the only nonnegative solution of (1) and (2), and we call $\pi$  the \textit{equilibrium} distribution. 

Since each move has an inverse, the probability $P(i,j)$ of going from a Hamiltonian cycle $i$ to another one $j$ via a single move is the same as the probability of going from $j$ to $i$. This means that the transition matrix $P$ is symmetric: $P(i,j) = P(j,i)$. We can verify that the uniform distribution satisfies equation (2). If $\pi$ is uniform (that is, $\pi(i) = c$ for all $i$), then it follows that
$$
\pi(j) = \sum_i c\,P(i,j) = c\sum_i P(i,j) = c.
$$
Since the uniform distribution $\pi$ satisfies (1) and (2), it is \textit{an} equilibrium distribution, and because the equilibrium distribution of an irreducible, aperiodic chain with finite state space is unique, it is \textit{the} equilibrium distribution. The harder condition to prove is irreducibility.

The authors in \cite{oberdorf2006secondary, jacobsen2008unbiased, deutsch1997long, mansfield1982monte} assume irreducibility but do not prove it. This dissertation proves that irreducibility holds for all rectangular grid graphs, thereby confirming that the corresponding Markov Chain Monte Carlo method indeed produces a uniform equilibrium distribution. We further extend this result to Hamiltonian paths by introducing the single-switch and backbite moves in addition to the double-switch move. For a more detailed discussion on Monte Carlo methods and reconfiguration of self-avoiding walks, see Chapter 9 in \cite{madras2013self}.

\null

\noindent \textbf{Related work} Nishat, Whitesides, and Srinivasan extended the result of \cite{nishat2017bend} to 1-complex Hamiltonian paths in rectangular grid graphs \cite{nishat2024hamiltonian, nishat20231, nishat2023reconfiguration}, and to 1-complex Hamiltonian cycles in L-shaped grid graphs \cite{nishat2019reconfiguring}. The authors define a 1-complex $s$,$t$ Hamiltonian path to be a 1-complex Hamiltonian path that begins and ends at diagonally opposite corners $s$ and $t$ of a rectangular grid graph. We note that Corollary 5.4.1 extends the results in \cite{nishat2024hamiltonian}, \cite{nishat20231}, and  \cite{nishat2023reconfiguration} to arbitrary $s$,$t$ Hamiltonian paths.

\null

\noindent We prove the existence of move sequences that reconfigure Hamiltonian cycles and paths in rectangular grid graphs (Theorems I.1 and I.2). The proofs rely on structural properties that Hamiltonian paths impose on grid graphs, which Chapter 1 analyzes in detail. Chapter 2 presents the reconfiguration algorithm for cycles (Theorem I.1\footnote{We prove Theorem 2.1 which is a slight generalization of Theorem 1.I.},), whose correctness depends on two algorithms: MLC and 1LC. Chapter 3 proves these algorithms exist and shows that the 1LC proof requires a technical lemma (Lemma 3.13) handling a difficult family of configurations. Chapter 4 analyzes this family and proves Lemma 3.13, while also establishing results needed for Chapter 5. Chapter 5 extends the reconfiguration result to Hamiltonian paths (Theorem I.2, restated there as Theorem 5.9), building on the cycle algorithms (Chapter 2), the structural analysis (Chapter 1), and the tools from Chapter 4.

\null

\newpage

\section{Structural properties of Hamiltonian paths and cycles on grid graphs}

A \index{grid graph|textbf}\textit{grid graph} is a subgraph of the integer grid $\mathbb{Z}^2$. A \index{lattice animal|textbf}\textit{lattice animal} is a finite connected subgraph of $\mathbb{Z}^2$. A \index{Hamiltonian path|textbf}\index{Hamiltonian cycle|textbf}\textit{Hamiltonian path (cycle)} of a graph $G$ is a path (cycle) that visits each vertex of the graph exactly once. Assume that $G$ has a cut vertex $v$. Then $G$ cannot have a Hamiltonian cycle. %\textcolor{red}{[Suppose, for contradiction, that $G$ has a Hamiltonian cycle $C$. Let the two neighbors of $v$ on $C$ be $u$ and $w$. Deleting $v$ from $C$ leaves a $u$--$w$ path $P$ that contains $V(G) \setminus \{v\}.$ In particular, $u$ and $w$ lie in the same connected component of $G - v$ (they are joined by $P$ inside $G - v$), and every vertex of $G - v$ lies on $P$, hence in that same component. But $v$ is a cut vertex, so $G - v$ is disconnected and has at least two components. The preceding paragraph shows all vertices of $G - v$ lie in a single component, a contradiction.]} 
Let $G_1$ and $G_2$ be the components of $G\setminus v$. Let $H_1$ be a Hamiltonian path of $G\setminus G_1$ and let  $H_2$  be a Hamiltonian path of $G\setminus G_2$ such that $H_1$ and $H_2$ have $v$ as an end-vertex. Then a Hamiltonian path $H$ of $G$ can be obtained by concatenating $H_1$ and $H_2$. Since $H_1$ and $H_2$ are smaller than $H$, they are easier to find and reconfigure. It follows that a graph that cannot be decomposed in this manner must be 2-connected. Thus, from here on, we will restrict our attention to 2-connected grid graphs.  %\textcolor{red}{[NOTE: Def of polyking is used in Lemmas 1.3.5, 1.3.6, 1.3.7]}

We will show that, for a subclass of 2-connected grid graphs called rectangular grid graphs, we can reconfigure a Hamiltonian path into any other Hamiltonian path, by making a sequence of moves. We define rectangular grid graphs further down and we define moves in Section 1.4. We remark that Section 1.4 requires only the definitions from the opening of this chapter and the definition of a polyomino graph from Section 1.1. Readers interested in understanding the moves before proceeding through the more dense technical Sections 1.1--1.3 may prefer to read Section 1.4 immediately.

Suppose $H$ and $K$ are two distinct Hamiltonian paths of $G$. We will show that there exists a sequence of moves $\mu_1,\ldots, \mu_r$ and a sequence of Hamiltonian paths $H_0, H_1, \dots, H_r$, where $H=H_0$ and $K=H_r$, such that for each $j \in \{1, \ldots, r\}$ $\mu_j$ reconfigures $H_{j-1}$ into $H_j$.

\null 

\noindent \textbf{Notation.} Let $G$ be a grid and let $H$ be a Hamiltonian path of $G$. We will call the unit square faces of the integer grid \index{box|textbf}\textit{boxes} (see Figure 1.2). We denote by $\textrm{Boxes}(G)$ the set of boxes of a grid graph $G$ that have all four of their edges in $G$. From here on we will assume that $G$ is finite. We will need some definitions to navigate $G$ and $H$. We will position $G$ in the first quadrant of the integer grid in such a way that the westernmost vertices of $G$ have x-coordinate zero and the southernmost vertices of $G$ have y-coordinate zero. We use the $x$ and $y$ coordinates to describe a rectangle in the graph and denote it $R(k_1,k_2; l_1,l_2)$. This rectangle corresponds to the rectangle on the plane which is the Cartesian product of closed intervals $[k_1,k_2] \times [l_1,l_2]$. We will denote a box of $G$ by $R(k,l)$, where $k$ and $l$ are the coordinates of the corner of the box that is closest to the origin. 

We denote a vertex $v$ by $v(k,l)$, where $k$ and $l$ are the vertex coordinates. We denote a horizontal edge $e$ by $e(k_1,k_2; l_1)$, where $k_1, k_2$ are the x-coordinates of the vertices of $e$, and $l_1$ is the y-coordinate of the vertices of $e$. Similarly, we write $e(k_1; l_1, l_2)$ for vertical edges. For a horizontal line consisting of $r$ consecutive edges, we write $e(k, k+r;l)$; we denote vertical lines analogously.  
It will be convenient to use the notation $\{u, v\}$ to denote edges of $G$, and the notation $(u, v)$ to denote directed edges of $G$. For a directed edge $e=(u,v)$, $u$ is said to be the \textit{tail} of $e$ and $v$ is said to be the \textit{head} of $e$.

Let $P(u,v)$ and $P(v,w)$ be paths. We denote the concatenation of $P(u,v)$ with $P(v,w)$ as $P(u,v), P(v,w)$ (a comma indicates concatenation).

\null

The rest of Chapter 1 contains definitions and technical results used in Chapters 2-5. A flowchart at the end of the chapter illustrates dependencies within Chapter 1 and connections to later chapters.  

Our reconfiguration strategy relies on controlling which edges belong to the Hamiltonian cycle by applying moves. To analyze when moves can add or remove specific edges while preserving the Hamiltonian property, we need to understand how $H$ decomposes the grid graph into components. Chapter 1 explores the structure that arises from this decomposition.

Consider a Hamiltonian path $H$ on the \index{m x n grid graph@$m \times n$ grid graph}\textit{$m \times n$ grid graph} $G$ in Figure 1.4. Section 1.2 describes how $H$ decomposes $G$ into \index{H-component@$H$-component}$H$-components. Since these components need not be rectangular, Section 1.1 first defines polyomino graphs as a generalization of rectangular grids.

Section 1.3 proves technical results about $H$-component structure needed in Chapters 3-5. Section 1.4 gives conditions for valid moves and describes their effects on edge orientations and path structure.  See the Flowchart 1 at the end of the chapter for detailed dependencies.

\null

\subsection{Polyomino graphs}

\noindent In this section we prove preliminary results about polyomino graphs. Corollary 1.1.5, a consequence of Jordan's Curve Theorem, is used throughout later proofs. Results 1.1.7-1.1.11 characterize structural properties of polyominoes needed in Section 1.3 and Chapters 3-5.

\null 

\noindent \textbf{Definitions.} \noindent Let $G$ be a planar graph with face set $F$, and let $F' \subseteq F$. Analogous to a vertex-induced subgraph, we define the \index{face-induced subgraph|textbf}\textit{face-induced subgraph} $G[F']$ to be the subgraph of $G$ consisting of all vertices and edges incident to the faces in $F'$.

We restrict our attention to a subclass of 2-connected, finite, box-induced subgraphs of the square lattice that we call \index{polyomino|textbf}\textit{polyomino graphs} (see Figure 1.2). We define polyomino graphs recursively as follows: 

\begin{enumerate}
    \item The polyomino graph of order 1 is the box-induced subgraph on a single box. This is the cycle $C_4$.
    
    \item Let $\mathcal{P}_k$ denote the set of polyomino graphs of order $k$. Let $\{X_1, \dots, X_k\}$ be a set of boxes in $\mathbb{Z}^2$, and let the polyomino graph $G \in \mathcal{P}_k$ be the box-induced subgraph on this set of boxes. We may obtain a polyomino graph $G' \in \mathcal{P}_{k+1}$ by taking the box-induced subgraph on the set $\{X_1, \dots, X_k, X\}$, where $X \in \mathbb{Z}^2 \setminus \{X_1, \dots, X_k\}$ shares an edge with some box in $\{X_1, \dots, X_k\}$.
\end{enumerate}

\noindent \textbf{Remark.} It is possible for two polyomino graphs $G \in \mathcal{P}_k$ and $G' \in \mathcal{P}_{k+r}$, for some $r > 0$, to be equal as subgraphs of the square lattice, since different sets of boxes can induce the same subgraph, as we shall now show. Consider a polyomino $G \in \mathcal{P}_k$, constructed recursively using the procedure above. Then $G$ is the box-induced subgraph on a set of boxes $\{X_1, \dots, X_k\}$ of $\mathbb{Z}^2$.

We define a \index{one-hole|textbf} \textit{one-hole} of $G$ to be a box $X \in \mathbb{Z}^2$ such that:

    1. $X \notin \{X_1, \dots, X_k\}$, and
    
    2. $X$ shares an edge with four distinct boxes in $\{X_1, \dots, X_k\}$. See Figure 1.1.

In this case, all four edges and all four vertices of $X$ are already contained in $G$, so $X \in \text{Boxes}(G) \setminus \{X_1, \dots, X_k\} $. We may now define $G' \in \mathcal{P}_{k+1}$ as the box-induced subgraph on $\{X_1, \dots, X_k, X\}$. Since the addition of $X$ introduces no new edges or vertices to the graph, we have $G = G'$.

\null 

\noindent It will be convenient to define here a generalization of a polyomino graph called a \index{polyking|textbf}\textit{polyking} graph. A polyking graph of order $k$ is defined to be a box-induced subgraph on $k$ boxes, constructed recursively in the same way as a polyomino graph, except that each new box is only required to share a vertex (rather than an edge) with an existing box.

\begin{center}
\begin{adjustbox}{trim=0cm 0cm 0cm 0cm}

\end{adjustbox}
\end{center}

\noindent Let $G$ be a polyomino graph. We say that $G$ has a \index{polyking junction|textbf}\textit{polyking junction} at a vertex $v(k,l)$ if there are exactly two boxes $Z$ and $Z'$ of $G$ incident on $v(k,l)$ and $V(Z) \cap V(Z') = \{v(k,l)\}$. 

\noindent From here on, we will use the term \textit{polyomino} to mean polyomino graph, and \textit{polyking} to mean polyking graph. 

\null 

\noindent Let $G$ be a polyomino, and let $X_1, X_2$ be two distinct boxes of $G$. If $X_1$ and $X_2$ share an edge of $G$, we say that $X_1$ and $X_2$ are adjacent. Define a \index{walk of boxes|textbf}\textit{walk of boxes in $G$} to be a sequence $X_1, \dots, X_r$ of boxes in $G$, not necessarily distinct, such that for all $j \in \{1, 2, \dots, r-1\}$, either $X_j$ is adjacent to $X_{j+1}$ or $X_j = X_{j+1}$. We denote such a walk by $W(X_1, X_r)$. For each $j \in \{1, \dots, r-1\}$, we call the edge of $G$ shared by $X_j$ and $X_{j+1}$ a \index{gluing edge|textbf}\textit{gluing edge of $W(X_1, X_r)$}, whenever $X_j$ and $X_{j+1}$ are distinct boxes. If for all $i, j \in \{1, 2, \dots, r\}$ with $i \neq j$, $X_i$ is distinct from $X_j$, we call the sequence a \index{path of boxes|textbf}\textit{path of boxes in $G$} and denote it by $P(X_1, X_r)$. The length of $P(X_1,X_r)$ is $r-1$. A \index{cycle of boxes|textbf}\textit{cycle of boxes in $G$} is a walk $X_1, \dots, X_r$ such that $X_1 = X_r$ and for all $i, j \in \{1, \dots, r-1\}$, $X_i \neq X_j$. Note that, by construction, for any two boxes $X$ and $Y$ in $G$, there is at least one path of boxes $P(X,Y)$ in $G$. We say that $G$ is \index{box-path connected|textbf}\textit{box-path connected}.

A graph $G$ is said to be \index{k-connected@$k$-connected|textbf}\textit{$k$-connected} if it remains connected after the removal of any  $k-1$ vertices.

\null

\noindent \textbf{Lemma 1.1.1.} Every face of a 2-connected planar graph is bounded by a cycle. (Proposition 2.1.5 in Mohar and Thomassen, Graphs on Surfaces, p.21.\cite{GraphsOnSurfaces})  $\square$  

\null 

\noindent \textbf{Lemma 1.1.2.} Polyominoes are 2-connected. 

\null

\noindent \textit{Proof.} Let $G$ be a polyomino. Suppose we remove the vertex $w$ from $G$. We want to show that $G$ is still connected.  Let $u,v \in G \setminus w$ be incident on the boxes $X_u$ and $X_v$ of $G$, respectively. Since polyominoes are box-path-connected, there is a path of boxes $P(X_u,X_v)$ in $G$. We show that we can find a path of vertices $P(u,v)$ in $G \setminus w$ by induction on the length of the shortest path of boxes $P(X_u,X_v)$.

\null

\noindent Base Case. Suppose $P(X_u,X_v)$ in $G$ has length zero. Then $P(X_u,X_v)=X_u=X_v$ is the cycle $C_4$, which is 2-connected.

\null

\noindent Inductive Case. Assume that if $P(X_u,X_v)$ in $G$ has length $k$ then we can find a path of vertices $P(u,v)$ in $G \setminus w$. Suppose $P(X_u,X_v)$ has length $k+1$, let the box $Y$ be the $(k+1)^{\text{st}}$ box of $P(X_u,X_v)$, and let $\{u_1,u_2\}$ be the gluing edge between $Y$ and $X_v$. Remove $w$ from $G$. Either $w$ is one of $u_1$ and $u_2$, or it is not. 

If the latter, by inductive hypothesis, there is a path $P(u,u_1)$ contained in $P(X_u,Y)$ and avoiding $w$. By base case, there is a path $P(u_1,v)$ contained in $X_v$ and avoiding $w$. Then, $P(u,u_1)$, $P(u_1,v)$ is a $(u,v)$-path in $G \setminus w$ contained in $P(X_u,X_v)$.

Assume the former. By symmetry, we may assume without loss of generality that $w=u_1$. By inductive hypothesis, there is a path $P(u,u_2)$ contained in $P(X_u,Y)$ and avoiding $w$. By base case, there is a path $P(u_2,v)$  contained in $X_v$ and avoiding $w$. Then, $P(u,u_2)$, $P(u_2,v)$ is a $(u,v)$-path in $G \setminus w$ contained in $P(X_u,X_v)$. $\square$.

\null 

\noindent \textbf{Observation 1.1.3.} Let $G$ be a polyomino. Then:

(a) Every edge of $G$ is incident on a box of $G$. 

(b) A box-path-connected set of boxes is a polyomino.

\null 

\noindent \textbf{Theorem 1.1.4 (Jordan’s Curve Theorem for polygons}\cite{GraphsOnSurfaces}).   
Let $Q$ be a simple polygon in the plane. Then the set $\mathbb{R}^2 \setminus Q$ consists of two disjoint subsets, called the ``interior'' (Int) and ``exterior'' (Ext), each of which has $Q$ as its boundary. Moreover, any two points in Int (or in Ext) can be joined by a polygonal path that does not intersect $Q$, and any polygonal path joining a point of Int to a point of Ext must intersect $Q$ $\square$.

% General JCT 
{
%\noindent \textbf{Theorem 1.1.4. Jordan's Curve Theorem (JCT).} A simple closed curve $Q$ divides the set of points of the plane not on $Q$ into two disjoint subsets Int (for ``Interior'') and Ext (for ``Exterior'') that have $Q$ as a common boundary. Furthermore,  any two points contained within Ext (Int) can be joined by a continuous curve that does not intersect $Q$, while any continuous curve joining a point of Int to a point of Ext must intersect $Q$. $\square$ 
}

\null 

\noindent We record here a useful consequence of Jordan's Curve Theorem for polygons (JCT).

\null

\noindent \textbf{Corollary 1.1.5} Let $Q$ be a simple polygon and let $p_1$ and $p_2$ be points in the plane not on $Q$. If the segment $[p_1,p_2]$ intersects $Q$ exactly once at a point $q$ that is not a vertex of $Q$, then one of $p_1$ and $p_2$ is in Ext and the other is in Int. $\square$

\null

\noindent \textbf{Definitions.} Let $G$ be a polyomino. We will denote the cycle that bounds the outer face $\text{Outer}(G)$ of a polyomino $G$ by $B_0(G)=B_0$. Define the  \index{enclosure (of a polyomino)|textbf}\textit{enclosure of $G$} to be the set of boxes contained in the region of the plane bounded by $B_0$, and denote it by $\textrm{Encl}(G)$. By JCT, $B_0$ divides the plane into two disjoint sets of boxes: $\text{Outer}(G)$ and  $\textrm{Encl}(G)$. 

We define the \index{boundary (of a polyomino)|textbf}\textit{boundary of $G$} to be the edge-induced subgraph on the set of edges that are incident on a box of $G$ and a box of $\mathbb{Z}^2$ that is not a box of $G$, and denote it by B$(G)$. We define a \index{boundary box|textbf}\textit{boundary box} of $G$ to be a box that is incident on B$(G)$ but that is not a box of $G$. Denote the set of all boundary boxes of $G$ by $\text{BBoxes}(G)$. See Figure 1.2. 

Consider the set of boxes $\textrm{Encl}(G) \setminus \text{Boxes}(G)$. Let $\{\mathcal{O}_1, \ldots, \mathcal{O}_q \}$ be the set of maximal box-path connected sets of boxes of $\textrm{Encl}(G) \setminus \text{Boxes}(G)$. We call each member of $\{\mathcal{O}_1, \ldots, \mathcal{O}_q \}$ a \index{hole|textbf}\textit{hole} of $G$ and we denote the $\textrm{Encl}(G) \setminus \text{Boxes}(G)$ by $\text{Holes}(G)$.

%\textcolor{red}{The proof of JCT in Courant and Robbins uses the ``clear'' fact that any two points on opposite sides of $Q$ but ``close enough" to each other have different parities. The corollary follows from this.} 

%\null

\begingroup 
\setlength{\intextsep}{0pt}
\setlength{\columnsep}{20pt}
\begin{center}
\begin{adjustbox}{trim=0cm 0cm 0cm 0cm}

\end{adjustbox}
\end{center}

%remark about the def of B(G) not holding for graphs more general than polyominoes, not needed.

%Note that an edge is in $\text{B}(G)$ if and only if it is incident on a box of $G$ and a boundary box of $G$. Note that this only holds for polyominoes. It need not hold for any two-connected grid graph. See Figure 1.3 (d). 

\noindent We define a \index{simply connected polyomino|textbf}\textit{simply connected polyomino $G$} to be a polyomino $G$ such that $\text{Boxes}(G)=\textrm{Encl}(G)$. Let $G$ be a simply connected polyomino such that $B(G)$ is an $m \times n$ vertices rectangle. That is, the horizontal side of $B(G)$ has length $m-1$ and the vertical side of $B(G)$ has length $n-1$. Then we call $G$ an \textit{ $m \times n$ rectangular grid graph without holes ($m \times n$ grid graph}). Note that this is the opposite of the usual convention for matrices. We have the following hierarchy of \index{lattice animal}lattice animals: 

\null 

\noindent \index{m x n grid graph@$m \times n$ grid graph|textbf}\textit{$m \times n$ grid graph}\ $m \times n$ grid graphs $\subset$ simply connected polyominoes $\subset$  polyominoes $\subset $ 2-connected grid graphs $\subset $ lattice animals. See Figure 1.3.

\begin{center}
\setlength{\intextsep}{0pt}
\setlength{\columnsep}{20pt}
\begin{adjustbox}{trim=0.5cm 0cm 0cm 0.5cm}

\end{adjustbox}
\end{center}

\noindent \textbf{Observation 1.1.6.} Let $G$ be a polyomino with holes $\{\mathcal{O}_1, \ldots, \mathcal{O}_q \}$. Then:

($a_1$) $\text{Boxes}(G) \cap \text{Outer}(G) =\emptyset$.

($a_2$) $\text{Boxes}(G)\subseteq \textrm{Encl}(G)$. 

($b_1$) For each $i\neq j$, $E(B_0(\mathcal{O}_i)) \cap E(B_0(\mathcal{O}_j)) = \emptyset$.

($b_2$) For each $i \in \{1, ..,q\}$, $E(B_0(G)) \cap E(B_0(\mathcal{O}_i)) = \emptyset$.

(c) For each $i \in \{1, ..,q\}$, $B_0(\mathcal{O}_i) \subset \textrm{Encl}(G)$. 

%(d) $B_0$ is a subgraph of $B(G)$. \textcolor{red}{[Where do I use this? (It used to be Observation 1.1.7)] }

(d)  For each $i \in \{1, ..,q\}$, $\mathcal{O}_i$ is a polyomino. $\square$

% PROOF FOR LEMMA 1.1.7 BEFORE IT BECAME OBVIOUS
{
%\textit{Proof.} Let $v \in B_0$. To show that $v \in B(G)$ we need to find two boxes $X$ and $Y$ both incident on $v$ and such that one box is a box of $G$ and the other is a boundary box of $G$. 

%Let $v'$ be adjacent to $v$ in $B_0$. Let $Y$ be the box in $\text{Outer}(G)$ that is incident on $\{v,v'\}$ and let $X$ be the box in $\text{Encl}(G)$ that is incident on $\{v,v'\}$. First we show that $X \notin \text{Boxes}(G)$ is impossible. By Corollary 1.1.6, $Y \notin \text{Boxes}(G)$. But if neither box is incident on $\{v,v'\}$ is a contradiction to Observation 1.1.3 (a). So we must have $X \in \text{Boxes}(G)$. Now $v \in V(Y)$  and $v\in V(G)$ but, by Observation 1.1.5, $Y$ is not a box of $G$. So $Y$ is a boundary box of $G$. $\square$
}

% OMITTED PROOF OF OBSERVATION 1.1.5 [[ FOR (a_1) in current version]]
{
%\textit{Proof.} BWOC let $X \in \text{Boxes}(G) \cap \text{Outer}(G)$. Either $X=\text{Outer}(G)$ or $X$ is properly contained in $\text{Outer}(G)$. The former contradicts that $\text{Outer}(G)$ is the unbounded face of $G$. The latter contradicts that $\text{Outer}(G)$ is connected since the region contained inside the box $X$ is not connected with the region contained in $\text{Outer}(G)\setminus X$. \textcolor{red}{[A face is maximally connected region of the plane].} Either way, we must have  that $\text{Boxes}(G) \cap \text{Outer}(G) =\emptyset$. $\square$
}

\null 

\noindent \textbf{Lemma 1.1.7.} Let $G$ be a polyomino with holes $\{\mathcal{O}_1, \ldots, \mathcal{O}_q \}$. Then:

  (a) An edge $e$ is in $B_0(G)$ iff $e$ is incident on a box of $\text{Outer}(G)$ and on a box of 
  
\hspace*{0.6cm} $\text{Boxes}(G)$. 

 (b) For each $i \in \{1, ..,q\}$, $E(B_0(\mathcal{O}_i)) \subset E(B(G))$. Furthermore, every edge  
 
\hspace*{0.6cm} $e \in E(B_0(\mathcal{O}_i))$ is incident on a box of $\mathcal{O}_i$ and on a box of $G$.

 (c) For each $i \in \{1, ..,q\}$, $\text{Encl}(\mathcal{O}_i) \subset \textrm{Encl}(G)$.

 (d) For every edge $e_X \in B_0(G)$ there is a box $X \in \text{Boxes(G)} \setminus \bigcup_{i=1}^q \text{Encl}(\mathcal{O}_i)$ such 
 
\hspace*{0.6cm} that $X$ is incident on $e_X$.

 (e) For each $i \in \{1, ..,q\}$, $\text{Boxes(G)} \cap\text{Encl}(\mathcal{O}_i) = \emptyset$.

 (f) The holes of $G$ are simply connected polyominoes.

 \null 

\noindent\textit{Proof of (a).} Let $e$ be an edge of $B_0$. Let $Y$ be the box contained in $\text{Outer}(G)$ that is incident on $e$, and let $X$ be the box contained in $\text{Encl}(G)$ that is incident on $e$. By Observation 1.1.3(a), $X$ must be a box in $\text{Boxes}(G)$.
Conversely, assume that $e$ is incident on a box of $\text{Outer}(G)$ and a box of $\text{Encl}(G)$. Then $e$ is contained in the cycle that bounds $\text{Outer}(G)$. That is, $e$ is an edge of $B_0$. End of proof for (a).

\null

\noindent\textit{Proof of (b).} Let $e \in E(B_0(\mathcal{O}_i))$. By Observation 1.1.6(d), $\mathcal{O}_i$ is a polyomino. By part (a), $e$ is incident on a box of $\mathcal{O}_i$ and a box of $\text{Outer}(\mathcal{O}_i)$. Let $X$ be the box of $\text{Outer}(\mathcal{O}_i)$ incident on $e$, and let $Y$ be the box of $\mathcal{O}_i$ incident on $e$. By Observation 1.1.6(b$_2$), $X$ cannot be in $\text{Outer}(G)$. Otherwise, since $Y \in \mathcal{O}_i \subset \text{Encl}(G)$ and $X \in \text{Outer}(G)$, it would follow that $e \in B_0(G)$. This would imply $e \in E(B_0(G)) \cap E(B_0(\mathcal{O}_i)) \neq \emptyset$, contradicting (b$_2$). Therefore, $X \in \text{Encl}(G)$. Note that if $X \in \text{Encl}(G) \setminus \text{Boxes}(G)$, this would contradict the maximality of $\mathcal{O}_i$ (since $X$ is adjacent to $\mathcal{O}_i$, we could simply include it in $\mathcal{O}_i$). Thus, $X$ must be a box of $G$. Therefore, by definition, $e \in E(B(G))$. End of proof for (b).

\null 

%\noindent\textit{Proof of (c).} Let $X \in \text{BBoxes}(\mathcal{O}_i)$. Then $X \notin \mathcal{O}_i$, and $X$ is incident on an edge $e$ of $B(\mathcal{O}_i)$. Let $Y$ be the other box incident on $e$. By part (b), $Y$ must belongs to $\mathcal{O}_i$. By part (b) again, $X$ must belong to $G$. End of proof for (c).

\noindent \textit{Proof of (c).} By Observation 1.1.3(b), $\text{Encl}(\mathcal{O}_i)$ is a polyomino. Let $X \in \text{Encl}(\mathcal{O}_i)$. Then $X \in \text{Encl}(G)$ or $X \in \text{Outer}(G)$. Assume, for contradiction, that $X \in \text{Outer}(G)$. Let $Y \in \text{Encl}(\mathcal{O}_i)$ be adjacent to a boundary box $Y'$ of $\mathcal{O}_i$. Consider the \index{path of boxes}path of boxes $P(X,Y)$. This path is either contained in $\text{Outer}(G)$ or not.

\null 

\noindent \textit{CASE 1:} $P(X,Y) \subset \text{Outer}(G)$. Let $e$ be the edge that $Y$ and $Y'$ share. Then $e \in E(B_0(\mathcal{O}_i))$. By Lemma 1.1.7(b), $e \in E(B(G))$. Since $Y \in \text{Outer}(G)$, by Lemma 1.1.7(a), $e \in E(B_0(G))$. But this contradicts Observation 1.1.6(b$_2$). End of Case 1.

\null 

\noindent \textit{CASE 2:} $P(X,Y) \not\subset \text{Outer}(G)$. Let $P(X,Y) = Z_1, \dots, Z_s$. Then there must be a box $Z_j$ in $P(X,Y)$ such that $P(Z_1, Z_j)$ is entirely contained in $\text{Outer}(G)$, while $P(Z_1, Z_{j+1})$ is not. Let $f$ be the edge that $Z_j$ and $Z_{j+1}$ share. Then, by definition, $f \in E(B_0(G))$, and by Lemma 1.1.7(a), $Z_{j+1}$ is a box of $G$. Since $Z_j \in \text{Encl}(\mathcal{O}_i)$ and $Z_{j+1} \in \text{Boxes}(G)$, it follows that $f \in E(B_0(\mathcal{O}_i))$. But this contradicts Observation 1.1.6(b$_2$).  End of proof for (c). 

\null 

\noindent \textit{Proof of (d).} Let $Z \in \text{Boxes}(G)$ be a box incident on an edge $e \in B_0(G)$, and let $Z'$ be the other box incident on $e$. If $Z$ is in $\text{Boxes}(G) \setminus \bigcup_{i=1}^q \text{Encl}(\mathcal{O}_i)$, we are done. We will now show that the alternative is impossible.

For contradiction, assume that there is some $i \in \{1, \dots, q \}$ such that $Z \in \text{Boxes}(G) \cap \text{Boxes}(\mathcal{O}_i)$. By definition, $Z' \in \text{Outer}(G)$. If there is some $i \in \{1, \dots, q\}$ such that $e \in B_0(\mathcal{O}_i)$, this contradicts Observation 1.1.6(b$_2$). Otherwise, if no such $i$ exists, then $Z'$ must also belong to $\text{Boxes}(\mathcal{O}_i)$. But then, by part (c), $Z' \in \text{Encl}(G)$, which contradicts the assumption that $Z' \in \text{Outer}(G)$. End of proof for (d). 

\null 

\noindent \textit{Proof of (e).} Assume, for contradiction, that there is some $i \in \{1, \dots, q\}$ such that $\text{Boxes}(G) \cap \text{Encl}(\mathcal{O}_i) \neq \emptyset$. Let $X$ be in $\text{Boxes}(G) \cap \text{Encl}(\mathcal{O}_i)$. By part (d), $\text{Boxes}(G) \setminus \bigcup_{i=1}^q \text{Encl}(\mathcal{O}_i) \neq \emptyset$. Let $Y \in \text{Boxes}(G) \setminus \bigcup_{i=1}^q \text{Encl}(\mathcal{O}_i)$. 

Since $Y \in \text{Boxes}(G) \setminus \bigcup_{i=1}^q \text{Encl}(\mathcal{O}_i)$, and $X \in \text{Boxes}(G) \cap \text{Encl}(\mathcal{O}_i)$, $X$ and $Y$ are on different sides of $B_0(\mathcal{O}_i)$, and they both belong to $\text{Boxes}(G)$. Then there is a path of boxes $P(X,Y)$ in $G$ with $Z_1=X$, and so that it contains $t$ boxes, with $t\geq 2$.  By the Jordan Curve Theorem, there is $j \in \{1, \dots, t-1\}$ such that $P(Z_1, Z_j) \subset \text{Encl}(\mathcal{O}_i)$, but $P(Z_1, Z_{j+1}) \not\subset \text{Encl}(\mathcal{O}_i)$. 

It follows that the gluing edge $e$ between $Z_j$ and $Z_{j+1}$ is in $B_0(\mathcal{O}_i)$. Then, by part (b), one of $Z_j$ and $Z_{j+1}$ is in $G$, and the other is in $\mathcal{O}_i$. But this contradicts the assumption that $P(X,Y)$ is entirely contained in $G$. End of proof for (e). 

\null 

\noindent \textit{Proof of (f).} We will check that for each $i \in  \{1, \ldots, q\}$, $\text{Encl}(\mathcal{O}_i) \subset \text{Boxes}(\mathcal{O}_i)$. If $\text{Encl}(\mathcal{O}_i)$ has no holes, we are done, so assume that $\text{Encl}(\mathcal{O}_i)$ has a hole $U$. By part (b), there is an edge $e \in B_0(U)$ that is incident on the boxes $Z \in U$ and $Z' \in \mathcal{O}_i$. By maximality of holes (of $G$, in this case), $Z$ does not belong to any hole of $G$. Then $Z \in \text{Boxes}(G)$ or $Z \in \text{Outer}(G)$. If the former, then, since $Z$ is also in $\text{Encl}(\mathcal{O}_i)$, this contradicts part (e). And if the latter, then, since $Z' \in \text{Boxes}(\mathcal{O}_i) \subset \text{Encl}(\mathcal{O}_i)$, $e$ belongs to $B_0(G)$. But then, since neither $Z$, nor $Z'$ belong to $G$, this contradicts part (a). End of proof for (f). $\square$

\null

\noindent \textbf{Lemma 1.1.8}. Let $G$ be a simply connected polyomino. Then the boundary of $G$ is $B_0$.

\noindent \textit{Proof.} Let $e\in B(G)$. Then $e$ is incident on a box $X$ of $G$ and on a boundary box $Y$ of $G$. By Observation 1.1.6($a_2$), $X \in \text{Boxes}(G)$ implies that $X \in \text{Encl}(G)$. Either $Y$ belongs to $\textrm{Encl}(G)$ or $Y$ belongs to $\text{Outer}(G)$. If $Y \in \textrm{Encl}(G)$, since $G$ is a simply connected polyomino, then $Y$ is a box of $G$. But then $Y$ is both a box of $G$ and a boundary box of $G$, which is a contradiction. It must be the case that $Y \in \text{Outer}(G)$. By Lemma 1.1.7(a), $e$ belongs to $B_0$. $\square$

\null

\noindent \textbf{Observation 1.1.9.} Let $G$ be a polyomino with holes $\{\mathcal{O}_1, \ldots, \mathcal{O}_q \}$, bounded by the cycles $B_1, \ldots, B_q$. It follows from Lemma 1.1.7(f) that $B(G)$ is the disjoint union of $B_0, B_1, \ldots, B_q$. For each $i \in \{0, \ldots,q \}$, we call $B_i$ a \index{boundary component|textbf}\textit{boundary component of $G$}.

\null 

% Remark 1.1.12, not used

{
%\textbf{Remark 1.1.12.} Let $G$ be a polyomino  with holes $\{\mathcal{O}_1, \ldots, \mathcal{O}_1 \}$.  Then, for each $i \in \{1, \ldots,q \}$, $\text{Boxes(B}_i(G))$ is box-path-connected in $\mathcal{O}_i$ and $\text{Boxes(B}_0(G))$ is box-path-connected in $\text{Outer}(G)$. $\square$

%\textcolor{red}{$\downarrow$  \_\_\_\_\_\_\_\_\_\_\_\_\_\_\_\_\_\_\_\_\_\_\_\_\_\_\_\_\_\_\_\_\_\_\_\_\_\_\_\_\_\_\_\_\_\_\_\_  ADDED AFTER JULY \_\_\_\_\_\_\_\_\_\_\_\_\_\_\_\_\_\_\_\_\_\_\_\_\_\_\_\_\_\_\_\_\_\_\_\_\_\_\_\_\_\_\_\_\_\_\_\_$\downarrow$  }

}

\noindent \textbf{Corollary 1.1.10}. $G$ is a simply connected polyomino if and only if $B(G)=B_0(G)$.

\null

\noindent \textit{Proof.} By Lemma 1.1.8 we only need to show that if $B(G)=B_0(G)$ then $G$ is a simply connected polyomino. We prove the contrapositive. Suppose that $G$ is not simply connected. Then the set $D=\textrm{Encl}(G) \setminus \text{Boxes}(G)$ is nonempty. Let $\{\mathcal{O}_1, \ldots, \mathcal{O}_q \}$
and $B_0, B_1, \ldots, B_q$ be as in Observation 1.1.9. Then $q \geq 1$. By Observation 1.1.9, $B_1$ and $B_0$ are disjoint. It follows that  $B(G)\neq B_0(G)$. $\square$

\null

\noindent \textbf{Lemma 1.1.11}. Simply connected polyominos have no \index{polyking junction}polyking junctions.

\null 

\noindent \textit{Proof.} Let $G$ be a simply connected polyomino. By Corollary 1.1.10 $B(G)=B_0(G)$. By Lemma 1.1.1, $B_0(G)$ is a cycle. For a contradiction, assume that $G$ has a polyking junction at a vertex $v$. Observe that this implies that all edges incident on $v$ are boundary edges of $G$. But then $B_0(G)$ has a vertex of degree 4, contradicting that it is a cycle. $\square$

%\textcolor{red}{$\uparrow$  \_\_\_\_\_\_\_\_\_\_\_\_\_\_\_\_\_\_\_\_\_\_\_\_\_\_\_\_\_\_\_\_\_\_\_\_\_\_\_\_\_\_\_\_\_\_\_\_  ADDED AFTER JULY \_\_\_\_\_\_\_\_\_\_\_\_\_\_\_\_\_\_\_\_\_\_\_\_\_\_\_\_\_\_\_\_\_\_\_\_\_\_\_\_\_\_\_\_\_\_\_\_$\uparrow$ }

\null

\noindent 

{
}

% Relocate or delete.

%\textcolor{blue}{\textbullet}  

\subsection{The $H$-decomposition of $G$ and the follow-the-wall construction}

\noindent This section describes the $H$-decomposition of $G$ in detail and proves that $H$-components have a tree-like structure. We introduce the follow-the-wall construction, a procedure for navigating \index{H-component@$H$-component}$H$-components that appears extensively in later proofs (Chapters 1, 3, 4, and 5).

\null 

\begingroup
\setlength{\intextsep}{0pt}
\setlength{\columnsep}{20pt}

\begin{adjustbox}{trim=0cm 0cm 0cm 0cm}

\end{adjustbox}

\null

\noindent \textbf{Definitions.} A \index{walk|textbf}\textit{walk} (of length $r$) in a graph $G$ is an alternating sequence $v_0e_1v_1e_2 \ldots e_rv_r$ of vertices and edges. Define a \index{lazy walk|textbf}\textit{lazy walk} to be a sequence of edges and vertices where every edge is in between two vertices and in between every two edges there is a vertex or multiple copies of a vertex. That is, a lazy walk is roughly a walk in which consecutive vertices can be the same, allowing the walk to remain at a vertex for one or more steps without traversing any edges.

Let $G$ be any polyomino and let $H$ be any subgraph in $G$. Let $X_1,X_2$ be two adjacent boxes of $G$. If $E(X_1)\cap E(X_2) \cap E(H)=\emptyset$, we say that $X_1$ and $X_2$ are \index{H-neighbour@$H$-neighbour|textbf}\textit{$H$-neighbours} or $X_1$ is \index{H-adjacent@$H$-adjacent|textbf}\textit{$H$-adjacent} to $X_2$. Define an \index{H-walk@$H$-walk|textbf}\textit{H-walk of boxes in G} ($H$-walk) to be a sequence $X_1, \ldots,X_r$ of boxes in $G$, not necessarily distinct, such that for all $j \in \{1,2, \ldots, r-1 \}$, $X_j$ is an $H$-neighbour of $X_{j+1}$ or $X_j = X_{j+1}$. 
%gluing edge extra definitions 
{
%and denote it by $W(X_1,X_r)$. For each $j \in \{1,\ldots,r-1\}$, whenever $X_j \neq X_{j+1}$, we call the edge of $G$ that $X_j$ and $X_{j+1}$ share a \textit{gluing edge} of $W(X_1,X_r)$, and denote it by $e(X_j,X_{j+1})$. We denote the set of gluing edges of $W(X_1,X_r)$ by $E_{\text{glue}}(W(X_1,X_r))$. 
} 
If for all $i,j \in \{1,2, \ldots, r \}$ with $i \neq j$, $X_i$ is distinct from $X_j$, we call the sequence an \index{H-path@$H$-path|textbf}\textit{H-path of boxes in G} and denote it by $P(X_1,X_r)$. If in addition, for all $i,j \in \{1,2, \ldots, r-1 \}$ with $ i \neq j$, $E(X_i)\cap E(X_j) \cap E(H)=\emptyset$, we say that $P(X_1,X_r)$ is a \index{non-self-adjacent|textbf}\textit{non-self-adjacent $H$-path of boxes in G}. Otherwise, if $E(X_i)\cap E(X_j) \cap E(H) \neq \emptyset$, we say that $P(X_1,X_r)$ is a \index{self-adjacent|textbf}\textit{self-adjacent $H$-path of boxes in G}. 
% self adjacent walks definition 
{%\textcolor{red}{[We define non-self-adjacent and self-adjacent $H$-walks of boxes analogously. (This may not be needed.)]}  
}
Let $r \geq 4$. Define an \index{H-cycle@$H$-cycle|textbf}\textit{H-cycle $C$ of boxes in G} ($H$-cycle) to be a set $X_1, X_2, \ldots, X_r=X_1$ of boxes in $G$ such that for each $j \in \{1,2, \ldots, r-1 \}$, $X_j$ is an $H$-neighbour of $X_{j+1}$ and the boxes $X_1, \ldots, X_{r-1}$ are distinct. We note that every box of $C$ has exactly two gluing edges. Proposition 1.2.1. will show that if $H$ is a Hamiltonian cycle of $G$, then there are no $H$-cycles of boxes in $G$.

By definition, $\text{BBoxes}(G) \cap \text{Boxes}(G) = \emptyset$. We use the notation $G_{-1}$ to denote the union of the graph on $\text{BBoxes}(G)$ and $G$. We extend the definitions of $H$-walks, $H$-cycles, $H$-paths, and self-adjacency to $G_{-1}$. See Figure 1.2, where $G_{-1}$ consists of the green boxes and blue boxes

% Pragraph removed
{
%We remark here that the definition of an $H$-walk in a polyomino $G$ coincides with the definition of a lazy walk in the plane dual graph of $G$ (with the added $H$-adjacency requirement). \textcolor{blue}{[Reference Diestel's ``Graph Theory book here for definition of plane dual graph?]} However, we will not work directly with the dual graph, as it adds no value to the discussion.
}

\null 

\begingroup
\setlength{\intextsep}{0pt}
\setlength{\columnsep}{20pt}
\begin{wrapfigure}[]{l}{0cm}
\begin{adjustbox}{trim=0cm 0cm 0cm 0.25cm}
\begin{tikzpicture}[scale=1.35]

\draw[gray,very thin, step=0.5cm, opacity=0.5] (0,0) grid (4,3);

\fill[green!50!white,opacity=0.5](0.5, 0.5) rectangle  (1,2.5); 
\fill[green!50!white,opacity=0.5](3, 0.5) rectangle  (3.5,2.5); 
\fill[green!50!white,opacity=0.5](1, 0.5) rectangle  (3,1); 
\fill[green!50!white,opacity=0.5](1, 2) rectangle  (3,2.5); 

\draw[green!50!black, line width=0.5mm] 
(0.75,0.75)--++(0,1.5)--++(2.5,0)--++(0,-1.5)--++(-2.5,0);

\draw [blue, line width=1mm] (2.5,2)--++(0,0.5);

\foreach \x in {0,...,4}
\draw[orange, line width=0.5mm] 
(1+0.5*\x,0.5)--++(0,0.5);

\foreach \x in {0,...,4}
\draw[orange, line width=0.5mm] 
(1+0.5*\x,2)--++(0,0.5);

\foreach \x in {0,...,2}
\draw[orange, line width=0.5mm] 
(0.5,1+0.5*\x)--++(0.5,0);

\foreach \x in {0,...,2}
\draw[orange, line width=0.5mm] 
(3,1+0.5*\x)--++(0.5,0);

\draw [blue, line width=1mm] plot [smooth, tension=0.75] coordinates {(1.5,2)(2.25,1.75)(2.5,2)};
\draw [blue, line width=1mm] plot [smooth, tension=0.75] coordinates {(2.5,2.5)(2.25,2.7)(1.75,2.75)(1,2.5)};

\draw[fill=blue, opacity=1] (1.5,2) circle [radius=0.05];
\draw[fill=blue, opacity=1] (1,2.5) circle [radius=0.05];

\node[right] at (2.45,2.35) [scale=1] {\small{$e$}};
\node[below] at (1.5,2) [scale=1] {\small{$a$}};
\node[above] at (1,2.5) [scale=1] {\small{$b$}};
\node at (1.25,2.25) [scale=1] {\small{$X_i$}};

\node[below, align=left, text width=6cm] at (2.15,-0.25)
{Fig. 1.5. $P(a,b)$ in blue, gluing edges  of $C$ in orange, $Q$ in dark  green, $C$ shaded light green. };

\end{tikzpicture}
\end{adjustbox}
\end{wrapfigure}

\noindent \textbf{Proposition 1.2.1.} Let $G$ be a polyomino and let $H$ be a Hamiltonian path of $G$. Then every $H$-cycle in $G_{-1}$ is contained in $\textrm{BBoxes}(G)$. 

\null

\noindent \textit{Proof.} 
For a contradiction, assume that there is an $H$-cycle $C$ in $G_{-1}$ with boxes $X_1, X_2, \ldots, X_r=X_1$ that has a box $X_i$ contained in $G$. Let $c_1, c_2, \ldots, c_r=c_1$ be the centers of the boxes of $C$. That is, if $X_j=R(k,l)$ then $c_j=(k+\frac{1}{2}, l+\frac{1}{2})$. Observe that for each $j \in \{1,2,\ldots,r-1\}$ $\{1,2,\ldots,r-1\}$, $[c_j,c_{j+1}]$ intersects the gluing edge of $X_j$ and $X_{j+1}$ and $[c_j,c_{j+1}]$ intersects no other edge of $G_{-1}$.

We will first show that the set of segments $[c_j,c_{j+1}]$ is a non-self-intersecting polygon $Q$. Since $c_1=c_r$, $Q$ is a polygon. For a contradiction, assume that $Q$ is self intersecting, so there are points $c_{j-1}$, $c_j$, $c_{j+1}$ and $c_i$ such that the segments $[c_{j-1},c_j]$, $[c_j,c_{j+1}]$ and $[c_j,c_i]$ are edges of $Q$. But then $X_j$ has three gluing edges, a contradiction.

By JCT, $Q$ divides the plane into two subsets Int and Ext such that any two points within a subset can be joined by a path that does not intersect $Q$ while a path joining a point of Int to a point of Ext must intersect $Q$. Let $V(\textrm{Int})$ be the vertices of $G_{-1}$ contained in Int and let $V(\textrm{Ext})$ be the vertices of $G_{-1}$  contained in Ext. Note that any box of $C$ contains at least one vertex in $V(\textrm{Int})$ and one vertex in $V(\textrm{Ext})$, so both sets are nonempty. In particular, this is true for $X_i$. For definiteness, assume that $a \in V(X_i)$ is contained in $V(\textrm{Int})$ and $b \in V(X_i)$ is contained in $V(\textrm{Ext})$. By JCT, $V(\textrm{Int}) \cup V(\textrm{Ext})=V(G)$ and $V(\textrm{Int}) \cap V(\textrm{Ext})=\emptyset$. %, otherwise we can find a path from Int to Ext that does not intersect $Q$.
Consider the subpath $P(a,b)$ of $H$. By JCT again, there is an edge $e \in P(a,b) \subset H$ intersecting $Q$ at some segment $[c_j, c_{j+1}]$. But then $e$ is a \index{gluing edge}gluing edge of $C$, so $e$ cannot belong to $H$. See Figure 1.5. $\square$

\endgroup 

\null

\noindent Let $G$ be a polyomino and let $H$ be a Hamiltonian path or cycle of $G$. Let $X_1,\ldots, X_r$ be a set of boxes in $G$ such that for any $i,j \in \{1,2, \ldots, r \}$, there is an $H$-path $P(X_i,X_j)$ between $X_i$ and $X_j$ contained in $G$. We say that $X_i$ and $X_j$ are \index{H-path-connected@$H$-path-connected|textbf}\textit{$H$-path-connected} in $G$ and the set of boxes $\{ X_1,\ldots, X_r \}$ is an \textit{$H$-path connected set of boxes in $G$}. If an $H$-path-connected set of boxes in $G$ contains no cycles of boxes we call it an \index{H-tree@$H$-tree|textbf}\textit{$H$-tree}. An \index{H-component@$H$-component|textbf}\textit{$H$-component} of $G$ is a maximal $H$-path connected set of boxes.

Let $J$ be an $H$-subtree of an $H$-component of $G$. We say that $J$ is non-self-adjacent if it contains no self-adjacent paths of boxes. Otherwise we say that $J$ is self-adjacent. We define a \index{neck (of an H-component)@neck (of an $H$-component)|textbf}\textit{neck of $J$} to be a box $N_J$ of $J$ that is incident on a boundary edge $e_J$ of $G$ such that $e_J \notin H$. We call $e_J$ a \index{neck edge|textbf}\textit{neck-edge} of $J$. Note that the other box incident on $e_J$ must be in $G_{-1}\setminus G$. See Figure 1.4.

%\textcolor{blue}{Any pair of boxes of $G_{-1}$ is $H$-path connected in $G_{-1}$. So, when we say that a pair of boxes is $H$-path connected, we mean that pair of boxes is $H$-path connected in $G$.}

\null 

\noindent \textbf{Corollary 1.2.2.}  Let $G$ be a polyomino and let $H$ be a Hamiltonian path of $G$. Then the $H$-path $P(X,Y)$ in $G$ is unique, if it exists. Furthermore, $H$ partitions the boxes of $G$ into $H$-path-connected $H$-components which are maximal \index{H-tree@$H$-tree}\textit{$H$-trees}. $\square$

\null 

\noindent The uniqueness of $H$-paths contained in $G$ is a key structural property of the decomposition that a Hamiltonian path $H$ imposes on the boxes of a polyomino $G$. We note that such an $H$-path exists if and only if $X$ and $Y$ belong to the same $H$-component. Most $H$-paths we consider in the results that follow lie entirely within a single $H$-component. We consider $H$-paths that contain boxes in $G_{-1}$ only a handful of times. This uniqueness property is used implicitly throughout the dissertation—appearing tacitly to avoid cluttering the text with excessive references.

\null

\noindent \textbf{Corollary 1.2.3}\footnote{This Corollary is only used once by Lemma 1.3.1.}  Let $G$ be a polyomino and let $H$ be a Hamiltonian path of $G$. Then every $H$-cycle in $\mathbb{Z}^2$ is contained in $\text{Boxes} (\mathbb{Z}^2) \setminus \text{Boxes}(G)$. $\square$

\null 

%We will call the partition that $H$ imparts of $G$ the \textit{$H$-decomposition of $G$}. 

%\textcolor{red}{ \_\_\_\_\_\_\_\_\_\_\_\_\_\_\_\_\_\_\_\_\_\_\_\_\_\_COROLLARY 1.2.3 is sent to SECTION 1 TEMP  \_\_\_\_\_\_\_\_\_\_\_\_\_\_\_\_\_\_\_\_\_\_\_\_ }

%DELETED from Corollary 1.4
%Consider the boxes $Y_1$ and $Y_s$ in $G_{-1}$ that are $H$-adjacent to $X_1$ and $X_r$, respectively. Note that there is an $H$-path $P(Y_s, Y_1)$ contained in $G_{-1} \setminus G_0$ and that $C=P(X_1,X_r), P(Y_1, Y_s)$ is an $H$-cycle of $G_{-1}$ \textcolor{blue}{(added the index -1 to $G_{-1}$)}. 

\noindent \textbf{Definitions.} Let $G$ be a graph. A \index{trail|textbf}\textit{trail} is a walk in $G$ where all edges are distinct. A trail where the first and last vertices coincide is called a \textit{closed trail} or a \index{circuit|textbf}\textit{circuit}. A \index{directed walk|textbf}\textit{directed walk} (of length $s$) is an alternating sequence $v_0e_1v_1e_2\ldots e_sv_s$ of vertices and directed edges such that for $j \in \{1,\ldots,s\}$, the directed edge $e_j$ has tail $v_{j-1}$ and head $v_j$. A \index{directed trail|textbf}\textit{directed trail} is a directed walk where all directed edges are distinct. Note that the edges $(u,v)$ and $(v,u)$ in a directed trail are distinct edges. Similarly, a directed trail $v_0e_1v_1e_2 \dots e_sv_s$ and its reversal $v_se_s \dots e_2v_2e_1v_1$ are distinct directed trails. We will use the notation $\overrightarrow{K}$ to denote directed trails. 

\begin{figure}[t]
\begingroup
\setlength{\intextsep}{0pt}
\setlength{\columnsep}{20pt}
\begin{center}
\begin{adjustbox}{trim=0cm 0cm 0cm 0cm}

\end{adjustbox}

\end{center}
\endgroup 
\end{figure}

\noindent Let the box $X$ be incident on a directed edge $(u,v)$ of the integer grid. Then $X$ and its vertices not incident on the edge $e_j$ are either on the \index{right (side of)|textbf}\textit{right} or on the \index{left (side of)|textbf}\textit{left} side of $(u,v)$. See Figure 1.6 for an illustration and Note 1.2.5 at the end of this section for a more precise definition of a box being on the right or left side of a directed edge. Note that if $X$ is on the right side of $(u,v)$, then the other box incident on $(u,v)$, say $X'$, is on the right side of $(v,u)$.

\noindent Let $H$ be a path in a polyomino $G$, that is not necessarily Hamiltonian. Consider the directed multigraph $H^*$ where each edge $\{u,v\}$ of $H$ is replaced by the directed edges $(u,v)$ and $(v,u)$. Define an \index{H star trail@$H^*$-trail|textbf}\textit{$H^*$-trail} to be a directed trail $\overrightarrow{K}=e_1, \ldots, e_s$ in $H^*$, such that if the edge $(u,v)=e_j$ and $(v,u)=e_{j+1}$, then $v$ is an end-vertex of $H$. We call a closed $H^*$-trail an \textit{$H^*$-circuit}.  So, $H^*$-trails are trails that can only ``turn around" at an end-vertex.

\null 

\noindent Now we will use an $H^*$-trail $\overrightarrow{K}$ to construct an $H$-walk of boxes in {$G_{-1}$} that we will call  \index{the right H-walk induced by K@the right $H$-walk induced by $\overrightarrow{K}$|textbf}\textit{the right $H$-walk induced by $\overrightarrow{K}$} and denote it by $W_{\textrm{right}}(X_1,X_r)$, where $X_1$ and $X_r$ are the end-boxes of $W_{\textrm{right}}(X_1,X_r)$. Roughly, $W_{\textrm{right}}(X_1,X_r)$ will be the \index{walk of boxes}walk of boxes that a ``walker'' would encounter as they followed along the side of $\overrightarrow{K}$ when starting on the right side of the first edge $e_1$ of $\overrightarrow{K}$. We will call this construction the \index{follow-the-wall (FTW)|textbf}\textit{follow-the-wall} construction (FTW). This is very similar to the well-known hand-on-the-wall maze-solving algorithm. 

Let $X$ be the box of $G_{-1}$ on the right of the edge $e_1$ of $\overrightarrow{K}$ and let $e_j$ be the $j^{\textrm{th}}$ edge of $\overrightarrow{K}$. Then $X=X_1$ is the first box of the $H$-walk. Let $e_j=(u,v)$, $e_{j+1}=(v,w)$ and let $X_i$ be on the right of the edge $e_j$. There are four possibilities for the position of $w$ with respect to $(u,v)$: $w$ is right of $(u,v)$, $w \neq u$ is colinear with $(u,v)$, $w$ is left of $(u,v)$, and $w=u$. For the last two cases define $e_j'=(v, v')$ to be the edge in $G_{-1} \setminus H$ that is colinear with $e_j$ and set $e_j''=(v', v)$. See Figure 1.6.

%We call $e_j'$ and $e_j''$ \textit{ghost} edges of $\overrightarrow{K}$.

\begin{itemize}
    \item [(i).] \textit{$w$ is right of $e_j$.} Then $X_{i+1}$ is on the right of the edge $e_{j+1}$. Note that in this case, the walk has a repeated box since $X_i=X_{i+1}.$
    \item [(ii).] \textit{$w \neq u$ is collinear with $e_j$.} Then $X_{i+1}$ is on the right of the edge $e_{j+1}$. 
    \item [(iii).] \textit{$w$ is left of $e_j$.} Then $X_{i+1}$ is on the right of the edge $e_j'$ and $X_{i+2}$ is on the right of the edge $e_{j+1}$.
    \item [(iv).] \textit{$w=u$.} Then $X_{i+1}$ is on the right of the edge $e_j'$ and $X_{i+2}$ is on the right of the edge $e_j''$ and $X_{i+3}$ is on the right of the edge $e_{j+1}$.
\end{itemize}

\noindent We say that the edge $e_{j+1}$ \textit{adds} to the $H$-walk $W_{\textrm{right}}(X_1,X_r)$ the box $X_{i+1}$, in cases (i) and (ii), boxes $X_{i+1}$ and $X_{i+2}$ in case (iii), and boxes $X_{i+1}$, $X_{i+2}$ and $X_{i+3}$ in case (iv). If $e_s$ adds more than one box to $W_{\textrm{right}}(X_1,X_r)$ then we will adopt the convention that $X_r$ is the last box added by the edge $e_s$ and that all the boxes added by each edge of $e_j$, $j\in \{1, \ldots, s \}$ are on the right side of $e_j$. Note that this convention is necessary for the box $X_{i+1}$ in Case (iv). We remark that the first edge $e_1$ can only add the single box $X_1$. The left $H$-walk $W_{\textrm{left}}(X_1,X_r)$ induced by $\overrightarrow{K}$ can be constructed analogously. 
% not needed unless I must refer to an H-walk without its inuducing trail
{
%Fix a $\text{side} \in \{\text{right},\text{left}\}$. Note that a box of $G$ may be visited multiple times, consecutively and non-consecutively, by the $H$-walk $W(X_1,X_r)_{\text{side}}$. 

%Let $Y_1, \ldots, Y_s$ be the distinct boxes of $G$ that are visited by $W(X_1,X_r)$, ordered by their first appearance in $W(X_1,X_r)$. Then the  $H$-walk $W(X_1,X_r)$ begins at with the first appearance of $Y_1$ in $W_{H^*}$ and ending at the last appearance of $Y_s$ in $W_{H^*}$.
}

It is straightforward to extend FTW to the case where $H$ is a cycle in $G$. We remark that in this case, no two directed edges of the $H^*$-trail may correspond to the same (undirected) edge of $H$ and that, since cycles have no end-vertices, Case (iv) never occurs. In fact, we can orient any path (or cycle) $H$ in $G$ to obtain a directed trail $\overrightarrow{K}$ and then use FTW to obtain an $H$-walk from $\overrightarrow{K}$. 
    
 Let $H$ be a path or cycle in $G$. Let $\mathcal{K}(H^*)$ be the set of all $H^*$-trails. We can view the FTW construction as a function $\Phi$ that assigns an $H$-walk to elements of $\mathcal{K}(H^*) \times \{\text{right, left}\}$. We will take a closer look at $H^*$-trails in the case where $H$ is a Hamiltonian path or a Hamiltonian 
cycle of $G$. 

\null 

\begingroup
\setlength{\intextsep}{0pt}
\setlength{\columnsep}{20pt}

\begin{adjustbox}{trim=0cm 0cm 0cm 0cm}

\end{adjustbox}

\null

\noindent Let $H=v_1, v_2,   \dots ,  v_r$ be a Hamiltonian path of $G$. Consider an Eulerian circuit $(v_j, v_{j+1})$, $(v_{j+1}, v_{j+2}),   \ldots,$  $(v_{r-1}, v_r)$, $(v_r,v_{r-1}), \ldots,$  $(v_2,v_1)$, $ (v_1,v_2), \ldots,$ $(v_{j-1}, v_j)$ of $H^*$. Note that there are $2r-2$ such distinct Eulerian circuits, one for each possible first edge. Observe that any subtrail of such a circuit is completely determined by its first and last edges. Therefore it will be fitting to use the notation $\overrightarrow{K}(e_s,e_t)$ to denote the unique subtrail starting at edge $e_s$ and ending at edge $e_t$. We will use the notation $\overrightarrow{K}_{H^*}(e_j)$ to denote an Eulerian circuit of $H^*$ starting at $e_j$ and abbreviate to $\overrightarrow{K}_{H^*}$ whenever the first edge $e_j$ is not relevant to our argument. Fix a side in $\{\text{right, left}\}$. We will use the notation $W_{\text{side},H^*}$ to denote $\Phi( \overrightarrow{K}_{H^*}, \text{side})=W_{\text{side},H^*}$.

%pagemarker
\noindent Consider an $H$-circuit $\overrightarrow{K}_{H^*}(e_{j+1})=\overrightarrow{K}_{H^*}$, starting at $e_{j+1}$ and ending at $e_j$. Fix a side in $ \{\text{right, left}\}$. If $e_j$ and $e_{j+1}$ are as in cases (iii) and (iv) of the description of FTW,  then $\Phi( \overrightarrow{K}_{H^*}, \text{side})=W_{\text{side},H^*}$ might miss at least one box of $G_{-1}$. We will call such Eulerian circuits of $H^*$ \index{flawed (Eulerian circuit)|textbf}\textit{flawed} and we avoid using them. We note that, for any path of length greater than two in a polyomino $G$ with  $G \notin \mathcal{P}_1$, it is possible to choose a starting edge for $\overrightarrow{K}_{H^*}$ so that $\overrightarrow{K}_{H^*}$ is not flawed. From here on, all Eulerian circuits of $H^*$ we consider will be assumed to be not flawed. See Figure 1.7.

\null 

\endgroup

%%%%%%%%%   MAY BE ABLE TO DO AWAY WITH THIS TOO  %%%%%%%%%
\noindent \textbf{Observation 1.2.4.} Let $G$ be a polyomino and let $H$ be a path or cycle in $G$. 

  (a) Let the $H^*$-trail $\overrightarrow{K}'$ be a subtrail of the $H^*$-trail $\overrightarrow{K}$ and fix a side in $\{\text{right, left}\}$.

  \hspace*{0.6cm}  Then $\Phi(\overrightarrow{K}', \text{side})$ is an $H$-subwalk of $\Phi(\overrightarrow{K}, \text{side})$.

  (b) Assume that $H$ is a Hamiltonian path of $G$ and fix a side in $\{\text{right, left}\}$. 
  
\hspace*{0.6cm} Then every box of $G_{-1}$ is added to $W_{\text{side},H^*}$ by an edge of $\overrightarrow{K}_{H^*}$. 

  (c) Let $\overrightarrow{K} = (v_s, v_{s+1}), \ldots, (v_{t-1}, v_t)$ and $\overrightarrow{K}' = (v_t, v_{t-1}), \ldots, (v_{s+1}, v_s)$ be $H$ trails. 
 
 \hspace*{0.6cm} Then $\text{Boxes}(\Phi(\overrightarrow{K}, \text{right})) = \text{Boxes}(\Phi(\overrightarrow{K}', \text{left}))$ and 

\hspace*{0.6cm} $\text{Boxes}(\Phi(\overrightarrow{K}, \text{left})) = \text{Boxes}(\Phi(\overrightarrow{K}', \text{right}))$. See Figure 1.8.

  (d) We may orient a subpath $P=v_s, v_{s+1}, \ldots, v_t$ of $H$ to obtain a directed $H^*$-trail 
 
 \hspace*{0.6cm} $\overrightarrow{K}_P = (v_s, v_{s+1}), \ldots, (v_{t-1}, v_t)$ of $H^*$. Similarly, an $H^*$-trail 

\hspace*{0.6cm} $\overrightarrow{K}= (v_s, v_{s+1}), \ldots, (v_{t-1}, v_t)$ not containing an end-vertex of $H$ determines an 

\hspace*{0.6cm} oriented subpath $P=v_s, v_{s+1}, \ldots, v_t$ of $H$.

\null

\begingroup
\setlength{\intextsep}{0pt}
\setlength{\columnsep}{20pt}
\begin{wrapfigure}[]{l}{0cm}
\begin{adjustbox}{trim=0cm 0cm 0cm 0cm}
\begin{tikzpicture}[scale=1]
\usetikzlibrary{decorations.markings}

\begin{scope}[xshift=0cm] 

\begin{scope}
[very thick,decoration={
    markings,
    mark=at position 0.85 with {\arrow{>}}}
    ]
    \draw[postaction={decorate}] (1,-1)--(-0.5,2);
    \draw[postaction={decorate}] (1,-1)--(2.25,-0.5);
    \draw[postaction={decorate}] (2.25,-0.5)--(1.5,1);
    \draw[postaction={decorate}] (1,-1)--(1.5,1);
    \draw[postaction={decorate}] (0,1)--++(1.5,0.75);
\end{scope}

\draw[fill=black, opacity=1] (-0.5,2) circle [radius=0.05];
\draw[fill=black,, opacity=1] (1,-1) circle [radius=0.05];
%\draw[fill=blue, opacity=1] (2.25,-0.5) circle [radius=0.05];
\draw[fill=black, opacity=1] (1.5,1) circle [radius=0.05];

\node[left] at  (-0.5,2) [scale=0.8]{\small${B}$};
\node[below] at  (1,-1) [scale=0.8]{\small${A}$};
\node[above] at  (1.5,1) [scale=0.8]{\small${P}$};

\node[right] at  (1.9,0.25) [scale=0.8]{\small{$\overrightarrow{d^{\parallel}}$}};

\node[below] at  (1.75,-0.7) [scale=0.8]{\small{$\overrightarrow{d^{\bot}}$}};

\node[left] at  (1.25,0) [scale=0.8]{\small{$\overrightarrow{d}$}};

\node[above] at  (0.75,1.4) [scale=0.8]{\small{$\overrightarrow{n}$}};

%\node at  (1.4 1.65) [scale=1.25]{U};

\node[below] at (1.25,-1.25) [scale=1]{\begin{tabular}{c} Fig. 1.8. $P$ on the right of $\overrightarrow{AB}$. \end{tabular}};;

\end{scope}

\end{tikzpicture}
\end{adjustbox}
\end{wrapfigure}

\noindent \textbf{Note 1.2.5.} Let $A=(x_1,y_1)$, $B=(x_2,y_2)$ and $P=(x,y)$ be points in the plane that are not collinear. We define $(x_2-x_1, y_2-y_1)$ to be the direction of the vector $\overrightarrow{AB}$. Then the direction of the normal $\overrightarrow{n}$ to $\overrightarrow{AB}$, obtained by rotating $\overrightarrow{AB}$ by $-\frac{\pi}{2}$, is $(y_2-y_1, x_1-x_2)$. We want to know whether the point $P$ is on the side of $\overrightarrow{AB}$ toward which $\overrightarrow{n}$ is pointing. Let $\overrightarrow{d}=\overrightarrow{AP}$. Let  $\overrightarrow{d^{\bot}}$ be the component of $\overrightarrow{d}$ that is perpendicular to $\overrightarrow{AB}$ and let $\overrightarrow{d^{\parallel}}$ be the component of $\overrightarrow{d}$ that is parallel to $\overrightarrow{AB}$. See Figure 1.8. Note that: 

\begin{align*}
\overrightarrow{d} \cdot \overrightarrow{n} &= (\overrightarrow{d^{\parallel}} + \overrightarrow{d^{\bot}}) \cdot \overrightarrow{n}= \overrightarrow{d^{\bot}} \cdot \overrightarrow{n} \\
&= (x-x_1, y-y_1) \cdot (y_2-y_1, x_1-x_2) \\
&= (x-x_1)(y_2-y_1)+(y-y_1)(x_1-x_2).
\end{align*}

\noindent We say that $P$ is on the \textit{right} of $\overrightarrow{AB}$ if $ \overrightarrow{d} \cdot \overrightarrow{n} > 0$ and we say that $P$ is on the \textit{left}  $\overrightarrow{AB}$ if $ \overrightarrow{d} \cdot \overrightarrow{n} < 0$.

\null

\noindent Let $e=(u,v)$ be an edge of a lattice animal $G$, where $u=v(k_1,l_1)$, $v=v(k_2,l_2)$. Let $X$ be a box of the square lattice that is incident on $(u,v)$. We say that $X$ is on the \index{right (side of)}\textit{right} of the edge $(u,v)$ if there is a vertex $w=v(k,l)$ in $V(X)\setminus V(e)$ such that $(k-k_1, l_2-l_1) \cdot (l-l_1, k_1-k_2)=1$ and we say that $X$ is on the \index{left (side of)}\textit{left} of the edge $(u,v)$ if there is a vertex $w=v(k,l)$ in $V(X)\setminus V(e)$ such that $(k-k_1, l_2-l_1) \cdot (l-l_1, k_1-k_2)=-1$. 

%\textcolor{red}{[$ \mathbf{u}\cdot\mathbf{v} = \|\mathbf{u}\|\,\|\mathbf{v}\|\cos\theta$, with $\cos(0)=1$ (when both $\mathbf{n}$ and $d^{\perp}$ point ``right''), and $\cos(180)=-1$ (when they point in opposite directions).]}

\endgroup

\null

\subsection{The structure of $H$-components}

\noindent In this section we prove several structural properties of $H$-components. Key results include: $H$-components in simply connected polyominoes have unique necks (Corollary 1.3.2), and each $H$-component has a unique main trail - a subtrail of $H^*$ that determines the component's structure (Proposition 1.3.9, Corollary 1.3.11). When $H$ is a cycle, Corollary 1.3.15 proves properties of $H$-components that lie in the exterior of $H$ (which we call cookies, and are central to the reconfiguration of Hamiltonian cycles in Chapter 2).

\null

\noindent \textbf{Lemma 1.3.1.} Let $G$ be a polyomino, let $H$ be a Hamiltonian path of $G$, and let $J$ be an $H$-component of $G$. Then $J$ has a neck. Furthermore, if the boundary of $G$ has $q$ boundary components then $J$ has at most $q$ \index{neck (of an H-component)@neck (of an $H$-component)}necks.

\null 

\begingroup 
\setlength{\intextsep}{0pt}
\setlength{\columnsep}{20pt}
\begin{wrapfigure}[]{r}{0cm}
\begin{adjustbox}{trim=0cm 0cm 0cm 0cm}
\begin{tikzpicture}[scale=1.5]
\usetikzlibrary{decorations.markings}
\begin{scope}[xshift=0cm] 

\draw[gray,very thin, step=0.5cm, opacity=0.5] (0,0) grid (2.5,2.5);

\fill[blue!40!white,opacity=0.5] plot [smooth, tension=0.75] coordinates {(2,1.5)(2,2) (1.85,2.25)(1.5,2.4)(1, 2.2)(0.75,1.5)(0.975,1)(1,0.5)  (1.15,0.65)(1.1,1.1)(1,1.25)(0.9,1.5)(0.95,1.75)(1,1.9)(1.25,2.2)(1.5,2.25)(1.6,2.25)(1.75,2.15)(1.85,2)(1.9,1.65)};

\draw[black, line width=0.5mm] (2,1.5)--++(0,0.5);
\draw[black, line width=0.5mm] (0.5,0.5)--++(0.5,0);
\draw[black, line width=0.5mm] (1,0.5)--++(0,0.5);

% X-path
\begin{scope}[very thick,decoration={
    markings,
    mark=at position 0.6 with {\arrow{>}}}
    ] 
    %green
    \draw[postaction={decorate}, blue, thick] (2,1.5) to [out=135,in=-135] (2,2); 
    \draw[postaction={decorate}, blue, thick] (1,1) to [out=-45,in=45] (1,0.5); 

\end{scope}

\draw [black, line width=0.5mm] plot [smooth, tension=0.75] coordinates {(2,2) (1.85,2.25)(1.5, 2.4)(1, 2.2)(0.75,1.5)(1,1)};

\node at  (0.75,0.75) [scale=0.8]{\small{$e_t$}};
\node at  (2.25,1.75) [scale=0.8]{\small{$e_z$}};
\node at  (1.3,0.75) [scale=0.8]{\small{$N_J$}};
\node at  (1.25,0.25) [scale=0.8]{\small{$Y$}};
\node at  (1.75,1.75) [scale=0.8]{\small{$Z$}};

\draw[black, line width=0.15mm] (1.2,0.45)--(1.2,0.55);
\draw[black, line width=0.15mm] (1.25,0.45)--(1.25,0.55);
\draw[black, line width=0.15mm] (1.3,0.45)--(1.3,0.55);

\node[below] at (1.25,0) [scale=1]{\small{\begin{tabular}{c} Fig. 1.9. $\Phi \big(\overrightarrow{K}(e_z, e_t) , \textrm{left}\big)$ \\  shaded in blue. \end{tabular}}};;

\end{scope}
\end{tikzpicture}
\end{adjustbox}
\end{wrapfigure}

\noindent \textit{Proof.} If $J$ only has one box, we're done, so assume that $J$ has more than one box. Let $Z$ be a box of $J$. Without loss of generality, we may assume that $Z$ is on the left of $e_z \in \overrightarrow{K}_H$. If no edge is incident on $Z$, we choose one of the four neighbours beside it. We claim that there exists a subtrail $\overrightarrow{K}(e_z, e_{t+1})$ of $\overrightarrow{K}_H$ such that $\Phi \big(\overrightarrow{K}(e_z, e_t) , \textrm{left}\big)$ is contained in $J$ but $\Phi \big(\overrightarrow{K}(e_z, e_{t+1}) , \textrm{left}\big)$ is not. For a contradiction, assume that for every subtrail of $\overrightarrow{K}_{H^*}$ starting at $e_z$, $\Phi \big(\overrightarrow{K}(e_z, e_j) , \textrm{left}\big)$ is contained in $J$, where $j \in \{z+1, \ldots, z\}$. But then $\Phi \big(\overrightarrow{K}(e_z),\textrm{left}\big)=\Phi \big(\overrightarrow{K}_{H^*},\textrm{left}\big)$ is contained in $J$, contradicting that 
$\Phi \big(\overrightarrow{K}_{H^*},\textrm{left}\big)$ contains the boxes of $G_{-1}\setminus G$. It follows that $e_{t+1}$ adds the first box $Y$ of $\Phi \big(\overrightarrow{K}(e_z, e_{t+1}) , \textrm{left}\big)$ that is not contained in $J$. Note that, by definition of $H$-component, $Y$ must belong to $G_{-1} \setminus G$. (Since $Y$ is $H$-adjacent to the box $X$ preceding it, but $Y$ does not belong to $J$, it must be the case that $Y$ is not in $G$). See Figure 1.9. Let $X$ be the box of $J$ preceding $Y$ in $\Phi \big(\overrightarrow{K}(e_z, e_{t+1}) , \textrm{left}\big)$. We have that $X$ and $Y$ are $H$-adjacent and share a boundary edge $e_J$ of $G$ that is not in $H$. By definition of neck of an $H$-component, $X=N_J$. 

Now we show that $J$ has at most $q$ necks. For a contradiction, assume that $J$ has $q'>q$ necks. By the pigeonhole principle there is a \index{boundary component}boundary component that has at least two necks $N_{J,1}$ and $N_{J,2}$ of $J$ incident on it. Let $N_{J,1}'$ and $N_{J,2}'$ be the boxes in BBoxes($G$) incident on $N_{J,1}$ and $N_{J,2}$, respectively. Note that Outer($G$) is $H$-path-connected and so is every hole of $G$. Then there is an $H$-cycle $N_{J,1}', \ldots, N_{J,2}'$ and $N_{J,2}, \ldots, N_{J,1}, N_{J,1}'$, containing the box $N_{J,1}$ of $G$, which contradicts Corollary 1.2.3. $\square$ 

%\textcolor{red}{Check that there is a t least one neck, by following a subtrail of $\overrightarrow{K}_{H^*}$. Then use the Def. of \textit{boundary components} and ``no $H$-cycles'' to check there are at most $q$-necks.}

\endgroup

\null 

\noindent \textbf{Corollary 1.3.2.} Let $G$ be a simply connected polyomino, let $H$ be a Hamiltonian path of $G$ and let $J$ be an $H$-component of $G$. Then $J$ has a unique neck. $\square$

%Remark 1.3.3 - Not Used

{
%\textbf{Remark 1.3.3} If $G$ is a simply connected polyomino then the boundary $Q$ of an $H$-component $J$ of $G$ consists of the neck edge $ \{v_1,v_s\}$ of $J$ and the subpath $P(v_1,v_s)$ of $H$. $\square$ \textcolor{red}{[Where do I use this?]}
}

\null

\noindent \textbf{Definitions.} Let $G$ be a polyomino and let $B_0$ be the cycle bounding its outer face. Recall that the enclosure of $G$ is the set of all boxes contained in the region of the plane bounded by $B_0$. 

Let $Q$ be a cycle of vertices in $G$ consisting of a path $P(v_1,v_s) \neq \{v_1,v_s\} $ and the edge $\{v_1,v_s\}$. $Q$ bounds a region $U$ of $\text{Encl}(G)$. Let $N_U$ be the box incident on $\{v_1,v_s\}$ that is contained in $U$. We call $N_U$ \index{neck (of a region)|textbf}\textit{the neck} of $U$ and we call $\{v_1,v_s\}$ the  \index{neck edge}\textit{neck edge} of $U$. 

Let $H=v_1, \ldots, v_s, \ldots, v_t, \ldots, v_r$ be any path of $G$ such that $\{v_s,v_t\} \in E(\textrm{Encl}(G)) \setminus E(H)$. Let $P_1=P(v_1,v_s)$, $P_2=P(v_s,v_t)$ and $P_3=P(v_t,v_r)$ be a partitioning of $H$ into subpaths. We denote the $H^*$-trails $(v_s,v_{s-1}), \ldots, (v_{s-1}, v_s)$ and $(v_t,v_{t+1}), \ldots, (v_{t+1}, v_t)$ by $\overrightarrow{K}_{P_1}$ and  $\overrightarrow{K}_{P_3}$, respectively.

% Remark 1.3.4 - Not used

{
%\textbf{Remark 1.3.4} There is  side $\in \{$left, right$\}$ such that $\textrm{Encl}(G) \cap \Phi(B_0, \textrm{side}) = \emptyset$ and $\textrm{Encl}(G) \cup \Phi(B_0, \textrm{side}) = G_{-1}$. \textcolor{red}{[Where do I use this?]} 

%\null 

%\textit{Proof.} By definition of FTW and by induction on the edges of $B_0$ we have that $\Phi(B_0, \textrm{side}) \subseteq G_{-1} \setminus \textrm{Encl}(G)$, for a side $\in \{$left, right$\}$. Then the conclusion follows. 

}

\null

\noindent \textbf{Lemma 1.3.3.} Let $G$ be a polyomino. Let $H=v_1, \ldots, v_s, \ldots, v_t, \ldots, v_r$ be a Hamiltonian path of $G$
\footnote{This result holds more generally when $H$ is not Hamiltonian, but we state it this way as the general case is not needed for our purposes.} such that $\{v_s,v_t\} \in E(\textrm{Encl}(G)) \setminus E(H)$. Let $P_1=P(v_1,v_s)$, $P_2=P(v_s,v_t)$ and $P_3=P(v_t,v_r)$ be a partitioning of $H$ into subpaths. Let $U$ be the region of the plane bounded by the polygon $Q$ consisting of the subpath $P_2$ of $H$ and the edge $\{v_s, v_t\}$. Let $\overrightarrow{K}_Q$ be the directed circuit obtained from orienting $Q$. %\textcolor{red}{[Remark somewhere that this means that $Q=B_0(U)$]}. 
Then:

 (a) $U \subseteq \textrm{Encl}(G)$

 (b) $\textrm{Boxes}(\Phi(\overrightarrow{K}_Q, \textrm{right}))\subset U$ and $\textrm{Boxes}(\Phi(\overrightarrow{K}_Q, \textrm{left})) \subset G_{-1} \setminus U$ iff 

\hspace*{0.6cm} $\Phi((v_t,v_s), \textrm{right})$ is a box of $U$.

 (c) For $i \in \{1,3\}$, $P_i$ is contained in $U$ or $P_i$ is contained in $(G_{-1} \setminus U)$. 

 (d) For a side in $\{\text{right}, \text{left}\}$ and $i \in \{1,3\}$, if $P_i \subset U$ then $\Phi(\overrightarrow{K}_{P_i}, \text{side})\subset U$ and, 

\hspace*{0.6cm}  if $P_i \subset G_{-1} \setminus U$ then $\Phi(\overrightarrow{K}_{P_i}, \text{side})\subset G_{-1} \setminus U$.

\null 

\noindent \textit{Proof.} Part (a) follows from the fact that the boundary of $U$ is contained in $\text{Encl}(G)$, and so we must have $U\subseteq \textrm{Encl}(G)$ as well. Part (b) follows by definition of FTW and induction on the edges of $Q$. Part (d) follows by definition of $\overrightarrow{K}_{P_i}$ and FTW and induction on the edges of $\overrightarrow{K}_{P_i}$. 

For part (c), we assume, for definiteness, that $v_1 \in U$. Since $H$ is a path, $V(P_1) \cap V(Q)= v_s$. Thus $(P_1 \setminus v_s) \subset U$. The proofs for the other inclusions are similar. $\square $

%removed. related to proof of part (a) in Lemma 1.3.3
{
%Since $\{v_s,v_t\}$ is a boundary edge of $U$, one of the boxes incident on it, say $Z$, belongs to $G_{-1} \setminus U$.  BWOC assume that there is a box $X$ of $U$ that is not contained in $\textrm{Encl}(G)$. By Lemma 1.3.6 there is  side $\in \{$left, right$\}$ such that $X \in \Phi(B_0, \textrm{side})$. Observe that the walk $\Phi(B_0, \textrm{side})$ determines a cycle of boxes in $G_{-1}$. Note that for any $Y \in \Phi(B_0, \textrm{side})$, we can draw a path contained in $\Phi(B_0, \textrm{side})$ that does not intersect $B_0$ and that joins the center $c_y$ of $Y$ to the center $c_x$ of $X$. By JCT, it follows that $\Phi(B_0, \textrm{side}) \subset U$. Note that this means that $\textrm{Encl}(G) \subset U$ as well. But now, by Remark (Z)gggg.2.c) we have that $G_{-1}$ is contained in $U$, contradicting that $Z$ belongs to $G_{-1} \setminus U$.  $\square$.
}

\null 

\noindent The next three Lemmas will only be used by Lemma 4.12 in Chapter 4.

\null

\noindent \textbf{Lemma 1.3.4.} If a graph has a closed walk with a non-repeated edge, then the graph has a cycle\cite{profnotes}. $\square$ 

%\textcolor{red}{[Proof. Let $W$ be a shortest closed walk with a non-repeated edge $e$. If $W$ is a cycle, we are done. Otherwise, there is a repeated vertex and $W$ is a union of two closed walks $W_1$ and $W_2$ that are shorter than $W$. One of them, say $W_1$, contains $e$, a non-repeatededge. This contradicts the minimality of $W$.]}

\null 

\noindent The following lemma is a proof of the fact that around each hole of a polyomino, there is a cycle of boxes. 

\null 

% partial dual graph
{
%Let $G$ be a polyomino and let $H$ be a lattice animal in $G$. Let $W=W(X_1,X_s)$ be an $H$-walk of boxes in $G_{-1}$. Mod out $W$ by consecutively repeated boxes to obtain the $H$-walk $W'$. That is, if the string of boxes $XY^kZ$ appears in $W$, substitute it by the string $XYZ$. Note that $W'$ is an $H$-walk that retains the order in which boxes of $G_{-1}$ appear in $W$.

%Define the \textit{partial dual graph} of $W$ to be the graph $W^{D}=(V^D, E^D)$ to be the graph with vertex set $V^{D}$ corresponding to the set of distinct boxes of $W$ and edge set $E^D$ corresponding to the set of distinct edges in $E_{\text{glue}}(W)$. We note that there is a one-to-one correspondence between cycles in $W^{D}$ and cycles in $W'$, and that every cycle in $W'$ has a corresponding cycle in $W$.
}

\noindent \textbf{Lemma 1.3.5.} Let $G$ be a polyomino and let $\mathcal{O}$ be a hole of $G$ with boundary $B(\mathcal{O})$. Let $\overrightarrow{K}$ be the circuit obtained by orienting the edges of $B(\mathcal{O})$ (Recall Lemma 1.1.7 (f)).
%Fix a side in $\{ \text{left}, \text{right} \}$. Assume that $\Phi(e_1, \text{side})$ is contained in $G_{-1} \setminus \mathcal{O}$.
Then it is possible to choose a starting edge $e_1$ for the $B(\mathcal{O})^*$-circuit $\overrightarrow{K}$ 
and a side in $\{ \text{right},\text{left}  \}$ so that the $B(\mathcal{O})$-walk $\Phi(\overrightarrow{K}, \text{side})$ is contained in $G_{-1} \setminus \mathcal{O}$ and it contains a cycle of boxes, or the $B(\mathcal{O})$-walk $\Phi(\overrightarrow{K}, \text{side}), \Phi(e_1, \text{side})$ is contained in $G_{-1} \setminus \mathcal{O}$ and it contains a cycle of boxes.

\null 

\noindent \textit{Proof.} By Lemma 1.3.3 (b), for any edge $e$ of $\overrightarrow{K}$, there is a side in $ \{ \text{right},\text{left}  \}$ such that $\Phi(e, \text{side}) \subset G_{-1} \setminus \mathcal{O}$ implies $\Phi(\overrightarrow{K}, \text{side}) \subset G_{-1} \setminus \mathcal{O}$. For definiteness, let $\Phi(\overrightarrow{K}, \text{left}) \subset G_{-1} \setminus \mathcal{O}$. Let $\overrightarrow{K}=e_1, \ldots, e_t$ be such that $e_1,e_t$ are as $e_{j+1}, e_j$, respectively, in Case (i) or Case (ii) of the definition of FTW. Let $W=W_{\text{left}}(X_1,X_s)=\Phi(\overrightarrow{K}, \text{left})$. If $e_1,e_t$ are as $e_{j+1}, e_j$, respectively, in Case (i) of the definition of FTW then $X_1=X_s$ (Figure 1.10 (a)); and if $e_1,e_t$ are as $e_{j+1}, e_j$, respectively, in Case (ii) of the definition of FTW then $X_1$ is adjacent to $X_s$ (Figure 1.10 (b)). If the latter, we may add $X_1$ to $W$ after $X_s$ to obtain a closed $H$-walk. From here on we will assume that $e_1,e_t$ are as $e_{j+1}$ and $e_j$, respectively, in Case (i) of FTW, so that $W$ is the closed walk starting and ending at $X_1$. The case where $e_1,e_t$ are as $e_{j+1}$ and $e_j$, respectively, in Case (ii) of FTW is very similar, so we omit it. It remains to check that $W$ contains a cycle.

\begin{center}
\begin{adjustbox}{trim=0cm 0cm 0cm 0cm}

\end{adjustbox}
\end{center}

\noindent Mod out $W$ by consecutively repeated boxes to obtain the $H$-walk $W'$. That is, if the string of boxes $XY^kZ$ appears in $W$, substitute it by the string $XYZ$. Note that $W'$ is an $H$-walk that retains the order in which boxes of $G_{-1}$ appear in $W$. Now we can view $W'$ as an alternating sequence of boxes and gluing edges starting and ending at the box $X_1$. If we can find a non-repeated gluing edge of $W'$, then Lemma 1.3.4 implies that $W$ contains a cycle of boxes of $G_{-1}$. 

Let $e_z$ be a horizontal northernmost edge in $\overrightarrow{K}$. Let $Z= \Phi (e_z, \text{left})$. Note that $Z\in W$, $Z'=Z+(1,0)$ is an $H$-neighbour of $Z$, and that no other edge of $\overrightarrow{K}$ is incident on $Z$. We claim that the vertical gluing edge $f=\{u,v\}$ between $Z$ and $Z'$, where $u\in \overrightarrow{K}$ and $v \notin \overrightarrow{K}$, is a non-repeated gluing edge of $W'$. See Figure 1.10 (c). For a contradiction, assume that there is a second occurrence of $f$ in $W'$. This means that the sequence of boxes $Z,Z'$ or the sequence of boxes $Z'Z$ occurs again in $W'$. By definition of FTW, and the fact that the edges of $\overrightarrow{K}$ are unique, such a sequence must be added to $W$ by edges of $\overrightarrow{K}$ incident on $v$, of which there are none. 

\null 

\noindent \textbf{Lemma 1.3.6.} Let $H$ be a Hamiltonian path of a polyomino $G$. Assume that $J$ is a \index{non-self-adjacent}non-self-adjacent $H$-subtree of an $H$-component of $G$, and that $J$ has no \index{polyking junction}polyking junctions. Then $J$ is a \index{simply connected polyomino}simply connected polyomino.

\null 

\noindent \textit{Proof.} Since $J$ is an $H$-subtree, $J$ is box-path-connected, so it is a polyomino. For a contradiction, assume that $J$ is not simply connected. Then $J$ has a hole $\mathcal{O}$ with boundary $B(\mathcal{O})$.  
Orient the edges of $B(\mathcal{O})$ to obtain a directed circuit $\overrightarrow{K}$. By Lemma 1.3.5, we may choose a side and a starting edge for $\overrightarrow{K}$ such that $\Phi(\overrightarrow{K}, \text{side}) \subseteq J_{-1} \setminus \mathcal{O}$ and such that $\Phi(\overrightarrow{K}, \text{side})$ contains a cycle $C$ of boxes in $J_{-1} \setminus \mathcal{O}$. For definiteness, we assume that $\Phi(\overrightarrow{K}, \text{left}) \subseteq J_{-1} \setminus \mathcal{O}$.

Suppose that $\Phi(\overrightarrow{K}, \text{left})$ is contained in $J \setminus \mathcal{O}$. Then, either at least one of the gluing edges of $C$ belongs to $H$, or no gluing edge of $C$ is in $H$. The former contradicts the assumption that $J$ is non-self-adjacent. The latter contradicts Proposition 1.2.1. It remains to show that $\text{Boxes}(\Phi(\overrightarrow{K},\text{left})) \subseteq J$.

\null

\noindent We use induction on the edges of $\overrightarrow{K}$ to show that $\text{Boxes} (\Phi(\overrightarrow{K}, \text{left})) \subseteq J$. Let $Z$ and $Z'$ be the boxes in $G_{-1}$ incident on an edge $e \in \overrightarrow{K}$. By Lemma 1.1.7 (b), one of $Z$ and $Z'$ belongs to $\mathcal{O}$ and the other belongs to $J$. Let this be observation $(*)$. Note that by $(*)$ and the assumption that $\Phi(\overrightarrow{K}, \text{left}) \subseteq J_{-1} \setminus \mathcal{O}$, the base case holds: $\Phi(e_1, \text{left})$ is a box of $J$. 

For the inductive case, assume that for $i \in \{2,\ldots,i_0\}$, $\text{Boxes}(\Phi(\overrightarrow{K}(e_1,e_{i_0}), \text{left})) \subset J$. For definiteness, assume that $e_{i_0}=(a,a+1;b)=(v(a,b), v(a+1;b))$ and let $X=R(a,b)$. Then $X\in J$ and $X+(0,-1) \in \mathcal{O}$. We want to check that $\text{Boxes}(\Phi(\overrightarrow{K}(e_1,e_{i_0+1}), \text{left})) \subset J$.

\begingroup 
\setlength{\intextsep}{0pt}
\setlength{\columnsep}{20pt}

\begin{adjustbox}{trim=0cm 0cm 0cm 0cm}

\end{adjustbox}

\noindent There are three possibilities: $e_{i_0+1}=e(a+1; b,b+1)$, $e_{i_0+1}=e(a+1, a+2; b)$, and $e_{i_0+1}=e(a+1; b-1,b)$. See Figure 1.11. If $e_{i_0+1}=e(a+1;b,b+1)$, then $\text{Boxes}(\Phi(\overrightarrow{K}(e_1,e_{i_0+1}), \text{left}))$$=\text{Boxes}(\Phi(\overrightarrow{K}(e_1,e_{i_0}), \text{left}))$ and we are done by the inductive hypothesis. If $e_{i_0+1}=e(a+1, a+2; b)$, then Observation 1.1.9 implies that $e(a+1; b,b+1) \notin B(J)$ and $e(a+1; b-1,b) \notin B(J)$. Then observation $(*)$ implies that one of $X+(1,0)$ and $X+(1,-1)$ belongs to $J$, and the other belongs to $\mathcal{O}$. The assumption that $J$ has no polyking junction implies that we must have $X+(1,0) \in J$ and $X+(1,-1) \in \mathcal{O}$.
{%It follows that $X$ and $X+(1,0)$ are box-path-connected and thus must both belong to $J$.
} If $e_{i_0+1}=e(a+1; b-1,b)$, then by $(*)$, $Z+(1,-1)$ belongs to $J$. Our assumption that $J$ has no polyking junctions implies that $Z+(1,0)$ must belong to $J$ as well. Thus, $(\Phi(\overrightarrow{K}, \text{left}))$ is contained in $J$. $\square$

\null 

\noindent \textbf{Proposition 1.3.7.} Let $G$, $H$, $P_i$ for $i \in \{1,2,3\}$, $\{v_s,v_t\}$, $Q$ and $U$ be as in Lemma 1.3.3. In addition, assume that $H$ is a Hamiltonian path of $G$, that $\{v_s,v_t\}$ is in $B(G)$, that $J$ is an $H$-component of $G$ that has $\{v_s,v_t\}$ as one of its neck-edges, and that the neck $N_U$ of $U$ coincides with a neck $N_J$ of $J$\footnote{This can be made more general by only assuming that $U \cap J \neq \emptyset$, but it makes the proof longer and we don't need it.}. Then $J \subseteq U$. Moreover, if $G$ is a simply connected polyomino, then $J=U$.

\null 

\noindent \textit{Proof.} For a contradiction, assume that there is a box $Z$ in $J \setminus U$. Then there is an $H$-path $P(X_1,X_r)$ in $J$, where $X_1=N_J=N_U$ and $X_r=Z$. Let $Y$ be the other box incident on $\{v_s,v_t\}$ and let $c_1, c_2, \ldots, c_r=c_1$ be the centers of the boxes of $P(X_1,X_r)$. Then for each $j \in \{1,2,\ldots,r-1\}$, $[c_j,c_{j+1}]$ intersects the gluing edge of $X_j$ and $X_{j+1}$ and $[c_j,c_{j+1}]$ intersects no other edge of $G$. By JCT, since $c_1 \in U$ and $c_r \notin U$, $P(c_1,c_r)$ intersects $Q$ at some edge $e$. Since $Y \notin G$, $e \neq \{v_s,v_t\}$. Then $e$ must be some other edge of $Q$. But all other edges of $Q$ belong to $H$, contradicting that $e$ is a gluing edge.

\null 

\noindent Assume that $G$ is a simply connected polyomino.  We will show that $U \subseteq J$. Since $G$ is a simply connected polyomino, we have that $N_U$ is the unique neck of $J$. For a contradiction, assume that there is a box $Z$ in $U \setminus J$. For definiteness, by Observation 1.2.4 (b), we may assume that $Z$ is on the left of an edge $e_z$ of $\overrightarrow{K}_{H^*}$. (If no edge of $H$ is incident on $Z$, as in Case (iii) of FTW, we replace $Z$ with the second box added by $e_z$). 

\begingroup 
\setlength{\intextsep}{10pt}
\setlength{\columnsep}{20pt}
\begin{wrapfigure}[]{r}{0cm}
\begin{adjustbox}{trim=0cm 0cm 0cm 0cm}
\begin{tikzpicture}[scale=2]
\usetikzlibrary{decorations.markings}
\begin{scope}[xshift=0cm] 

\draw[gray,very thin, step=0.5cm, opacity=0.5] (0,0) grid (2,1.5);

\fill[blue!40!white,opacity=0.5] plot [smooth, tension=0.2] coordinates {(0,0.5)(0.5,0.5)(0.65,0.6)(0.85,0.9)(1,1)(1,1.5)(0.5,1)(0,1)};

\fill[blue!40!white,opacity=0.5] (1,1) rectangle (2,1.5);

\fill[green!40!white,opacity=0.5] (1.5,0.5) rectangle (2,1);

 \draw[orange, line width=0.5mm] (1.5,1)--++(0.5,0);

\begin{scope}
[very thick,decoration={
    markings,
    mark=at position 0.6 with {\arrow{>}}}
    ]
    
    \draw[postaction={decorate}, blue, line width=0.5mm] (0,0.5)--++(0.5,0);
    
    \draw[postaction={decorate}, blue, line width=0.5mm] (1,1)--++(0.5,0);
    
    \draw[postaction={decorate}, blue, line width=0.5mm] (1.5,1)--++(0,-0.5);
    
\end{scope}

\draw [blue, line width=0.5mm] plot [smooth, tension=0.75] coordinates {(0.5,0.5) (0.65,0.6)(0.85,0.9)(1, 1)};

\node[below] at  (0.2,0.5) [scale=0.8]{\small{$e_z$}};

\node[below] at  (1.25,1) [scale=0.8]{\small{$e_{y-1}$}};

\node[left] at  (1.55,0.65) [scale=0.8]{\small{$e_y$}};

\node[right] at  (2,1) [scale=0.8]{\small{$f$}};

\node at  (1.75,1.25) [scale=1]{\small{$Y'$}};
\node at  (1.75,0.75) [scale=1]{\small{$Y$}};
\node at  (0.25,0.75) [scale=1]{\small{$Z$}};

\node[below] at (1,0) [scale=1]{\begin{tabular}{c} Fig. 1.12. $\Phi \big(\overrightarrow{K}(e_z, e_t) , \textrm{left}\big)$ \\  shaded in blue; $f$ in orange. \end{tabular}};;

\end{scope}
\end{tikzpicture}
\end{adjustbox}
\end{wrapfigure}

\noindent Let $\overrightarrow{K}(e_z,e_{y-1})$ be a subtrail of $\overrightarrow{K}_{H^*}$ such that $\Phi(\overrightarrow{K}(e_z,e_{y-1}),\textrm{left})$ is contained in $U\setminus J$ but $\Phi(\overrightarrow{K}(e_z,e_y),\textrm{left})$ is not. Let $Y$ be the first box not in $U \setminus J$ added by $e_y$ and let $Y'$ be the box preceding $Y$ in $\Phi(\overrightarrow{K}(e_z,e_y),\textrm{left})$. See Figure 1.12. This means that $Y' \in U \setminus J$. Let $f$ be the edge that $Y$ and $Y'$ share.  Since $J \subset U$, either $Y \in J$ or $Y \in G_{-1} \setminus U$.

Suppose that $Y \in J$. Since $G$ is simply connected, by Corollary 1.3.2, $J$ has a unique neck. Since $Y'$ is not in $J$, $f$ must be the neck-edge of $J$. By assumption, $Y=N_J=N_U$. But then, by definition of the neck of $U$, $Y' \notin U$, contradicting our assumption that $Y' \in U \setminus J$.

Suppose that $Y \in G_{-1} \setminus U$. Since $Y' \in U$, $f$ must be a boundary edge of $U$, so $f \in Q$. Since $\Phi(\overrightarrow{K}(e_z,e_y),\textrm{left})$ is an $H$-walk, $f$ is a gluing edge of $\Phi(\overrightarrow{K}(e_z,e_y),\textrm{left})$ and so $f \notin H$. It follows that $f=\{v_s,v_t\}$. Since $Y' \in U$ and $Y \in  G_{-1} \setminus U$, we have  $Y'=N_U=N_J \in J$. But this contradicts our assumption that $Y' \in U \setminus J$. $\square$ 

\endgroup

\null 

\noindent \textbf{Corollary 1.3.8.} Let $G$ be a simply connected polyomino and let $H$, $P_i$ for $i \in \{1,2,3\}$, $\{v_s,v_t\}$, $Q$, $U$ and $J$ be as in Proposition 1.3.7, and let $\overrightarrow{K}_{P_i}$ be as in Lemma 1.3.3. Then:

 (a) Exactly one of the following is true:

\hspace*{0.6cm} (I) $\textrm{Boxes}(\Phi(P_2, \textrm{right}))\subseteq J$ and $\textrm{Boxes}(\Phi(P_2, \textrm{left})) \subseteq G_{-1} \setminus J$, and 

\hspace*{0.6cm} (II) $\textrm{Boxes}(\Phi(P_2, \textrm{left}))\subseteq J$ and $\textrm{Boxes}(\Phi(P_2, \textrm{right})) \subseteq G_{-1} \setminus J$.

(b) For a side in $\{\text{right}, \text{left}\}$ and $i \in \{1,3\}$, $\text{Boxes}(\Phi(\overrightarrow{K}_{P_i}, \text{side}))\subseteq J$ or 

\hspace*{0.6cm} $\text{Boxes}(\Phi(\overrightarrow{K}_{P_i}, \text{side}))\subseteq G_{-1} \setminus J$. Furthermore, $\text{Boxes}(\Phi(\overrightarrow{K}_{P_i}, \text{side}))\subseteq J$ 

\hspace*{0.6cm}  if and only if $P_i$ is contained in $J$ and $\text{Boxes}(\Phi(\overrightarrow{K}_{P_i}, \text{side}))\subseteq G_{-1} \setminus J$ 

\hspace*{0.6cm}  if and only if $P_i$ is contained in $G_{-1} \setminus J$. $\square$
%c) $J$ is a simply connected polyomino with boundary $Q=B(J)$. \textcolor{red}{[follows from Corollary 1.1.10]}

\null 

%\noindent \textit{Proof of c).} By JCT, $U$ is a simply connected polyomino. By the `moreover' part of Proposition 1.3.7, $U=J$, so $J$ is also a simply connected polyomino and $B(J)=B(U)=Q$. 

%\noindent Lemma 1.3.3 allows us to later prove Lemmas A.2, 4.9, and 4.17 in the more general case where the underlying graph $G$ is a polyomino. We note here that while we can soften this assumption to $G$ being a simply connected polyomino in Lemmas A.2 and 4.17, we cannot do the same for for Lemma 4.9.

\noindent \textbf{Proposition 1.3.9.} Let $G$ be a simply connected polyomino and let $H$, $P_i$ for $i \in \{1,2,3\}$, $\{v_s,v_t\}$, $Q$, and $J$ be as in Corollary 1.3.8. Then there are exactly two maximal and unique elements $(\overrightarrow{K}_J, \text{right})$ and $(\overrightarrow{K}'_J, \text{left})$ of $\mathcal{K} \times \{\text{right, left}\}$ such that $\text{Boxes}(\Phi(\overrightarrow{K}_J, \text{right}))=J$ and $\text{Boxes}(\Phi(\overrightarrow{K}'_J, \text{left}))=J$.

\null

\noindent \textit{Proof.} We may assume that $J$ has more than one box. Let $\overrightarrow{K}_{P_1}$ and $\overrightarrow{K}_{P_3}$ be as in Lemma 1.3.3. By Observation 1.2.4 (d), $P_2$ determines the trail $\overrightarrow{K}_{P_2}=(v_s,v_{s+1}), \ldots,(v_{t-1},v_t)$. Let $\overrightarrow{K}'_{P_2}$ be the trail $=(v_t,v_{t-1}), \ldots,(v_{s+1},v_s)$. Note that  $\overrightarrow{K}_{H^*}=$ $\overrightarrow{K}_{P_1}$, $\overrightarrow{K}_{P_2}$, $\overrightarrow{K}_{P_3}$, $\overrightarrow{K}'_{P_2}$ is an Eulerian circuit of $H^*$. We will construct $\overrightarrow{K}_{J}$ and $\overrightarrow{K}'_{J}$ explicitly and then check maximality and uniqueness.

Note that $J$ may contain both end-vertices of $H$, exactly one end-vertex or neither end-vertex, so there are three cases to consider.

\null 

\noindent \textit{CASE 1. $J$ contains no end-vertices.} By Corollary 1.3.8 (a) $\textrm{Boxes} (\Phi (P_2, \textrm{right}))\subseteq J$ or $\textrm{Boxes} (\Phi (P_2, \textrm{left}))$ $\subseteq J$. For definiteness, assume that $\textrm{Boxes} (\Phi (P_2, \textrm{right}))\subseteq J$. Then $\textrm{Boxes} (\Phi (P_2, \textrm{left}))\subseteq G_{-1}\setminus J$.

First we will show that $\overrightarrow{K}_J=\overrightarrow{K}_{P_2}$ satisfies $\text{Boxes}(\Phi(\overrightarrow{K}_J, \textrm{right}))=J$. We already have that $\textrm{Boxes} (\Phi (\overrightarrow{K}_{P_2}, \textrm{right}))\subseteq J$. We need to check that $J \subset \textrm{Boxes} (\Phi (\overrightarrow{K}_{P_2}, \textrm{right}))$. Let $Z$ be a box of $J$. Note that at least one vertex of $Z$, say $v_z$, is neither an end-vertex of $H$ nor incident on the neck of $J$. Either at least one of the two edges of $Z$ incident on $v_z$ is in $H$, or neither is.

\null 

\noindent \textit{CASE 1.1: At least one of the two edges of $Z$ incident on $v_z$, say $e_z$, is in $H$.} If we can show that $e_z \in \overrightarrow{K}_{P_2}$ then, by Observation 1.2.4 (a), we have that $Z=\Phi (e_z, \text{right})$ is contained in $\Phi(\overrightarrow{K}_{P_2}, \text{right})$ and we're done. Now, Corollary 1.3.8 (b) implies that $e_z \notin \overrightarrow{K}_{P_1}$ and $e_z \notin \overrightarrow{K}_{P_3}$. BWOC assume that $e_z\in \overrightarrow{K}'_{P_2}$. Then, by Observation 1.2.4 (a), $Z=\Phi(e_z, \text{right})$ is contained in $\textrm{Boxes} (\Phi (\overrightarrow{K}'_{P_2}, \textrm{right}))$. By Observation 1.2.4 (c), $\textrm{Boxes} (\Phi (\overrightarrow{K}'_{P_2}, \textrm{right}))$= $\textrm{Boxes} (\Phi (\overrightarrow{K}_{P_2}, \textrm{left}))$ and $\textrm{Boxes} (\Phi (\overrightarrow{K}_{P_2}, \textrm{left}))$ is contained in $G_{-1} \setminus J$. But this contradicts that $Z\in J$. Thus we must have that $e_z\in \overrightarrow{K}_{P_2}$. End of Case 1.1.

\begingroup 
\setlength{\intextsep}{0pt}
\setlength{\columnsep}{20pt}
\begin{wrapfigure}[]{l}{0cm}
\begin{adjustbox}{trim=0cm 0.25cm 0cm 0.5cm}
\begin{tikzpicture}[scale=2]
\usetikzlibrary{decorations.markings}
\begin{scope}[xshift=0cm] 

\draw[gray,very thin, step=0.5cm, opacity=0.5] (0,0) grid (1,1);

\fill[blue!40!white,opacity=0.5] (0,0) rectangle (1,0.5);
\fill[blue!40!white,opacity=0.5] (0.5,0.5) rectangle (1,1);

\draw[fill=blue] (0.5,0.5) circle [radius=0.05];

\begin{scope}
[very thick,decoration={
    markings,
    mark=at position 0.6 with {\arrow{>}}}
    ]
    
    \draw[postaction={decorate}, blue, line width=0.5mm] (0,0.5)--++(0.5,0);
    
    \draw[postaction={decorate}, blue, line width=0.5mm] (0.5,0.5)--++(0,0.5);

\end{scope}

%black lines
{

\draw[black, line width=0.15mm] (0.45,0.2)--++(0.1,0);
\draw[black, line width=0.15mm] (0.45,0.25)--++(0.1,0);
\draw[black, line width=0.15mm] (0.45,0.3)--++(0.1,0);

\draw[black, line width=0.15mm] (0.7,0.45)--++(0,0.1);
\draw[black, line width=0.15mm] (0.75,0.45)--++(0,0.1);
\draw[black, line width=0.15mm] (0.8,0.45)--++(0,0.1);
}

\node[above] at  (0.1,0.475) [scale=0.8]{\small{$e_z$}};
\node[right] at  (0.475,0.9) [scale=0.8]{\small{$e_{z+1}$}};
\node at  (0.4,0.6) [scale=0.8]{\small{$v_z$}};

\node at  (0.25,0.25) [scale=0.8]{\small{$Z'$}};
\node at  (0.75,0.25) [scale=0.8]{\small{$Z$}};
\node at  (0.75,0.75) [scale=0.8]{\small{$Z''$}};

\node[right] at (1,0) [scale=1]
{\tiny{$b$}};
\node[right] at (1,0.5) [scale=1]
{\tiny{$+1$}};
\node[right] at (1,1) [scale=1]
{\tiny{$+2$}};

\node[above] at (0, 1) [scale=1]
{\tiny{-1}};
\node[above] at (0.5,1) [scale=1]
{\tiny{$a$}};

\node[below] at (0.6,0) [scale=1]{\begin{tabular}{c} Fig. 1.13. Case 1.2. \end{tabular}};;

\end{scope}
\end{tikzpicture}
\end{adjustbox}
\end{wrapfigure}

\null 
 
\noindent \textit{CASE 1.2: Neither of the two edges of $Z$ incident on $v_z$ is in $H$.} For definiteness, let $Z=R(a,b)$ and $v_z=v(a,b+1)$. Let $Z'=Z+(-1,0)$ and $Z''=Z+(0,1)$. Then we have that $e(a-1,a;b+1) \in H$, $e(a;b+1,b+2) \in H$, $Z' \in J$ and $Z'' \in J$. Let $(v(a-1,b+1),v(a,b+1))=e_z$ and $(v(a,b+1),v(a,b+2))=e_{z+1}$. By Case 1.1, $e_z \in \overrightarrow{K}_{P_2}$ and $e_{z+1} \in \overrightarrow{K}_{P_2}$. Then $Z \in \text{Boxes}(\Phi (\overrightarrow{K}(e_z, e_{z+1}), \text{right})) \subset \text{Boxes}(\Phi (\overrightarrow{K}_{P_2}), \text{right})$. End of Case 1.2. See Figure  1.13.

%pagemarker
\noindent Cases 1.1 and 1.2 showed that $\text{Boxes}(\Phi(\overrightarrow{K}_{P_2}, \textrm{right}))=J$ and that the edges of $\overrightarrow{K}_J$ cannot belong to $\overrightarrow{K}_{P_1}$, $\overrightarrow{K}'_{P_2}$ or $\overrightarrow{K}_{P_3}$. It follows that $\overrightarrow{K}_J$ cannot be extended and thus is maximal. 

To see that $(K_J, \text{right})$ is unique, assume toward a contradiction, that there exists an element $(\overrightarrow{K}_J^{\dag}, \text{right})$ in  $\mathcal{K} \times \{\text{right}\}$, distinct from $\overrightarrow{K}_J$, such that $\text{Boxes}(\Phi(\overrightarrow{K}_J^{\dag}, \text{right}))=J$. Then there is an edge $e \in \overrightarrow{K}_J^{\dag} \setminus \overrightarrow{K}_J$. Since $\overrightarrow{K}_J = \overrightarrow{K}_{P_2}$, $e \in \overrightarrow{K}_{P_1}$ or $e \in \overrightarrow{K}_{P_3}$ or  $e \in \overrightarrow{K}'_{P_2}$. By Corollary 1.3.8(b), $e \notin \overrightarrow{K}_{P_1}$ and $e \notin \overrightarrow{K}_{P_3}$. By Cases 1.1 and 1.2, $e \notin \overrightarrow{K}'_{P_2}$. But then $e$ must belong to  $\overrightarrow{K}_{P_2}=\overrightarrow{K}_J$, contradicting that $e \in \overrightarrow{K}_J^{\dag} \setminus \overrightarrow{K}_J$. Thus, $(\overrightarrow{K}_J, \text{right})$ is unique.

By Observation 1.2.4 (c) we can see that  $\textrm{Boxes} (\Phi(\overrightarrow{K}'_{P_2},\textrm{left}))=J$. The proof that $\overrightarrow{K}'_{P_2}=\overrightarrow{K}'_J$ is maximal and unique is the same as the proof for $\overrightarrow{K}_J$, so we omit it. End of Case 1. 

\null

\noindent \textit{CASE 2. $J$ contains exactly one end-vertex.} For definiteness, assume that $J$ contains $v_1$. By Corollary 1.3.8 (a) $\textrm{Boxes} (\Phi (P_2, \textrm{right}))$ $\subseteq J$ or $\textrm{Boxes} (\Phi (P_2, \textrm{left}))\subseteq J$. For definiteness, assume that $\textrm{Boxes} (\Phi (P_2, \textrm{right}))$ $\subseteq J$. We will show that $\overrightarrow{K}_J =\overrightarrow{K}_{P_1}, \overrightarrow{K}_{P_2}.$ By Corollary 1.3.8 (b), $\textrm{Boxes} (\Phi (\overrightarrow{K}_{P_1}, \textrm{right}))\subseteq J$. Then we have that $\textrm{Boxes} (\Phi (\overrightarrow{K}_{P_1},\overrightarrow{K}_{P_2}, \textrm{right}))\subseteq J$. We need to check that $J \subseteq \textrm{Boxes} (\Phi (\overrightarrow{K}_{P_1},\overrightarrow{K}_{P_2}, \textrm{right}))$.

Let $Z$ be a box of $J$ and let $v_z$ be as in Case 1. By Case 1, it is sufficient to check the case where at least one of the two edges of $Z$ incident on $v_z$, say $e_z$, is in $H$. By Corollary 1.3.8(b), $e_z \notin \overrightarrow{K}_{P_3}$ and by Case 1, $e_z \notin \overrightarrow{K}'_{P_2}$. Then $e_z \in \overrightarrow{K}_{P_1},\overrightarrow{K}_{P_2}$. 

Let $\overrightarrow{K}'_J=\overrightarrow{K}'_{P_2}, \overrightarrow{K}_{P_1}$. By Observation 1.2.4 (c) we can see that $\textrm{Boxes} (\Phi(\overrightarrow{K}'_J,\textrm{left}))=J$ as well. Proofs for maximality and uniqueness are similar to those in Case 1, so we omit them. End of Case 2.

\null

\noindent \textit{CASE 3. $J$ contains both end-vertices.} As in previous cases, we may assume, for definiteness, that $\textrm{Boxes} (\Phi (P_2, \textrm{right}))\subseteq J$. Let $\overrightarrow{K}_J=\overrightarrow{K}_{P_1}, \overrightarrow{K}_{P_2}, \overrightarrow{K}_{P_3}$ and $\overrightarrow{K}'_J=\overrightarrow{K}_{P_3}, \overrightarrow{K}'_{P_2}, \overrightarrow{K}_{P_1}$. Using the same arguments as above we can see that the elements of $\mathcal{K} \times \{\text{right, left}\}$ that satisfy the conclusion are $(\overrightarrow{K}_J, \text{right})$ and $(\overrightarrow{K}'_J, \text{left})$.  End of Case 3. $\square$

\null

\noindent \textbf{Corollary 1.3.10.} Let $G$ be a simply connected polyomino, let $H$ be a Hamiltonian path of $G$, and let $J$ be an $H$-component of $G$. If some edge $e_z \in \overrightarrow{K}_{H^*}$ adds a box $Z$ in $J$ to $W_{H^*,\textrm{right}}$, then $e_z \in \overrightarrow{K}_J$; and if $e_z \in \overrightarrow{K}_{H^*}$ adds a box $Z$ in $J$ to $W_{H^*,\textrm{left}}$, then $e_z \in \overrightarrow{K}'_J$. $\square$

\null

\noindent \textbf{Definitions.} Let $G$ be a simply connected polyomino, and let $J$ be an $H$-component of $G$. We call the elements $(\overrightarrow{K}_J, \text{right})$ and $(\overrightarrow{K}'_J, \text{left})$ of $\mathcal{K} \times \{\text{right, left}\}$  \textit{the right main trail of $J$} and \textit{the left main trail of $J$}, respectively. We call $\Phi((\overrightarrow{K}_J, \text{right}))$ and  $\Phi((\overrightarrow{K}'_J, \text{left}))$ and \textit{the right main walk of $J$} and  the \textit{left main walk of $J$}, respectively.

\null

\noindent \textbf{Corollary 1.3.11} Let $G$ be a simply connected polyomino and let $H$, $P_i$ for $i \in \{1,2,3\}$, $\{v_s,v_t\}$, $Q$, and $J$ be as in Corollary 1.3.8. For $i \in \{1,2,3\}$, let $\overrightarrow{K}_{P_i}$ and $\overrightarrow{K}_{P_i}'$ be as in Proposition 1.3.9. %Let $\{u,v\}$ be the neck-edge of $J$. 
Then:

(a) If $J$ contains no end-vertices,  then $\overrightarrow{K}_J=\overrightarrow{K}_{P_2}$. 

(b) If $J$ contains exactly one end-vertex, then $\overrightarrow{K}_J=\overrightarrow{K}_{P_1},\overrightarrow{K}_{P_2}$, or $\overrightarrow{K}_J=\overrightarrow{K}_{P_2},\overrightarrow{K}_{P_3}$. 

(c) If $J$ contains both end-vertices, then $\overrightarrow{K}_J=\overrightarrow{K}_{P_1},\overrightarrow{K}_{P_2},\overrightarrow{K}_{P_3} $. 

(d) $J$ is a simply connected polyomino with boundary $B(J)=P_2, \{v_s,v_t\}$. $\square $

\null 

\begingroup 
\setlength{\intextsep}{0pt}
\setlength{\columnsep}{20pt}
\begin{wrapfigure}[]{l}{0cm}
\begin{adjustbox}{trim=0cm 0.5cm 0cm 0cm}

\end{adjustbox}
\end{wrapfigure}

\noindent \textbf{Lemma 1.3.12.} Let $H$ be a Hamiltonian path of a simply connected polyomino $G$ and let $J$ be an $H$-component of $G$. Then $J$ is \index{self-adjacent}\index{non-self-adjacent}self-adjacent if and only if an edge of a main trail of $J$ has an end-vertex of $H$ incident on it in $G \setminus B(G)$

\null 

\noindent \textit{Proof.} Note that if the condition is satisfied, then the four boxes incident on the end-vertex belong to $J$ and are $H$-path connected, and so $J$ is self-adjacent. It remains to prove the converse. 

\noindent Suppose that $J$ is a self-adjacent $H$-component of $G$. We will show that some edge of the main left trail $\overrightarrow{K}'_J$ of $J$ has an end-vertex of $H$ incident on it in $G \setminus B(G)$.  Since $J$ is self-adjacent, there are $H$-path-connected boxes $Z$ and $Z'$ in $J$ sharing an edge $e=\{v_j,v_{j+1}\}$ of $H$. Without loss of generality, we may assume that $Z$ is on the left of $(v_j,v_{j+1})$. Then $(v_j,v_{j+1})$ adds $Z$ to $\Phi((\overrightarrow{K}'_J, \text{left}))$. See Figure 1.14. By Corollary 1.3.10, $(v_j,v_{j+1}) \in \overrightarrow{K}'_J$. Similarly, $(v_{j+1},v_j) \in \overrightarrow{K}'_J$. The maximality of $\overrightarrow{K}'_J$ implies that exactly one of $\overrightarrow{K}\big((v_j,v_{j+1}),(v_{j+1},v_j) \big)$ and $\overrightarrow{K}\big((v_{j+1},v_j),(v_j,v_{j+1}) \big)$ is a subtrail of $\overrightarrow{K}'_J$. For definiteness, assume that $\overrightarrow{K}\big((v_j,v_{j+1}),$ $(v_{j+1},v_j) \big)$ is a subtrail of $\overrightarrow{K}'_J$. Note that the end-vertex $v_r$ of $H$ is incident on the edge $(v_{r-1},v_r)$ of $\overrightarrow{K}\big((v_j,v_{j+1}),(v_{j+1},v_j) \big)$ and so $v_r$ is incident on the edge $(v_{r-1},v_r)$ of $\overrightarrow{K}'_J$.

To see that $v_r \in G \setminus B(G)$, BWOC assume that $v_r \in B(G)$. By definition of FTW, (see Figure. 1.6 (iv)) all four boxes incident on $v_r$ must belong to $J$. But then at least one box incident on $v_r$ belongs to $G_{-1} \setminus G$, and so $(G_{-1} \setminus G) \cap J \neq \emptyset$, contradicting the definition of $J$. $\square$

\endgroup 

%\noindent \textbf{Lemma 1.3.13.} Let $H$ be a Hamiltonian path of a simply connected polyomino $G$ and let $J$ be an $H$-component of $G$. Assume that there are no end-vertices of $H$ incident on $J \cap (G \setminus B(G))$. Let $P_J$ be the subpath of $H$ determined by a main trail of $J$ and let $e'$ be the neck edge of $J$. Then $B(J) \setminus e'=P_J$.

%\null 

% \noindent \textit{Proof.} We will make use of Observation 1.2.4 (d) and Corollary 1.3.8 (c) implicitly and repeatedly. Since there are no end-vertices of $H$ incident on $J \cap (G \setminus B(G))$, by Lemma 1.3.12, $J$ is a non-self-adjacent $H$-component of $G$. 

% Let the edge $e \in P_J$, adjacent to boxes $X$ and $Y$ of $G_{-1}$, correspond to the directed edge $\overrightarrow{e}$ of the right main trail $\overrightarrow{K}_J$ of $J$. For definiteness, assume that $\overrightarrow{e}$ adds the box $X$ to $J$. Since $J$ is non-self-adjacent, $Y$ must belong to $G_{-1} \setminus J$. Then $e \in B(J)\setminus e'$.

% Conversely, let $e=\{u,v\} \in B(J)\setminus e'$. Then by Corollary 1.3.8 (a), $e$ is incident on a box $Z$ of $J$ and either $Z=\Phi ((u,v), \text{right})$ or $Z=\Phi ((v,u), \text{right})$. For definiteness, assume that $\Phi ((u,v), \text{right})=Z$. By Corollary 1.3.10, $(u,v)$ belongs to the right main trail $\overrightarrow{K}_J$ of $J$. It follows that $e \in P_J$ $\square$

\null 

%We will first show that we can reconfigure Hamiltonian cycles and e-cycles into one another. Then we will show that we can reconfigure paths into e-cycles. `

\noindent \noindent \textbf{Hamiltonian e-cycles, cycles, and cookies.} Let $G$ be a polyomino, let $H$ be a Hamiltonian cycle of $G$ and let $e$ be an edge of $H$ in the boundary of $G$. We call the path $H \setminus e$ a Hamiltonian \index{e-cycle|textbf}\textit{e-cycle} of $G$ with $H$-components $J_0, \ldots, J_s$, where $J_0$ is the $H$-component that contains both end-vertices of $H$. By Lemma 1.3.12, $J_0$ and all other $H$-components of $G$ are non-self-adjacent. 
We define $\text{int}(H)$ to be the $H$-component $J_0$ and ext$H$ to be the complement of $J_0$ in $G$. We will call the $H$-components $J_1, \ldots, J_s$ \index{cookie|textbf}\textit{cookies}\footnote{This definition is similar to the one given in ~\cite{nishat2017bend}.}. If a cookie $J$ consists of exactly one box, we call that $J$ a \index{small cookie|textbf}\textit{small} cookie. Otherwise we call $J$ a \index{small cookie|textbf}\textit{large} cookie See Figure 1.15.

\begin{center}

\begin{adjustbox}{trim=0cm 0cm 0cm 0cm}
\begin{tikzpicture}[scale=1.5]

\begin{scope}[xshift=0cm]
{
\draw[gray,very thin, step=0.5cm, opacity=0.5] (0,0) grid (3.5,3.5);

\fill[orange!75!white,opacity=0.5](0.5,3) rectangle  (1,3.5); 

\fill[orange!75!white,opacity=0.5](3,0.5) rectangle  (3.5,1); 
\fill[orange!75!white,opacity=0.5](3,1.5) rectangle  (3.5,2);

\fill[blue!50!white,opacity=0.5](1,0) rectangle  (1.5,1); 
\fill[blue!50!white,opacity=0.5](0.5,0.5) rectangle (1,2.5); 

\fill[blue!50!white,opacity=0.5](2,0) rectangle  (2.5,1); 
\fill[blue!50!white,opacity=0.5](1.5,1.5) rectangle (2,3); 

\fill[blue!50!white,opacity=0.5](2,1.5) rectangle (2.5,2); 
\fill[blue!50!white,opacity=0.5](2,2.5) rectangle (3.5,3);

\draw[blue, line width=0.65mm] (0.5,0)--++(-0.5,0)--++(0,3.5)--++(0.5,0)--++(0,-0.5)--++(0.5,0)--++(0,0.5)--++(2.5,0)--++(0,-0.5)--++(-2,0)--++(0,-1.5)--++(1,0)--++(0,0.5)--++(-0.5,0)--++(0,0.5)--++(1.5,0)--++(0,-0.5)--++(-0.5,0)--++(0,-0.5)--++(0.5,0)--++(0,-0.5)--++(-0.5,0)--++(0,-0.5)--++(0.5,0)--++(0,-0.5)--++(-1,0)--++(0,1)--++(-0.5,0)--++(0,-1)--++(-0.5,0)--++(0,1)--++(-0.5,0)--++(0,1.5)--++(-0.5,0)--++(0,-2)--++(0.5,0)--++(0,-0.5);

\draw[fill=blue, opacity=1] (0.5,0) circle [radius=0.05];
\draw[fill=blue, opacity=1] (1,0) circle [radius=0.05];

\node[right, align=left, text width=6cm] at (4,1) 
{Fig. 1.15. An $8 \times 8$ grid graph  with a Hamiltonian e-cycle in   blue.  Large cookies  shaded blue.  Small cookies shaded orange.  $J_0=\textrm{int}(H)$ in white.};

}
\end{scope}

\end{tikzpicture}
\end{adjustbox}
\end{center}

%pagemarker

\noindent Let $H$ be a Hamiltonian path of $G$. A box of $G$ with vertices $a,b,c,d$ is \index{switchable|textbf}\textit{switchable} in $H$ if it has exactly two edges in $H$ and those edges are parallel to each other.

\null

\noindent \textbf{Lemma 1.3.14.} Let $G$ be a simply connected polyomino, let $H$ be an $e$-cycle of $G$ and let $J_0, \ldots, J_s$ be the $H$-components of $G$. Then every edge of $H$ is incident on a box of $J_0$ and a box of $G_{-1} \setminus J_0$.

\null

\noindent \textit{Proof.} Let $H=v_1, \ldots, v_r$. We adopt here the partitioning of the Hamiltonian path into three subpaths and all relevant notation from Corollary 1.3.8. Note that $\{v_s,v_t\}$ in Corollary 1.3.8 corresponds to $\{v_1,v_r\}$ here. Then $P_1=v_1$, $P_3=v_r$ and $P_2=H$. For definiteness, by Corollary 1.3.8 (a), we may assume that $\textrm{Boxes} (\Phi (H, \textrm{right}))\subseteq J_0$ and $\textrm{Boxes} (\Phi (P_2, \textrm{left})) \subseteq G_{-1} \setminus J_0$. Then every edge of $H$ is incident on a box of $J_0$ and a box of $G_{-1} \setminus J_0$. $\square$

\null 

\noindent \textbf{Corollary 1.3.15.} Let $G$, $H$ and $J_0, \ldots, J_s$ be as in Lemma 1.3.14. Then:

 (a) Boxes of distinct cookies\index{cookie} cannot be adjacent to one another.

 (b) A large cookie has exactly one box incident on a boundary edge of $G$, namely 
 
\hspace*{0.6cm} its neck. Furthermore, if $G$ is an $m \times n$ grid graph, then the neck of each large  

\hspace*{0.6cm} cookie is switchable.

\null 

\noindent \textit{Proof of (a).} BWOC assume there are boxes $X\in J_i$, $Y\in J_j$ with $1 \leq i,j \leq s$, $i\neq j$ such that $X$ and $Y$ are adjacent sharing an edge $e$ of $H$. Then, by Lemma 1.3.14, at least one of $X$ and $Y$ must belong to $J_0$, contradicting that $i,j \geq 1$.

\null

\noindent \textit{Proof of (b).} First we show that a large cookie has exactly one box incident on a boundary edge of $G$, namely its neck. Let $J_i$ be a cookie. By Corollary 1.3.2, $J_i$ has a neck $N_{J_i}$ and $N_{J_i}$ is incident on $B(G)$. Suppose that there is another box $X$ of $J_i$ that is incident on $B(G)$. Note that $X \in G_{-1} \setminus J_0$. Let $e$ be the boundary edge of $X$ and let $Y$ be the box in $G_{-1} \setminus G$ that is incident on $e$. Then either $e\in H$ or $e\notin H$. Note that $e \notin H$ contradicts Corollary 1.3.2 so we only need to check the case where $e\in H$. Suppose that $e\in H$. $X \in G_{-1} \setminus J_0$ and Lemma 1.3.14 imply that $Y$ is in $J_0 \subset G$. But then $Y \in G$ and $Y\in G_{-1} \setminus G$ is a contradiction.

\null 

\noindent Now we show that in an $m \times n$ grid graph, the neck of each large cookie is a switchable box. Let $X=R(k,l)$ be the neck of a cookie $J_i$. Let $v(k,l)=a$, $v(k+1,l)=b$, $v(k+1,l+1)=c$ and $v(k,l+1)=d$. For definiteness, assume that $\{a,b\}$ is the neck edge of $J_i$. Since $i \neq 0$, neither $a$ nor $b$ is an end-vertex of $H$. This, together with the assumption that $G$ is an $m \times n$ grid graph, implies that $\{a,d\} \in H$ and $\{b,c\} \in H$. Since $J_i$ is a large cookie,  $\{d,c\} \notin H$. Thus $X$ is switchable. $\square$

\null

\noindent \textbf{Remark.} All the definitions made for the case where $H$ is a Hamiltonian e-cycle of a simply connected polyomino $G$, as well as Lemma 1.3.14 and Corollary 1.3.15, translate immediately to the case where $H$ is a Hamiltonian cycle of $G$. Let $H'$ be a Hamiltonian e-cycle of $G$. We may just add the edge $e$ incident on the end-vertices of the $H'$ to obtain the cycle $H$. $J_0$ no longer has a neck, and all other $H$-components and the properties we have observed remain unchanged.

\subsection{Moves}

Let $G$ be a polyomino and let $H$ be a Hamiltonian path of $G$. We will define here three types of moves we can apply to $H$ in order to obtain a new Hamiltonian cycle $H'$: switches, double-switches and backbites. In Chapter 2 we describe an algorithm that can reconfigure any two Hamiltonian cycles $H$ and $H'$ of an $m \times n$ grid graph into one another by using only double-switch moves. In Chapter 5 we describe an algorithm that can reconfigure any two Hamiltonian paths $H$ and $H'$ of an $m \times n$ grid graph into one another by using switch, double-switch and backbite moves.

\null

\noindent \textbf{Switch and double-switch moves.} Let $H$ be a Hamiltonian path of a polyomino $G$. Let $abcd$ be a switchable box with edges $ab$ and $cd$ in $H$. We define a \index{switch move|textbf}\textit{switch move} on the box $abcd$ in $H$ as follows: remove edges $ab$ and $cd$ and add edges $bc$ and $ad$. Let $X \in G$ be a switchable box in $H$. We denote a switch move by $\textrm{Sw}(X)$.

A \index{double-switch move|textbf}\textit{double-switch move} is a pair of switch operations where we first switch $X$ and then find a switchable $Y$ and switch it. We denote a double-switch move by $X \mapsto Y$. See Figure 1.14. If after a move, be it switch or double switch, we get a new Hamiltonian path, then we call the move a \index{valid move|textbf}\textit{valid move}. From here on we will often say ``move" to mean ``valid move''. If $X$ and $Y$ share an edge and $X \mapsto Y$ is a valid move, we call $X \mapsto Y$ a \textit{flip} move. 

It will be useful to consider $H$ as a directed path, so we assign an orientation  $v_1, \ldots, v_r$ to $H$. Let $X$ be a switchable box in $H$ with edges $e_1=(v_s,v_{s+1})$ and $e_2=(v_t,v_{t+1})$. We say $e_1$ and $e_2$ are \index{parallel|textbf}\textit{parallel} if $v_s$ is adjacent to $v_t$ in $G$, and \index{anti-parallel|textbf}\textit{anti-parallel} if $v_s$ is adjacent to $v_{t+1}$. Similarly, we call the box $X$ a \textit{parallel (anti-parallel)} box if its edges are parallel (anti-parallel). 

%pagemarker
\noindent Let $P(v_x,v_y)$ be a subpath of $H$ and assume that $\mu$ is a valid move that gives a Hamiltonian path $H'=v_1', \ldots, v_r'$, where $v_1'=v_1$ and $v_r'=v_r$. If there is a subpath $P(v_x', v_y')$ of $H'$ such that $P(v_x', v_y')=P(v_x,v_y)$, then we say that $\mu$ \index{fix (a subpath of a Hamiltonian path)|textbf}\textit{fixes} $P(v_x,v_y)$; and if there is a subpath $P(v_x', v_y')$ of $H'$ such that $P(v_x', v_y')=P(v_y,v_x)$, then we say that $\mu$ \index{reverse (a subpath of a Hamiltonian path)|textbf}\textit{reverses} $P(v_x,v_y)$.

Let $H=v_1, \ldots, v_r$.
We will show below that if we switch a switchable box of $H$ with anti-parallel edges, we get a path $H_p$ and a cycle $H_c$. We define a \index{H p H c port@$(H_p,H_c)$-port|textbf}$(H_p,H_c)$\textit{-port} to be a switchable box in $H_p \cup H_c$ that has one edge in $H_p$ and the other in $H_c$. See Figure 1.16 

\null 

\begingroup
\setlength{\intextsep}{0pt}
\setlength{\columnsep}{20pt}
\begin{center}
\begin{adjustbox}{trim=0cm 0cm 0cm 0cm}

\end{adjustbox}
\end{center}

\noindent \textbf{Lemma 1.4.1.} Let $G$ be a polyomino. Let  $H=v_1, \ldots, v_r$ be a Hamiltonian path of $G$ and let $s+1<t$. Let $X$ be a switchable box of $H$ with edges $e_1=(v_s,v_{s+1})$ and $e_2=(v_t,v_{t+1})$.  Let $P_1=P(v_1, v_s)$, $P_2=P(v_{s+1},v_t)$ and $P_3=P(v_{t+1},v_r)$. 

(i) \ \ If $e_1$ and $e_2$ are parallel, $\textrm{Sw}(X)$ is a valid move that fixes $P_1$ and $P_3$ and 

\hspace*{0.6cm} \ reverses $P_2$.

(ii) \ If $e_1$ and $e_2$ are anti-parallel $\textrm{Sw}(X)$ splits $H$ into a cycle $H_c$ and a path $H_p$.

  (iii) Suppose $e_1$ and $e_2$ are anti-parallel and assume we apply $\textrm{Sw}(X)$. If $Y$ is an 

\hspace*{0.6cm} $(H_p,H_c)$-port then $X \mapsto Y$ is a valid double-switch move.

\null 

\noindent \textit{Proof.} Let $X$ be a switchable box of $H$. Without loss of generality we may assume that $s+1 < t$. If we remove $e_1$ and $e_2$, $H$ splits into three disjoint sub-paths: $P_1=P(v_1, v_s)$, $P_2=P(v_{s+1},v_t)$ and $P_3=P(v_{t+1},v_r)$. 

\null 

\noindent \textit{Proof of (i).} Suppose that $e_1$ and $e_2$ are parallel. Then $v_s$ is adjacent to $v_t$, and $v_{s+1}$ is adjacent to $v_{t+1}$. Now, adding $e_1'=(v_s,v_t)$ and $e_2'=(v_{s+1},v_{t+1})$ gives a new Hamiltonian path $H'$ on $v_1, \ldots, v_s, v_t, \ldots, v_{s+1}, v_{t+1}, \ldots, v_r$. See Figure 1.17. End of proof for (i).

\null 



\noindent \textit{Proof of (ii).} Suppose that $e_1$ and $e_2$ are anti-parallel. Then we have that $v_s$ is adjacent to $v_{t+1}$ and $v_{s+1}$ is adjacent to $v_t$. Now adding $e_1'=(v_s,v_{t+1})$ gives path $H_p$ on vertices 
$v_1, \ldots, v_s, v_{t+1}, \ldots, v_r$ that joins $P_1$ with $P_3$; and adding $e_2'=(v_{s+1},v_t)$ joins the endpoints of $P_2$ to give a cycle $H_c$ on $v_{s+1}, v_{s+2}, \ldots, v_t, v_{s+1}$. See Figure 1.18. End of proof for (ii).

The proof of (iii) is essentially the same as the proof for (ii). $\square$ 

\null

\noindent \textbf{Corollary 1.4.2.} Let $H=v_1,\ldots, v_r$ be a Hamiltonian path of a polyomino $G$. Let $X$ be an anti-parallel switchable box with edges $e_1(X)$ and $e_2(X)$, and let $Y$ be a parallel switchable box that is edge-disjoint from $X$, on edges $(v_s, v_{s+1})$ and $(v_t, v_{t+1})$. Let $P_1=P(v_1,v_s)$,  $P_2=P(v_{s+1},v_t)$ and $P_3=P(v_{t+1},v_r)$, and assume that $s+1<t$. Then:

  (i) \ If $e_1(X) \in P_i$ and $e_2(X) \in P_j$ where $(i,j) \in \{(1,1), (2,2),(3,3), (1,3), (3,1)\}$, then 
 
  \hspace*{0.6cm} after Sw($Y$), $X$ remains anti-parallel.

  (ii) If $e_1(X) \in P_i$ and $e_2(X) \in P_j$ where $(i,j) \in \{(1,2), (2,1),(2,3), (3,2)\}$, then after  
  
\hspace*{0.6cm} Sw($Y$), $X$ becomes parallel. $\square$

\null

\begingroup
\setlength{\intextsep}{0pt}
\setlength{\columnsep}{20pt}

\begin{adjustbox}{trim=0cm 0cm 0cm 0cm}

\end{adjustbox}

\noindent \textbf{Corollary 1.4.3.}  Let $H$ be a Hamiltonian cycle or e-cycle of an $m\times n$ grid graph $G$, and let $X$ be a switchable box with edges $e_1$ and $e_2$. By Lemma 1.5.1 in the next section, all switchable boxes in a Hamiltonian cycle or e-cycle are anti-parallel. Thus, $e_1$ and $e_2$ are anti-parallel.

 (i) \ $\textrm{Sw}(X)$ splits $H$ into the cycles $H_1$ and $H_2$.

 (ii) Assume we apply $\textrm{Sw}(X)$. If $Y$ is an $(H_1,H_2)$-port then $X \mapsto Y$ is a valid 
 
\hspace*{0.6cm} double-switch move. See Figure 1.19. $\square$

\endgroup 

\null

\noindent \textbf{Backbite moves.} A \index{backbite move|textbf}\textit{backbite move} consists of adding an edge $e$ of $G \setminus H$ incident on an endpoint of $H$ and simultaneously removing an edge $e'$ of $H$ adjacent to $e$ such that the resulting graph $H'$ is a Hamiltonian path. We will use the notation $e' \mapsto e$ to indicate that we are removing the edge $e' \in H$ and adding the edge $e \in G \setminus H$. 

Here is a more detailed description of a backbite move. Let $H=v_1, v_2, \ldots, v_r$ be a Hamiltonian path of $G$. First we choose a vertex $v_s \in G$ adjacent to $v_1$ such that $\{v_1, v_s\} \in G \setminus H$. We will add to $H$ the edge $\{v_1, v_s\}$ and remove from $H$ the edge $(v_{s-1},v_s)$. This is a valid backbite move that gives a Hamiltonian path $H'$ distinct from $H$. See Figure 1.20.



\noindent \textbf{Lemma 1.4.4.} Let $H=v_1, v_2, \ldots, v_r$ be a Hamiltonian path of a graph $G$  and assume that $\{v_1,v_s\} \in G \setminus H$. Then $(v_s,v_{s-1}) \mapsto \{v_1 ,v_s\}$ is a valid backbite move and $(v_s,v_{s+1}) \mapsto \{v_1,v_s\}$ is not. Furthermore, $(v_s,v_{s-1}) \mapsto \{v_1 ,v_s\}$ reverses $P(v_1, v_{s-1})$ and fixes $P(v_s, v_r)$, and the end-vertex in the resulting Hamiltonian path is located on the vertex of $G$ where $v_{s-1}$ was located in $H$.

%Similarly, $(v_t,v_{s+1}) \mapsto \{v_t,v_r\}$ is a valid backbite move and $(v_{t-1}, v_t) \mapsto \{v_t,v_r\}$ is not. 

\null

\noindent \textit{Proof.} Suppose we add the edge $\{v_1,v_s\}$ to $H$. Note that if we remove the edge $\{v_s,v_{s+1}\}$ we get a cycle $H_c$ consisting of $P(v_1,v_s)$ and $\{v_1,v_s\}$ and a path $H_p=P(v_{s+1},v_r)$ on the vertices $v_{s+1}, v_{s+2}, \ldots, v_r$. Thus $(v_s,v_{s+1}) \mapsto \{v_1,v_s\}$ is not a valid backbite move. If we remove the edge $(v_{s-1},v_s)$ we get a new Hamiltonian path $H'=v_{s-1}, v_{s-2}, \ldots, v_1, v_s, v_{s+1}, \ldots,$ $v_r$, where we see that $(v_s,v_{s-1}) \mapsto \{v_1 ,v_s\}$ fixes $P(v_s, v_r)$ and reverses $P(v_1, v_{s-1})$. See Figure 1.20. $\square$ 

\null

\noindent \textbf{Corollary 1.4.5.} Let $H=v_1,\ldots, v_r$ be a Hamiltonian path of a polyomino $G$. Let $X$ be an anti-parallel switchable box with edges $e_1(X)$ and $e_2(X)$, such that neither $\{v_1,v_s\}$ nor $\{v_{s-1},v_s\}$ is an edge of $X$. Let $P_1=P(v_1,v_{s-1})$ and  $P_2=P(v_s,v_r)$. Then:

  (i) If $e_1(X) \in P_i$ and $e_2(X) \in P_j$ where $(i,j) \in \{(1,1), (2,2)\}$, then after 
  
\hspace*{0.6cm} $(v_s,v_{s-1}) \mapsto \{v_1,v_s\}$, $X$ remains anti-parallel.

 (ii) If $e_1(X) \in P_i$ and $e_2(X) \in P_j$ where $(i,j) \in \{(1,2), (2,1)\}$, then after 
 
\hspace*{0.6cm} $(v_s,v_{s-1}) \mapsto \{v_1,v_s\}$, $X$ becomes parallel. $\square$

%pagemarker
\noindent Let $e \mapsto f$ be a backbite move. Note that if the end-vertices of $H$ are not adjacent then $e \mapsto f$ is fully characterized by the edge $f$ that we choose to add to $H$; and if the end-vertices of $H$ are adjacent then $e \mapsto f$ is fully characterized by the edge $f$ that we choose to add to $H$ and the end-vertex to which the edge $e$ is added. Whenever the edge that a backbite move removes is irrelevant to the argument, we will use the notation $bb_{v_1}(north)$ to indicate that the backbite move we're applying is the one that adds to $H$ the edge of $G$ on $\{v_1,v_{\textrm{north}}\}$, where $v_{\textrm{north}}$ is the vertex of $G$ north of $v_1$. The notation for backbite moves in other directions is defined analogously.

% DELETED stuff about reversing subpaths

{%\textbf{Reversed subpaths} Let $H$ be the Hamiltonian path of $G$. Fix an orientation $v_1, \ldots v_r$ of $H$. Let $P(v_{s_1}, v_{s_2})$ and $P(v_{t_1}, v_{t_2})$ be subpaths of $H$. Let $\mu$ be a  backbite move or a switch move on $H$. Applying $\mu$ gives a new Hamiltonian path $H'$. If there is a subpath $P(u,w)$ of $H'$ such that $P(u,w)=P(v_{s_2},v_{s_1})$, then we say that $\mu$ \textit{reverses} the subpath $P(v_s, v_t)$ of $H$; and if there is a subpath $P(x,y)$ of $H'$ such that $P(x,y)=P(v_{t_1},v_{t_2})$, then we say that $\mu$ \textit{fixes} the subpath $P(v_{t_1}, v_{t_2})$ of $H$. 

%Let $H$, $e_1$,  $e_2$ and $X$ be as in (i) of  Then Sw($X$) reverses the subpath $P(v_t,v_{s+1})$ and \textit{fixes} the subpaths $P(v_1,v_s)$ and $P(v_{t+1},v_r)$.

%Let Let $H$ be as in Observation 1.4.4 Then the backbite  move $(v_s,v_{s-1}) \mapsto \{v_1 ,v_s\}$ reverses the subpath $P(v_1, v_{s-1})$ and fixes the subpath $P(v_s,v_r)$. \textcolor{red}{I need to bring $H$-components into this $\ldots$}
}

\null

\noindent \textbf{Definitions.} We define the \index{inverse of a move|textbf}\textit{inverse} of a move $\mu$ to be the move $\mu^{-1}$, such that after applying $\mu$ and $\mu^{-1}$, we get back the original Hamiltonian path. We call a move $X \mapsto X$ a \index{trivial move|textbf}\textit{trivial} move. We define a \index{cascade|textbf}\textit{cascade} to be a sequence of moves $\mu_1, \ldots, \mu_r$ such that for $0\leq j\leq r-1$:

1) $\mu_1$ is valid, 

2) if $\mu_1, \ldots,\mu_j$ have been applied then $\mu_{j+1}$ is valid, and

3) the sequence does not create any new \index{cookie}cookies.

\null 

\noindent Let $H$ be a Hamiltonian e-cycle of an $m \times n$ grid graph $G$ and let $J$ be a cookie of $H$ with neck $N_J$. Consider a cascade $\mu_1, \ldots, \mu_r$ where $\mu_r$ is the nontrivial move $Z \mapsto N_J$. We say that the cascade $\mu_1, \ldots, \mu_r$ \index{collect (a cookie, or a leaf)|textbf}\textit{collects} the cookie $J$. Note that switches, double-switches and backbites are all invertible moves. For non-adjacent boxes $X$ and $Y$, the moves $X \mapsto Y$ and $Y \mapsto X$ yield the same result. When $X$ and $Y$ are adjacent with $X$ switchable and $Y$ a leaf (i.e. $X \mapsto Y$ is a flip move), $X$ must be switched first before $Y$ becomes switchable, so the order matters. See Figure I.3.

%\noindent \textcolor{red}{ Section 2 contains algorithms we can use to reconfigure one Hamiltonian cycle (e-cycle) into another. Proofs of existence are in Section 3. In Section 4 we prove a lemma used in Section 3, and some other results that will be needed for Section 5. In Section 5 we extend the result to show that we can reconfigure a Hamiltonian path into another.}

\subsection{Supporting lemmas}

\noindent The Lemmas in this section are used in Chapters 4 and 5.

\null 

\noindent \textbf{Lemma 1.5.1.} Let $H=v_1, \ldots, v_r$ be a Hamiltonian cycle or e-cycle of a polyomino $G$. Then all switchable boxes of $G$ are anti-parallel.

\begingroup
\setlength{\intextsep}{0pt}
\setlength{\columnsep}{20pt}
\begin{wrapfigure}[]{l}{0cm}
\begin{adjustbox}{trim=0cm 0cm 0cm 0cm}
\begin{tikzpicture}[scale=2]
\begin{scope}[xshift=0cm] 
\draw[gray,very thin, step=0.5cm, opacity=0.5] (0,0) grid (1.5,1.5);
\fill[green!40!white,opacity=0.5] plot [smooth, tension=0.75] coordinates {(0.5,1) (1,1.5)(1.5,1)(1,0.5)(0.5,0.5)};
\begin{scope}
[very thick,decoration={
    markings,
    mark=at position 0.6 with {\arrow{>}}}
    ]
    \draw[postaction={decorate}, blue, line width=0.5mm] (0.5,0.5)--++(0,0.5);
    \draw[postaction={decorate}, blue, line width=0.5mm] (1,0.5)--++(0,0.5);
\end{scope}
\draw[red, dotted, line width=0.5mm] (0,0.5)--++(1,0.5);
%for the shaded region
\draw [blue, line width=0.5mm] plot [smooth, tension=0.75] coordinates  {(0.5,1) (1,1.5)(1.5,1)(1,0.5)};
\draw[fill=blue, opacity=1] (0.5, 0.5) circle [radius=0.035];
\draw[fill=blue, opacity=1] (0.5, 1) circle [radius=0.035];
\draw[fill=blue, opacity=1] (1, 0.5) circle [radius=0.035];
\draw[fill=blue, opacity=1] (1, 1) circle [radius=0.035];
\draw[opacity=1] (0,0.5) circle [radius=0.05];
\draw[opacity=1] (0.5,0) circle [radius=0.05];
\node[below] at  (0,0.5) [scale=0.8]{\small${v_{s-1}}$};
\node[left] at  (0.5,0.5) [scale=0.8]{\small${v_s}$};
\node[left] at  (0.5,1) [scale=0.8]{\small${v_{s+1}}$};
\node[left] at  (1,0.5) [scale=0.8]{\small${v_t}$};
\node[left] at  (1,1) [scale=0.8]{\small${v_{t+1}}$};
\node at  (1.25, 1) [scale=1.5]{U};
\node at  (0.75,0.75) [scale=1]{X};
\node[above] at  (0.5, 1.5) [scale=1]{\tiny{$k$}};
\node[above] at  (1, 1.5) [scale=1]{\tiny{$+1$}};
\node[above] at  (0, 1.5) [scale=1]{\tiny{$-1$}};
\node[right] at  (1.5, 0.5) [scale=1]{\tiny{$\ell$}};
\node[right] at  (1.5, 1) [scale=1]{\tiny{$+1$}};
\node[right] at  (1.5, 0) [scale=1]{\tiny{$-1$}};

\node[below] at (0.75,0) [scale=1]{\begin{tabular}{c} Fig. 1.21.  Case 1: $X\in U$. \end{tabular}};
 
\end{scope}
\end{tikzpicture}
\end{adjustbox}
\end{wrapfigure}

\null

\noindent \textit{Proof.} Assume that $H$ is a Hamiltonian cycle or e-cycle of $G$. For a contradiction, assume that there is a parallel box $X$ with edges $(v_s, v_{s+1})$ and $(v_t, v_{t+1})$, and without loss of generality, that $s <t+1$. For definiteness, let $X=R(k,l)$, $v_s=v(k,l)$, and $v_{s+1}=v(k,l+1)$. Then $v_t=v(k+1,l)$ and $v_{t+1}=v(k+1,l+1)$. Let $Q$ be the cycle $P(v_s,v_t), \{v_t,v_s\}$ and let $U$ be the region of $G$ bounded by $Q$. Either $X \in U$ or $X \notin U$.

\endgroup 

\null 

\noindent  \textit{CASE 1: $X \in U$.} Then, by Lemma 1.3.3 (b) $v_{t+1} \in U \setminus Q$. By Lemma 1.3.3 (c), $P(v_{t+1},v_r)$ is contained in $U \setminus Q$ as well, so $v_r \in U \setminus Q$. Now, $v_{s-1}=v(k-1,l)$, or $v_{s-1}=v(k,l-1)$, or $v_s=v_1$. 
Note that if $v_s=v_1$, then all neighbours of $v_s$ must belong to $Q$. But then none of them can be $v_r$ which contradicts that  $H$ is a Hamiltonian cycle or e-cycle. 

Then, either $v_{s-1}=v(k-1,l)$, or $v_{s-1}=v(k,l-1)$. See Figure 1.21. Without loss of generality, assume that $v_{s-1}=v(k-1,l)$. Since the segment $[v(k-1,l),v_{t+1}]$ (in red in Figure 1.21) satisfies the hypotheses of Corollary 1.1.5, $v(k-1,l)$ must belong to $G \setminus (U \cup Q)$. By JCT, $P(v_1, v(k-1,l))$ is contained in $G \setminus (U \cup Q)$ as well, so $v_1 \in G \setminus (U \cup Q)$. But then JCT implies that the path $v_{t+1}, \ldots, v_r,v_1$ intersects $Q$, contradicting that $H$ is a Hamiltonian cycle or e-cycle. End of Case 1.

\endgroup 

\null 

\noindent \textit{CASE 2: $X \notin U$.} This is similar to Case 1, so we omit the proof. $\square$

\null 

\noindent \textbf{Corollary 1.5.2.} Let $G$ be an $m \times n$ grid graph and let $H$ be a Hamiltonian path of $G$. If both end-vertices of $H$ are in $B_0$, then every switchable box is anti-parallel. $\square$

\begingroup
\setlength{\intextsep}{0pt}
\setlength{\columnsep}{20pt}
\begin{wrapfigure}[]{r}{0cm}
\begin{adjustbox}{trim=0cm 0cm 0cm 0.25cm}
\begin{tikzpicture}[scale=1.75]
\usetikzlibrary{decorations.markings}

\begin{scope}[xshift=0cm] 

\draw[gray,very thin, step=0.5cm, opacity=0.5] (0,0) grid (2.5,1.5);

\fill[green!40!white,opacity=0.5] (0,0.5)--++(0,1)--++(2.5,0)--++(0,-1);
\fill[orange!40!white,opacity=0.5] (0.5,0)--++(0,0.5)--++(1.5,0)--++(0,-0.5);

\draw[red, dotted, line width=0.5mm] (0.5,0.5)--++(0,-0.5)--++(1.5,0)--++(0,0.5);

%for the shaded region
\draw [blue, line width=0.5mm] plot [smooth, tension=0.75] coordinates  {(0.5,0.5) (0,1)(0.5,1.5)(1,1.25)(0.85,0.85)(1.25,0.5)(1.5,1)(1.4,1.4)(1.75,1.475)(2.5,1)(2,0.5)};

\draw[fill=blue, opacity=1] (0.5, 0.5) circle [radius=0.035];
\draw[fill=blue, opacity=1] (2,0.5) circle [radius=0.035];

\node[above] at  (0.5,0.5) [scale=0.8]{\small${u}$};
\node[above] at  (2,0.5) [scale=0.8]{\small${v}$};

\node[above] at  (0.5, 1.5) [scale=1]{\tiny{$a$}};
\node[above] at  (2, 1.5) [scale=1]{\tiny{$b$}};

\node[right] at  (2.5, 0) [scale=1]{\tiny{$-1$}};
\node[right] at  (2.5, 0.5) [scale=1]{\tiny{$0$}};

\node[below] at (1.25,0) [scale=1]{\begin{tabular}{c} Fig. 1.22.  Case 1: $G$ shaded \\  green, $G' \setminus G$ shaded orange, \\ $H$ in blue, $H' \setminus H$ in red. \end{tabular}};;

\end{scope}

\end{tikzpicture}
\end{adjustbox}
\end{wrapfigure}

\null 

\noindent \textit{Proof.}  Let $u$ and $v$ be the end-vertices of $H$. Then $u$ and $v$ are on the same side of $B_0$, or on adjacent sides, or on opposite sides, so there are three cases to check.

\null

\noindent \textit{CASE 1. $u$ and $v$ are on the same side of $B_0$.} Without loss of generality assume that $u=v(a,0)$ and $v=v(b,0)$ are on the southern side of $B_0$. Let $G'$ be the simply connected polyomino obtained by adding to $G$ the boxes $R(a,-1), R(a,-2), \ldots, R(b-1,-1)$, and let $H'$ be the Hamiltonian cycle of $G'$ obtained by adding to $H$ the edges $e(a;0,-1), e(a,a+1;-1), e(a+1,a+2;-1), \ldots, e(b-1,b;-1), e(b;-1,0)$. See Figure 1.22. By Lemma 1.5.1, $G'$ has no boxes that are parallel in $H'$. It follows that $G$ has no boxes that are parallel in $H$. End of Case 1.

\null 

\noindent  The same argument works for the other two cases, so we omit the proof. $\square $

\endgroup 

\null

\noindent \textbf{Corollary 1.5.3.} Let $G$ be an $m \times n$ grid graph, let $H$ be a Hamiltonian path of $G$, and let $J$ be an $H$-component of $G$ with switchable neck $N_J$ such that no end-vertices of $H$ are incident on $J \setminus B(J)$. Then the box $N_J$ is anti-parallel.

\null

\noindent \textit{Proof.} By Corollary 1.3.11, $J$ is a simply connected polyomino, and $\overrightarrow{K}_J$ is a Hamiltonian e-cycle of $J$. Then, by Lemma 1.5.1, $N_J$ is anti-parallel. $\square$

\null

\noindent \textbf{Lemma 1.5.4.} Let $H=v_1, \ldots, v_r$ be a path in a polyomino $G$ such that the box $X$ with edges $(v_s,v_{s+1})$ and $(v_t,v_{t+1})$ in $H$ is anti-parallel. Let $s>1$, $t<r$, and $s+1 <t$. Let $Q$ be the cycle $P(v_s,v_{t+1}), \{v_{t+1},v_s\}$ and let $U$ be the region of $G$ bounded by $Q$. Then, either both end-vertices of $H$ are incident on $U \setminus Q$, or neither is.

\begingroup
\setlength{\intextsep}{0pt}
\setlength{\columnsep}{20pt}
\begin{wrapfigure}[]{l}{0cm}
\begin{adjustbox}{trim=0cm 0.5cm 0cm 0cm}

\end{adjustbox}
\end{wrapfigure}

\null 

\noindent \textit{Proof.}  For definiteness, let $X=R(k,l)$, $v_s=v(k,l)$, and $v_{s+1}=v(k,l+1)$. Then $v_t=v(k+1,l+1)$ and $v_{t+1}=v(k+1,l)$. Let $c=(k+\frac{1}{2}, l+\frac{1}{2})$ be the center of the box $X$. Either $X \in U$ or $X \notin U$.

\null 

\noindent \textit{CASE 1: $X \in U$.} Note that $v_{s-1}=v(k-1,l)$ or $v_{s-1}=v(k,l-1)$. In either case, the segment $[v_{s-1},c]$ satisfies the premises of Corollary 1.1.5, and so $v_{s-1}$ must belong to $G \setminus  (U\cup Q)$. By JCT, $v_1$ must be in $G \setminus (U\cup Q)$. The same argument can be used to show that $v_{t+2}$, and consequently, $v_r$, must also belong to $G \setminus (U\cup Q)$. End of Case 1. See Figure 1.23.

\endgroup 

\null 

\noindent \textit{CASE 2: $X \notin U$.} A similar argument to the one in Case 1 can be used to show that in this case, both end-vertices must be incident on $U \setminus Q$. End of Case 2. $\square$

\null

\noindent \textbf{Corollary 1.5.5.} Let $G$ be an $m \times n$ grid graph, let $H$ be a Hamiltonian path of $G$, and let $J$ be an $H$-component of $G$ with neck $N_J$. If $N_J$ is anti-parallel, then either both end-vertices of $H$ are incident on $J \setminus B(J)$, or neither is. $\square$

\null 

\subsection{Summary}

This chapter proved structural properties of rectangular grid graphs, their Hamiltonian paths and cycles, and their $H$-components that are used throughout the dissertation. Key definitions include the Follow-The-Wall construction, $H$-components, cookies, and necks. Key results include Corollary 1.2.2, Corollary 1.3.11, and the move characterization results in Section 1.4. See Flowchart 1 on the next page for an illustration of how Chapter 1 results are used by Chapters 3--5.

\null 

\hspace*{0.5cm}
\vspace*{-1.5cm}
\makebox[\textwidth][l]{%
  \hspace{-2cm}%
  \includegraphics[scale=0.85]{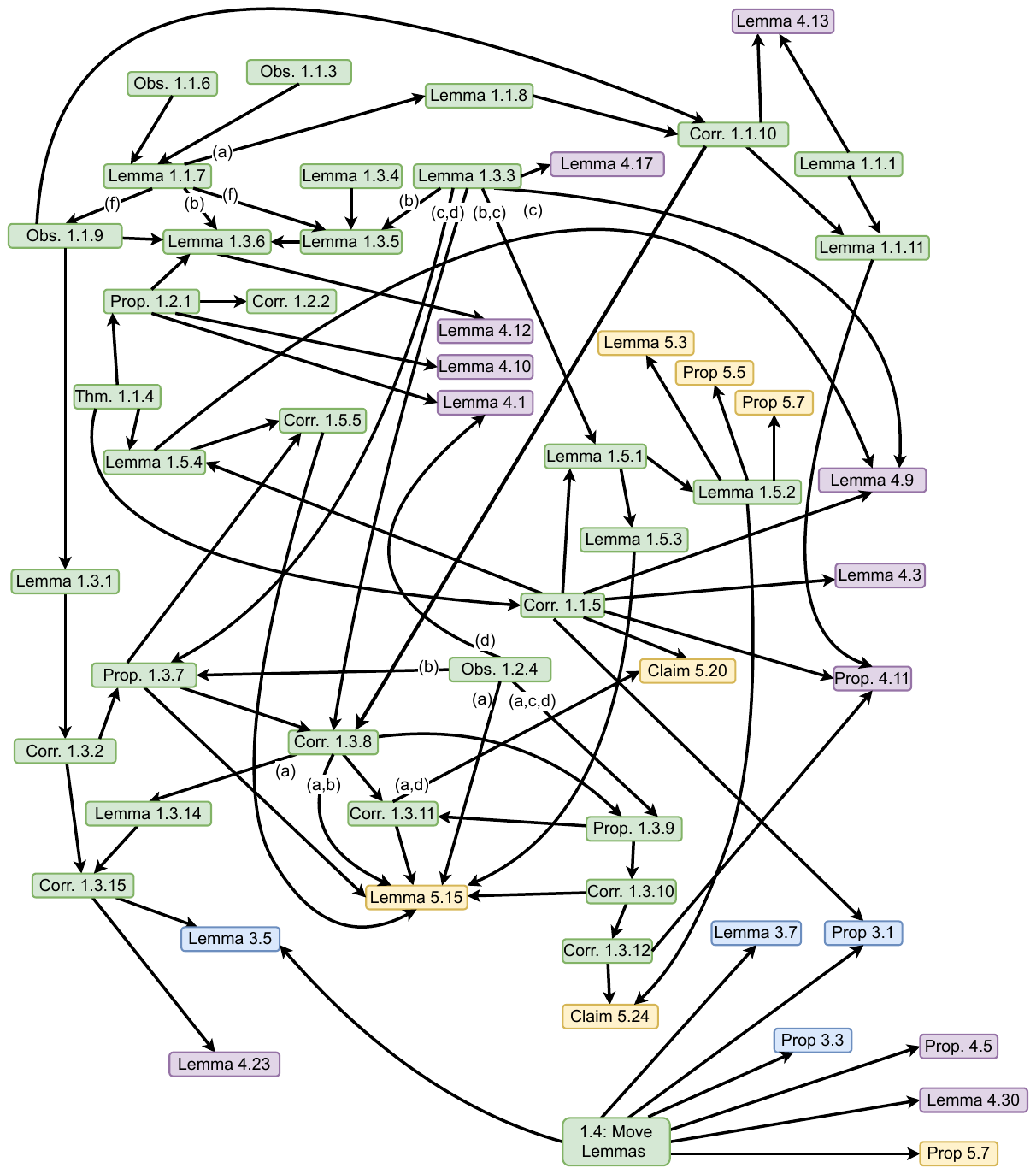}%
}

\vspace*{-1cm}
\noindent \textbf{Flowchart 1.} Dependencies within Chapter 1 and connections to later chapters. Arrows indicate prerequisite relationships between results.

\newpage

\section{Reconfiguration algorithm for Hamiltonian cycles}

\null

\noindent In this chapter we give the algorithms required for the reconfiguration of Hamiltonian cycles on an $m \times n$ grid graph $G$. The $3 \times n$ and $4 \times n$ cases were done by Nishat in \cite{nishat2020reconfiguration}, so from here on, we will assume that $m,n \geq 5$. Furthermore we assume that $m$ and $n$ are not both odd, since in that case $G$ does not have a Hamiltonian cycle. See Note 2.3 at the end of the chapter.

We start with the statement of the main theorem and a few definitions, and follow by an overview of the general reconfiguration strategy.

\null 

\noindent \textbf{Definitions.} Let $G$ be an $m \times n$ grid graph. Denote by $G_s$ the induced subgraph of $G$ on all the vertices with distance $s$ or greater from the boundary of $G$. Denote by $R_s$ the rectangular induced subgraph on vertices of $G$ with distance $s$ from the boundary. Then $R_s$ is the boundary of $G_s$ and the edges of $R_0$ are the boundary edges of $G$. We define $(R_i,R_{i+1})$ to be the set of all the boxes of $G$ adjacent to both $R_i$ and $R_{i+1}$. See Figure 2.1.

\null

\begin{adjustbox}{trim=0 0 0 0}
\begin{tikzpicture}[scale=1.25]

\draw[gray,very thin, step=0.5cm, opacity=0.5] (0,0) grid (4,2.5);

\fill[green!10!white,opacity=0.5]
(0,0) rectangle  (4,2.5);

\usetikzlibrary{patterns} 

\draw[pattern=checkerboard, opacity=0.20] (0.5,0.5) rectangle  (3.5,2);

\draw[pattern=dots, opacity=1] (1,1) rectangle  (3,1.5);

\draw[green!50!black, line width=0.5mm] 
(0,0)--++(0,2.5)--++(4,0)--++(0,-2.5)--++(-4,0);

\draw[orange, line width=0.5mm] 
(0.5,0.5)--++(0,1.5)--++(3,0)--++(0,-1.5)--++(-3,0);

\draw[blue, line width=0.5mm] 
(1,1)--++(0,0.5)--++(2,0)--++(0,-0.5)--++(-2,0);

%\fill[blue!30!white,opacity=0.35](0.5,0.5) rectangle  (3.5,2.5);

\foreach \x in {1,...,4}
\node[left] at (0,0+0.5*\x) [scale=1]
{\small{\x}};

\foreach \x in {1,...,7}
\node[below] at (0+0.5*\x, 0) [scale=1]
{\small{\x}};

\node at (0.5,2.6) [scale=1]{\begin{tabular}{c} $R_0$ \end{tabular}};;

\node at (1.25,2.1) [scale=1]{\begin{tabular}{c} $R_1$ \end{tabular}};;

\node at (2,1.6) [scale=1]{\begin{tabular}{c} $R_2$ \end{tabular}};;

\node[right, align=left, text width=10cm] at (4.5,1.25)
{Fig. 2.1. $R_0$ in green, $G=G_0$ is a $9 \times 6$ grid graph on all light green boxes; $R_1$ in orange, $G_1$ is the rectangle $R(1,7;1,4)$ on all checkered boxes; $R_2$ in blue, $G_2$ is the rectangle $R(2,6;2,3)$ on all dotted boxes. $(R_1,R_2)$ consists of all checkered but not dotted boxes.};

\end{tikzpicture}
\end{adjustbox}

\null

\noindent Theorem 2.1 below is the main result of this dissertation. Its proof will require Chapters 3 and 4.

\null 

\noindent \textbf{Theorem 2.1.} Let $G$ be an $m\times n$ grid graph with $n \geq m$. Let $H$ and $K$ be two \index{Hamiltonian cycle}Hamiltonian cycles or Hamiltonian \index{e-cycle}e-cycles of $G$ with the same edge $e$. Then there is a sequence of at most $n^2m$ valid double-switch moves that reconfigures $H$ into $K$.

\null

\subsection{Canonical forms}

\noindent \textbf{Overview of reconfiguration strategy.} Let $G$ be an $m\times n$ grid graph. We denote the set of all Hamiltonian cycles and e-cycles on an $G$ by $\mathcal{H}(m,n)$. We shall define a subset $\mathcal{H}_{\textrm{can}}(m,n)$ of Hamiltonian cycles of $G$, which we call canonical forms. Then we check that reconfiguration between canonical forms is easy. The problem then reduces to reconfiguring an arbitrary Hamiltonian cycle into a canonical form. Roughly, \index{canonical form}canonical forms have the following structure: the subpath $H_i=G_i \cap H$ of $H$ on each of the nested rectangles $G_0 = G, G_1, G_2, \ldots$ down to the rectangle before the central rectangle is a Hamiltonian e-cycle of that rectangle with exactly one \index{large cookie|textbf}large cookie and no small cookies. See Figure 2.2.

\null

\begin{adjustbox}{trim=0 0 0 0}
\begin{tikzpicture}[scale=1.25]

\draw[gray,very thin, step=0.5cm, opacity=0.5] (0,0) grid (4,2.5);

\draw[green!70!white, line width=1.25mm] (2,1.5)--++(1,0)--++(0,-0.5)--++(-2,0)--++(0,0.5)--++(0.5,0);

\draw[blue, line width=0.5mm] 
(0,0)--++(4,0)--++(0,1)--++(-0.5,0)--++(0,-0.5)--++(-3,0)--++(0,1.5)--++(1,0)--++(0,-0.5)--++(-0.5,0)--++(0,-0.5)--++(2,0)--++(0,0.5)--++(-1,0)--++(0,0.5)--++(1.5,0)--++(0,-0.5)--++(0.5,0)--++(0,1)--++(-4,0)--++(0,-2.5);

\draw[orange, dashed, line width=0.5mm] (3.5,1)--++(0,-0.5)--++(-3,0)--++(0,1.5)--++(1,0)--++(0,-0.5)--++(-0.5,0)--++(0,-0.5)--++(2,0)--++(0,0.5)--++(-1,0)--++(0,0.5)--++(1.5,0)--++(0,-0.5);

\node at (0.5,2.6) [scale=1]{\begin{tabular}{c} $H_0$ \end{tabular}};;

\node at (1.25,2.1) [scale=1]{\begin{tabular}{c} $H_1$ \end{tabular}};;

\node at (2.5,1.6) [scale=1]{\begin{tabular}{c} $H_2$ \end{tabular}};;

%\draw[orange, line width=0.5mm] 
%(0.5,0.5)--++(0,2)--++(3,0)--++(0,-2)--++(-3,0);
%\draw[blue, line width=0.5mm] 
%(1,1)--++(0,1)--++(2,0)--++(0,-1)--++(-2,0);

%\fill[blue!30!white,opacity=0.35](0.5,0.5) rectangle  (3.5,2.5);

\foreach \x in {1,...,5}
\node[left] at (0,0+0.5*\x) [scale=1]
{\small{\x}};

\foreach \x in {1,...,7}
\node[below] at (0+0.5*\x, 0) [scale=1]
{\small{\x}};

\node[right, align=left, text width=10cm] at (4.5,1.25)
{Fig. 2.2. $H=H_0$ in blue is a Hamiltonian e-cycle of a $9 \times 6$ grid graph $G$. The subpath $H_1$ of $H$ in $G_1$ is in orange; it is a Hamiltonian e-cycle of $G_1$. The subpath $H_2$ of $H$ in $G_2$ is highlighted in green; it is a Hamiltonian e-cycle of $G_2$. $H_0$ and $H_1$ have exactly one cookie.};

\end{tikzpicture}
\end{adjustbox}

\null

{%\noindent The $3 \times n$ and $4 \times n$ cases were done by Nishat in \cite{nishat2020reconfiguration}, so from here on, we will assume that $m,n \geq 5$, and that $m$ and $n$ are not both odd. First we will describe canonical forms for Hamiltonian cycles and e-cycles. Then we show that we can reconfigure any two canonical forms into one another. Then we show that any Hamiltonian cycle (e-cycle) can be reconfigured into a canonical form. Observing that double-switch moves are invertible completes the proof. That is, $X \mapsto Y$ followed by $Y \mapsto X$ results in no net change. More specifically, suppose we want to reconfigure a Hamiltonian cycle (e-cycle) $H$ into a Hamiltonian cycle (e-cycle) $K$. Let $\mu_1, \ldots, \mu_k$ and $\nu_1,  \ldots, \nu_s$ be the sequences of moves that reconfigure $H$ and $K$ into the canonical forms $H_{can}$ and $K_{can}$, respectively. Then $\nu_s, \nu_{s-1}, \ldots,\nu_1$ reconfigures $K_{can}$ into $K$. Let $\eta_1, \ldots, \eta_t$ be the sequence of moves that reconfigures $H_{can}$ into $K_{can}$. Then the sequence of moves $\mu_1, \ldots, \mu_k, \eta_1, \ldots, \eta_t, \nu_s, \nu_{s-1}, \ldots,\nu_1 $ reconfigures $H$ into $K$. 
}

\noindent \textbf{Description of \index{canonical form|textbf}canonical forms.} A Hamiltonian cycle $H$ of an $m \times n$  grid graph $G$ belongs to $\mathcal{H}_{\textrm{can}}(m,n)$  if and only if $H$ can be constructed by the ``Canonical Form Builder"  algorithm described below.

Let $t=\Big\lfloor \frac{\min(m,n)-4}{2}  \Big\rfloor$. Let $k_1=|m-n|+2$ and $k_2=|m-n|+3$. If $\min(m,n)$ is even, let $D$ be the Hamiltonian cycle of the $2 \times k_1$ grid graph $G_{t+1}$. If $\min(m,n)$ is odd, let $D$ be any Hamiltonian cycle of the $3 \times k_2$ grid graph $G_{t+1}$. Let $U = D \cup  \bigcup_{i=0}^t R_i$. See Figure 2.3.

\begin{center}

\end{center}

\noindent Now we can state the Canonical Form Builder algorithm \index{CFB Algorithm|textbf} (CFB) that takes as inputs $m$ and $n$ and  outputs an element of $\mathcal{H}_{\textrm{can}}(m,n)$.

\begin{itemize}
    
    \item [Step 1.] Set $i=0$. Switch one of the $2(m-3)+2(n-3)$ switchable boxes of $(R_0,R_1)$ of the graph $U$. This switch removes some edge, say $e_1$, from $E(R_1)$. If $t=0$, stop. If $t>0$, go to Step 2.
    
    \item [Step 2.] Increase $i$ by $1$. Switch one of the switchable boxes of $(R_i,R_{i+1}) \setminus e_i$.

    \item [Step 3.] If $i \leq t$, go to Step 2. If $i=t+1$, then stop.

    We have arrived at a canonical form $H$. Record the switched boxes $X_0, \ldots, X_t$ in a list \textit{List(H)}. So $\textrm{List(H)}=(X_0, \ldots, X_t)$  consists of the faces of $G$ that were chosen to be switched to make $U$ into a canonical form, listed in order.

\end{itemize}

\noindent We observe that the CFB algorithm above works just as well for e-cycles if we remove $e$ from $U$.

\null 

\noindent \textbf{Reconfiguration between canonical forms.} Let $H,K \in \mathcal{H}_{\textrm{can}}(m,n)$. Let $\textrm{List(H)}=(X_0,\ldots,X_t)$ and $\textrm{List(K)}=(Y_0,\ldots, Y_t)$ be the switched boxes of $H$ and $K$ respectively. We will reconfigure $H$ into $K$, so the algorithm will run on $H$. 

Let $D(H)=H \cap G_{t+1}$ and note that $D(H)$ is an e-cycle of $G_{t+1}$. Using the result of Nishat in \cite{nishat2020reconfiguration}, $D(H)$ can be reconfigured into $D(K)$ by a sequence of valid moves. 

\begin{itemize}
    \item [Step 1.] Set $i=0$. If $t=0$, go to Step 3. If $t>0$, go to Step 2.

    \item [Step 2.] If $Y_i$ is switchable after switching $X_i$, switch both $X_i$ and $Y_i$. 
    
    If $Y_i$ is not switchable after switching $X_i$, switch $X_{i+1}$ and any other switchable box in $(R_{i+1}, R_{i+2})$, say $X_{i+1}'$, such that $X_{i+1} \mapsto X_{i+1}'$ is valid. We remark that the only case where $Y_i$ would not be switchable after switching $X_i$ occurs when $Y_i$ is adjacent to $X_{i+1}$. Note that there are many possible choices for $X_{i+1}'$. Now $Y_i$ is switchable. Switch both $X_i$ and $Y_i$. Update $\textrm{List(H)}$ by setting the $(i+1)^{st}$ slot to  $X_{i+1}'$. Increase $i$ by $1$. 

    \item [Step 3.] If $i<t$, go to Step 2. 
    
    If $i=t$ and $\min(m,n)$ is even, switch $X_i$ and $Y_i$, and then stop.

    If $i=t$ and $\min(m,n)$ is odd, go to Step 3.1.

     \item [Step 3.1.] Switch $X_t$ and any one of the four switchable boxes, say $X_t'$, located on the short sides of $D$. Run NRI's algorithm to reconfigure $D(H)$ into $D(K)$. Switch $X_t'$ and $Y_t$. Stop.

\end{itemize}

\subsection{The RtCF algorithm}

\noindent \textbf{Reconfiguration of a cycle into a canonical form \index{RtCF Algorithm|textbf}(RtCF).} Let $H=H_0$ be a Hamiltonian cycle of an $m \times n$ grid graph $G=G_0$. The RtCF algorithm takes as input a Hamiltonian e-cycle and outputs a canonical e-cycle. It works iteratively to reduce the number of cookies over each of the subgraphs $G_0, G_1, ... G_t$. The existence of required moves is guaranteed by the MLC and 1LC algorithms, stated in Proposition 2.2 below, and proved in Chapter 3. Once it reduces the number of cookies of $G=G_0$ to one, it moves on to $G_1$. Now, consider the resulting Hamiltonian cycle $H_0'$. The subpath $H_1=H_0' \cap G_1$ is a Hamiltonian e-cycle of $G_1$, reducing the problem to a smaller grid. This process continues until reaching the central rectangle, at which point $H$ has been reconfigured into canonical form. See Figure 2.4 on page 42 for an illustration of a full execution of the RtCF algorithm.

\null

\noindent \textbf{Proposition 2.2.} Let $H \in \mathcal{H}(m,n)$. \index{small cookie}\index{large cookie}\index{cookie}

  (a) If $H$ has more than one large cookie, then there is a cascade that reduces the 

\hspace*{0.6cm} number of large cookies of $H$ by one. This is the \index{MLC Algorithm|textbf}ManyLargeCookies (MLC) 

\hspace*{0.6cm} algorithm.

  (b) If $H$ has exactly one large cookie and at least one small cookie, then there is a 

\hspace*{0.6cm}  cascade that reduces the number of small cookies of $H$ by one and such that it 

\hspace*{0.6cm}  does not increase the number of large cookies. This is the OneLargeCookie 

\hspace*{0.6cm} \index{1LC Algorithm|textbf}(1LC) algorithm.

\null

\noindent Now we can describe the RtCF algorithm. Without loss of generality, assume $H$ is a Hamiltonian e-cycle of $G$. Suppose $H=H_0$ has $c_{0;\textrm{large}}$ large cookies and $c_{0;\textrm{small}}$ small cookies. We run MLC $c_{0;\textrm{large}}-1$ times and then run 1LC $c_{0;\textrm{small}}$ times to reconfigure $H_0$ into $H_0'$, where $H_0'$ has exactly one (necessarily large) cookie $C_0$. We define  $H_1=(G_1 \cap H_0')$ and observe that $H_1$ is a Hamiltonian e-cycle of $G_1$. This is the first iteration of (RtCF). Now we describe the $j^{\textrm{th}}$ iteration. We run MLC $c_{j;\textrm{large}}-1$ times and then run 1LC $c_{j;\textrm{small}}$ times to reconfigure $H_{j}$ into $H_{j}'$, where $H_{j}'$ has exactly one (necessarily large) cookie $C_j$. The RtCF algorithm stops when $j=t$, having completed $t+1$ iterations. We give a summary of the algorithm below.

\begin{itemize}
    \item [Step 1.] Set $j=0$. Run MLC $c_{0;\textrm{large}}-1$ times and then 1LC $c_{0;\textrm{small}}$ times on $H_0$ to reconfigure $H_0$ into $H_0'$.
    
    \item [Step 2.] Increase $j$ by $1$. Set $H_j=G_j \cap H'_{j-1}$ and note that $H_j$ is a Hamiltonian e-cycle in $G_j$. Run MLC $c_{j;\textrm{large}}-1$ times and then 1LC $c_{j;\textrm{small}}$ times on $H_{j}$ to reconfigure $H_{j}$ into $H_{j}'$. 

    \item [Step 3.] If $j<t$, go to Step 2. If $j=t$, stop. 

\end{itemize}

\noindent \textit{Proof of the RtCF algorithm}. Let $N_j$ be the neck of the only cookie $C_j$ of $H_j'$ in $G_j$. Define  $e_1(N_j)=N_j \cap R_j$, $e_2(N_j)=N_j \cap R_{j+1}$ and $\{ e_3(N_j),e_4(N_j)\} = N_j \cap H_j'$. We observe that when the RtCF algorithm stops, we have reconfigured $H$ into

$$H_c=D(H) \cup \bigcup_{j=0}^t \big(R_j \cup e_3(N_{j}) \cup e_4(N_{j})\big) \setminus \big(e_1(N_{j}) \cup e_2(N_{j})\big).$$

\noindent Now we can see that $H_c$ is an element of  $\mathcal{H}_{\textrm{can}}(m,n)$  by setting $X_{j}=N_j$ for $j=0,1,\ldots,t$ and running CFB. 

\null

\noindent \textbf{Bound for Theorem 2.1.} Recall that $n \geq m$. Note that it takes at most $2m$ \footnote{If $\min(m,n)$ is odd, add another $(n-m)/2$ moves to reconfigure $D$.} moves to reconfigure canonical forms into one another. Now we count the moves required for RtCF to terminate. Observe that for each $j \in \{0, \ldots, t-1\}$, $H_j$ has at most $2 \big(\frac{n-2j}{2}+ \frac{m-2j}{2}\big)=n+m-4j$ cookies. This is the number of iterations of MLC or 1LC required for each $j$. It will follow from the proofs in Sections 3 and 4 that each application of MLC or 1LC in $H_j$ requires at most $\frac{1}{2}n+m-3j+3$ moves. So, RtCF requires at most:

\begin{align*}
S&=\sum_{j=1}^{\lfloor \frac{1}{2}(m-2) \rfloor}\Big(n+m-4j\Big) \Big(\frac{n}{2}+m-3j+2\Big) \\
&= \sum_{j=1}^ {\lfloor \frac{1}{2}(m-2) \rfloor} 12j^2 + (-5n - 7m - 8)j + \frac{n^2}{2} + \frac{3nm}{2} + 2n + m^2 + 2m \\
&\leq  \frac{m^3}{2} - \frac{3m^2}{2} + m + (-5n - 7m - 8)\left( \frac{m^2}{8} - \frac{m}{4} \right)  + \frac{ \left( \frac{n^2}{2} + \frac{3nm}{2} + 2n + m^2 + 2m \right)(m - 2) }{2}  \\
&= \frac{n^2m}{4}+\frac{nm^2}{8}+\frac{3nm}{4}+\frac{m^3}{8}-\frac{3m^2}{4}-\frac{n^2}{2}+m-2n\\
&= \frac{n^2m}{2} + \frac{nm}{4} \big( 3-\frac{3m}{n}-\frac{2n}{m} \big)+m-2n. 
\end{align*}

\null

\noindent Let $x=\frac{m}{n}$. Then $\frac{3m}{n}+\frac{2n}{m}=3x+\frac{2}{x}$. Using calculus, we find that it attains a minimum of $2\sqrt{6}$ at $x=\frac{\sqrt{6}}{3}$. Then $\big( 3-\frac{3m}{n}-\frac{2n}{m} \big)$ can be at most $3-2\sqrt{6} \leq -1$. It follows that RtCF requires at most $\frac{n^2m}{2} -\frac{nm}{4}+m-2n$ to terminate. For a complete reconfiguration we need to run RtCF once for each e-cycle and reconfigure the resulting canonical forms. So, we need at most $2 \big( \frac{n^2m}{2}-\frac{nm}{4}+m-2n \big) +m=n^2m-\frac{nm}{2}-4n+3m < n^2m$ moves for a complete reconfiguration. We remark that this is a worst case scenario and conjecture that the typical number of moves required is of the order of $n^2$.

\null 

\noindent \textbf{Note 2.3.} An $m \times n$ rectangular grid graph does not admit a Hamiltonian cycle when both $m$ and $n$ are odd. To see this, observe that the grid graph is bipartite: color a vertex $(i,j)$ red if $i+j$ is even and blue if $i+j$ is odd. Any cycle in a bipartite graph must alternate between the two color classes, and therefore has even length. However, when both $m$ and $n$ are odd, the grid has $mn$ vertices, which is an odd number. Thus, no Hamiltonian cycle can exist.

\null 

\subsection{Summary}

This chapter defined canonical forms and the RtCF algorithm for transforming Hamiltonian cycles into canonical form. We observed that reconfiguration between canonical forms is straightforward, reducing the general reconfiguration problem to reaching a canonical form. The main result is Theorem 2.1. Proving that RtCF can always find its required moves is handled by the MLC and 1LC algorithms, in Chapter 3.

\null

\begin{center}

\end{center}

%We give the proofs of existence for the MLC and 1LC algorithms in the next chapter.

\section{Existence of the MLC and 1LC Algorithms}

Recall that RtCF algorithm (Chapter 2) presupposes the existence of the moves required for its execution. The proofs of existence were deferred to the ManyLargeCookies (MLC) and OneLargeCookie (1LC) algorithms, which we prove in Sections 3.1 and 3.2, respectively. These algorithms ensure that whenever RtCF requires a particular move, either the move is immediately available, or else there exists a cascade after which the required move becomes available. Importantly, such cascades do not undo the progress already made: RtCF does not regress. The restrictions in the definition for cascades at the end of Section 1.4 were designed precisely for this reason. 

Consider an iteration of RtCF on rectangle $G_i$ with Hamiltonian e-cycle $H_i$. At this stage, $H_i$ either has more than one large cookie, exactly one large cookie with at least one small cookie, or exactly one large cookie with no small cookies. In the last case, $H_i$ is already in the desired form for this iteration, and RtCF proceeds to $G_{i+1}$. The \index{MLC Algorithm} MLC algorithm handles the first case by finding the required cascades to collect large cookies when multiple large cookies are present. The \index{1LC Algorithm}1LC algorithm handles the second case by finding the required cascades to collect small cookies when exactly one large cookie remains.

Why do we need two separate algorithms for what appears to be the same task? This is because small cookies can be harder to collect than large ones. A second large cookie $J$ always has a switchable neck $N_J$; to collect $J$ we need only find another switchable box $X$ such that $N_J \mapsto X$ is a valid move, or a cascade delivering such a switchable box. In Section 3.1, we show that it takes at most two moves to accomplish this (Proposition 3.8). Small cookies, by contrast, consist of a single non-switchable box. To collect a small cookie, either the box $Y$ adjacent to it in $(R_i, R_{i+1})$ must be switchable, or we must find a cascade that makes $Y$ switchable. The latter task can require much longer cascades, and it is more difficult to deal with. It requires Lemma 3.7, all of Section 3.2, most of Chapter 4, and several results from Chapter 1. Furthermore, the assumption that exactly one large cookie is present significantly shortens and simplifies the proofs of Proposition 3.10 and Lemmas 3.11--3.15 in Section 3.2, by precluding the possibility of several tedious cases.

\null 

\noindent \textbf{Notation.} In this chapter, we prove results for a Hamiltonian cycle $H$ of an $m \times n$ grid graph $G$. These results apply verbatim to $H_i$ on $G_i$ at each iteration of RtCF, so we omit subscripts to simplify notation.

\subsection{Existence of the MLC Algorithm}

Let $G$ be an $m \times n$ grid graph, let $H$ be a Hamiltonian path or cycle of $G$, and let $W$ be a switchable box in $H$. Let $X$ and $Y$ be the boxes adjacent to $W$ that are not its $H$-neighbours, and assume that $X$ and $Y$ belong to the same $H$-component. By Corollary 1.2.2, the $H$-path $P(X,Y)$ is unique. We call $P(X,Y)$ \index{the looping H-path of a box@the looping $H$-path of a box|textbf}\textit{the looping $H$-path of $W$}. See Figure 3.1 for an illustration of the looping $H$-path of a switchable box $W$ in a Hamiltonian cycle of a $6 \times 6$ grid graph.

\begingroup
\setlength{\intextsep}{0pt}
\setlength{\columnsep}{20pt}
\begin{wrapfigure}[]{r}{0cm}
\begin{adjustbox}{trim=0cm 0cm 0 0cm}
\begin{tikzpicture}[scale=1.5]
\begin{scope}[xshift=0cm, yshift=0cm]

\draw[gray,very thin, step=0.5cm, opacity=0.5] (0,0) grid (2.5,2.5);

\fill[orange!50!white, opacity=0.5] (1,0)--++(0.5,0)--++(0,1.5)--++(-0.5,0);
\fill[orange!50!white, opacity=0.5] (2,0)--++(0.5,0)--++(0,1.5)--++(-0.5,0);
\fill[orange!50!white, opacity=0.5] (1.5,0)--++(0.5,0)--++(0,0.5)--++(-0.5,0);

\draw[blue, line width=0.5mm] (0,0)--++(2.5,0)--++(0,2.5)--++(-1,0)--++(0,-0.5)--++(0.5,0)--++(0,-1.5)--++(-0.5,0)--++(0,1)--++(-0.5,0)--++(0,-1)--++(-0.5,0)--++(0,1.5)--++(0.5,0)--++(0,0.5)--++(-1,0)--++(0,-2.5);

\node at (1.75,1.25) [scale=0.8] {W};
\node at (1.25,1.25) [scale=0.8] {X};
\node at (2.25,1.25) [scale=0.8] {Y};

\node[below, align=center, text width=5cm] at (1.25,0)
{Fig. 3.1. The looping  $H$-path of $W$ shaded orange. };

\end{scope}
\end{tikzpicture}
\end{adjustbox}
\end{wrapfigure}

We remark that all boxes $Z$ for which $W \mapsto Z$ is a valid move (that is, in the terminology of Section 1.4, all $(H_1, H_2)$-ports that arise once $W$ is switched) lie within the looping H-path $P(X,Y)$.

\null

\noindent \textbf{Outline of the MLC algorithm.} Let $H$ be a Hamiltonian cycle of an $m \times n$ grid graph $G$ with multiple large cookies. We first identify a large cookie $J$ with switchable neck $N_J$. Consider what happens if $N_J$ is switched: this would produce two cycles, $H_1$ and $H_2$. First we observe that there must be some edge $\{v_1,v_2\}$ in $R_2$ (recall the nested rectangles from Chapter 2) with $v_1 \in H_1$ and $v_2 \in H_2$ (Lemma 3.7). The proximity of $\{v_1,v_2\}$ to the boundary constrains the possible configurations of edges in its vicinity. We analyze those configurations (Lemma 3.5) and show that either an $(H_1,H_2)$-port already exists near $\{v_1, v_2\}$, or a single-move cascade on the original $H$ yields a Hamiltonian cycle $H'$ where such a port exists after switching $N_J$.

\endgroup

\null 

\noindent \textbf{Proposition 3.1.} Let $G$ be an $m \times n$ grid graph, let $H$ be a Hamiltonian cycle of $G$, and let $P(X,Y)$ be the looping $H$-path of a switchable box $W$. Let $H'$ be the graph consisting of the cycles $H_1$ and $H_2$ obtained after switching $W$. Then a box $Z$ of $G$ belongs to the $H'$-cycle $P(X,Y), W,X$  if and only if $Z$ is incident on a vertex of $H_1$ and on a vertex of $H_2$.

\null 

\noindent \textit{Proof.} Orient $H$ as a directed cycle $H=v_1, \ldots, v_r,v_1$. By Corollary 1.4.3, $W$ is anti-parallel. Let the edges of $W$ in $H$ be  $\{v_x, v_{x+1}\}$ and  $\{v_{y-1}, v_y\}$. For definiteness, assume that $X$ is adjacent to $\{v_x,v_{x+1}\}$, $Y$ is adjacent to $\{v_{y-1},v_y\}$ and that $W$ is on the right of $\{v_x,v_{x+1}\}$. Then we have that $\Phi((v_x,v_{x+1}),\text{left})=X$ and $\Phi((v_{y-1},v_y),\text{left})=Y$. Define $\overrightarrow{K}_1$ and $\overrightarrow{K}_2$ to be the subtrails $\overrightarrow{K}((v_x,v_{x+1})$, $(v_{y-1}, v_y))$ and $\overrightarrow{K}((v_{y-1},v_y), (v_x,v_{x+1}))$ of $\overrightarrow{K}_H$, respectively. By Lemma 1.4.1 (ii), switching $W$ splits $H$ into two cycles $H_1$ and $H_2$, with  $V(H_1)=V(\overrightarrow{K}_1) \setminus \{v_x,v_y\}$ and $V(H_2)=V(\overrightarrow{K}_2) \setminus \{v_{x+1},v_{y-1}\}$. See Figure 3.2.

\null

\noindent \textit{Proof of $(\implies)$.} Since  $P(X,Y)$ is unique, any $H$-walk of boxes between $X$ and $Y$ contains $P(X,Y)$. In particular, $\Phi(\overrightarrow{K}_1,\text{left})$ contains $P(X,Y)$ and $\Phi(\overrightarrow{K}_2,\text{left})$ contains $P(X,Y)$. Let $Z$ be a box of $P(X,Y)$. Then $Z$ is added to $\Phi(\overrightarrow{K}_1,\text{left})$ by an edge of $\overrightarrow{K}_1$ and $Z$ is added to $\Phi(\overrightarrow{K}_2,\text{left})$ by an edge of $\overrightarrow{K}_2$. By definition of FTW, $Z$ is incident on a vertex of $\overrightarrow{K}_1$ and a vertex of $\overrightarrow{K}_2$. Since for $i \in \{1,2\}$, $V(\overrightarrow{K}_i) \supset V(H_i)$, we have that $Z$ is incident on a vertex of $H_1$ and a vertex of $H_2$. Lastly, note that $W$ is incident on $v_{x+1} \in H_1$ and $v_x \in H_2$. End of proof for $(\implies)$.

\endgroup 

\null 

\begingroup
\setlength{\intextsep}{0pt}
\setlength{\columnsep}{20pt}
\begin{wrapfigure}[]{l}{0cm}
\begin{adjustbox}{trim=0cm 0cm 0cm 0cm}

\end{adjustbox}
\end{wrapfigure}

\noindent \textit{Proof of $(\impliedby)$.} Suppose we switch $W$ and obtain the graph $H'$ consisting of the cycles $H_1$ and $H_2$. Observe that $P(X,Y), W, X$ is the only $H'$-cycle in $G$. We will say that a box $Z$ of $G$ satisfies $(*)$ if $Z$ is incident on a vertex in $H_1$ and $H_2$. We will show that if a box $Z$ of $G$ satisfies $(*)$ then it must belong to an $H'$-cycle of boxes that satisfy $(*)$. Then, since there is only one $H'$-cycle in $G$, $Z=W$ or $Z \in P(X,Y)$.

Let $Z$ be a box in $G$ that satisfies $(*)$. For definiteness, assume that $Z=R(k,l)$. We will show that $Z$ has exactly two neighbours in $G$ that satisfy $(*)$ and that $Z$ is $H'$-adjacent to those two neighbours. Since $Z$ satisfies $(*)$, either $Z$ has two vertices in $H_1$ and two vertices in $H_2$ or $Z$ has one vertex in one of $H_1$ and $H_2$ and three vertices in the other.

\null 

\begingroup 
\setlength{\intextsep}{0pt}
\setlength{\columnsep}{20pt}
\begin{wrapfigure}[]{r}{0cm}
\begin{adjustbox}{trim=0cm 0cm 0cm 0cm}
\begin{tikzpicture}[scale=1.75]

\begin{scope}[xshift=0cm]{
\draw[gray,very thin, step=0.5cm, opacity=0.5] (0,0) grid (1.5,1);

\draw [blue, line width=0.5mm] plot [smooth, tension=0.75] coordinates {(0.5,0)(0.25,0.5)(0.5,0.75)(1,0.5)};

\draw [orange, line width=0.5mm] plot [smooth, tension=0.75] coordinates {(0.5,0.5)(1,0.75)(1.25,0.5)(1,0)};

\draw[blue, line width=0.5mm] (0.5,0)--++(0.5,0.5);

\draw[fill=orange] (0.5,0.5) circle [radius=0.05];
\draw[fill=orange] (1,0) circle [radius=0.05];
\draw[fill=blue] (0.5,0) circle [radius=0.05];
\draw[fill=blue] (1,0.5) circle [radius=0.05];

% nodes
{

\node[left] at (0,0.5) [scale=1]{\tiny{+1}};
\node[left] at (0,0) [scale=1]
{\tiny{$\ell$}};
%\node[left] at (0,1) [scale=1]{\tiny{1}};

\node[above] at (0.5,1) [scale=1]
{\tiny{$k$}};
\node[above] at (1, 1) [scale=1]{\tiny{+1}};

\node at (0.85,0.15) [scale=0.8] {\small{$Z$}};
}

\node[below, align=center, text width=2cm] at (0.75,-0.15)
{Fig. 3.3.};

} \end{scope}

\end{tikzpicture}
\end{adjustbox}
\end{wrapfigure}

\noindent \textit{CASE 1: $Z$ has two vertices in $H_1$ and two vertices in $H_2$.} First we will check that the pair of vertices belonging to $H_i$, $i\in \{1,2\}$ must be adjacent in $Z$. For contradiction, assume that $v(k,l)$ and $v(k+1,l+1)$ belong to $H_1$ and $v(k+1,l)$ and $v(k,l+1)$ belong to $H_2$. Let $Q$ be the closed polygon consisting of the subpath $P(v(k,l), v(k+1,l+1))$ of $H_1$ and the segment $[v(k,l), v(k+1,l+1)]$. Then, by Corollary 1.1.5, $v(k+1,l)$ and $v(k,l+1)$ are on different sides of $H_1$. It follows that the subpath $P(v(k+1,l),v(k,l+1))$ intersects $Q$. Since $P(v(k+1,l),v(k,l+1))$ does not intersect $Q$ at the segment $[v(k,l), v(k+1,l+1)]$, it must intersect $Q$ at some vertex in $P(v(k,l), v(k+1,l+1))$. But this contradicts that $H_1$ and $H_2$ are disjoint. It follows that the pair of vertices belonging to $H_i$, $i\in \{1,2\}$ must be adjacent in $Z$.

For definiteness, assume that $v(k,l)$ and $v(k+1,l)$ belong to $H_1$ and that $v(k+1,l+1)$ and $v(k,l+1)$ belong to $H_2$. Since $H_1$ and $H_2$ are disjoint $e(k;l,l+1) \notin H'$ and $e(k+1;l,l+1) \notin H'$ so $Z+(-1,0)$ and $Z+(1,0)$ are $H'$-adjacent to $Z$. Since $v(k,l) \in H_1 \cap V(Z+(-1,0))$ and $v(k,l+1) \in H_2 \cap V(Z+(-1,0))$, $Z+(-1,0)$ satisfies $(*)$. Similarly, $Z+(1,0)$ satisfies $(*)$. It remains to check that $Z+(0,1)$ and  $Z+(0,-1)$ do not satisfy $(*)$.

\begingroup 
\setlength{\intextsep}{0pt}
\setlength{\columnsep}{20pt}
\begin{wrapfigure}[]{l}{0cm}
\begin{adjustbox}{trim=0cm 0cm 0cm 0cm}
\begin{tikzpicture}[scale=1.75]

\begin{scope}[xshift=0cm]{
\draw[gray,very thin, step=0.5cm, opacity=0.5] (0,0) grid (1.5,1);

\fill[green!50!white, opacity=0.5] (0,0) rectangle (1.5,0.5);

\draw[orange, line width=0.5mm] (0,0.5)--++(1,0);

\draw[fill=blue] (1,0) circle [radius=0.05];
\draw[fill=blue] (0.5,0) circle [radius=0.05];
\draw[fill=blue] (0.5,1) circle [radius=0.05];

\draw[fill=orange] (0,0.5) circle [radius=0.05];
\draw[fill=orange] (0.5,0.5) circle [radius=0.05];
\draw[fill=orange] (1,0.5) circle [radius=0.05];

%black lines
{

\draw[black, line width=0.15mm] (0.45,0.2)--++(0.1,0);
\draw[black, line width=0.15mm] (0.45,0.25)--++(0.1,0);
\draw[black, line width=0.15mm] (0.45,0.3)--++(0.1,0);

\draw[black, line width=0.15mm] (0.95,0.2)--++(0.1,0);
\draw[black, line width=0.15mm] (0.95,0.25)--++(0.1,0);
\draw[black, line width=0.15mm] (0.95,0.3)--++(0.1,0);

\draw[black, line width=0.15mm] (0.45,0.7)--++(0.1,0);
\draw[black, line width=0.15mm] (0.45,0.75)--++(0.1,0);
\draw[black, line width=0.15mm] (0.45,0.8)--++(0.1,0);

}

% nodes
{

\node[left] at (0,0.5) [scale=1]{\tiny{+1}};
\node[left] at (0,0) [scale=1]
{\tiny{$\ell$}};
%\node[left] at (0,1) [scale=1]{\tiny{1}};

\node[above] at (0,1) [scale=1]{\tiny{-1}};
\node[above] at (0.5,1) [scale=1]
{\tiny{$k$}};
\node[above] at (1, 1) [scale=1]{\tiny{+1}};

\node at (0.75,0.25) [scale=0.8] {\small{$Z$}};
}

\node[below, align=center, text width=2cm] at (0.75,-0.15)
{Fig. 3.4. };

} \end{scope}

\end{tikzpicture}
\end{adjustbox}
\end{wrapfigure}

For contradiction, assume that one of $Z+(0,1)$ and $Z+(0,-1)$  satisfies $(*)$. By symmetry we may assume WLOG that $Z+(0,1)$ satisfies $(*)$. Then at least one of $v(k,l+2)$ and $v(k+1,l+2)$ belongs to $H_1$. By symmetry we may assume WLOG that $v(k,l+2) \in H_1$. It follows that $e(k;l+1,l+2) \notin H'$ and that $v(k-1,l+1) \in H_2$. Then we must have $e(k-1,k;l+1) \in H_2$ and $e(k,k+1;l+1) \in H_2$. But then by Corollary 1.1.5, one of $v(k,l)$ and $v(k,l+2)$ lies inside the region bounded by $H_2$ and the other lies outside it. It follows that the subpath $P(v(k,l), v(k,l+2))$ of $H_1$ intersects $H_2$, contradicting that $H_1$ and $H_2$ are disjoint. End of Case 1.

\endgroup 

\null 

\noindent \textit{CASE 2: $Z$ has one vertex in one of $H_1$ and $H_2$ and three vertices in the other.} For definiteness, assume that $v(k,l)$, $v(k,l+1)$ and $v(k+1,l+1)$ belong to $H_1$ and that $v(k+1,l)$ belongs to $H_2$. Then $e(k,k+1;l) \notin H$ and $e(k+1;l,l+1) \notin H$, so $Z+(1,0)$ and $Z+(0,-1)$ are $H'$-neighbours of $Z$. Since $v(k,l) \in H_1 \cap V(Z+(0,-1))$ and $v(k+1,l) \in H_2 \cap V(Z+(0,-1))$, $Z+(0,-1)$ satisfies $(*)$. Similarly, $Z+(1,0)$ satisfies $(*)$. It remains to check that $Z+(0,1)$ and $Z+(0,-1)$ do not satisfy $(*)$. 

\null 

\begingroup 
\setlength{\intextsep}{0pt}
\setlength{\columnsep}{20pt}
\begin{wrapfigure}[]{r}{0cm}
\begin{adjustbox}{trim=0cm 0cm 0cm 0cm}

\end{adjustbox}
\end{wrapfigure}

\noindent For a contradiction, assume that one of $Z+(-1,0)$ and $Z+(0,1)$  satisfies $(*)$. By symmetry we may assume WLOG that $Z+(0,1)$ satisfies $(*)$. Then one of $v(k,l+2)$ and $v(k+1,l+2)$ belongs to $H_2$. Note that if $v(k+1,l+2) \in H_2$ we run into the same contradiction as in Case 1, so we only need to check the case where $v(k,l+2) \in H_2$. Now, either $e(k,k+1;l+1) \in H_1$, or $e(k,k+1;l+1) \notin H'$. 

\null 

\noindent \textit{CASE 2.1: $e(k,k+1;l+1) \in H_1$.} Then the segment $[v(k,l+2), v(k+1,l)]$ intersects $H_1$ at the point $(k \frac{1}{2}, l+1)$. by Corollary 1.1.5, the vertices $e(k,l+2)$ and $v(k+1,l)$ are on different sides of $H_1$, and we run into the same contradiction as in Case 1 again. See Figure 3.5 (a).

\null 

\noindent \textit{CASE 2.2: $e(k,k+1;l+1) \notin H_1$.} Consider the polygon $Q$ consisting of the segment $[v(k,l+2), v(k+1,l)]$ and the subpath $P(v(k,l+2), v(k+1,l))$ of $H_2$. By Corollary 1.1.5, the vertices $v(k,l+1)$ and $v(k+1,l+1)$ are on different sides of $Q$. By JCT the subpath $P(v(k,l+1),v(k+1,l+1))$ of $H_1$ intersects $Q$. Since $P(v(k,l+1),v(k+1,l+1))$ does not intersect $Q$ at the segment $[v(k,l+2), v(k+1,l)]$, it must do so at some vertex of $P(v(k,l+2), v(k+1,l))$, contradicting that $H_1$ and $H_2$ are disjoint. See Figure 3.5 (b). End of Case 2.2. End of Case 2. $\square$

\null

\noindent \textbf{Corollary 3.2 (i).} Let $G$ be an $m \times n$ grid graph, let $H$ be a Hamiltonian cycle of $G$, and let $W$ be a switchable box in $H$. Let $H'$ be the graph consisting of the cycles $H_1$ and $H_2$ obtained after switching $W$. Let $a,b$ and $c$ be colinear vertices such that $b$ is adjacent to $a$ and $c$. Then, for $i \in \{1,2\}$, If $a$ and $c$ belong to $H_i$, so must $b$. See Figure 3.4. $\square$

\null

\noindent \textbf{Corollary 3.2 (ii).} Let $G$ be an $m \times n$ grid graph, let $H$ be a Hamiltonian cycle of $G$, and let $W$ be a switchable box in $H$. Let $H'$ be the graph consisting of the cycles $H_1$ and $H_2$ obtained after switching $W$. Let $Z$ be a box on vertices $a,b,c$, and $d$ such that $a$ and $b$ belong to $H_1$, and $c$ and $d$ belong to $H_2$. Then $a$ is adjacent to $b$, and $c$ is adjacent to $d$. See Figure 3.3. $\square$

\null 

\noindent \textbf{Proposition 3.3.} Let $H$ be a Hamiltonian cycle of an $m \times n$ grid graph $G$, let $W$ be a switchable box in $H$ and let $P(X,Y)$ be the looping $H$-path of $W$. If $P(X,Y)$ has a switchable box $Z$, then $Z \mapsto W$ is a valid move.

\begingroup
\setlength{\intextsep}{0pt}
\setlength{\columnsep}{20pt}
\begin{wrapfigure}[]{l}{0cm}
\begin{adjustbox}{trim=0cm 0cm 0cm 0cm}

\end{adjustbox}
\end{wrapfigure}

\noindent \textit{Proof.} Let $H=v_1, v_2 \ldots, v_r, v_1$. By Corollary 1.4.3, $W$ is anti-parallel. Let $W$, $P(X,Y)$, $\{v_x, v_{x+1}\}$ and $\{v_{y-1}, v_y\}$ be as in Proposition 3.1. Let $P_1=P(v_x,v_y)$ and let $P_2=P(v_y,v_x)$. By Proposition 3.1, every box of $P(X,Y)$ is incident on a vertex of $P_1$ and a vertex of $P_2$.

Let $Z$ be a switchable box of $P(X,Y)$. Let $(v_s,v_{s+1})$ and $(v_t, v_{t+1})$ be the edges of $Z$ in $H$. For definiteness, assume $s+1<t$. Proposition 3.1 implies that exactly one of $(v_s,v_{s+1})$ and $(v_t, v_{t+1})$ is in $P_1$ and the other is in $P_2$. WLOG assume that $(v_s,v_{s+1})$ is in $P_1$ and that $(v_t, v_{t+1})$ is in $P_2$. Then we can partition the edges of $H$ as follows: $P(v_1,v_s)$, $(v_s,v_{s+1})$, $P(v_{s+1},v_t)$, $(v_t,v_{t+1})$, $P(v_{t+1},v_r)$, $\{v_r,v_1\}$ where $1<x<s<y<t<r$. 

We check that $Z\mapsto W$ is a valid move. After removing the edges $(v_s,v_{s+1})$ and $(v_t,v_{t+1})$ we are left with two paths: $P(v_{t+1},v_s)$ and $ P(v_{s+1},v_t)$. Note that adding the edge $\{v_s, v_{t+1}\}$ gives a cycle $H_1$ consisting of the path $P(v_s,v_{t+1})$ and the edge $\{v_s,v_{t+1}\}$, and adding the edge $\{v_{s+1}, v_t\}$ gives a cycle $H_2$ consisting of the path $P(v_{s+1},v_t)$ and the edge $\{v_{s+1}, v_t\}$. Now $1<x<s$ implies that $(v_x,v_{x+1}) \in H_1$ and $s<y<t$ implies that $(v_{y-1}, v_y) \in H_2$. It follows that $W$ is now an $(H_1,H_2)$-port \index{H p H c port@$(H_p,H_c)$-port}. By Lemma 1.4.1 (iii), $Z \mapsto W$ is a valid move. $\square$ 

\endgroup

\null

\noindent \textbf{Remark.} Proposition 3.3 is the primary tool for collecting cookies and characterizes when double-switch moves are valid. To collect a large cookie with switchable neck $W$, we search for a switchable box $Z$ in the looping $H$-path of $W$. If such a $Z$ exists, we can collect the cookie immediately. If not, this imposes strong structural constraints on the $H$-path. This dichotomy—either we find a switchable box and succeed, or we find none and exploit the resulting structure (and eventually find one anyway)—drives the case analyses of Lemma 3.5 below as well as much of Section 3.2 and Chapters 4 and 5. Beyond collecting large cookies, the proposition is essential for collecting small cookies (which requires finding a switchable box adjacent to the small cookie) and for characterizing valid double-switch moves more generally. In Chapter 4 we prove analogous results that are used in the Hamiltonian path reconfiguration algorithms of Chapter 5. Proposition 4.5 extends Proposition 3.3 and Proposition 4.4 extends Proposition 3.1.

\null 

\noindent Let $H$ be a Hamiltonian cycle of an $m \times n$ grid graph $G$. Recall the definitions of the subgraphs $G_i$ and $R_i$ of $G$ on page 37 in Chapter 2.

\null

\noindent \textbf{Observation 3.4.}  Let $X \mapsto Y$ be a valid move. If $X \in \text{ext}(H)$ and $Y \in \text{int}(H)$ then:

(i) $X \mapsto Y$ increases the total number of cookies if and only if $X \in G_1$ and 

\hspace{0.6cm} \ \ $Y\in G_0 \setminus G_1$.

(ii) $X \mapsto Y$ increases the total number of large cookies, leaving the total number of 

\hspace{0.6cm}  cookies unchanged, if and only if $X \in G_1$, $Y \in G_1 \setminus G_2$ and $Y$ is adjacent 

\hspace{0.6cm}  to a small cookie. 

(iii) $X \mapsto Y$ decreases the total number of cookies if and only if $X \in G_0 \setminus G_1$, $Y\in G_1.$ 

%\textcolor{red}{[There is a case where the move should be $Y \mapsto X$, if X is a small cookie. ]}. 

\null

\noindent \textbf{Lemma 3.5.} Let $H$ be a Hamiltonian cycle of an $m \times n$ grid graph $G$ and let $Z$ be a switchable box in $\text{ext}(H) \cap ((G_0 \setminus G_1) \cup G_3)$. Assume that switching $Z$ splits $H$ into the cycles $H_1$ and $H_2$ that are such that there is $v_1 \in H_1 \cap R_2$ and $v_2 \in H_2 \cap R_2$ with $v_1$ adjacent to $v_2$. Then either $Z \mapsto Z'$ is a cascade, or there is a cascade $\mu, Z \mapsto Z'$ (of length two), with $Z \mapsto Z'$ nontrivial in either case.

\null

\noindent \textbf{Remark.} For the MLC algorithm (Proposition 3.8) in this section, we only need the case where $Z \in \text{ext}(H) \cap (G_0 \setminus G_1)$, since cookie necks reside in $G_0 \setminus G_1$. However, the more general statement (including $Z \in G_3$) will be needed for the 1LC algorithm in Section 3.2, specifically in the proof of Lemma 3.12. We prove the general result here to avoid duplication. 

In the figures illustrating the proof, vertices and edges of $H_1$ and $H_2$ are in blue and orange, respectively, and boxes of $\text{int}(H)$ are shaded in green.

\begingroup
\setlength{\intextsep}{0pt}
\setlength{\columnsep}{10pt}
\begin{wrapfigure}[]{l}{0cm}
\begin{adjustbox}{trim=0cm 0cm 0cm 0cm}

\end{adjustbox}
\end{wrapfigure}

\noindent \textit{Proof.} We will use the assumption that $Z \in  (G_0 \setminus G_1) \cup G_3$ repeatedly and implicitly throughout the proof. Switch $Z$ to obtain $H'$ consisting of the disjoint cycles $H_1$ and $H_2$. Note that now, if a vertex belongs to $H_i$ for $i\in \{1,2\}$, both edges incident on it must also belong to $H_i$.

For definiteness, let $v_1 \in H_1 \cap R_2$ be the vertex $v(2,l)$ for some $l \in \{ 2, \ldots, n-3\}$ and let $v_2=v(2,l-1) \in H_2 \cap R_2$. Then $e(2;l-1,l) \notin H$, and by Corollary 3.2 (i), $v(2,l+1) \in H_1$ as well. By Proposition 3.1, $R(1,l-1)$ and $R(2,l-1)$ belong to $\text{int}(H)$. Now, either $v(1,l) \in H_1$, or $v(1,l) \in H_2$. See Figure 3.7. 

\null 

\noindent \textit{CASE 1: $v(1,l) \in H_2$.} Then $e(1,2;l) \notin H$.  Corollary 3.2(i), $v(3,l) \in H_1$. It follows that $e(2;l,l+1) \in H'$ and $e(2,3;l) \in H'$. Now, by Corollary 3.2 (ii), $v(1,l-1) \in H_2$ and by Corollary 3.2 (i), $v(0,l) \in H_2$. At this point we must either have, $e(1;l,l+1) \in H'$ or $e(1;l,l+1) \notin H'$. See Figure 3.7.

\endgroup 

\null

\begingroup
\setlength{\intextsep}{0pt}
\setlength{\columnsep}{20pt}
\begin{wrapfigure}[]{r}{0cm}
\begin{adjustbox}{trim=0cm 0cm 0cm 0cm}

\end{adjustbox}
\end{wrapfigure}

\noindent \textit{CASE 1.1: $e(1;l,l+1) \notin H'$.} Then $e(1;l-1,l) \in H_2$ and $e(0,1;l) \in H_2$. Since  $Z \in \text{ext}(H)$, by Proposition 3.1, $R(1,l-1) \in \text{int}(H)$. Then $R(0,l-1)$ must be a small cookie of $G$, so we must have that $e(0,1;l-1) \in H_2$ and that $e(1,2;l-1) \notin H'$. It follows that $e(2,3;l-1) \in H'$ and $e(3;l-1,l) \notin H'$. Now note that $R(2,l-1)$ is an $(H_1,H_2)$-port. Then, by Lemma 1.4.1 (iii) and Observation 3.4, $Z \mapsto R(2,l-1)$ is valid move that does not create new cookies. So, $Z \mapsto R(2,l-1)$ is the cascade we seek. End of Case 1.1.

\null 

\noindent \textit{CASE 1.2: $e(1;l,l+1) \in H'$.} Proposition 3.1, and the assumption that $Z \in \text{ext}(H)$ imply that $R(1,l) \in \text{int}(H)$. Then, Corollary 1.3.15 (b) implies that $R(0,l)$ is a small cookie of $G$. Note that if $e(2,3;l-1) \in H'$, then we're back to Case 1.1, so we may assume that $e(2,3;l-1) \notin H'$. It follows that $e(1,2; l-1) \in H_2$ and that $R(0,l-1)$ is not a small cookie of $H$. Then, by Observation 3.4 and Proposition 3.3, $R(0,l) \mapsto R(1,l)$, $Z \mapsto R(1,l-1)$ is the cascade we seek. End of Case 1.2. End of Case 1.

\endgroup 

\null

\begingroup
\setlength{\intextsep}{0pt}
\setlength{\columnsep}{10pt}
\begin{wrapfigure}[]{l}{0cm}
\begin{adjustbox}{trim=0cm 0cm 0cm 0cm}

\end{adjustbox}
\end{wrapfigure}

\noindent \textit{CASE 2: $v(1,l) \in H_1$.} Then $e(2,3;l) \in H_1$ or $e(2,3;l) \notin H'$. See Figure 3.9.

\null

\noindent \textit{CASE 2.1: $e(2,3;l) \in H_1$.} Note that if $e(2,3;l-1) \in H_2$, then we're back to essentially the same scenario as Case 1.1, so we may assume that $e(2,3;l-1) \notin H_2$. Then $e(1,2;l-1) \in H_2$. It follows that $e(1;l-1,l) \notin H'$ and so $R(0,l-1)$ is not a small cookie of $H$. Now, either $e(1,2;l) \in H_1$ (Figure 3.10 (a)), or $e(1,2;l) \notin H'$ (Figure 3.10 (b)).

\endgroup

\begingroup
\setlength{\intextsep}{0pt}
\setlength{\columnsep}{10pt}
\begin{wrapfigure}[]{r}{0cm}
\begin{adjustbox}{trim=0cm 0.5cm 0cm 0.25cm}

\end{adjustbox}
\end{wrapfigure}

\null

\noindent \textit{CASE 2.1 (a):  $e(1,2;l) \in H_1$.} By Observation 3.4 and Proposition 3.3,  $Z \mapsto R(1,l-1)$ is the cascade we seek. End of Case 2.1(a). 

\null

\noindent \textit{CASE 2.1 (b):  $e(1,2;l) \notin H'$.} Then we have that $e(2;l,l+1) \in H_1$, that $e(1;l,l+1) \in H_1$ and that $R(1,l) \in \text{int}(H)$. It follows that $R(0,l) \in \text{ext}(H)$, so $R(0,l)$ must be a small cookie. Then, after $R(0,l) \mapsto R(1,l)$, we are back to Case 2.1 (a). End of Case 2.1(b). End of Case 2.1.

\endgroup 

\begingroup
\setlength{\intextsep}{0pt}
\setlength{\columnsep}{10pt}
\begin{wrapfigure}[]{l}{0cm}
\begin{adjustbox}{trim=0cm 0cm 0cm 0cm}
\begin{tikzpicture}[scale=1.5]

\begin{scope}[xshift=0cm]
{
\draw[gray,very thin, step=0.5cm, opacity=0.5] (0,0) grid (1.5,1.5);

\fill[green!50!white, opacity=0.5] (0,0.5)--++(1.5,0)--++(0,0.5)--++(-1.5,0);
\fill[green!50!white, opacity=0.5] (0.5,0)--++(0.5,0)--++(0,0.5)--++(-0.5,0);

\draw[blue, line width =0.5mm] (0.5,1)--++(0.5,0);

\draw[orange, line width =0.5mm] (1,0)--++(0,0.5);
\draw[green!50!black, line width =0.5mm] (0,0)--++(0.5,0)--++(0,0.5)--++(-0.5,0);

\draw[yellow, line width =0.2mm] (0,0)--++(0,1.5);

\draw[fill=blue] (1,1) circle [radius=0.05];
\draw[fill=blue] (1,1.5) circle [radius=0.05];
\draw[fill=blue] (0.5,1) circle [radius=0.05];

\draw[fill=orange] (1,0.5) circle [radius=0.05];
\draw[fill=orange] (1,0) circle [radius=0.05];

% black lines
{

\draw[black, line width=0.15mm] (0.95, 0.7)--++(0.1,0);
\draw[black, line width=0.15mm] (0.95, 0.75)--++(0.1,0);
\draw[black, line width=0.15mm] (0.95, 0.8)--++(0.1,0);

\draw[black, line width=0.15mm] (1.2, 0.95)--++(0,0.1);
\draw[black, line width=0.15mm] (1.25, 0.95)--++(0,0.1);
\draw[black, line width=0.15mm] (1.3, 0.95)--++(0,0.1);

\draw[black, line width=0.15mm] (0.7, 0.45)--++(0,0.1);
\draw[black, line width=0.15mm] (0.75, 0.45)--++(0,0.1);
\draw[black, line width=0.15mm] (0.8, 0.45)--++(0,0.1);

\draw[black, line width=0.15mm] (0.45, 0.7)--++(0.1,0);
\draw[black, line width=0.15mm] (0.45, 0.75)--++(0.1,0);
\draw[black, line width=0.15mm] (0.45, 0.8)--++(0.1,0);

}

\node[left] at (0,0.5) [scale=1] {\tiny{-1}};
\node[left] at (-0.05,1) [scale=1] {\tiny{$\ell$}};
\node[left] at (0,1.5) [scale=1] {\tiny{+1}};

\node[below] at (0.5,0) [scale=1] {\tiny{1}};
\node[below] at (1,0) [scale=1] {\tiny{2}};
\node[below] at (1.5,0) [scale=1] {\tiny{3}};

\node[below, align=left, text width=3cm] at (0.75,-0.25) 
{Fig. 3.11. Case 2.2.};

}
\end{scope}

\end{tikzpicture}
\end{adjustbox}
\end{wrapfigure}

\null

\noindent \textit{CASE 2.2: $e(2,3;l) \notin H'$.} Then we have that $e(1,2;l) \in H_1$. Note that if $e(1,2;l-1) \in H_2$, then we're back to Case 2.1, so we may assume that $e(1,2;l-1) \notin H'$. It follows that $e(2;l-2,l-1) \in H_2$. Note that if $e(1; l-1,l) \in H_1$, then $R(0,l-1) \in \text{ext}(H)$ and $R(0,l) \in \text{ext}(H)$, contradicting Corollary 1.3.15 (b), so we may assume that $e(1; l-1,l) \notin H'$. Then we must have that $e(1;l-2,l-1) \in H'$. This implies that $R(0,l-2)$ is a small cookie of $H$. Then after $R(0,l-2) \mapsto R(1,l-2)$, we're back to Case 2.1(a). End of Case 2.2. End of Case 2. $\square$.

\endgroup 

\null

\noindent \textbf{Observation 3.6.} Let $H$ be a Hamiltonian cycle of an $m \times n$ grid graph $G$, with $m,n \geq 5$ and let $J$ be a large cookie of $G$. Then $J \cap R_2 \neq \emptyset$.

\null

\noindent \textbf{Lemma 3.7.} Let $H$ be a Hamiltonian cycle of an $m \times n$ grid graph $G$, with $m,n \geq 5$ and assume that $G$ has at least two large cookies $J_1$ and $J_2$. Then switching the neck $N_{J_1}$ of $J_1$ splits $H$ into two cycles $H_1$ and $H_2$ such that there is $v_1 \in H_1 \cap R_2$ and $v_2 \in H_2 \cap R_2$ with $v_1$ adjacent to $v_2$.

\null

\noindent \textit{Proof.} Orient $H$. Let $\{v_x, v_y\}$ be the boundary edge of the neck $N_{J_1}$ of $J_1$. Define $\overrightarrow{K}_1$ and $\overrightarrow{K}_2$ to be the subtrails $\overrightarrow{K}((v_x,v_{x+1})$, $(v_{y-1}, v_y))$ and $\overrightarrow{K}((v_y,v_{y+1}), (v_{x-1},v_x))$ of $\overrightarrow{K}_H$, respectively. By Lemma 1.4.1 (ii), switching $N_{J_1}$ gives two cycles $H_1$ and $H_2$, with  $V(H_1)=V(\overrightarrow{K}_1 \setminus \{(v_x, v_y)\})$ and $V(H_2)=V(\overrightarrow{K}_2)$. By Corollary 1.3.11 (a) and (d), $V(J_1) = V(\overrightarrow{K}_1)$. By Observation 3.6, $V(J_1) \cap R_2 \neq \emptyset$. Then $V(H_1) = R_2 \neq \emptyset$. Since $V(J_1) = V(\overrightarrow{K}_1)$, we have that $V(J_2) \subseteq V(\overrightarrow{K}_2)=V(H_2)$. By Observation 3.6, $V(J_2) \cap R_2 \neq \emptyset$. It follows that $ V(H_2) \cap R_2 \neq \emptyset$. 

We have shown that $V(H_1) = R_2 \neq \emptyset$ and that $ V(H_2) \cap R_2 \neq \emptyset$. It remains to check there is $v_1 \in H_1 \cap R_2$ and $v_2 \in H_2 \cap R_2$ with $v_1$ adjacent to $v_2$. Let $v_1\in H_1 \cap R_2$ and let $R_2 = w_1, \ldots, w_s$, with $v_1=w_1$. Sweep $R_2$ starting at $w_1$. If there is $i \in \{1,\ldots,s-2\} $ such that $v_{i} \in H_1$ and $v_{i+1} \in H_2$, we are done. If there is no such $i$ then $R_2 \cap H_2 =\emptyset$, contradicting that $R_2 \cup H_2 \neq \emptyset$. $\square$

\null

\noindent \textbf{Proposition 3.8. \index{MLC Algorithm}(MLC Algorithm.)} Let $H$ be a Hamiltonian cycle of an $m \times n$ grid graph $G$. Assume that $G$ has more than one large cookie. Then there is a cascade of length at most two that reduces the number of large cookies of $G$ by one. 

\null 

\noindent \textit{Proof.} Let $J$ be a large cookie of $G$ with neck $N_J$. Switching $N_J$ splits $H$ into the cycles $H_1$ and $H_2$. By Lemma 3.7 there is $v_1 \in H_1 \cap R_2$ and $v_2 \in H_2\cap R_2$ with $v_1$ adjacent to $v_2$. By Lemma 3.5 there is a cascade of length at most two, whose last move is $N_J \mapsto N_J'$ with $N_J \neq N_J'$. By Observation 3.4, this cascade decreases the number of large cookies of $G$ by one. $\square$

\begin{algorithm}[H]
\caption{ManyLargeCookies}\label{alg:mlc}
\KwIn{Hamiltonian cycle $H$ on $m \times n$ rectangular grid graph.}
\KwOut{Hamiltonian cycle $H'$ with all large cookies collected.}
\BlankLine
\textbf{Assumptions and Notation:} 
Without loss of generality, assume $m \leq n$. 
For a large cookie $J$ with neck $N_J$, let $(v_1, v_2)$ denote a pair of adjacent vertices belonging to cycles $H_1$ and $H_2$ (arising from switching $N_J$) that are incident on $R_2$ (guaranteed by Lemma~3.7).
Let $\mathbf{B}$ denote the $4 \times 4$ rectangle centered at $v_1$ (or $v_2$).
\BlankLine

Identify all large cookies \tcp{$O(m+n)$}
Let numLargeCookie $\gets$ number of large cookies\;
\While{numLargeCookie $\geq 2$}{
    Let $J$ be any remaining large cookie\;
    Let $N_J \gets$ neck of $J$ \tcp{$O(1)$}
    Let $(v_1, v_2) \gets$ vertex pair for $J$ \tcp{Lemma~3.7}
    Let $\mathbf{B} \gets$ $4 \times 4$ box at $v_1$\;
    \BlankLine
    \tcp{Try direct move}
    \ForEach{box $Z'$ in region $\mathbf{B}$}{
        \If{move $N_J \to Z'$ is valid}{
            Execute($N_J \to Z'$)\;
            numLargeCookie $\gets$ numLargeCookie $-1$\;
            \textbf{break}\;
        }
    }
    \BlankLine
    \tcp{Direct move failed - find prep move}
    \ForEach{box $Z''$ in region $\mathbf{B}$}{
        \ForEach{box $Z'$ in region $\mathbf{B}$}{
            \If{both $Z'' \to Z'$ and $N_J \to Z''$ valid}{
                Execute($Z'' \to Z'$)  \tcp{Prep move} 
                Execute($N_J \to Z''$) \tcp{direct move} 
                numLargeCookie $\gets$ numLargeCookie $-1$\;
                \textbf{break}\;
            }
        }
    }
}
\Return $H'$\;
%\BlankLine
%\textbf{Time complexity:} $O(mn)$ total
\end{algorithm}

\BlankLine
\noindent \textbf{Time complexity:} $O(m+n)$ total. The initial scan takes $O(m+n)$ time. Each cookie collection takes $O(1)$ time (search is localized to region $\mathbf{B}$), and there are at most $O(m+n)$ cookies.

\subsection{Existence of the 1LC Algorithm}

\noindent In this section we prove the 1LC algorithm (Proposition 3.10). We need some definitions before we can state it.

\null

\noindent \textbf{Definitions}.  Let $G$ be an $m \times n$ grid graph and $H$ be a Hamiltonian cycle of $G$. We call a subpath of $H$ on the edges $e(k;l,l+1), e(k,k+1;l+1)$ and $e(k+1;l,l+1)$ a \index{leaf|textbf}\textit{northern leaf}. We will often say that $R(k,l)$ is a northern leaf to mean that $e(k;l,l+1), e(k,k+1;l+1)$ and $e(k+1;l,l+1)$ belong to $H$. Southern, eastern and western leaves are defined analogously. We call the subgraph of $H$ on the edges $e(k-1,k;l), e(k;l,l+1)$, $e(k+1,k+2;l)$, and $e(k+1;l,l+1)$ a \index{A-type@$A$-type|textbf}\textit{northern $A$-type}. Suppose $H$ has a northern $A$-type. We call the subgraph $A \cup e(k,k+1;l+1)$ of $H$ a \index{A0-type@$A_0$-type|textbf}\textit{northern $A_0$-type}, and we call the subgraph $A \cup e(k;l+1,l+2) \cup e(k+1;l+1,l+2)$ of $H$ a \index{A1-type@$A_1$-type|textbf}\textit{northern $A_1$-type}. We make analogous definitions for eastern, southern and western $A$-types. See Figure 3.12.

\setlength{\intextsep}{0pt}
\setlength{\columnsep}{20pt}
\begin{adjustbox}{trim=0cm 0cm 0cm 0cm}

\end{adjustbox}

\noindent Let $R(k,l-1)$ be a northern leaf. If $H$ has a northern $A$-type on $e(k-1,k;l+1), e(k;l+1,l+2)$ and $e(k+1,k+2;l+1), e(k+1;l+1,l+2)$ then we say that $A$-type \index{follow (a leaf, or an edge)|textbf}\textit{follows} the northern leaf $R(k,l-1)$ \textit{northward}. We call the boxes $R(k,l+1)$ and $R(k,l+2)$ the \index{middle box (of an A1 type)@middle box of an $A_1$-type|textbf}\textit{middle-boxes} of the $A_1$-type. We call the box $R(k,l+1)$ the \index{switchable middle-box of an A1 type@switchable middle-box of an $A_1$-type|textbf}\textit{switchable middle-box} of the $A_1$-type. Analogous definitions apply for other compass directions.

Let $e(k,k+1;l)$ be an edge in $H$, let $A$ be a northern $A_0$-type in $H$ on the edges $e(k-1,k;l+1)$, $e(k;l+1,l+2)$, $e(k,k+1,l+2)$,  $e(k+1;l+1,l+2)$, and $e(k+1,k+2;l+1)$, and let $j \in \big\{1,2, \ldots \big\lfloor \frac{n-1-l}{2} \big\rfloor \big\}$. We define a \index{j-stack (of A0s)@$j$-stack of $A_0$-types|textbf}\textit{northern j-stack of $A_0$s starting at $A$} to be a subgraph \textit{$stack(j; A_0)$} of $H$, where $\textrm{stack}(j; A_0)= \bigcup_{i=0}^{j-1} \Big(A_{0} +(0,2i) \Big)$. If $j=\big\lfloor \frac{n-1-l}{2}   \big\rfloor$, we call the $j$-stack a \index{full j-stack (of A0s)@full $j$-stack of $A_0$-types|textbf}\textit{full j-stack of $A_0$'s}. We will often say that the $j$-stack of $A_0$'s is \textit{following} the edge $e(k,k+1;l)$ whenever this edge plays an important role. Eastern, southern, and western $j$-stacks are defined analogously. See Figure 3.14(a) on page 53.

We extend the notion of \index{collect (a cookie, or a leaf)}\textit{collecting} to leaves. Let $R(k,l)$ be leaf of $H$. WLOG assume that $R(k,l)$ is northern. We say that the cascade $\mu_1, \ldots, \mu_r$ \textit{collects} $R(k,l)$ if $\mu_r$ is the move $R(k,l+1) \mapsto R(k,l)$. 

We denote the set of northern and southern small cookies by $\text{SmallCookies}\{N,S\}$ and the set of eastern and western small cookies by $\text{SmallCookies}\{E,W\}$. Let $C$ be an easternmost or westernmost northern or southern small cookie in $\text{SmallCookies}\{N,S\}$. Then we say that $C$ is an \index{outermost small cookie|textbf}\textit{outermost} small cookie in $\text{SmallCookies}\{N,S\}$. Outermost small cookies in $\text{SmallCookies}\{E,W\}$ are defined analogously. 

\begingroup
\setlength{\intextsep}{0pt}
\setlength{\columnsep}{20pt}
\begin{wrapfigure}[]{r}{0cm}
\begin{adjustbox}{trim=0cm 0cm 0cm 0cm}
\begin{tikzpicture}[scale=1.5]

% TOP %
%%%%%%%%%%%%  Initial Configuration %%%%%%%%%%%%% 
\begin{scope} 
{
\draw[gray,very thin, step=0.5cm, opacity=0.5] (0,0) grid (1.5,2.5);

% colors
{
\fill[blue!50!white, opacity=0.5] (0,1) rectangle (0.5,1.5);
\fill[blue!50!white, opacity=0.5] (1,1) rectangle (1.5,1.5);
\fill[blue!50!white, opacity=0.5] (0,1.5) rectangle (0.5,2);
\fill[blue!50!white, opacity=0.5] (1,1.5) rectangle (1.5,2);

\draw[blue, line width=0.5mm] (0,1)--++(0.5,0)--++(0,0.5); 
\draw[blue, line width=0.5mm] (1,1.5)--++(0,-0.5)--++(0.5,0); 

\fill[blue!40!white, opacity=0.4](0.5,0)--(0.5,0.5)--(1,0.5)--(1,0);
\draw[blue, line width=0.5mm] (0.5,0)--(0.5,0.5)--(1,0.5)--(1,0);
}
% labellings 
{
\node[left] at (0,0) [scale=1]
{\tiny{$\ell$}};

\foreach \x in {1, ...,4}
\node[left] at (0,0.5*\x) [scale=1]
{\tiny{+\x}};

\node at (0.5,-0.15) [scale=1]
{\tiny{k}};

\node at (1,-0.15) [scale=1]
{\tiny{+1}};

\node at (1.5,-0.15) [scale=1]
{\tiny{+2}};
}
%black lines
{
\draw[black, line width=0.15mm] (1.2,1.45)--(1.2,1.55);
\draw[black, line width=0.15mm] (1.25,1.45)--(1.25,1.55);
\draw[black, line width=0.15mm] (1.3,1.45)--(1.3,1.55);

\draw[black, line width=0.15mm] (0.2,1.45)--(0.2,1.55);
\draw[black, line width=0.15mm] (0.25,1.45)--(0.25,1.55);
\draw[black, line width=0.15mm] (0.3,1.45)--(0.3,1.55);

\draw[black, line width=0.15mm] (0.7,0.95)--(0.7,1.05);
\draw[black, line width=0.15mm] (0.75,0.95)--(0.75,1.05);
\draw[black, line width=0.15mm] (0.8,0.95)--(0.8,1.05);
}

\node[below, align=center, text width=3cm] at (0.75,-0.25) 
{Fig. 3.13.};

}
\end{scope}

\end{tikzpicture} 
\end{adjustbox}
\end{wrapfigure}

\null 

\noindent Given a leaf  $L$, we want to show that there is a cascade that collects $L$. For definiteness, assume that $L$ is a northern leaf $R(k,l)$. Furthermore, assume that $L$ is not incident on the northern side of the boundary, so $l+1 <n-1$. If $e(k,k+1; l+2) \in H$, then $L+(0,1) \mapsto L$ is the cascade we seek, so we only need to consider the case where $e(k,k+1; l+1) \notin H$. Then we must have $e(k-1,k;l+2) \in H$, $e(k+1,k+2; l+2) \in H$, $e(k;l+2,l+3) \in H$ and $e(k+1; l+2,l+3) \in H$. Note that if $e(k-1,k;l+3) \in H$, then $L +(0,2) \mapsto L+(-1,2)$ followed by $L+(0,1) \mapsto L$ is the cascade we seek, so we consider the case where $e(k-1,k;l+3) \notin H$ and, by symmetry, where $e(k+1,k+2; l+3) \notin H$. 
Now, we either have $e(k;l+3,l+4) \in H$  and $e(k+1;l+3,l+4) \in H$ or  $e(k,k+1;l+3) \in H$. That is, $L$ is followed northward by an $A_0$-type or by an $A_1$-type. See Figure 3.13. From this point onward, we will omit the compass direction when it does not introduce ambiguity. We coalesce this paragraph into the following lemma:

\endgroup

\null 

\noindent \textbf{Lemma 3.9.} Let $L=R(k,l)$ be a northern leaf with $l+1<n-1$. If $L$ is not followed by an $A_0$-type or by an $A_1$-type, then there is a cascade of length at most two that collects $L$.

\null 

\noindent \textbf{Proposition 3.10. \index{1LC Algorithm}(1LC Algorithm).} Let $G$ be an $m \times n$ grid graph, let $H$ be a Hamiltonian cycle of $G$. If $H$ has exactly one large cookie and at least one small cookie, then there is a cascade of length at most $\frac{1}{2}\max(m,n)+\min(m,n)+2$ moves that reduces the number of small cookies of $H$ by one and such that it does not increase the number of large cookies.

\null 

\noindent The proof of Proposition 3.10 requires the following two lemmas.

\noindent \textbf{Lemma 3.11.} Let $G$ be an $m \times n$ grid graph, let $H$ be a Hamiltonian cycle of $G$, and let $C \in \text{SmallCookies}\{N,S\}$ be an easternmost small cookie. Assume that $G$ has only one large cookie. Then there cannot be a full $j$-stack of $A_0s$ starting at the $A_0$-type that contains $C$. 

\null

\noindent \textbf{Lemma 3.12.} Let $G$ be an $m \times n$ grid graph, let $H$ be a Hamiltonian cycle of $G$, and let $C \in \text{SmallCookies}\{N,S\}$ be an easternmost small cookie. Assume that there is a $j$-stack of $A_0$ starting at the $A_0$-type containing $C$. Let $L$ be the leaf in the top ($j^{\text{th}}$) $A_0$ of the stack. Assume that $L$ is followed by an $A_1$-type and that $G$ has only one large cookie. Then there is a cascade of at most $\min(m,n)+3$ moves that collects $L$. 

\null

\noindent \textit{Proof of Proposition 3.10.} Since there is at least one small cookie, one of $\text{SmallCookies}{N,S}$ or $\text{SmallCookies}{E,W}$ must be nonempty. Without loss of generality we may assume that $\text{SmallCookies}\{N,S\}$ is nonempty. Let $C \in \text{SmallCookies}\{N,S\}$ be an easternmost small cookie. For definiteness, assume that $C$ is a small northern cookie on the southern boundary.  Let $Q(j)$ be the statement ``There is a $j$-stack of $A_0$'s starting at the $A_0$-type containing $C$''. Note that $C$ is contained in an $A_0$-type, so $Q(1)$ is true. By Lemma 3.11,  there is a $j_0 \in \Big\{2,3, \ldots\Big\lfloor \frac{n}{2} \Big\rfloor \Big\}$ such that for each $j \in \{2.., j_0-1\}$, $Q(j)$ is true for each $j < j_0$ but $Q(j_0)$ is not true.

%By Lemma 3.9, we may assume that $C$ is followed by an $A_0$-type or an $A_1$-type. If $C$ is followed by an $A_1$-type, by Lemma 3.12, there is a cascade that collects $C$, so we may assume that $C$ is followed by an $A_0$-type. 

\begingroup
\setlength{\intextsep}{0pt}
\setlength{\columnsep}{20pt}
\begin{center}

\end{center}

\noindent For $j \in \{1.., j_0-1\}$ let $L_j$ be the northern leaf of the $j^{\text{th}}$ $A_0$-type in the stack. Note that $L_1=C$. Lemma 3.9 implies that $L_{j_0-1}$ is either followed by an $A_1$-type or there is a cascade that collects $L_{j_0-1}$. If $L_{j_0-1}$ is followed by an $A_1$-type, then by Lemma 3.12, we can find a cascade that collects it, so we only need to check the case in which there is a cascade $\mu_1, \ldots, \mu_s$ that collects $L_{j_0-1}$. Note that $\mu_s$ must be the move $L_{j_0-1}+(1,0) \mapsto L_{j_0-1}$.
Then $\mu_1, \ldots, \mu_s, L_{j_0-2}+(0,1) \mapsto L_{j_0-2} , \ldots, L_1+(0,1) \mapsto L_1$ is a cascade that collects $C$. Note that $j \leq \frac{n}{2}$, and that by Lemma 3.12, there are at most $\min(m,n)+3$ moves required to collect $L$. After that, we need at most another $j-1$ flips to collect $C$, so $C$ can be collected after at most, $\frac{1}{2}\max(m,n)+\min(m,n)+2$ moves to collect $C$. See Figure 3.14 for an illustration with $j_0-1=3$. $\square$

\null 

\noindent We prove Lemma 3.11 next. The proof of Lemma 3.12 takes up the remainder of the section.

\null 

\noindent \textit{Proof of Lemma 3.11.} Assume for a contradiction that $C$ is in a full stack of $A_0$'s starting at the $A_0$ that contains $C$. For definiteness, assume that $C=R(k,0)$ is a small northern cookie on the southern boundary. First we check that $m-1 >k+2$. If $m-1=k+2$, then we must have $e(k+2;0,1) \in H$ and $e(k+2;1,2) \in H$. But then $H$ misses $v(k+2,3)$ (in green in Figure 3.15 (a)). There fore, we must have that $m-1 >k+2$.

The number $j$ of $A_0$'s in the full stack is even or an odd so there are two cases to check. Note that for each odd $i \in \{1, 2, \ldots, j \}$, the $i^{\text{th}}$ $A_0$ belongs to $\text{ext}(H)$ and for each even $i \in \{1, 2, \ldots, j\}$, the $i^{\text{th}}$ $A_0$ belongs to $\text{int}(H)$.

\begingroup
\setlength{\intextsep}{0pt}
\setlength{\columnsep}{20pt}
\begin{wrapfigure}[]{r}{0cm}
\begin{adjustbox}{trim= 0cm 0cm 0cm 0cm}
 %%% CASE 3.1 %%%

\end{adjustbox}
\end{wrapfigure}

\null

\noindent \textit{CASE 1: $j$ is even.} Note that the top leaf of the stack is in $\text{int}(H)$. Now, $n-1$ is either even or odd. 

\null 

\noindent \textit{CASE 1.1: $n-1$ is even.} We have that $R(k,n-3) \in \text{int}(H)$. But then we must have $R(k,n-2) \in \text{ext}(H)$ and $R(k+1,n-2) \in \text{ext}(H)$, contradicting Lemma 1.14. End of Case 1.1. See Figure 3.15 (a).

\null 

\noindent \textit{CASE 1.2: $n-1$ is odd.} Then we must have that $R(k+1,n-2)$ is a small southern cookie. But this contradicts our assumption that $C$
is the easternmost small cookie in $\text{SmallCookies}\{N,S\}$.  End of Case 1.2. End of Case 1. See Figure 3.15 (b).

\endgroup

\null 

\noindent \textit{CASE 2: $j$ is odd.} Note that the top leaf of the stack is in $\text{ext}(H)$. Again, $n-1$ is either even or odd.

\setlength{\intextsep}{0pt}
\setlength{\columnsep}{20pt}
\begin{center}
\begin{adjustbox}{trim= 0cm 0cm 0cm 0cm}
 %%% CASE 3.1 %%%

\end{adjustbox}
\end{center}

\noindent \textit{CASE 2.1: $n-1$ is odd.} We have that  $R(k,n-2) \in \text{ext}(H)$. But then, the fact that $e(k,k+1;n-1) \in H$ implies that $R(k,n-2)$ is not a cookie neck, contradicting Lemma 1.14. End of Case 2.1. See Figure 3.16 (a).

\null

\noindent \textit{CASE 2.2: $n-1$ is even.} We have that $e(k+1,k+2;0) \in H$. Then, either $e(k+2, k+3;0) \in H$, or $e(k+2;0,1) \in H$.

\null

\noindent \textit{CASE 2.2(a): $e(k+2, k+3;0) \in H$.} Then we must have $e(k+2;1,2) \in H$. Note that for $i \in \{ 1,3, \ldots, n-2\}$, $e(k+2;i,i+1)\in H$ implies $e(k+2;i+2,i+3)\in H$. Then, for $i \in \{ 1,3, \ldots, n-2\}$, we have that $e(k+2;i,i+1)\in H$. Note that we must also have $e(k+2,k+3;n-1) \in H$. Then $R(k+2,n-2)$ must be a southern small cookie,  contradicting the eastmost assumption. See Figure 3.16 (b) End of Case 2.2(a).

\null

\noindent \textit{CASE 2.2(b): $e(k+2;0,1) \in H$.} Note that if $e(k+2,k+3;1) \in H$, then $R(k+2,0)$ must be a small cookie, contradicting the eastmost assumption. Then $e(k+2,k+3;1) \notin H$. But then we have $e(k+2;1,2) \in H$, and we are back to Case 2.2(a). See Figure 3.16 (c). End of Case 2.2(b). End of Case 2.2. End of Case 2. $\square$

\endgroup

\null

\noindent It remains to prove Lemma 3.12. Its proof will require Lemmas 3.13-3.16.

\null

\noindent \textbf{Lemma 3.13.} Let $G$ be an $m \times n$ grid graph, and let $H$ be a Hamiltonian cycle of $G$. Let $C$ be a small cookie of $G$.  Assume that $G$ has only one large cookie, and that there is a $j$-stack of $A_0$ starting at the $A_0$-type containing $C$. Let $L$ be the leaf in the top ($j^{\text{th}}$) $A_0$ of the stack, and assume that $L$ is followed by an $A_1$-type. Let $X$ and $Y$ be the boxes adjacent to the switchable middle-box of the $A_1$-type that are not its $H$-neighbours. If $P(X,Y)$ has no switchable boxes, then either:

(i)  there is a cascade of length at most $\min(m,n)$, which avoids the stack of $A_0$'s, and 

\hspace{0.6 cm} after which $P(X,Y)$ gains a switchable box, or 

(ii) there is a cascade of length at most $\min(m,n)+1$, that collects $L$ and avoids the 

\hspace{0.6 cm} stack of $A_0$'s.

\null 

\noindent We postpone the proof of Lemma~3.13 to Chapter~4. For now, we assume Lemma~3.13 and use it to prove Lemma~3.12. Recall the setup: we have a $j$-stack of $A_0$-types followed by an $A_1$-type with middle-box $X'$. A consequence of the assumption that there is only one large cookie is that $X'$ must determine an $H$-path. That is, if $X$ and $Y$ are the boxes adjacent to $X'$ that are not its $H$-neighbors, then $X$ and $Y$ lie in the same $H$-component, so the $H$-path $P(X,Y)$ is well-defined. Lemma~3.13 guarantees that either $P(X,Y)$ already has a switchable box or there exists a cascade that produces one. 
 
Lemma~3.14 shows that $X'$ must lie in $G_2$. Lemma~3.15 handles the case $X' \in G_2 \setminus G_3$. For the remaining case, $X' \in G_3$, we use the existence of a switchable box in $P(X,Y)$ (guaranteed by Lemma~3.13) to prove Lemma~3.16. This lemma plays a role analogous to Lemma~3.7 in Section~3.1: it shows that switching $X'$ yields two cycles with adjacent vertices on $R_2$. We then can apply the part of Lemma~3.5 covering the case $Z \in G_3$ to find either a switchable box or a cascade of length at most two that delivers the required move. 

\null

\noindent \textbf{Lemma 3.14.} Let $G$ be an $m \times n$ grid graph, let $H$ be a Hamiltonian cycle of $G$, and let $C \in \text{SmallCookies}\{N,S\}$ be an easternmost small cookie. Assume that $G$ has only one large cookie, and that there is a $j$-stack of $A_0$ starting at the $A_0$-type containing $C$. Let $L$ be the leaf contained in the top ($j^{\text{th}}$) $A_0$ of the stack. Assume that $L$ is followed by an $A_1$-type with looping $H$-path $P(X,Y)$. Let $X'$ be the box of $G$ that shares edges with $X$ and $Y$. Then $X' \in G_2$.

\null 

\setlength{\intextsep}{0pt}
\setlength{\columnsep}{20pt}
\begin{wrapfigure}[]{l}{0cm}
\begin{adjustbox}{trim=0cm 0.5cm 0cm 0.25cm} 

\end{adjustbox}
\end{wrapfigure}

\noindent \textit{Proof.} For definiteness, assume that $L$ is the northern leaf $R(k,l-2)$, and that $X=R(k-1,l)$. Then $X'=R(k,l)$ and $Y=R(k+1,l)$. Note that $l-2 \geq 0$ and $l+2 \leq n-1$. Either $P(X,Y)$ is contained in $\text{ext}(H)$, or $P(X,Y)$ is contained in $\text{int}(H)$, so there are two cases to check.

\null
\null 

\noindent \textit{CASE 1}: $P(X,Y) \subset \text{ext}(H)$. By Lemma 1.14, we must have that $m-1 >k+2$ and $k-1>0$. To see that $n-1 > l+2$, assume for a contradiction that $n-1 = l+2$. By Lemma 1.14, $X+(1,0)$ and $Y+(1,0)$ are cookie necks. But this contradicts the assumption that there is only one large cookie in $G$. See Figure 3.17 (a). End of Case 1. 

\null

\noindent \textit{CASE 2}: $P(X,Y) \subset \text{int}(H)$. By Lemma 1.14, we must have that $m-1 >k+2$ and $k-1>0$. To see that $n-1 > l+2$, assume for a contradiction that $n-1 = l+2$. Lemma 1.14 implies that $X'+(0,1)$
is the neck of the large cookie of $G$. But now $X'+(2,1)$ must be a small cookie of $G$, contradicting the easternmost assumption. See Figure 3.17 (b).

\null

\endgroup 

\noindent \textbf{Lemma 3.15.} Let $G$ be an $m \times n$ grid graph, let $H$ be a Hamiltonian cycle of $G$, and let $C \in \text{SmallCookies}\{N,S\}$ be an easternmost small cookie. Assume that $G$ has only one large cookie, and that there is a $j$-stack of $A_0$ starting at the $A_0$-type containing $C$. Let $L$ be the leaf in the top ($j^{\text{th}}$) $A_0$ of the stack. Assume that $L \in \text{int}(H)$ and that $L$ is followed by an $A_1$-type with looping $H$-path $P(X,Y)$. Let $X'$ be the box of $G$ that shares edges with $X$ and $Y$. If $X'$ is not in $G_3$ then, either $X' \mapsto W$ is a cascade, or there is a cascade $\mu, X' \mapsto W$, of length two, with $X' \mapsto W$ nontrivial in either case.

\null

\begingroup
\setlength{\intextsep}{0pt}
\setlength{\columnsep}{20pt}
\begin{wrapfigure}[]{r}{0cm}
\begin{adjustbox}{trim= 0cm 0cm 0cm 0cm}
 %%% CASE 3.1 %%%
\begin{tikzpicture}[scale=1.5]

\begin{scope}[xshift=0cm, yshift=0cm ]

\draw[gray,very thin, step=0.5cm, opacity=0.5] (0,0) grid (2,1.5);

{
\draw[blue, line width=0.5mm] (0.5,0)--(0.5,0.5)--(1,0.5)--(1,0)--++(0.5,0);

\draw[blue, line width=0.5mm] (0,1)--++(0.5,0)--++(0,0.5)--++(0.5,0)--++(0,-0.5)--++(0.5,0);

\draw[blue, thick] (2,1) circle [radius=0.05];

\draw[blue, line width=0.5mm] (2,0)--++(0,0.5);

\draw[blue, line width=0.5mm] (1.5,0)--++(0.5,0);
\draw[blue, line width=0.5mm] (1.5,1)--++(0,-0.5)--++(0.5,0);

}

{

\node[below] at (0.5,0) [scale=1]
{\tiny{k}};
\node[below] at (1,0) [scale=1]
{\tiny{+1}};
\node[below] at (1.5,0) [scale=1]
{\tiny{+2}};
\node[below] at (2,0) [scale=1]
{\tiny{+3}};

\node at (0.75,0.25) [scale=0.8]
{\small{$C$}};

}

\node[below] at (1,0) [scale=0.75] {\small };

\node[below, align=center, text width=3cm] at (1,-0.25) {Fig 3.18. $m-1=k+3$.};

\end{scope}

\end{tikzpicture}
\end{adjustbox}
\end{wrapfigure}

\noindent \textit{Proof.} Suppose that $X'$ is not in $G_3$. By Lemma 3.14, $X' \in G_2 \setminus G_3$. For definiteness, assume that $L$ is a northern leaf, and let $X'=R(k,l)$. The assumption that $L \in \text{int}(H)$ implies that $l-2>0$.  See Figure 3.18.

Now we check that $m-1 > k+3$. By Lemma 3.14, $m-1 > k+2$. For a contradiction, assume that $m-1=k+3$. Note that we must have $e(k+1,k+2;0) \in H$, $e(k+2,k+3;0) \in H$, and $e(k+3;0,1) \in H$. This implies that we must have $e(k+2; 1,2)\in H$ and $e(k+2,k+3; 1)\in H$. But now $H$ misses $v(k+3;2)$. It follows that we must have $m-1 > k+3$. See Figure 3.18. By symmetry, $0 < k-2$. It follows that $l+3=n-1$. 

\noindent The same argument used in Case 2.2 of Lemma 3.11 (see Figure 3.16 (b) and (c)) shows that we have $e(k+2;l+1,l+2)\in H$ and $e(k+2,k+3;l+1)\in H$. Now either $e(k,k+1;l+2) \in H$ or $e(k,k+1;l+2) \notin H$. See Figure 3.19. 

\null 

\noindent \textit{CASE 1: $e(k,k+1;l+2) \in H$. } By Lemma 1.16, $X' \mapsto X' +(1,1)$ is a valid move, and by Observation 3.4, $X' \mapsto X' +(1,1)$ creates no new cookies. End of Case 1.

\null

\setlength{\intextsep}{0pt}
\setlength{\columnsep}{40pt}
\begin{wrapfigure}[]{l}{0cm}
\begin{adjustbox}{trim=0cm 0cm 0cm 0cm} 

\end{adjustbox}
\end{wrapfigure}

\noindent \textit{CASE 2: $e(k,k+1;l+2) \notin H$. } Lemma 1.14 implies that $e(k+1,k+2;l+2)$ cannot be in $H$ either. It follows that $X'+(0,2)$ must be the neck of the large cookie. The assumption that there is only one large cookie implies that $e(k+2;l+2,l+3) \notin H$. Then we must have $e(k+2, k+3;l+2) \in H$. Then, by Lemma 1.16, $X'+(1,1) \mapsto X'+(2,1)$, $X' \mapsto X'+(1,0)$ is the cascade we seek. $\square$

\endgroup 

\null 

\noindent \textbf{Lemma 3.16.}  Let $H$ be a Hamiltonian cycle of an $m \times n$ grid graph $G$. Let $X' \in G_3 \cap \text{ext}(H)$ be a switchable box, and let $P(X,Y)$ be the looping $H$-path of $X'$. Assume that $P(X,Y)$ has a switchable box in $G_0 \setminus G_2$. Then switching $X'$ splits $H$ into two cycles $H_1$ and $H_2$ such that there is $v_1 \in H_1 \cap R_2$ and $v_2 \in H_2 \cap R_2$ with $v_1$ adjacent to $v_2$.

\null

\noindent \textit{Proof.} Let $Z$ be a switchable box of $P(X,Y)$ in $G_0 \setminus G_2$. Orient $H$. Let $(v_x, v_{x+1})$ and $(v_{y-1}, v_y)$ be the edges of $X'$ in $H$. Define $\overrightarrow{K}_1$ and $\overrightarrow{K}_2$ to be the subtrails $\overrightarrow{K}((v_x,v_{x+1})$, $(v_{y-1}, v_y))$ and $\overrightarrow{K}((v_y,v_{y+1}),$ $(v_{x-1},v_x))$ of $\overrightarrow{K}_H$, respectively. By Lemma 1.16 (i), switching $X'$ gives two cycles $H_1$ and $H_2$, with  $V(H_1)=V(\overrightarrow{K}_1 \setminus \{v_x, v_y\})$ and $V(H_2)=V(\overrightarrow{K}_2)$. By Proposition 3.1, $Z$ has a vertex in $H_1$ and another in $H_2$, and the same holds for $X'$. Since $X' \in G_3$ and $Z \in G_0 \setminus G_2$, by JCT, $H_1 \cap R_2 \neq \emptyset$. Similarly, $H_2 \cap R_2 \neq \emptyset$. Now the argument in the last paragraph of Lemma 3.7 shows that there must be $v_1 \in H_1 \cap R_2$ and $v_2 \in H_2 \cap R_2$ with $v_1$ adjacent to $v_2$. $\square$

\null 

\noindent \textit{Proof of Lemma 3.12.} For definiteness, assume that $L$ is the northern leaf $R(k,l)$. Let $P(X,Y)$ be the looping $H$-path following $L$, with $X=R(k-1,l+2)$ and $Y=R(k+1,l+2)$. By Lemma 3.13, either there is a cascade that collects $L$, or a cascade after which $P(X,Y)$ gains a switchable box, with each cascade having length at most $\min(m,n)+1$, and both avoiding the $j$-stack of $A_0$'s starting at $C$. If the former, we are done, so may assume that $P(X,Y)$ has a switchable box $Z$. Let $J$ be the large cookie of $G$ and let $N_J$ be the neck of $J$. Note that $N_J$ cannot be a box of $P(X,Y)$. If it were, since the neck has exactly one $H$-neighbour in $G$, $N_J=X$ or $N_J=Y$, but neither box is switchable . Now, $P(X,Y)$ is either contained in $\text{ext}(H)$ or $\text{int}(H)$. 

\null

\noindent \textit{CASE 1}: $P(X,Y) \subseteq \text{ext}(H)$. Then $X' \subset \text{int}(H)$.  By Lemma 3.14, $X' \in G_2$. By Proposition 3.3, $Z \mapsto X'$ is a valid move. By Observation 3.4, $Z \mapsto X'$ does not create additional cookies. Then $Z \mapsto X'$, $L+(0,1) \mapsto L$ is a cascade that collects $L$. End of Case 1. 

\null

\noindent \textit{CASE 2:} $P(X,Y) \subseteq  \text{int}(H)$. Then $X' \subset \text{ext}(H)$.  If $Z \subset G_2$, by Proposition 3.3, $Z \mapsto X'$ is a valid move, and by Observation 3.4, $Z \mapsto X'$ does not create additional cookies. Then $Z \mapsto X'$, $L \mapsto L+(0,1)$ is a cascade that collects $L$.

Suppose then that $Z \subset G_0 \setminus G_2$. By Lemma 3.15, we only need to check the case where $X' \in G_3$. Note that switching $X'$ splits $H$ into two cycles $H_1$ and $H_2$. By Lemma 3.16 there is $v_1 \in H_1 \cap R_2$ and $v_2 \in H_2 \cap R_2$ with $v_1$ adjacent to $v_2$. By Lemma 3.5 there is a cascade $\mu, X' \mapsto W$, or a cascade $X' \mapsto W$ with $X' \neq W$. Note that here, $X'$ plays the role that $Z$ played in Lemma 3.5. Then $\mu, X' \mapsto W, L+(0,1) \mapsto L$ or $X' \mapsto W, L+(0,1) \mapsto L$ is a cascade that collects $L$. End of Case 2.

\null 

\noindent We have just shown that if $P(X,Y)$ has a switchable box, then the cascade required to collect $L$ has length at most three. By Lemma 3.13, the cascade after which $P(X,Y)$ gains a switchable box has length at most $\min(m,n)$. Thus, at most $\min(m,n)+3$ moves are required to collect $L$. $\square$

\null

\subsection{Summary}

In this chapter we proved the MLC and 1LC algorithms. The proof of the MLC algorithm is fully contained here, while the proof of the 1LC algorithm depends on Lemma~3.13, whose proof is given in Chapter~4.

Proposition 3.3 characterizes when double-switch moves are valid and serves as the primary tool for both algorithms. 

The MLC algorithm handles the case where $H$ has multiple large cookies. To collect a large cookie $J$ with switchable neck $N_J$, we look for a switchable box $Z$ in the looping $H$-path of $N_J$. Proposition~3.8 shows that either such a $Z$ already exists, or a there is single preparatory move that produces one.

The 1LC algorithm handles the case where $H$ has exactly one large cookie and at least one small cookie. It collects outermost small cookies. Suppose that $C$ is an outermost small cookie. Either $C$ can be collected immediately by a single move, or $C$ is followed by a $j$-stack of $A_0$-types and an $A_1$-type with switchable middle-box $X'$. If the latter, let $P(X,Y)$ be the $H$-path determined by $X'$. If $P(X,Y)$ contains a switchable box $Z$, then $C$ can be collected by either switching $X'$ directly (if $Z \in G_2$) or by using Lemma~3.5 to find a cascade of length at most two that enables switching $X'$ (if $Z \in G_0 \setminus G_2$). In both cases, a cascade of flips then collects $C$. The existence of such a switchable box $Z$ is guaranteed by Lemma~3.13, whose proof takes up much of Chapter~4. 

\newpage

\section{Fat Paths}

In this chapter we prove Lemma 3.13 and some related results that will be needed for the reconfiguration of Hamiltonian paths in Chapter 5. Recall the setup of Lemma 3.13: we have a Hamiltonian cycle $H$ in an $m \times n$ grid graph $G$ with a $j$-stack of $A_0$-types followed by an $A_1$-type with middle-box $W$. We must show that if the $H$-path $P(X,Y)$ determined by $W$ contains no switchable boxes, then there exists a cascade that produces one. The absence of switchable boxes in $P(X,Y)$ constrains the possible edge configurations of $H$ in the vicinity of $P(X,Y)$, making it highly structured. We will be able to use this structure to find the cascade we seek.

We begin with the necessary definitions. Among these, \textit{turn} and \textit{looping fat path} are used most extensively. We then sketch the proof of Lemma 3.13 and give an outline of the chapter.

\null

\begingroup
\setlength{\intextsep}{0pt}
\setlength{\columnsep}{20pt}
\begin{wrapfigure}[]{r}{0cm}
\begin{adjustbox}{trim=0cm 0cm 0cm 0cm}
\begin{tikzpicture}[scale=1.5]

\begin{scope}[xshift=0cm, yshift=0cm]
{
\draw[gray,very thin, step=0.5cm, opacity=0.5] (0,0) grid (1.5,3.5);

\draw[blue, line width=0.5mm] (0.5,1.5)--++(0.5,0); 

\draw[blue, line width=0.5mm] (0,1)--++(0.5,0)--++(0,-0.5)--++(0.5,0)--++(0,0.5)--++(0.5,0);

\draw[blue, line width=0.5mm] (0,2)--++(0.5,0)--++(0,1); 
\draw[blue, line width=0.5mm] (1.5,2)--++(-0.5,0)--++(0,1);

\begin{scope}
[very thick,decoration={
    markings,
    mark=at position 1 with {\arrow{>}}}
    ]   
    \draw[postaction={decorate}, green!50!black, line width=0.25mm] (0.75,1.5)--++(0,2.25);
\end{scope}

\draw[fill=orange, opacity=1] (0.5, 2) circle [radius=0.05];

\draw[fill=orange, opacity=1] (1, 2) circle [radius=0.05];

% labellings 
{

\node[left] at (0,1.5) [scale=1]
{\tiny{$\ell$}};
\node[left] at (0,2) [scale=1]
{\tiny{+1}};
\node[left] at (0,2.5) [scale=1]
{\tiny{+2}};
{\tiny{+1}};
\node[left] at (0,3) [scale=1]
{\tiny{+3}};

\node at (0.5,-0.15) [scale=1]
{\tiny{$k$}};
\node at (1,-0.15) [scale=1]
{\tiny{+1}};

}

\node[below] at (0.5,2) [scale=1]
{\small{$v_x$}};
\node[below] at (1,2) [scale=1]
{\small{$v_y$}};

\node[right, align=center, text width=3cm] at (1.5,1.75) {Fig. 4.1. Edge followed by a  northern $A_1$-type and a southern $A_0$-type.};

}
\end{scope}

\end{tikzpicture}
\end{adjustbox}
\end{wrapfigure}

\noindent \textbf{Definitions}. Let $G$ be a simply connected polyomino and let $H$ be a Hamiltonian path or cycle of $G$. 
Let $e(k,k+1;l) \in H$. If $H$ has a northern $A$-type on $e(k-1,k;l+1), e(k;l+1,l+2)$ and $e(k+1,k+2;l+1), e(k+1;l+1,l+2)$ then we say that $A$-type \index{follow (a leaf, or an edge)}\textit{follows} the edge $e(k,k+1;l)$ \textit{northward}. Similarly, if  $H$ has a southern $A$-type on $e(k-1,k;l-1), e(k;l-1,l-2)$ and $e(k+1,k+2;l-1), e(k+1;l-1,l-2)$ then we say that $A$-type \textit{follows} the edge $e(k,k+1;l+1)$ \textit{southward}. We make analogous definitions for the other directions. See Figure 4.1.

\null 

% FRESHLY COMMENTED OUT

%\noindent In Chapter 5 we give an algorithm that can be used to reconfigure Hamiltonian paths into one another. The proofs of these algorithms and the 1LC algorithm  essentially rest on being able to find cascades that collect a leaf followed by an $A_0$-type or by an $A_1$-type. Proving that such a cascade exists for a leaf followed by an $A_1$-type is the more involved portion. In this section, we shall lay out some of the elements that are used for those proofs.

\null

\noindent Let \( e(k,k+1;l) \in H \) be followed by a northern \( A_1 \)-type and let $PBR_{\text{North}}(e(k,k+1;l))$ denote the northward \index{perpendicular bisector ray (PBR)|textbf}perpendicular bisector ray $PBR_{\text{North}}(e(k,k+1;l))$ of the edge $e(k,k+1;l)$ originating from the point $v(k+\frac{1}{2},l)$ (green in Figure 4.1).

%%%% ALREADY DEFINED IN SECTION 3.2 DEFINITIONS %%%%%
%We refer to the boxes $R(k,l+1)$ and $R(k,l+2)$ as \textit{the middle-boxes of the $A_1$-type}. Note that the middle-box $R(k,l+1)$ is switchable. We refer to this fact in Lemma 4.1 below.

% right and left boxes of the A_1-type and repeated paragraph
{
%We call the boxes $R(k-1,l+1)$, $R(k-1,l+2)$, $R(k+1,l+1)$ and $R(k+1,l+2)$ the \textit{boxes of the $A_1$-type}
%We define the boxes $R(k-1,l+1)$ and $R(k-1,l+2)$ as \textit{the left boxes}, and $R(k+1,l+1)$ and $R(k+1,l+2)$ as \textit{the right boxes} of the $A_1$-type.

%We refer to the boxes $R(k-1,l+1)$, $R(k-1,l+2)$, $R(k+1,l+1)$, and $R(k+1,l+2)$ as \textit{the boxes of the $A_1$-type}, and the boxes $R(k,l+1)$ and $R(k,l+2)$ as \textit{the middle-boxes of the $A_1$-type}. More specifically, we define $R(k-1,l+1)$ and $R(k-1,l+2)$ as \textit{the left boxes}, and $R(k+1,l+1)$ and $R(k+1,l+2)$ as \textit{the right boxes} of the $A_1$-type. Note that the middle-box $R(k,l+1)$ is switchable, as referenced in Lemma 4.1 below.

}

We call the vertices  $v(k,l+1)$ and  $v(k+1,l+1)$ \index{corner vertex of an A1 type@corner vertex of an $A_1$-type|textbf}\textit{the corner vertices of $A_1$}. More specifically, we call the vertex $v(k,l+1)$ of $A_1$, located on the left of $PBR_{\text{North}}(e(k,k+1;l))$ the \index{left corner of an A1 type@left corner of an $A_1$-type|textbf}\textit{left corner of $A_1$} and denote it by $v_{\text{left}}(A_1)$. Similarly, we call the vertex of $v(k+1,l+1)$ of $A_1$ located on the right of $PBR_{\text{North}}(e(k,k+1;l))$ the \index{right corner of an A1 type@right corner of an $A_1$-type|textbf}\textit{right corner of $A_1$}, and denote it by $v_{\text{right}}(A_1)$ (orange in Figure 4.1). 

We call the edges $e(k;l+1,l+2)$ and $e(k;l+2,l+3)$ the  \index{left colinear edges of an A1 type@left colinear edges of an $A_1$-type|textbf}\textit{left colinear edges of the $A_1$-type}. Similarly, the edges $e(k+1;l+1,l+2)$ and $e(k+1;l+2,l+3)$ are called the  \index{right colinear edges of an A1 type@right colinear edges of an $A_1$-type|textbf}\textit{right colinear edges of the $A_1$-type}. We refer to these in Proposition 4.14.
We make analogous definitions for the other directions. See Figure 4.1. Note that if the $A_1$-type is southern, then the left corner of the $A_1$-type is on the right of the image in our page; and similarly if the $A_1$-type is western (eastern), then the left corner of the $A_1$-type is south (north) of the perpendicular bisector ray. 

Suppose that the edge $e(k,k+1,l)$ is followed by an $A_1$-type southward. Define an  \index{A1-partitioning of a Hamiltonian path@$A_1$-partitioning of a Hamiltonian path|textbf}\textit{$A_1$-partitioning of $H$} to be the partition of $H$ into the sub-paths $P_0=P(v_x,v_y)$, $P_1=P(v_1,v_x)$ and $P_2=P(v_y,v_r)$, where $v_x = v_{\text{left}}(A_1)$ and $v_y=v_{\text{right}}(A_1)$ or $v_x = v_{\text{right}}(A_1)$ and $v_y=v_{\text{left}}(A_1)$.

Let $W=R(k,l-2)$, $X=R(k+1,l-2)$, and $Y=R(k-1,l-2)$, and assume that the looping $H$-path $P(X,Y)$ of $W$ is contained in an $H$-component of $G$. Then we say that $P(X,Y)$ is a \index{the looping H-path of a box@the looping $H$-path of a box}\textit{southern looping $H$-path} and that $P(X,Y)$ \index{follow (a leaf, or an edge)}\textit{follows} $e(k,k+1,l)$ southward. Similarly, if $R(k,l)$ is a southern leaf, we say that $P(X,Y)$ \textit{follows} $R(k,l)$ southward. We make analogous definitions for western, northern and eastern looping $H$-paths. From here on, assume that every looping $H$-path follows an edge of $H$ and that its end-boxes are incident on an $A_1$-type, as described above.

% The purpose of this section is to prove the following: given any looping $H$-path of a switchable box $W$ following some edge of $H$, there is a cascade, after the execution of which, there is a nontrivial valid move $W \mapsto W'$

\noindent Let $J=\{X_1, X_2, \ldots, X_r\}$ be a collection of boxes of a polyomino $G$ and let $H$ be a Hamiltonian path or cycle of $G$. We will use the notation $G\langle J \rangle$ to denote the subgraph of $G$ with vertex set $V(G\langle J \rangle)=V(J)$ and edge set $E(G\langle J \rangle)=E(J)\cap E(H)$. The boxes of $G\langle J \rangle$ are the boxes of $J$. We call $G\langle J \rangle$ the \textit{subgraph of $G$ induced by $J$}.

\null

\begingroup
\setlength{\intextsep}{0pt}
\setlength{\columnsep}{20pt}
\begin{wrapfigure}[]{r}{0cm}
\begin{adjustbox}{trim=0cm 0cm 0 0cm}
\begin{tikzpicture}[scale=1.5]

\begin{scope}[xshift=0]
\draw[gray,very thin, step=0.5cm, opacity=0.5] (0,0) grid (2.5,2.5);
%%%%%%%%%%%%    FP edges and sFP boxes   %%%%%%%%%%%%  

\fill[blue!50!white, opacity=0.5] (0.0,1.5)--++(0.5,0)--++(0,0.5)--++(0.5,0)--++(0,-1)--++(0.5,0)--++(0,1)--++(0.5,0)--++(0,-0.5)--++(0.5,0)--++(0,-0.5)--++(-0.5,0)--++(0,-0.5)--++(-0.5,0)--++(0,-0.5)--++(-0.5,0)--++(0,0.5)--++(-0.5,0)--++(0,0.5)--++(-0.5,0);

% inner left
\draw[blue, line width=0.5mm] (1.0,2.5)--++(0.5,0);

% inner center
\draw[blue, line width=0.5mm] (0.0,1.5)--++(0.5,0)--++(0,0.5)--++(0.5,0)--++(0,-1)--++(0.5,0)--++(0,1)--++(0.5,0)--++(0,-0.5)--++(0.5,0)--++(0,-0.5)--++(-0.5,0)--++(0,-0.5)--++(-0.5,0)--++(0,-0.5);

\draw[blue, line width=0.5mm] (0.0,1.0)--++(0.5,0)--++(0,-0.5)--++(0.5,0)--++(0,-0.5);

% red traced path
\draw[red, line width=0.25mm, dotted](0.75,1.75)--++(0,-1)--++(1,0)--++(0,1);

\node[above] at (1.55,1.95) [scale=1]
{\small{$v_x$}};
\node[above] at (0.95,1.95) [scale=1]
{\small{$v_y$}};

\node[above] at (1.25,2.45) [scale=1]
{\small{$e$}};

\node[above] at (1.25,1.95) [scale=1]
{\small{$e'$}};

\draw[fill=blue, opacity=1] (1, 2) circle [radius=0.05];

\draw[fill=blue, opacity=1] (1.5, 2) circle [radius=0.05];

% labeling

\node[below, align=center, text width=6cm] at (1.25,0) {Fig. 4.2. A southern looping fat path $G\langle N[P(X,Y] \rangle$. $N[P(X,Y]$ shaded in light blue. $P(X,Y)$ traced in red. };

{

\node at  (1.75, 1.75) [scale=0.8]{X};
\node at  (0.75, 1.75) [scale=0.8]{Y};
\node at  (1.25, 1.75) [scale=0.8]{W};

%\node[left] at (-0.5,2.5) [scale=1]{\tiny{$\ell$}};
%\node[above] at (1.0,3.0) [scale=1]{\tiny{$k$}};
}

\end{scope}

\end{tikzpicture}
\end{adjustbox}
\end{wrapfigure}

\noindent Let $e$ be an edge of $H$ followed by an $A_1$-type with switchable middle-box $W$. Let $X$ and $Y$ be the boxes adjacent to $W$ that are not its $H$-neighbours, and assume that $X$ and $Y$ belong to the same $H$-component. We will call the set of boxes consisting of the boxes of the $H$-path $P(X,Y)$ together with all the $H$-neighbours of $P(X,Y)$  \index{H-neighbourhood@$H$-neighbourhood|textbf}\textit{the H-neighbourhood of $P(X,Y)$} and denote it by $N[P(X,Y)]$. Consider the subgraph $F=G\langle N[P(X,Y)] \rangle$ of $G$ induced by $N[P(X,Y)]$. We call $F$ a \index{looping fat path|textbf}\textit{looping fat path} if:

1. no end-vertex of $H$ is incident on $P(X,Y)$, 

2. $W$ is anti-parallel, and

3. $P(X,Y)$ has no switchable boxes.

\noindent If furthermore, $F$ is such that:

4. There is no switch or flip or transpose (recall definitions from Introduction) move after which $W$ is switched, or $P(X,Y)$ gains a switchable   box, or $W\mapsto X$ or $W\mapsto Y$ is valid,

\noindent then we say that $F$ is a \index{sturdy looping fat path|textbf}\textit{sturdy looping fat path}. We call these conditions the \textit{fat path conditions} (FPC), and say FPC-j when referring to condition j, for $j\in \{1,\ldots,4\}$. The role of FPC-1 in handling end-vertices is addressed in Sections 4.3 ( Remark after Observation 4.20) and later sections. Roughly, our goal will be to find a cascade after which the edge $e'=\{v_x,v_y\}$ (Figure 4.2.) is in the resulting Hamiltonian path. If any of the conditions FPC-j fail, we can find such a cascade immediately. In that sense, a sturdy looping fat path represents an unyielding configuration -- the kind of subgraph where there are no obvious good moves, and we have to look further into structure to obtain the required cascade.

\endgroup

\null

\noindent Recall the definition of an outermost small cookie from the beginning of Section 3.2. Suppose $G$ has only one large cookie, and that there is a $j$-stack of $A_0$'s starting at the $A_0$-type containing an outermost southern small cookie $C$. Let $L$ be the leaf in the top ($j^{\text{th}}$) $A_0$ of the stack, and assume that $L$ is followed by an $A_1$-type with looping $H$-path $P(X,Y)$. Let $F = G\langle N[P(X,Y)] \rangle$, and assume that $F$ is a fat path. We say that $F$ is a \textit{southern} looping fat path \index{anchored|textbf}\textit{anchored} at the outermost southern small cookie $C$. Analogous definitions apply for northern, eastern, and western looping fat paths.

From here on, whenever $H$ is a cycle and we are discussing a looping fat path $F$, we will assume that $G$ has exactly one large cookie and that $F$ is anchored at an \index{outermost small cookie}outermost small cookie.

Assume that $H$ is a Hamiltonian path of $G$, and that $F = G\langle N[P(X,Y)] \rangle$ is a looping fat path such that $P(X,Y)$ follows an edge of $H$ southward. Then we say that $F$ is a \textit{southern} looping fat path. Analogous definitions apply for northern, eastern, and western looping fat paths.

\null 

\noindent Let $H$ be a Hamiltonian path of an $m \times n$ grid graph $G$. In Chapter 5 we will define an appropriate canonical form and show that $H$ can be reconfigured into such a form. If $H$ has adjacent end-vertices, the reconfiguration process shortens and simplifies. The more involved cases, requiring some of the tools we build in this chapter, arise when $H$ does not have adjacent end-vertices. Therefore, from here on, we will often assume that the Hamiltonian path we are considering does not have adjacent end-vertices. We call this assumption the \index{non-adjacency assumption (NAA)|textbf}non-adjacency assumption (NAA).

\null 

\noindent We define below a subgraph of $G$ consisting of the union of translations of two adjacent and perpendicular edges of $G$. Let $r \in \mathbb{N}$ and let \index{stairs|textbf}\textit{the stairs from (k,l) to (k+r,l-r) east} be denoted by $S_{\rightarrow}(k,l;k+r,l-r)$ and be defined as:

\setlength{\abovedisplayskip}{0pt}
\setlength{\belowdisplayskip}{10pt}

\[S_{\rightarrow}(k,l;k+r,l-r) = \bigcup_{j=0}^{r-1} \Big( e(k,k+1;l)+(j,-j) \Big )\cup \Big( e(k+1;l-1,l)+(j,-j) \Big).\]

\noindent We define $d(S)=r$ to be the \textit{length} of $S_{\rightarrow}(k,l;k+r,l-r)$. We say that $S_{\rightarrow}(k,l;k+r,l-r)$ starts at $v(k,l)$ and ends at $v(k+r,l-r)$. The subscripted arrow indicates the direction from $v(k,l)$ of the first edge of the subgraph. By choosing an ``up", ``down", ``left" or ``right" arrow for direction and a sign for the third and fourth arguments of $S_{\square}(k,l;k\pm r,l \pm r)$ we may describe any of the eight possible steps subgraphs starting at the vertex $v(k,l)$. See Figure 4.3 (a).

\null

\begingroup
\setlength{\intextsep}{0pt}
\setlength{\columnsep}{20pt}
\begin{wrapfigure}[]{r}{0cm}
\begin{adjustbox}{trim=0cm 0cm 0cm 1cm}
\begin{tikzpicture}[scale=1.5]
\begin{scope}[xshift=0cm, yshift=0cm]

\draw[gray,very thin, step=0.5cm, opacity=0.5] (0,0) grid (1.5,1.5);

%%%%%%%%%% Integer Coordinates %%%%%%%%
{
\node[right] at (1.5,1.5) [scale=1]{\tiny{$\ell$}};

\node[above] at (0,1.5) [scale=1]{\tiny{$k$}};
}
%%%%%%%%%%%%     DFP-0   %%%%%%%%%%%%  
{
\draw[blue, line width=0.5mm] (0,1.5)--++(0.5,0)--++(0,-0.5)--++(0.5,0)--++(0,-0.5)--++(0.5,0)--++(0,-0.5);

\fill[blue!40!white, opacity=0.5] (0,1.5)--++(0.5,0)--++(0,-0.5)--++(0.5,0)--++(0,-0.5)--++(0.5,0)--++(0,-0.5)--++(-1,0)--++(0,0.5)--++(-0.5,0);
}

\node[below, align=center, text width=7cm] at (2,0) {Fig. 4.3. (a) $S_{\rightarrow}(k,l;$  $k+3,l-3)$. (b) A half-open northeastern turn.};

\node[above] at (0.75,1.5) [scale=1.25] {(a)};

\end{scope}

\begin{scope}[xshift=2.75cm, yshift=0cm]

\draw[gray,very thin, step=0.5cm, opacity=0.5] (0,0) grid (1.5,1.5);

%%%%%%%%%% Integer Coordinates %%%%%%%%
{
\node[left] at (0,1.5) [scale=1]{\tiny{$\ell$}};
\node[above] at (0.0,1.5) [scale=1]{\tiny{$k$}};

\node[left] at (0,0.0) [scale=1]{\tiny{$\ell'$}};
\node[above] at (1.5,1.5) [scale=1]{\tiny{$k'$}};
}

%%%%%%%%%%%%     DFP-0   %%%%%%%%%%%%  
{
\draw[blue, line width=0.5mm] (0.5,1.5)--++(0,-0.5)--++(0.5,0)--++(0,-0.5)--++(0.5,0)--++(0,-0.5)--++(-0.5,0);
\draw[blue, line width=0.5mm] (0.0,1.0)--++(0,0.5);

\fill[blue!40!white, opacity=0.5] 
  (0.0,1.5)--++(0.5,0)--++(0,-0.5)--++(0.5,0)--++(0,-0.5)--++(0.5,0)--++(0,-0.5)--++(-1,0)--++(0,0.5)--++(-0.5,0);
}

\node[above] at (0.75,1.5) [scale=1.25] {(b)};

\end{scope}

\end{tikzpicture}
\end{adjustbox}
\end{wrapfigure}

\noindent  Let $H$ be a Hamiltonian path or cycle of a polyomino $G$. Let $T$ be the subgraph of $H$ on the edges $S_{\downarrow}(k+1,l;k',l'+1)$, $e(k;l-1,l)$ and $e(k'-1,k';l')$, where $k'=k+d(T)$, $l'=l-d(T)$, where $d(T)=d(S)+1$ is \textit{the length of T} and $d(T) \geq 2$. We call $T$ a \index{turn|textbf}\textit{northeastern turn}. If both $e(k,k+1;l)$ and $e(k';l',l'+1)$ belong to $G \setminus H$, call $T$ an \textit{open northeastern turn}. If exactly one of $e(k,k+1;l)$ and $e(k';l',l'+1)$ is in $H$, then $T$ is a \textit{half-open northeastern turn}. See Figure 4.3 (b). If both $e(k,k+1;l)\in H$ and $e(k';l',l'+1) \in H$, then $T$ is a \textit{closed northeastern turn}. For any northeastern turn $T$, we say that $R(k,l-1)$ is the \index{leaf of a turn|textbf}\textit{northern leaf} of $T$ and $R(k'-1,l')$ is the \textit{eastern leaf} of $T$. If $e(k,k+1;l)\notin H$ we call $R(k,l-1)$ an \textit{open northern leaf of T} and if $e(k,k+1;l)\in H$ we call $R(k,l-1)$ a \textit{closed northern leaf of T}. We note that the two leaves of a turn will determine its ``leaf prefix": If a turn has a northern leaf and an eastern leaf then the turn is a northeastern turn.

We will say that a looping fat path $F$ \textit{has a turn} (open, half-open or closed) to mean that  there exists some turn $T$ of $H$ such that $E(F) \supset E(T)$.

\null

\noindent \textbf{Sketch of proof of Lemma 3.13.} Let $H$ be a Hamiltonian cycle of an $m \times n$ grid graph $G$. Assume that $P(X,Y)$ is a looping $H$-path with no switchable boxes, following a leaf $L$. It follows that $P(X,Y)$ is contained in a \index{looping fat path}looping fat path $F$. In Section 4.2, we show that every looping fat path has a turn. In Sections 4.3-4.6 we show that given a turn, we can find a cascade we call a \index{weakening}\textit{weakening} (precise definition in  Section 4.3) that collects one of its leaves. It then follows that after this cascade, either $P(X,Y)$ gains a switchable box, or we can extend the cascade by a single move to collect $L$.

\null

\noindent The rest of the chapter is organized as follows. Section 4.1 proves structural properties of fat paths, which the later sections build on. Section 4.2 shows that every fat path contains a turn (Proposition 4.18). Section 4.3 defines weakenings and proves some basic properties. Sections 4.4-4.6 prove that turns have weakenings: Section 4.4 handles the case where $H$ is a Hamiltonian cycle (Proposition 4.24), while Sections 4.5 and 4.6 handle the case where $H$ is a Hamiltonian path (Propositions 4.27 and 4.28). The proof of Lemma 3.13 appears at the end of Section 4.6. Before we proceed to Section 4.1 we prove results 4.1-4.5, which generalize Propositions 3.1 and 3.3 to Hamiltonian paths in polyominoes.

\null 

\noindent \textbf{Lemma 4.1.} Let $G$ be a polyomino and let $H=v_1, \ldots, v_r$ be a Hamiltonian path of $G$. Let $P(X,Y)$ be the looping $H$-path of a switchable box $W$. Assume that $P(X,Y)$ is contained in an $H$-component of $G$ and that the edges of $W$ are anti-parallel. Then the $A_1$-partitioning of $H$ is such that every box of $P(X,Y)$ is incident on a vertex of $P_0$.

\null 
 
\noindent \textit{Proof.} For definiteness, assume that  $P(X,Y)$ is a southern looping $H$ path following the edge $e(k,k+1,l)$ southward, that $v_x=v_{\text{left}}(A_1)$, and $v_y=v_{\text{right}}(A_1)$. Then $v_x=v(k+1,l-1)$ and $v_y=v(k,l-1)$. Now,  $v_{x+1}=v(k+1,l-2)$ or $v_{x+1}=v(k+2,l-1)$.

Assume that $v_{x+1}=v(k+1,l-2)$. Then $v_{y-1}=v(k,l-2)$. Then $\Phi((v_x,v_{x+1}), \text{left})=X$ and $\Phi((v_{y-1},v_y), \text{left})=Y$. By Observation 1.2.4 (d), we can identify $\overrightarrow{K}((v_x,v_{x+1}),(v_{y-1},v_y))$ with $P_0$. Since the end-boxes $X$ and $Y$ of $P(X,Y)$ are contained in $G$, it follows from Proposition 1.2.1 that $P(X,Y)$ is not contained in any cycle of boxes in $G_{-1}$. Thus, $P(X,Y)$ is unique. It follows that any $H$-walk of boxes between $X$ and $Y$ contains $P(X,Y)$. In particular, $\Phi(P_0,\text{left})$ contains $P(X,Y)$. Then, by definition of FTW, any box of $P(X,Y)$ is incident on a vertex of $P_0$. 

The case where $v_{x+1}=v(k+2,l-1)$ uses a similar argument, so we omit the proof. $\square$

\null 

\noindent \textbf{Corollary 4.2.} Let $G$, $H$, $P(X,Y)$, $P_0$, $P_1$, and $P_2$ be as in Lemma 4.1. Then:

(a) For every box $Z$ of $P(X,Y)$ there is an edge $e_z \in P_0$, such that $e_z$ adds 

%and such that $Z$ is on the left side of $e_z$. 

\hspace{0.6cm} $Z$ to $\Phi(P_0, \text{left})$. 

(b) The middle-boxes of the $A_1$-type following $e(k,k+1,l)$ are not in $P(X,Y)$.

\null 

\noindent \textit{Proof.} Part (a) follows from the definition of FTW. It remains to prove part (b). Either $v_{x+1}=v(k+1,l-2)$ or $v_{x+1}=v(k+2,l-1)$. Assume that $v_{x+1}=v(k+1,l-2)$. For a contradiction, assume that $X+(-1,0)$ is in $P(X,Y)$. Then $X+(-1,0)$ must be added to $\Phi(P_0, \text{left})$ by $(v_{x+1},v_x)$ or by $(v_y,v_{y-1})$, but neither edge belongs to $P_0$. Similarly if $X+(-1,-1)$ were in $P(X,Y)$, it would have to be added to $\Phi(P_0, \text{left})$ by $(v_{x+2},v_{x+1})$ or by $(v_{y-1},v_{y-2})$, or, if $\{v_{y-2}, v_{x+2}\} \in H$, by $(v_{y-2}, v_{x+2})$ but none of those edge belong to $P_0$.

\noindent The case where $v_{x+1}=v(k+2,l-1)$ uses a similar argument, so we omit the proof. $\square$.

\null

\noindent \textbf{Lemma 4.3.}  Let $G$ be a polyomino and let $H=v_1, \ldots, v_r$ be a Hamiltonian path of $G$. Let $P(X,Y)$ be a looping $H$-path. Let $P_0$, $P_1$ and $P_2$ be the partitioning of $H$ given in Lemma 4.1. Then, no edge of $P_0$ has both of its incident boxes belonging to $P(X,Y)$.

\null

\noindent \textit{Proof.} For contradiction, assume that there is an edge $e$ of $P_0$ with both its incident boxes $Z$ and $Z'$ belonging to $P(X,Y)$. Let $Q$ be the cycle consisting of the edges of $P_0$ and the edge $\{v_x,v_y\}$ of $G \setminus H$. Consider the $H$-subpath $P(Z_1,Z_s)$ of $P(X,Y)$, where $Z_1=Z$ and $Z_s=Z'$. Let $c_j$ be the center of the box $Z_j$, for each $j \in  \{1, \ldots, s\}$. Note that each segment $[c_j,c_{j+1}]$ intersects one gluing edge of $P(Z_1,Z_s)$ and no other edge or vertex of $G$, and that $[c_1,c_2], \ldots, [c_{s-1}, c_s]$ is a path $P(c_1,c_s)$. By Corollary 1.1.5, $c_1$ and $c_s$ are on different sides of $Q$. By JCT, $P(c_1,c_s)$ intersects $Q$. Since the edges of $P(c_1,c_s)$ only intersects gluing edges of $P(Z_1,Z_s)$, some edge $\{c_i,c_{i+1}\}$ of $P(c_1,c_s)$, must intersect $Q$ at $\{v_x,v_y\}$. Now, $Z_i$ and $Z_{i+1}$ belong to $P(X,Y)$, and, since they are incident on $\{v_x,v_y\}$, at least one of them must be a middle-box of $F$, which contradicts Corollary 4.2 (b). $\square$

\null

\noindent \textbf{Proposition 4.4.} Let $G$ be a polyomino and let $H=v_1, \ldots, v_r$ be a Hamiltonian path of $G$. Let $P(X,Y)$ be a looping $H$-path in $G$ and assume that the switchable box between $X$ and $Y$ is anti-parallel. Then the partitioning of $H$ into the subpaths $P_0=P(v_x,v_y)$, $P_1=P(v_1,v_x)$ and $P_2=P(v_y,v_r)$, where $v_x$ and $v_y$ are the corners of the $A_1$-type, is such that every box of $P(X,Y)$ is incident on a vertex of $P_0$ and a vertex of $P_1$ or $P_2$.

\null 

\noindent \textit{Proof.} Let $v_x$ be incident on $X$ and $v_y$ be incident on $Y$. By Lemma 4.1, every box of $P(X,Y)$ has a vertex incident on $P_0$. To finish the proof, we proceed by induction. Note that, since $v_x \in V(P_0) \cap V(P_1)$ and $v_y \in V(P_0) \cap V(P_2)$,  $X$ and $Y$ satisfy the conclusion. Let $P(X,Y)$ be the path $X=X_1, \ldots, X_t=Y$. Assume that the box $X_j$ is incident on a vertex in $P_0$ and a vertex of $P_1$ or $P_2$. Let $\{v_p, v_q\}$ be the gluing edge between $X_j$ and $X_{j+1}$. By Lemma 4.1, at least one vertex of $X_{j+1}$ belongs to $V(P_0)$. For a contradiction, assume that no other vertex of $X_{j+1}$ belongs to $V(P_0)$.

Without loss of generality we may assume that $p<q$. Note that if $x=p$ or $y=q$ then $X_{j+1}$ satisfies the conclusion of the proposition, so we may assume that $x<p$ and $q<y$. Then we have that $1<x<p<q<y<r$. Let $U$ be the region of $G$ bounded by the polygon $Q$ consisting of the subpath $P(v_p,v_q)$ of $P_0$ and the edge $\{v_p, v_q\}$ of $G \setminus H$. Note that $X_j$ and $X_{j+1}$ are on distinct sides of $Q$. Assume for now that this implies that $X$ and $Y$ are on distinct sides of $Q$ as well. Then $v_x$ and $v_y$ must be on distinct sides of $Q$ ($v_x,v_y$ cannot be in $Q$ since $x$ is too small, and $y$ is too large. Therefore, if $X$ and $Y$ are on distinct sides of $Q$, so must $v_x$ and $v_y$). But this contradicts the fact that $v_x$ and $v_y$ are adjacent. It remains to show that $X$ and $Y$ are on distinct sides of $Q$. Either $X_{j+1} \in U$ or $X_{j+1} \notin U$.

\null 

\noindent \textit{CASE 1: $X_{j+1} \in U$.} We check that $Y \in U$. For a contradiction, assume $Y \notin U$. Consider the $H$-path $P(X_{j+1},Y)$. For each $i \in \{j+1, j+2,\ldots,t\}$, let $c_i$ be the center of the boxes $X_i$. Observe that for each $i \in \{j+1, j+2,\ldots,t\}$, the segment $[c_i,c_{i+1}]$ intersects the gluing edge of $X_i$ and $X_{i+1}$ and $[c_i,c_{i+1}]$ intersects no other edge of $G$. By JCT, any path between a point in $U$ and a point in $G \setminus U$ intersects $Q$. Since $X_{j+1} \in U$ and $Y\notin U$, we have that $c_{j+1} \in U$ and $c_t \notin U$. In particular, the path $P(c_{j+1}, c_t)$ intersects $Q$. By the ordering of $P(X,Y)$,  $X_{j} \notin P(X_{j+1},Y)$. Then the segment $[c_j,c_{j+1}]$ is not in $P(c_{j+1}, c_t)$, so  $P(c_{j+1}, c_t)$ does not intersect $Q$ at the edge $\{v_p, v_q\}$. It follows that $P(c_{j+1}, c_t)$ intersects $Q$ at some other edge $e$ of $P(v_p,v_q)$. But then $e$ is not a gluing edge of $P(X_{j+1},Y)$, contradicting our observation above. So we must have that $Y \in U$. 

A similar argument can be used to show that $X$ must belong to $G \setminus U$. End of Case 1.

\null 

\noindent \textit{CASE 2:  $X_{j+1} \notin U$.} A similar argument can be used to show that $X$ and $Y$ are on distinct sides of $Q$. End of Case 2. $\square$. 

\null

\noindent \textbf{Proposition 4.5.} Let $H$ be a Hamiltonian path of a polyomino $G$ with end-vertices $u$ and $v$, let $W$ be a switchable box in $H$, and let $P(X,Y)$ be the looping $H$-path of $W$, with $v_x=v_{\text{left}}(A_1)$ and $v_y=v_{\text{right}}(A_1)$. Assume that $P(X,Y)$ is contained in an $H$-component of $G$, and that $W$ is anti-parallel. If $P(X,Y)$ has a switchable box $Z$ then either $\textrm{Sw}(Z)$, $\textrm{Sw}(W)$ is a cascade, or $Z \mapsto W$ is a valid move.

\null

\noindent \textit{Proof.} Without loss of generality, let $x<y$. By Lemma 4.1, we may partition $H$ into the subpaths $P_0=P(v_x,v_y)$, $P_1=P(v_1,v_x)$ and $P_2=P(v_y,v_r)$ such that every box of $P(X,Y)$ is incident on a vertex of $P_0$ and a vertex of $P_1$ or $P_2$. Then the edges of $W$ in $H$ are either $(v_x, v_{x+1})$ and $(v_{y-1}, v_y)$, or  $(v_{x-1}, v_x)$ and $(v_y, v_{y+1})$.

\null 

\noindent \textit{CASE 1: The edges of $W$ in $H$ are $(v_x, v_{x+1})$ and $(v_{y-1}, v_y)$. } Suppose that $P(X,Y)$ has a switchable box $Z$. By Corollary 4.2 (b), $Z \neq W$. Since $Z$ is switchable, $Z$ is not any of the boxes incident on $v_x$ or $v_y$. Let $(v_s,v_{s+1})$ and $(v_t, v_{t+1})$ be the edges of $Z$ in $H$. Then exactly one of $(v_s,v_{s+1})$ and $(v_t, v_{t+1})$ is in $P_0$ and the other is in $P_1$ or $P_2$. WLOG assume that $(v_s,v_{s+1})$ is in $P_0$ and that $(v_t, v_{t+1})$ is in $P_2$. Then we can partition the edges of $H$ as follows: $P(v_1,v_s)$, $(v_s,v_{s+1})$, $P(v_{s+1},v_t)$, $(v_t,v_{t+1})$, $P(v_{t+1},v_r)$, where $1<x<s<y<t<r$. Now, $(v_s,v_{s+1})$ and $(v_t, v_{t+1})$ are either parallel or anti-parallel.

\null

\begin{center}
\begin{adjustbox}{trim=0cm 0cm 0cm 0.5cm}

\end{adjustbox}
\end{center}

\noindent \textit{CASE 1.1: $Z$ is anti-parallel.} We show that $Z\mapsto W$ is a valid move. After removing the edges $(v_s,v_{s+1})$ and $(v_t,v_{t+1})$ we are left with three paths: $P(v_1,v_s), P(v_{s+1},v_t)$ and $  P(v_{t+1},v_r)$. Note that adding the edge $\{v_{s+1}, v_t\}$ gives a cycle $H_c$ consisting of the path $P(v_{s+1},v_t)$ and the edge $\{v_{s+1}, v_t\}$, and adding the edge $\{v_s, v_{t+1}\}$ gives a path $H_p=P(v_1,v_s), \{v_s,v_{t+1}\}, P(v_{t+1},v_r)$. Now $1<x<s$ implies that $e_x \in H_p$ and $s<y<t$  implies that $e_y \in H_c$. Thus we have that $W$ is now an $(H_c,H_p)$-port \index{H p H c port@$(H_p,H_c)$-port}. By Lemma 1.4.1 (iii), $Z \mapsto W$ is a valid move. See Figure 4.4 (a). End of Case 1.1.

\endgroup

\null

\noindent \textit{CASE 1.2: $Z$ is parallel.} We show that $\textrm{Sw}(Z)$, $\textrm{Sw}(W)$ is a valid sequence of moves. Since $Z$ is parallel, by Lemma 1.4.1(i), Sw$(Z)$ is a valid move. After the removal of the edges $(v_s,v_{s+1})$ and $(v_t,v_{t+1})$ we are left with three paths: $P(v_1,v_s), P(v_{s+1},v_t)$ and $  P(v_{t+1},v_r)$. Adding $\{v_s,v_t\}$, $\{v_{s+1}, v_{t+1}\}$ gives a new Hamiltonian path $H'=v_1, \ldots, v_x, \ldots, v_s,v_t, \ldots, v_y,$ $ \ldots, v_{s+1}, v_{t+1}$, $ \ldots, v_r$. Note that now the edges $(v_x,v_{x+1})$ and $(v_y,v_{y+1})$ of $H'$ are parallel so $\textrm{Sw}(W)$ is a valid move. See Figure 4.4 (b). End of Case 1.2. End of Case 1.

\null

\noindent \textit{CASE 2: The edges of $W$ in $H$ are $(v_{x-1}, v_x)$ and $(v_y, v_{y+1})$.} This follows by an argument similar to that of Case 1, so we omit the proof. End of Case 2. $\square$.

\endgroup

\subsection{Properties of looping fat paths}

Section 4.1 derives and formalizes the structural properties of looping fat paths. There are three important results here: Lemmas 4.6 and 4.13, and Proposition 4.14. Lemma 4.6 gives some preliminary properties. Lemma 4.13 shows that the boundary $B(F)$ of a looping fat path $F$ is a Hamiltonian cycle of $F$. It builds on the facts that looping fat paths are non-self-adjacent (Proposition 4.7) and have no polyking junctions (Proposition 4.11). Proposition 4.14 shows that fat paths have no consecutive colinear edges other than those belonging to the $A_1$-type configuration. These three results are used in Section 4.2, where we scan looping fat paths for turns.

\null

\noindent \textbf{Lemma 4.6.} Let $F=G\langle N[P(X,Y)] \rangle$ be a \index{looping fat path}looping fat path. Suppose that $W=R(k,l)$ in $P$ has an $H$-neighbour $Z$ southward in $N[P] \setminus P$. Then:

(a) $Z$ has exactly one $H$-neighbour in $P$ and $W$ has no other $H$-neighbour in $N[P] \setminus P$.

(b) If $W$ is not an end-box of $P$, then the $H$-neighbours of $W$ in $P$ are $W+(-1,0)$ 

\hspace{0.6 cm} and $W+(1,0)$.  Furthermore, $S_{\rightarrow}(k-1,l;k,l-1) \in H$, $S_{\uparrow}(k+1,l-1;k+2,l) \in H$, 

\hspace{0.6 cm} and $e(k,k+1;l+1)\in H$.

(c) If $W$ is an end-box of $P$, then $e(k,k+1;l+1) \in H$ and exactly one of $e(k;l,l+1)$ 

\hspace{0.6 cm} and  $e(k+1;l,l+1)$ belong to $H$.

(d) $Z$ is a leaf or $Z$ is a switchable box in $H$.

\noindent Analogous statements apply when $Z$ is west, north or east of $W$.

\null 

\noindent \textit{Proof of (a).} Note that if $Z$ has more than one $H$-neighbour in $P$ then we can make an $H$-cycle, which contradicts Proposition 1.2.1. To see that $W$ has no other $H$-neighbour in $N[P] \setminus P$, assume, BWOC, that $W$ has at least two $H$-neighbours in $N[P] \setminus P$. 

If $W$ is an end-box of $P$, then, by definition of $A_1$, $W$ has at most two $H$-neighbours, and at least one of them must belong to $P$, contradicting our assumption that  $W$ has at least two $H$-neighbours in $N[P] \setminus P$.  

If $W$ is not an end-box, then $W$ must have four $H$ neighbours: two in $N[P] \setminus P$ and two in $P$. By definition of a looping fat path, at least one neighbour of $W$, say $W'$, is not an end-box. But then $W'$ must be switchable, which conflicts with the definition of a looping fat path. End of proof for (a).

\null 

\begingroup 
\setlength{\intextsep}{0pt}
\setlength{\columnsep}{18pt}
\begin{wrapfigure}[]{l}{0cm}
\begin{adjustbox}{trim=0cm 0.5cm 0cm 0.25cm}

\end{adjustbox}
\end{wrapfigure}

\noindent \textit{Proof of (b).} First, we show that the $H$-neighbours of $W$ are $W+(1,0)$ and $W+(-1,0)$. BWOC, assume that $W+(1,0)$ is not an $H$-neighbour of $W$ (See fig). Then the $H$-neighbours of $W$ in $P$ must be $W+(-1,0)$ and $W+(0,1)$. It follows that $S_{\rightarrow}(k-1,l;k,l-1)\in H$ and $S_{\rightarrow}(k-1,l+1;k,l+2)\in H$. Note that, by definition of $A_1$ and looping fat paths, $W+(-1,0)$ is not an end-box of $P$. But then $W+(-1,0)$ is a switchable box of $P$, which

\endgroup

\noindent conflicts with the definition of a looping fat path. Therefore the $H$-neighbours of $W$ in $P$ are $W+(-1,0)$ and $W+(1,0)$. It follows that $S_{\rightarrow}(k-1,l;k,l-1) \in H$ and $S_{\uparrow}(k+1,l-1;k+2,l) \in H$, and by part (a), $(k,k+1;l+1)\in H$. See Figure 4.5. End of proof for (b).

\begingroup 
\setlength{\intextsep}{0pt}
\setlength{\columnsep}{18pt}
\begin{wrapfigure}[]{r}{0cm}
\begin{adjustbox}{trim=0cm 0cm 0cm 0cm}
\begin{tikzpicture}[scale=1.5]
\begin{scope}[xshift=0cm]{
\draw[gray,very thin, step=0.5cm, opacity=0.5] (0,0) grid (1.5,1.5);

\draw[blue, line width=0.5mm] (0.5,0.5)--++(0,0.5)--++(0.5,0);

\node[right] at (1.5,0) [scale=1]
{\tiny{-1}};
\node[right] at (1.5,0.5) [scale=1]
{\tiny{$\ell$}};
\node[right] at (1.5,1) [scale=1]
{\tiny{+1}};

\node[above] at (0.5, 1.5) [scale=1]
{\tiny{$k$}};
\node[above] at (1, 1.5) [scale=1]
{\tiny{+1}};

\node at (0.75,0.25) [scale=0.8] {\small{Z}};
\node at (0.75,0.75) [scale=0.8] {\small{W}};

%black lines
{
\draw[black, line width=0.15mm] (0.70,0.45)--++(0,0.1);
\draw[black, line width=0.15mm] (0.75,0.45)--++(0,0.1);
\draw[black, line width=0.15mm] (0.8,0.45)--++(0,0.1);

\draw[black, line width=0.15mm] (0.95,0.70)--++(0.1,0);
\draw[black, line width=0.15mm] (0.95,0.75)--++(0.1,0);
\draw[black, line width=0.15mm] (0.95,0.80)--++(0.1,0);

}

\node[below, align=center, text width=3cm] at (0.75, 0) { Fig. 4.6. };

} 

\end{scope}

\end{tikzpicture}
\end{adjustbox}
\end{wrapfigure}

\null

\noindent \textit{Proof of (c).} By part (a) and the assumption that $W$ is an end-box of $P$, $W$ has exactly one $H$-neighbour in $P$ and no other $H$-neighbours in $N[P]\setminus P$. It follows that $W$ has exactly two edges in $H$ and two edges not in $H$. BWOC, assume that the other edge of $W$ not in $H$ is $e(k,k+1;l+1)$. But then $e(k;l-1,l) \in H$ and $e(k+1;l-1,l) \in H$, and $W$ is switchable, which conflicts with the definition of a looping fat path. It follows that $e(k,k+1;l+1) \in H$. See Figure 4.6.

That exactly one of $e(k;l,l+1)$ and $e(k+1;l,l+1)$ belong to $H$, follows from the fact that $W$ has exactly one $H$-neighbour in $P$ and exactly one $H$-neighbour, $Z$, in $N[P]\setminus P$. End of proof for (c).

\setlength{\intextsep}{0pt}
\setlength{\columnsep}{20pt}
\begin{center}
\begin{adjustbox}{trim=0cm 0cm 0cm 0cm}

\end{adjustbox}
\end{center}

\noindent \textit{Proof of (d).} $W$ is either an end-box of $P$ or it is not.

 \null

\noindent \textit{CASE 1: W is an end-box of P.} By part (c), we may assume WLOG that $e(k;l,l+1)\in H$ and $e(k+1;l,l+1)\notin H$. Then $F$ is eastern or southern. Suppose that $F$ is eastern. Then $e(k+1,k+2;l+1) \in H$. It follows that $S_{\uparrow}(k+1,l-1;k+2,l)\in H$. But then $W+(1,0) \in P$ is switchable, which conflicts with the definition of a looping fat path. So $F$ must be southern. Then $e(k;l-1,l) \in H$. It follows that $S_{\uparrow}(k+1,l-1;k+2,l)\in H$. Now, either $e(k,k+1;l-1) \in H$ or  $e(k,k+1;l-1)\notin H$. Either way, (d) is satisfied. See Figure 4.7 (a) and (b).

\null

\noindent \textit{CASE 2: W is not an end-box of P.} By part (b), the $H$-neighbours of $W$ in $P$ are $W+(-1,0)$ and $W+(1,0)$ and we have that $S_{\rightarrow}(k-1,l;k,l-1) \in H$, $S_{\uparrow}(k+1,l-1;k+2,l) \in H$. Then, either $e(k,k+1;l-1) \in H$ or  $e(k,k+1;l-1) \notin H$. Either way, (d) is satisfied.  See Figure 4.7 (c).$\square$

\null

\noindent Let $F=G\langle N[P(X,Y)] \rangle$ be a looping fat path following a leaf $L$. We will often write ``the middle-boxes of $F$'' to refer to the middle-boxes of the $A_1$-type following $L$. 

\null 

\noindent \textbf{Proposition 4.7.} Let $G$ be a polyomino and let $H$ be a Hamiltonian path of $G$. Then every \index{sturdy looping fat path}sturdy looping fat path of $G$ is non-self-adjacent.

\null

\noindent We will need the following lemmas and observations before we prove Proposition 4.7.

% Already added these definitions in first page
{
%Let $F=G\langle N[P(X,Y)] \rangle$ be a southern looping fat path following a leaf $L=R(k,l)$. Then we call the end-box of $P$ that is located on the half-plane on the right of the directed line determined by the edge $e(k;l,l+1)$, \textit{the X-end-box of F}; and we call the other end-box of $F$ \textit{the Y-end-box of F}.
}

\null

\noindent \textbf{Observation 4.8.} Let $F=G\langle N[P(X,Y)] \rangle$ be a southern looping fat path. Let $P_0$, $P_1$ and $P_2$ be the partitioning of $H$ given in Lemma 4.1. Then $(v_x, v_{x+1})$ and $(v_{y-1}, v_y)$ are either anti-parallel and incident on the box $X+(-1,0)$ (as in Figure 4.8 (a)) or collinear with the edge $\{v_x,v_y \}$ of $G\setminus H$ (as in Figure 4.8 (b)). 
%Furthermore, if $(v_x, v_{x+1})$ and $(v_{y-1}, v_y)$ are either anti-parallel then $v(v_{x-1},v_x) \in P_1$ and $v(v_y,v_{y+1}) \in P_3$; and if $(v_x, v_{x+1})$ and $(v_{y-1}, v_y)$ are collinear with the edge $\{v_x,v_y \}$ of $G\setminus H$, 
Analogous observations apply to eastern, northern and western looping fat paths. $\square$

\begingroup
\setlength{\intextsep}{0pt}
\setlength{\columnsep}{20pt}
\begin{center}
\begin{adjustbox}{trim=0cm 0.25cm 0cm 0cm}

\end{adjustbox}
\end{center}

\noindent \textbf{Lemma 4.9.} Let $G$, $H$, $P(X,Y)$, $P_0$, $P_1$ and $P_2$ be  as in Lemma 4.1. Let $Q$ be the cycle $P_0, \{v_x, v_y\}$ in $G$. Let the vertices $v$ and $v'$  belong to $V((P_1 \setminus \{v_x\}) \cup (P_2 \setminus \{v_y\}))$. Then the segment $[v,v']$ does not intersect $Q$ exactly once at point $p$ of an edge $\{u,u'\}$ of $Q$ with $p\neq u$ and $p\neq u'$. 

\null

\noindent \textit{Proof.} Let $U$ be the region bounded by $Q$. Let the vertices $v$ and $v'$  belong to $V((P_1 \setminus \{v_x\}) \cup (P_2 \setminus \{v_y\}))$. Then $v \in V(U)$ or $v \in V(G) \setminus V(U)$. For definiteness, assume that  $v \in V(U)$. For a contradiction, assume that the segment $[v,v']$ that intersects $Q$ exactly once at a point $p$ of an edge $\{u,u'\}$ of $Q$ with $p\neq u$ and $p\neq u'$. Now, Corollary 1.1.5 implies that $v' \in V(G) \setminus V(U)$. Lemma 1.3.3 (c) implies that $v_1 \in V(U)$. By Lemma 1.5.4 in Chapter 1, $v_r \in V(U)$ as well. By Lemma 1.3.3 (c) again, $v' \in V(U)$ contradicting that $v' \in V(G) \setminus V(U)$ $\square$.

\null

\noindent\textbf{Lemma 4.10.} Let $G$ be a polyomino and let $H=v_1, \ldots, v_r$ be a Hamiltonian path of $G$. Let $F=G\langle N[P(X,Y)] \rangle$ be a \index{sturdy looping fat path}sturdy looping fat path. Let $P_0$, $P_1$ and $P_2$ be the partitioning of $H$ given in Lemma 4.1. Then, no edge of $P_1$ or $P_2$ has both of its incident boxes belonging to $F$.

\null 

\begingroup

\noindent \textbf{Remark.} The non-self-adjacency of sturdy looping fat paths is a natural consequence of their local structure but the proof requires lengthy and detailed case analysis. In the proofs of Lemma 4.10 and Proposition 4.7 we will use repeatedly and implicitly Proposition 4.5 (a), Lemma 4.6 (a), and the facts that a box of $P$ that is not an end-box has exactly two $H$-neighbours in $P$ while an end-box of $P$ has exactly one $H$-neighbour in $P$. 

\setlength{\intextsep}{0pt}
\setlength{\columnsep}{20pt}
\begin{wrapfigure}[]{r}{0cm}
\begin{adjustbox}{trim=0cm 0cm 0cm 0cm}
\begin{tikzpicture}[scale=1.75]

\begin{scope}[xshift=0cm]{
\draw[gray,very thin, step=0.5cm, opacity=0.5] (0,0) grid (1,1);

\fill[blue!40!white, opacity=0.5] (0,0.5) rectangle (0.5,1);

\fill[green!40!white, opacity=0.5] (0.5,0.5) rectangle (1,1);

\draw[blue, line width=0.5mm] (0,0.5)--++(0,0.5)--++(1,0);

%black lines
{

\draw[black, line width=0.15mm] (0.2,0.45)--++(0,0.1);
\draw[black, line width=0.15mm] (0.25,0.45)--++(0,0.1);
\draw[black, line width=0.15mm] (0.3,0.45)--++(0,0.1);

}

% nodes
{
\node[right] at (1,0) [scale=1]
{\tiny{-1}};
\node[right] at (1,0.5) [scale=1]
{\tiny{$b$}};
%\node[right] at (1,1) [scale=1]
{\tiny{+1}};

\node[below] at (0,0) [scale=1]
{\tiny{$a$}};
\node[below] at (0.5,0) [scale=1]
{\tiny{+1}};

\node at (0.25,0.75) [scale=0.8] {\small{$Z$}};
}

\node[above] at (0.5,1) [scale=1.25]{(a)};

\node[below, align=center, text width=3cm] at (1.25, -0.25) { Fig. 4.9. };

} \end{scope}

\begin{scope}[xshift=1.75cm]{
\draw[gray,very thin, step=0.5cm, opacity=0.5] (0,0) grid (1,1);

\fill[blue!40!white, opacity=0.5] (0,0) rectangle (0.5,1);

\fill[green!40!white, opacity=0.5] (0.5,0.5) rectangle (1,1);

\draw[blue, line width=0.5mm] (0,0.5)--++(0,0.5)--++(1,0);

\draw[blue, line width=0.5mm] (0.5,0)--++(0,0.5)--++(0.5,0);
\draw[blue, line width=0.5mm] (0,0)--++(0,0.5);
%black lines
{

\draw[black, line width=0.15mm] (0.2,0.45)--++(0,0.1);
\draw[black, line width=0.15mm] (0.25,0.45)--++(0,0.1);
\draw[black, line width=0.15mm] (0.3,0.45)--++(0,0.1);

}

% nodes
{
\node[right] at (1,0) [scale=1]
{\tiny{-1}};
\node[right] at (1,0.5) [scale=1]
{\tiny{$b$}};
\node[right] at (1,1) [scale=1]
{\tiny{+1}};

\node[below] at (0,0) [scale=1]
{\tiny{$a$}};
\node[below] at (0.5,0) [scale=1]
{\tiny{+1}};

\node at (0.25,0.75) [scale=0.8] {\small{$Z$}};
}

\node[above] at (0.5,1) [scale=1.25]{(b)};

} \end{scope}

\end{tikzpicture}
\end{adjustbox}
\end{wrapfigure}

\noindent Lastly, the following type of scenario occurs several times. Let $F$ be a looping fat path. Suppose $Z=R(a,b)$ is an end-box of $P$ and $Z+(1,0)$ is its $H$-neighbour in $N[P] \setminus P$. Assume that $e(a;b,b+1)\in H$, $e(a,a+1;b) \notin H$, $e(a,a+1;b+1) \in H$ $e(a+1,a+2;b+1) \in H$. See Figure 4.9 (a). Then the $F$ can be eastern or southern. But if $F$ is southern then $e(a;b-1,b) \in H$ and $S_{\uparrow}(a+1,b-1;a+2,b) \in H$, and then $Z+(0,-1)$ is in $P$ and it is switchable. Hence, $F$ must be eastern. See Figure 4.9 (b). The cases where the $H$-neighbour of $Z$ is north, west, or south of $Z$ are analogous. We denote the argument in this paragraph ($**$) for reference. 

\endgroup 

\null 

\noindent \textit{Proof of Lemma 4.10.} Let $F$ be a sturdy looping fat path. Assume that some edge $e_j=(v_j,v_{j+1})$ of $P_1$ or $P_2$ has both its incident boxes in $F$. Without loss of generality assume that $e_j$ is an edge of $P_2$. We will use induction to show that now we can find a box of $P$ that has an end-vertex of $H$ incident on it, which contradicts the assumption that $F$ is a sturdy looping fat path.

Let $e_j$ be an edge of $P_2$. Then $e_{j+1}$ is one edge closer to $v_r$ then $e_j$.  We denote this remark by $(*)$ for reference.

\null

\begingroup 
\setlength{\intextsep}{0pt}
\setlength{\columnsep}{20pt}
\begin{wrapfigure}[]{r}{0cm}
\begin{adjustbox}{trim=0cm 0cm 0cm 0.25cm}
\begin{tikzpicture}[scale=2]

\begin{scope}[xshift=0cm]{
\draw[gray,very thin, step=0.5cm, opacity=0.5] (0,0) grid (1,1);

%f the dotted $H$-cycle
{
\draw [blue!40!white, line width=2mm, dotted, opacity=0.5] plot [smooth, tension=0.75] coordinates {(0.25,0.25) (-0.25,0.35)(-0.25,0.65)(0.25,0.75)(0.75,0.65)(0.75, 0.35)(0.25,0.25)};

}

\fill[green!50!white, opacity=0.5] (0,0) rectangle (0.5,1);

\draw[blue, line width=0.5mm] (0,0.5)--++(0.5,0);

\draw[fill=blue] (00.5,0.5) circle [radius=0.05];

%black lines
{

}

% nodes
{

\node[left] at (0,0.5) [scale=1]
{\tiny{$b$}};

\node[above] at (0,1) [scale=1]{\tiny{$a$}};
\node[above] at (0.5,1) [scale=1]{\tiny{+1}};

\node at (0.25,0.75) [scale=0.8] {\small{$W_1$}};
\node at (0.25,0.25) [scale=0.8] {\small{$W_2$}};
}

\node[below, align=center, text width=3cm] at (0.5, 0) { Fig. 4.10. };

} \end{scope}

\end{tikzpicture}
\end{adjustbox}
\end{wrapfigure}

\noindent Let $Q(j)$ be the statement ``Both boxes incident on $e_j$ belong to $F$''. Assume that $Q(j)$ implies $Q(j+1)$. BWOC assume $Q(j_0)$ for some edge $e_{j_0}$ of $P_2$. By induction, we have $Q(r-1)$. For definiteness let $W_1=R(a,b)$ and let $W_2=R(a,b-1)$ be the boxes of $F$ incident on $e_{r-1}$, and let $v_r=v(a+1,b)$. If any of the boxes incident to $v_1$ belong to $P$, then this conflicts with the definition of $F$. Then we must have that $W_1$ and $W_2$ belong to $N[P] \setminus P$. By Lemma 4.6 (a), $W_1+(1,0)$ and $W_2+(1,0)$ do not belong to $F$. Let $W_1'$ and $W_2'$ be the $H$-neighbours of $W_1$ and $W_2$ in $P$, respectively. Then there is an $H$-cycle $W_2, W_2', \ldots, W_1',W_1, W_1+(1,0), W_2+(1,0), W_2$, contradicting Proposition 1.2.1. See Figure 4.10.

\endgroup

\null

\noindent It remains to show that $Q(j)$ implies $Q(j+1)$. For definiteness, let $e_j=(v_j,v_{j+1})=(v(a,b), v(a+1,b))$.
We remark that assuming $e_{(j)}$ to be horizontal and fixing its orientation comes at the expense of allowing $F$ to have any direction; however, this trade-off is reasonable as it simplifies and shortens the proof. Assume $Q(j)$ is true, and let $Z_1$ and $Z_2$ be the boxes of $F$ incident on $e_j$. There are three possibilities: both $Z_1$ and $Z_2$ belong to $P$, or both $Z_1$ and $Z_2$ belong to $N[P]\setminus P$, or one of $Z_1$ and $Z_2$ belongs to $P$ and the other belongs to $N[P]\setminus P$. WLOG let $Z_1$ be on the left of $e_j$ and let $Z_2$ be on the right of $e_j$.

\null

\noindent \textit{CASE 1: Both $Z_1$ and $Z_2$ belong to $P$.} There are three possibilities: $e_{j+1}=e(a+1;b,b+1)$, $e_{j+1}=e(a+1,a+2;b)$, or $e_{j+1}=e(a+1;b-1,b)$. By symmetry of the first and third, we only need to check the first two.

\null 

\begingroup 
\setlength{\intextsep}{0pt}
\setlength{\columnsep}{20pt}
\begin{wrapfigure}[]{l}{0cm}
\begin{adjustbox}{trim=0cm 0.25cm 0cm 0cm}

\end{adjustbox}
\end{wrapfigure}

\noindent \textit{CASE 1.1: $e_{j+1}=e(a+1;b,b+1)$.} Then $e(a+1;b-1,b) \notin H$ and $e(a+1,a+2;b) \notin H$. Lemma 4.6 (d) implies that $Z_2+(1,0) \in P$. Then $Z_1+(1,0) \in P$ or $Z_1+(1,0) \in N[P] \setminus P$. Either way, $Q(j+1)$ holds. End of Case 1.1.

\null 

\noindent\textit{CASE 1.2: $e_{j+1}=e(a+1,a+2;b)$.} Then $e(a+1;b-1,b) \notin H$ and $e(a+1,b;b+1) \notin H$. Then $Z_1+(1,0)$ and $Z_2+(1,0)$ either belong to $P$ or to $N[P]\setminus P$. Either way $Q(j+1)$ holds. End of Case 1.2. End of Case 1.

\endgroup

\setlength{\intextsep}{0pt}
\setlength{\columnsep}{20pt}
\begin{center}
\begin{adjustbox}{trim=0cm 0cm 0cm 0.5cm}

\end{adjustbox}
\end{center}

\noindent \textit{CASE 2: One of $Z_1$ and $Z_2$ belongs to $P$ and the other belongs to $N[P]\setminus P$.} For definiteness, let $Z_1 \in P$ and $Z_2 \in N[P]\setminus P$. Then $e_{j+1}=e(a+1;b,b+1)$, $e_{j+1}=e(a+1,a+2;b)$ or $e_{j+1}=e(a+1,b-1;b)$. Note that the case where $e_{j+1}=e(a+1,b-1;b)$ is the same as Case 1.1, so we only need to check the cases where $e_{j+1}=e(a+1;b,b+1)$ and where $e_{j+1}=e(a+1,a+2;b)$.

\null 

\noindent \textit{CASE 2.1: $e_{j+1}=e(a+1;b,b+1)$.} Then $e(a+1;b-1,b) \notin H$. By Lemma 4.6 (d), $e(a,a+1;b-1)\in H$. By Lemma 4.6 (a), $Z_2+(1,0)$ either belongs to $P$ or it is not a box of $F$. If $Z_2+(1,0) \in P$, then $Z_2+(1,1) \in P$ or $Z_2+(1,1) \in N[P] \setminus P$. Either way, $Q(j+1)$ holds. We check that $Z_2+(1,0) \notin P$ is impossible.

Toward a contradiction, assume that $Z_2+(1,0)\notin P$. Then we have that $Z_2+(-1,0) \in P$, and Lemma 4.6 (b) and (c) imply that $e(a-1;b-1,b) \in H$. Now, $Z_2+(-1,0)$ is an end-box of $P$ or it is not. See Figure 4.12 (a).

\null 

\noindent \textit{CASE 2.1(a): $Z_2+(-1,0)$ is not an end-box of $P$.} By Lemma 4.6 (b), $Z_2+(-1,-1)\in P$, $Z_2+(-1,1)\in P$, $S_{\downarrow}(a,b+1;a+1,b) \in H$, $S_{\uparrow}(a,b-2;a+1,b-1) \in H$ and $e(a-1;b-1,b) \in H$. Then, after $Z_2 \mapsto Z_1$, $Z_2+(-1,0)\in P$ is switchable, contradicting FPC-4. See Figure 4.12 (b). End of Case 2.1(a). 

\endgroup 

\null

\noindent \textit{CASE 2.1(b): $Z_2+(-1,0)$ is an end-box of $P$.} By Lemma 4.6(c) exactly one of $e(a-1,a;b-1)$ and $e(a-1,a;b)$ belong to $H$.

\null 

\noindent \textit{CASE 2.1($b_1$): $e(a-1,a;b-1)\in H$ and $e(a-1,a;b) \notin H$.} Then we have that $e(a;b,b+1) \in H$. But this implies that $Z_1$ is an end-box of $P$ as well, which is impossible. See Figure 4.12 (c). End of Case 2.1($b_1$).

% NOT NEEDED
{
%Now, note that $F$ could be northern or eastern. 

%\textit{CASE 2.1($a_1$): $F$ is northern.} Then we must have $e(a-1;b,b+1)\in H$, $e(a;b,b+1)\in H$ and $e(a-1,a;b+1)\notin H$. But now $Z_1+(-1,0) \in P$ and switchable. End of Case 2.1($a_1$).

%\textit{CASE 2.1($a_2$): $F$ is eastern.} Then we must have that $e(a;b,b+1) \in H$ and that $X'=Z_2+(-1,-1)$. Then $Z_1 \mapsto Z_2$, $Z_2+(-1,0) \mapsto X'$ is a short weakening of $P(X,Y)$. End of Case 2.1($a_2$). End of Case 2.1(a).
}

\null 

\begingroup 
\setlength{\intextsep}{0pt}
\setlength{\columnsep}{20pt}
\begin{wrapfigure}[]{l}{0cm}
\begin{adjustbox}{trim=0cm 0cm 0cm 0cm}

\end{adjustbox}
\end{wrapfigure}

\noindent \textit{CASE 2.1($b_2$):  $e(a-1,a;b) \in H$ and  $e(a-1,a;b-1)\notin H$.} Then $e(a;b,b+1) \notin H$ and $e(a-1;b,b+1) \notin H$. Note that if $e(a,a+1;b+1) \in H$, then $Z_1$ is the other end-box of $P$, which does not agree with the assumption that $Z_2+(-1,0)$ is an end-box of $P$. So we only need to check the case where $e(a,a+1;b+1) \notin H$. In this case, we must have $S_{\rightarrow}(a-1,b+1;a,b+2) \in H$. Now, the fact that $Z_1$ is not an end-box of $P$ implies that both $H$-neighbours $Z_1+(-1,0)$ and $Z_1+(0,1)$ of $Z_1$ are boxes of $P$. But then $Z_1+(-1,0)$ is switchable. See Figure 4.13. End of Case 2.1$(b_2)$. End of Case 2.1(b). End of Case 2.1. 

\null 

\noindent \textit{CASE 2.2: $e_{j+1}=e(a+1,a+2;b)$.} Then $e(a+1;b-1,b) \notin H$, $e(a+1;b,b+1) \notin H$ and $Z_1+(1,0)$ must be a box of $F$. Now, $Z_2+(1,0)$ either belongs to $P$ or it is not a box of $F$. If $Z_2+(1,0) \in P$, then $Q(j+1)$ holds. We will show that $Z_2+(1,0) \notin P$ is impossible. 

\begingroup 
\setlength{\intextsep}{0pt}
\setlength{\columnsep}{20pt}
\begin{wrapfigure}[]{r}{0cm}
\begin{adjustbox}{trim=0cm 0.5cm 0cm 0.25cm}

\end{adjustbox}
\end{wrapfigure}

\noindent Toward a contradiction, assume that $Z_2+(1,0)$ is not a box of $F$. Then we have that $Z_2+(-1,0) \in P$, and, by Lemma 4.6 (d), $e(a,a+1;b-1) \in H$. Lemma 4.6 (b) and (c) imply that $e(a-1;b-1,b) \in H$. See Figure 4.14. Now, either $Z_2+(-1,0)$ is an end-box of $P$ or it is not.
 
\null

\begingroup 
\setlength{\intextsep}{0pt}
\setlength{\columnsep}{20pt}
\begin{wrapfigure}[]{l}{0cm}
\begin{adjustbox}{trim=0cm 0cm 0cm 0.5cm}
%
\end{adjustbox}
\end{wrapfigure}

\noindent \textit{CASE 2.2(a): $Z_2+(-1,0)$ is not an end-box of $P$.} By Lemma 4.6 (b) we have that $Z_2+(-1,-1) \in P$, $Z_2+(-1,1) \in P$, $e(a;b-2,b-1)\in H$ and $e(a;b,b+1)\in H$. See Figure 4.15 (a). Now, either $Z_1$ is an end-box of $P$ or it is not.

\null

\noindent \textit{CASE 2.2($a_1$): $Z_1$ is not an end-box of $P$.} Then $e(a,a+1;b+1) \notin H$. It follows that $S_{\downarrow}(a+1,b+2;a+2,b+1) \in H$. Note that if $e(a+2;b,b+1) \notin H$, then $Z_1+(1,0)$ is switchable, so we only need to check the case where $e(a+2;b,b+1) \in H$. This implies that $Z_1+(1,0)$ is an end-box of $P$. It follows that $F$ is eastern and that $v(a+2,b)$ is a corner of $A_1$. But now, the $A_1$-partitioning of $H$ implies that $e(a+2;b,b+1)=e_{j+2} \in P_0$, contradicting $(*)$ in the beginning of the lemma. See Figure 4.15 (b). End of Case 2.2($a_1$).  

\null 

\begingroup 
\setlength{\intextsep}{0pt}
\setlength{\columnsep}{20pt}
\begin{wrapfigure}[]{r}{0cm}
\begin{adjustbox}{trim=0cm 0cm 0cm 0.75cm}

\end{adjustbox}
\end{wrapfigure}

\noindent \textit{CASE 2.2($a_2$): $Z_1$ is an end-box of $P$.} Then $e(a;b+1,b+2)\in H$ or $e(a;b+1,b+2)\notin H$. If $e(a;b+1,b+2)\in H$, then we must have  $S_{\downarrow}(a+1,b+2;a+2,b+1) \in H$. But this implies that at least one of $Z_1+(0,1)$ and $Z_1+(1,0)$ belongs to $P$ and is switchable. So we only need to check the case where $e(a;b+1,b+2)\notin H$. In this case, we observe that $F$ must be western. This implies that one of $v(a,b)$ and $v(a,b+1)$ is a corner of $A_1$. If the former, then $e(a;b,b+1)\in P_0$, which contradicts Lemma 4.3. See Figure 4.16 (a). If the latter then the other end-box of $P$ is $Z_1+(0,2)$ and $F$ follows the edge $e(a-1;b+1,b+2)$. This means that $e(a,a+1;b+1) \in H$, and that $e(a+1,a+2;b+1) \in H$. But then $Z_1+(1,0) \in H$ is switchable. See Figure 4.16 (b). End of Case 2.2($a_2$).

\null 

\noindent \textit{CASE 2.2(b): $Z_2+(-1,0)$ is an end-box of $P$.} By Lemma 4.6 (c), exactly one of $e(a-1,a;b-1)$ and $e(a-1,a;b)$ belongs to $H$.

\null 

\noindent \textit{CASE 2.2($b_1$): $e(a-1,a;b-1) \in H$ and $e(a-1,a;b)\notin H$.} Then $e(a;b,b+1) \in H$. Since $Z_1$ cannot be the other end-box of $P$, its $H$-neighbours in $P$ must be $Z_1+(0,1)$ and $Z_1+(1,0)$, so $e(a,a+1;b+1) \notin H$. Then we must have that $e(a+1,a+2;b+1) \in H$. Since $Z_1+(1,0)$ cannot be the other end-box of $P$ either, it must be switchable, contradicting FPC-3. See Figure 4.17 (a). End of Case 2.2($b_1$).

\begingroup 
\setlength{\intextsep}{0pt}
\setlength{\columnsep}{20pt}
\begin{wrapfigure}[]{l}{0cm}
\begin{adjustbox}{trim=0cm 0cm 0cm 0cm}

\end{adjustbox}
\end{wrapfigure}

\null

\noindent \null  \textit{CASE 2.2($b_2$): $e(a-1,a;b-1) \notin H$ and $e(a-1,a;b)\in H$.} Then $e(a;b,b+1)\notin H$. Note that $Z_1$ and $Z_1+(1,0)$ are not end-boxes. If $e(a,a+1;b+1) \in H$, then $Z_1$ is switchable, so we only need to check the case where $e(a,a+1;b+1) \notin H$. Then we have that $S_{\downarrow}(a+1,b+2;a+2,b+1) \in H$. Lemma 4.6(b) implies $Z_1+(1,0) \in P$. Now, either $e(a+2;b,b+1)\in H$ or $e(a+2;b,b+1)\notin H$. The former implies that $Z_1+(1,0)$ is an end-box of $P$, which is impossible; and if the latter, then $Z_1+(1,0)$ is switchable, contradicting FPC-3. See Figure 4.17 (b). End of Case 2.2($b_2$). End of Case 2.2(b). End of Case 2.2. End of Case 2.

\null

\noindent \textit{CASE 3: both $Z_1$ and $Z_2$ belong to $N[P]\setminus P$.} There are three possibilities: $e_{j+1}=e(a+1;b,b+1)$, $e_{j+1}=e(a+1,a+2;b)$ or $e_{j+1}=e(a+1;b-1,b)$. By symmetry of the first and third, we only need to check the first two.

\null 

\begingroup 
\setlength{\intextsep}{0pt}
\setlength{\columnsep}{20pt}
\begin{wrapfigure}[]{r}{0cm}
\begin{adjustbox}{trim=0cm 0cm 0cm 0cm}
\begin{tikzpicture}[scale=1.5]

\begin{scope}[xshift=0cm]{
\draw[gray,very thin, step=0.5cm, opacity=0.5] (0,0) grid (1.5,1);

\fill[blue!40!white, opacity=0.5] (0,0) rectangle (0.5,0.5);

\fill[green!40!white, opacity=0.5] (0.5,0.5) rectangle (1,1);
\fill[green!40!white, opacity=0.5] (0.5,0) rectangle (1,0.5);

\begin{scope}
[very thick,decoration={
    markings,
    mark=at position 0.6 with {\arrow{>}}}
    ]
    
    \draw[postaction={decorate}, blue, line width=0.5mm] (0.5,0.5)--++(0.5,0);
    
    \draw[postaction={decorate}, blue, line width=0.5mm] (1,0.5)--++(0,0.5);
    
\end{scope}

\draw[blue, line width=0.5mm] (0,0)--++(0,0.5);
\draw[blue, line width=0.5mm] (0.5,0)--++(0.5,0);

%\draw[fill=orange] (0,1.5) circle [radius=0.05];
%\draw[red!75!black, line width=0.5mm] (0,1.5)--++(1.5,-1);

%black lines
{

\draw[black, line width=0.15mm] (0.45,0.2)--++(0.1,0);
\draw[black, line width=0.15mm] (0.45,0.25)--++(0.1,0);
\draw[black, line width=0.15mm] (0.45,0.3)--++(0.1,0);

\draw[black, line width=0.15mm] (0.95,0.2)--++(0.1,0);
\draw[black, line width=0.15mm] (0.95,0.25)--++(0.1,0);
\draw[black, line width=0.15mm] (0.95,0.3)--++(0.1,0);

\draw[black, line width=0.15mm] (1.2,0.45)--++(0,0.1);
\draw[black, line width=0.15mm] (1.25,0.45)--++(0,0.1);
\draw[black, line width=0.15mm] (1.3,0.45)--++(0,0.1);

}

% nodes
{
\node[left] at (0,0) [scale=1]
{\tiny{-1}};
\node[left] at (0,0.5) [scale=1]
{\tiny{$b$}};
\node[left] at (0,1) [scale=1]
{\tiny{+1}};

\node[below] at (0,0) [scale=1]
{\tiny{-1}};
\node[below] at (0.5,0) [scale=1]
{\tiny{$a$}};
\node[below] at (1, 0) [scale=1]
{\tiny{+1}};

\node at (0.75,0.25) [scale=0.8] {\small{$Z_2$}};
\node at (0.75,0.75) [scale=0.8] {\small{$Z_1$}};
}

\node[below, align=center, text width=3cm] at (0.75, -0.25) { Fig. 4.18.};

} \end{scope}

\end{tikzpicture}
\end{adjustbox}
\end{wrapfigure}

\noindent \textit{CASE 3.1: $e_{(j+1)}=e(a+1;b,b+1)$.} Then $e(a+1;b-1,b)\notin H$. By Lemma 4.6 (d), $e(a,a+1;b-1)\in H$. Now the $H$-neighbour of $Z_2$ in $P$ is $Z_2+(-1,0)$ or $Z_2+(1,0)$. If the latter, then we're done by Case 1.1. We will show the former is impossible.

Assume that the $H$-neighbour of $Z_2$ in $P$ is $Z_2+(-1,0)$. It follows that $e(a;b-1,b) \notin H$ and, by Lemma 4.6 (b) and (c), that $e(a-1;b-1,b) \in H$. Now, $Z_2+(-1,0)$ is an end-box of $P$ or it is not. Note that if $Z_2+(-1,0)$ is not an end-box of $P$, then, as in Case 2.1(a), after $Z_2\mapsto Z_1$, $Z_2+(-1,0)\in P$ is switchable, contradicting FPC-4. Therefore,
we only need to check the case where $Z_2+(-1,0)$ is an end-box of $P$. See Figure 4.18. By Lemma 4.6(c), exactly one of $e(a-1,a;b-1)$ and $e(a-1,a;b)$ belongs to $H$.

\endgroup

\begingroup 
\setlength{\intextsep}{0pt}
\setlength{\columnsep}{20pt}
\begin{wrapfigure}[]{l}{0cm}
\begin{adjustbox}{trim=0cm 0.5cm 0cm 0.25cm}

\end{adjustbox}
\end{wrapfigure}

\noindent \textit{CASE 3.1(a): $e(a-1,a;b)\in H$ and $e(a-1,a;b-1)\notin H$.} Then $e(a;b,b+1) \notin H$. By Lemma 4.6 (d), it follows that $e(a,a+1;b+1)\in H$. Then we have that $Z_1+(1,0) \in P$. Note that Lemma 4.6(b) implies that $Z_1+(1,0)$ is an end-box of $P$, but this disagrees with the assumption that $Z_2+(-1,0)$ is the other end-box of $P$. See Figure 4.19 (a). End of Case 3.1(a). 

\null 

\noindent \textit{CASE 3.1(b): $e(a-1,a;b-1)\in H$ and $e(a-1,a;b)\notin H$.} Then $e(a;b,b+1)\in H$. Note that $F$ can be northern or eastern. If $F$ is northern, then $e(a-1;b,b+1)\in H$, and the box $Z_2+(-1,1) \in P$ is switchable; and if $F$ is eastern, then the switchable middle-box of $F$ is $W=Z_2+(-1,-1)$. But then, after $Z_2 \mapsto Z_1$, $W \mapsto Z_2+(-1,0)$ is valid, contradicting FPC-4. See Figure 4.19 (b). End of Case 3.1(b). End of Case 3.1 

\endgroup

\null 

\begingroup 
\setlength{\intextsep}{0pt}
\setlength{\columnsep}{20pt}
\begin{wrapfigure}[]{r}{0cm}
\begin{adjustbox}{trim=0cm 0cm 0cm 0cm}
\begin{tikzpicture}[scale=1.5]

\begin{scope}[xshift=0cm]{
\draw[gray,very thin, step=0.5cm, opacity=0.5] (0,0) grid (1.5,1);

\fill[blue!40!white, opacity=0.5] (0,0) rectangle (0.5,0.5);

\fill[green!40!white, opacity=0.5] (0.5,0.5) rectangle (1,1);
\fill[green!40!white, opacity=0.5] (0.5,0) rectangle (1,0.5);

\begin{scope}
[very thick,decoration={
    markings,
    mark=at position 0.6 with {\arrow{>}}}
    ]
    
    \draw[postaction={decorate}, blue, line width=0.5mm] (0.5,0.5)--++(0.5,0);
    
    \draw[postaction={decorate}, blue, line width=0.5mm] (1,0.5)--++(0.5,0);
    
\end{scope}

\draw[blue, line width=0.5mm] (0,0)--++(0,0.5);
\draw[blue, line width=0.5mm] (0.5,0)--++(0.5,0);

%\draw[fill=orange] (0,1.5) circle [radius=0.05];
%\draw[red!75!black, line width=0.5mm] (0,1.5)--++(1.5,-1);

%black lines
{

\draw[black, line width=0.15mm] (0.45,0.2)--++(0.1,0);
\draw[black, line width=0.15mm] (0.45,0.25)--++(0.1,0);
\draw[black, line width=0.15mm] (0.45,0.3)--++(0.1,0);

\draw[black, line width=0.15mm] (0.95,0.2)--++(0.1,0);
\draw[black, line width=0.15mm] (0.95,0.25)--++(0.1,0);
\draw[black, line width=0.15mm] (0.95,0.3)--++(0.1,0);

\draw[black, line width=0.15mm] (0.95,0.7)--++(0.1,0);
\draw[black, line width=0.15mm] (0.95,0.75)--++(0.1,0);
\draw[black, line width=0.15mm] (0.95,0.8)--++(0.1,0);

}

% nodes
{
\node[left] at (0,0) [scale=1]
{\tiny{-1}};
\node[left] at (0,0.5) [scale=1]
{\tiny{$b$}};
\node[left] at (0,1) [scale=1]
{\tiny{+1}};

\node[below] at (0,0) [scale=1]
{\tiny{-1}};
\node[below] at (0.5,0) [scale=1]
{\tiny{$a$}};
\node[below] at (1, 0) [scale=1]
{\tiny{+1}};

\node at (0.75,0.25) [scale=0.8] {\small{$Z_2$}};
\node at (0.75,0.75) [scale=0.8] {\small{$Z_1$}};
}

\node[below, align=center, text width=4cm] at (0.75, -0.25) { Fig. 4.20. Case 3.2.};

} \end{scope}

\end{tikzpicture}
\end{adjustbox}
\end{wrapfigure}

\noindent \textit{CASE 3.2: $e_{(j+1)}=e(a+1,a+2;b)$.} If both $Z_1+(1,0)$ and $Z_2+(1,0)$ belong to $F$ we're done, so we assume that at least one of $Z_1+(1,0)$ and $Z_2+(1,0)$ is not a box of $F$. By symmetry, we may assume WLOG that $Z_2+(1,0)$ is not a box of $F$. By Lemma 4.6(d) we have that $e(a,a+1;b-1)\in H$, $Z_2+(-1,0) \in P$ and $e(a;b-1,b)\notin H$. By Lemma 4.6 (b) and (c), $e(a-1;b-1,b)\in H$. See Figure 4.20. Now, $Z_2+(-1,0)$ is an end-box of $P$ or it is not.

\null

\begingroup 
\setlength{\intextsep}{0pt}
\setlength{\columnsep}{20pt}
\begin{wrapfigure}[]{l}{0cm}
\begin{adjustbox}{trim=0cm 0cm 0cm 0cm}

\end{adjustbox}
\end{wrapfigure}

\noindent \textit{CASE 3.2(a): $Z_2+(-1,0)$ is not an end-box of $P$.}  Then, by Lemma 4.6 (b), $Z_1+(-1,0)\in P$ and $e(a;b,b+1) \in H$. By Lemma 4.6 (d) we have that $e(a,a+1;b+1) \in H$. It follows that the $H$-neighbour of $Z_1$ in $P$ is $Z_1+(1,0)$. Lemma 4.6 (b) and (c) imply that $e(a+2;b,b+1) \in H$. Then it must be the case that $Z_1+(1,0)$ is an end-box of $P$. Observe that this implies
that $v(a+2,b)$ is a corner vertex of $A_1$. Now, the $A_1$-partitioning of $H$ implies that $e(a+2;b,b+1)\in P_0$, while our assumption that $e(a,a+1;b)\in E(P_1) \cup E(P_2)$ implies that $e(a+2;b,b+1)=e_{(j+2)}$. But then $e_{(j+2)} \in P_0$, contradicting $(*)$. See Figure 4.21. End of Case 3.2(a).

\null 

\noindent \textit{CASE 3.2(b): $Z_2+(-1,0)$ is an end-box of $P$.}  By Lemma 4.6(c), exactly one of $e(a-1,a;b-1)$ and $e(a-1,a;b)$ belongs to $H$.

\null 
%pagemarker
\noindent \textit{CASE 3.2($b_1$): $e(a-1,a;b-1) \in H$ and $e(a-1,a;b) \notin H$.} It follows that $Z_1+(-1,0) \in P$ and $e(a;b,b+1)\in H$. As in Case 3.2(a), $Z_1+(1,0)$ must be an end-box of $P$. But this conflicts with our assumption that $Z_2+(-1,0)$ is an end-box of $P$. End of Case 3.2($b_1$)

% Probably not needed
{
%By Lemma 4.6 (c), $e(a,a+1;b+1) \in H$. It follows that $Z_1+(1,0)$ is the other end-box of $P$, but this contradicts our assumption that $Z_2+(-1,0)$ is an end-box of $P$. End of Case 3.2($b_1$). \textcolor{red}{[See fig for Case 3.2 (a).]}
}

\null

\begingroup 
\setlength{\intextsep}{0pt}
\setlength{\columnsep}{20pt}
\begin{wrapfigure}[]{l}{0cm}
\begin{adjustbox}{trim=0cm 0cm 0cm 0cm}

\end{adjustbox}
\end{wrapfigure}

\noindent \textit{CASE 3.2($b_2$): $e(a-1,a;b) \in H$ and $e(a-1,a;b-1) \notin H$.} Then $e(a;b,b+1)\notin H$ and, by Lemma 4.6 (d), $e(a,a+1;b+1)\in H$. Lemma 4.6 (b) and (c) imply that $e(a+2;b,b+1) \in H$. Then it must be the case that $Z_1+(1,0)$ is an end-box of $P$, but again, this conflicts with the assumption that  $Z_2+(-1,0)$ is an end-box of $P$. See Figure 4.22. End of Case 3.2($b_2$). End of Case 3.2(b). End of Case 3.2. End of Case 3. $\square$

\endgroup

\null

\noindent \textit{Proof of Proposition 4.7.} Let $F=G\langle N[P(X,Y)] \rangle$ be a sturdy looping fat path of $G$. Let $P_0$, $P_1$ and $P_2$ be the partitioning of $H$ given in Lemma 4.1. For a contradiction, assume that there is an edge $e$ of $H$ such that both its incident boxes $Z_1$ and $Z_2$ belong to $F$. By Lemma 4.10, we may assume that  $e$ is an edge of $P_0$. By Lemma 4.3, $Z_1$ and $Z_2$ cannot both be boxes of $P(X,Y)=P$. Then, either both $Z_1$ and $Z_2$ are boxes of $N[P] \setminus P$, or exactly one of $Z_1$ and $Z_2$ is a box of $P$, and the other is a box of $N[P] \setminus P$. For definiteness, assume that $e$ is the horizontal edge $e(a,a+1;b)$, that $Z_1$ is north of $e(a,a+1;b)$ and that $Z_2$ is south of $e(a,a+1;b)$,

\null

\noindent  \textit{CASE 1: $Z_1$ and $Z_2$ belong to $N[P] \setminus P$.} Let $Z_1'$ and $Z_2'$ be the $H$-neighbours of $Z_1$ and $Z_2$ in $P$, respectively. Note that it is not possible for both $Z_1'$ and $Z_2'$ to be end-boxes of $P$. Then, either exactly one of $Z_1$ and $Z_2$ is adjacent to an end-box of $P$ or neither is.

\null 

\noindent \textit{CASE 1.1: Neither $Z_1$ nor $Z_2$ is adjacent to an end-box of $P$.} There are three possibilities: $Z_2'=Z_2+(0,-1)$; $Z_2'=Z_2+(-1,0)$; and $Z_2'=Z_2+(1,0)$. Since the second and third are symmetric, we only need to check the first two.

\null

\noindent \textit{CASE 1.1(a): $Z_2'=Z_2+(0,-1)$.} By Lemma 4.6 (b), $Z_2'+(-1,-1)\in P$, $Z_2'+(1,-1)\in P$, $S_{\rightarrow}(a-1,b-1;a,b) \in H$, $S_{\downarrow}(a+1,b;a+2,b-1) \in H$ and $e(a,a+1;b-2)\in H$. It follows that $e(a;b,b+1)\notin H$ and $e(a+1;b,b+1)\notin H$. By Lemma 4.6(d), $e(a,a+1;b+1) \in H$. Then, after $Z_1 \mapsto Z_2$, the box $Z_2+(0,-1)$ is in $P$ and is switchable, contradicting FPC-4. See Figure 4.23 (a). End of Case 1.1(a).

\begingroup 
\setlength{\intextsep}{0pt}
\setlength{\columnsep}{20pt}
\begin{wrapfigure}[]{r}{0cm}
\begin{adjustbox}{trim=0cm 0cm 0cm 0cm}

\end{adjustbox}
\end{wrapfigure}

\null

\noindent \textit{CASE 1.1(b): $Z_2'=Z_2+(-1,0)$.} By Lemma 4.6 (b), $Z_2'+(0,-1)\in P$, $Z_2'+(0,1)\in P$, $S_{\downarrow}(a,b+1;a+1,b) \in H$, $S_{\uparrow}(a,b-2;a+1,b-1) \in H$ and $e(a-1;b-1,b) \in H$. Now $Z_1'=Z_1+(0,1)$ or $Z_1'=Z_1+(1,0)$. See Figure 4.23 (b).

\null

\noindent \textit{CASE 1.1($b_1$): $Z_1'=Z_1+(0,1)$.} By Lemma 4.6 (d), $e(a+1;b,b+1) \in H$. Then, after $Z_2 \mapsto Z_1$, $Z_2+(-1,0)\in P$ is switchable, contradicting FPC-4. See Figure 4.23 (c). End of Case 1.1($b_1$).

\endgroup 

\null

\noindent \textit{CASE 1.1($b_2$): $Z_1'=Z_1+(1,0)$.}  By Lemma 4.6 (b), $Z_1'+(0,-1)\in P$, $Z_1'+(0,1)\in P$, $S_{\rightarrow}(a,b+1;a+1,b+2) \in H$, $e(a+1;b-1,b) \in H$ and $e(a+2;b,b+1) \in H$. Note if $Z_1+(2,0) \in P$ then,  $e(a+2;b,b+1)\in E(P_1) \cup E(P_2)$ contradicts Lemma 4.10; and $e(a+2;b,b+1)\in P_0$ contradicts Lemma 4.3. Therefore we may assume that $Z_1+(2,0) \notin P$. If $e(a+2;b-1,b) \in H$ then $Z_2+(1,0)\in P$ is switchable, so we may also assume that $e(a+2;b-1,b) \notin H$. It follows that $S_{\uparrow}(a+2,b-2;a+3,b-1) \in H$ and $e(a+2,a+3;b) \in H$. By symmetry, we may assume that $e(a-1;b,b+1) \notin H$, $e(a-2,a-1;b) \in H$ and $S_{\rightarrow}(a-2,b+1;a-1,b+2) \in H$.

\begingroup 
\setlength{\intextsep}{0pt}
\setlength{\columnsep}{20pt}
\begin{wrapfigure}[]{r}{0cm}
\begin{adjustbox}{trim=0cm 0cm 0cm 0cm}

\end{adjustbox}
\end{wrapfigure}

Note that, by definition of $A_1$, none of the vertices  of $Z_1$ and $Z_2$ can be corners of $A_1$. If one of $v(a+2,b)$ and $v(a+2,b+1)$ is a corner of $A_1$, then $F$ has to be eastern. But this implies that $Z_1+(2,0)\in P$, contradicting our assumption in the previous paragraph. So we may assume that $v(a+2,b)$ and $v(a+2,b+1)$ are not corners of $A_1$. By symmetry $v(a-1,b-1)$ and  $v(a-1,b)$ are not corners of $A_1$ either. By Proposition 4.4, $v(a+2,b)$ and $v(a+2,b+1)$ must belong to $(V(P_1) \cup V(P_2)) \setminus \{v_x,v_y\}$. By symmetry $v(a-1,b-1)$ and $v(a-1,b)$ must belong to $(V(P_1) \cup V(P_2)) \setminus \{v_x,v_y\}$ as well. But now, the existence of the segment $[v(a-1,b-1), v(a+2,b+1)]$ contradicts Lemma 4.9. See Figure 4.24. End of Case 1.1$(b_2)$. End of Case 1.1(b). End of Case 1.1.

\endgroup 

\null

\noindent \textit{CASE 1.2: Exactly one of $Z_1$ and $Z_2$ is adjacent to an end-box of $P$.} WLOG we may assume that $Z_1$ is adjacent to the end-box $Z_1'$. Then $Z_1'=Z_1+(0,1)$, $Z_1'=Z_1+(-1,0)$ or $Z_1'=Z_1+(1,0)$. Since the second and third are symmetric, we only need to check the first two.

\begingroup 
\setlength{\intextsep}{0pt}
\setlength{\columnsep}{20pt}
\begin{wrapfigure}[]{l}{0cm}
\begin{adjustbox}{trim=0cm 0cm 0cm 0cm}

\end{adjustbox}
\end{wrapfigure}

\null

\noindent \textit{CASE 1.2(a): $Z_1'=Z_1+(0,1)$.} By Lemma 4.6 (d), $e(a;b,b+1) \in H$ and $e(a+1; b,b+1) \in H$. Then $e(a; b-1,b) \notin H$ and $e(a+1; b-1,b) \notin H$.  By Lemma 4.6 (d), $e(a,a+1;b-1) \in H$. Now, $Z_2'=Z_2+(-1,0)$ or $Z_2'=Z_2+(1,0)$. By symmetry, we may assume WLOG that $Z_2'=Z_2+(-1,0)$. But now we're back to Case 1.1$(b_1)$. See Figure 4.25 (a). End of Case 1.2(a).

\null 

\noindent \textit{CASE 1.2(b): $Z_1'=Z_1+(-1,0)$.} Then $e(a;b,b+1) \notin H$, and by Lemma 4.6 (d), $e(a,a+1;b+1) \in H$. Now there are three possibilities for the $H$-neighbour $Z_2'$ of $Z_2$:  $Z_2'=Z_2+(-1,0)$, $Z_2'=Z_2+(0,-1)$, or $Z_2'=Z_2+(1,0)$.

\endgroup 

\null

\noindent \textit{CASE 1.2$(b_1$): $Z_2'=Z_2+(-1,0)$.} Lemma 4.6 (b) implies that $e(a;b,b+1)\in H$. But this contradicts the assumption of Case 1.2(b). See Figure 4.25 (b). End of Case 1.2$(b_1$).

\null 

\noindent \textit{CASE 1.2$(b_2$): $Z_2'=Z_2+(0,-1)$.} We note that this is the same as Case 1.1(a). See Figure 4.23 (a). End of Case 1.2 $(b_2$).

\null 

\begingroup 
\setlength{\intextsep}{0pt}
\setlength{\columnsep}{20pt}
\begin{wrapfigure}[]{r}{0cm}
\begin{adjustbox}{trim=0cm 0cm 0cm 0.25cm}

\end{adjustbox}
\end{wrapfigure}

\noindent \textit{CASE 1.2$(b_3$): $Z_2'=Z_2+(1,0)$.} By Lemma 4.6 (b), $Z_2'+(0,1) \in P$, $Z_2'+(0,-1) \in P$, $e(a+1;b,b+1) \in H$, $S_{\rightarrow}(a,b-1;a+1,b-2)\in H$ and $e(a+2;b-1,b) \in H$. By Lemma 4.6(c) we have that $e(a-1;b,b+1)\in H$ and that exactly one of $e(a-1,a;b)$ and $e(a-1,a;b+1)$ belongs to $H$. See Figure 4.26.

\endgroup 

\null 

\begingroup 
\setlength{\intextsep}{0pt}
\setlength{\columnsep}{20pt}
\begin{wrapfigure}[]{l}{0cm}
\begin{adjustbox}{trim=0cm 0cm 0cm 0cm}
%
\end{adjustbox}
\end{wrapfigure}

\noindent \textit{CASE 1.2$(b_3).(i)$: $e(a-1,a;b) \in H$.} By ($(**)$), $F$ is eastern. Note that the corner vertices of $A_1$ must be $v(a-1,b)$ and $v(a-1,b-1)$. It follows that $V(Z_2')$ contains no corners of $A_1$ and that $v(a-1,b+2) \in (V(P_1) \cup V(P_2))$. Our assumption that $e(a,a+1;b)\in P_0$ and Observation 4.8, now imply that $v(a+1,b-1)$ belongs to $V(P_0) \setminus (V(P_1) \cup V(P_2))$. By Proposition 4.4, we have that $v(a+2;b-1)$ and $v(a+2;b)$ belong to $ (V(P_1) \cup V(P_2))\setminus V(P_0)$.  But then, the existence of the segment $[v(a-1,b+1), v(a+2,b-1)]$ contradicts Lemma 4.9. See Figure 4.27. End of Case 1.2$(b_3).(i)$. 

\endgroup

\null

\begingroup 
\setlength{\intextsep}{0pt}
\setlength{\columnsep}{20pt}
\begin{wrapfigure}[]{r}{0cm}
\begin{adjustbox}{trim=0cm 0cm 0cm 0cm}

\end{adjustbox}
\end{wrapfigure}

\noindent \textit{CASE 1.2$(b_3).(ii)$: $e(a-1,a;b+1) \in H$.} Then we have that $e(a;b-1,b) \in H$, and, by ($(**)$), that $F$ is eastern. {
%Suppose $F$ is southern. \textcolor{red}{[Add fig. reference.]} Then we have that $v(a-1,b+1)=v_x$ and that $e(a-1;b-1,b) \in H$. But then $Z_1'+(0,-1)$ is a switchable box of $P$, contradicting Proposition 4.5 (a).
} Using the same arguments as in Case 1.2$(b_3).(i)$, we have that $e(a;b-1,b) \in (P_0)$ and that $v(a-1,b)$ and $v(a+2,b-1)$ belong to $(V(P_1) \cup V(P_2)) \setminus V(P_0)$. But then the existence of the segment $[v(a-1,b), v(a+2,b-1)]$ contradicts Lemma 4.9. See Figure 4.28. End of Case 1.2$(b_3).(ii)$. End of Case 1.2$(b_3)$. End of Case 1.2(b). End of Case 1.2. End of Case 1.

\endgroup

\begingroup 
\setlength{\intextsep}{0pt}
\setlength{\columnsep}{15pt}
\begin{wrapfigure}[]{l}{0cm}
\begin{adjustbox}{trim=0cm 0.75cm 0cm 0.5cm}
\begin{tikzpicture}[scale=1.5]

\begin{scope}[xshift=0cm]{
\draw[gray,very thin, step=0.5cm, opacity=0.5] (0,0) grid (1,1.5);

\fill[blue!40!white, opacity=0.5] (0,0) rectangle (0.5,1.5);

\fill[blue!40!white, opacity=0.5] (0.5,1) rectangle (1,1.5);

\fill[green!40!white, opacity=0.5] (0.5,0.5) rectangle (1,1);

\draw[blue, line width=0.5mm] (0.5,1.5)--++(0,-0.5)--++(0.5,0);

\draw[blue, line width=0.5mm] (0.5,0)--++(0,0.5)--++(0.5,0);

\draw[blue, line width=0.5mm] (1,1)--++(0,0.5);

\draw[blue, line width=0.5mm] (0,0.5)--++(0,0.5);

%black lines
{

\draw[black, line width=0.15mm] (0.45,0.7)--++(0.1,0);
\draw[black, line width=0.15mm] (0.45,0.75)--++(0.1,0);
\draw[black, line width=0.15mm] (0.45,0.8)--++(0.1,0);

\draw[black, line width=0.15mm] (0.95,0.7)--++(0.1,0);
\draw[black, line width=0.15mm] (0.95,0.75)--++(0.1,0);
\draw[black, line width=0.15mm] (0.95,0.8)--++(0.1,0);

\draw[black, line width=0.15mm] (0.7,1.45)--++(0,0.1);
\draw[black, line width=0.15mm] (0.75,1.45)--++(0,0.1);
\draw[black, line width=0.15mm] (0.8,1.45)--++(0,0.1);

}

% nodes
{
\node[left] at (0,0) [scale=1]
{\tiny{-2}};
\node[left] at (0,0.5) [scale=1]
{\tiny{-1}};
\node[left] at (0,1) [scale=1]
{\tiny{$b$}};
\node[left] at (0,1.5) [scale=1]
{\tiny{+1}};

\node[above] at (0.5, 1.5) [scale=1]
{\tiny{$a$}};
\node[above] at (1, 1.5) [scale=1]
{\tiny{+1}};
%\node[above] at (1.5, 2) [scale=1]
%{\tiny{+2}};

\node at (0.75,1.25) [scale=0.8] {\small{$Z_1$}};
\node at (0.75,0.75) [scale=0.8] {\small{$Z_2$}};
\node at (0.25,0.75) [scale=0.8] {\small{$Z_2'$}};
}

 \node[below, align=center, text width=2.5cm] at (0.5, -0.1) { Fig. 4.29.  Case 2.1(a).};

} \end{scope}

\end{tikzpicture}
\end{adjustbox}
\end{wrapfigure}

\null 

\noindent \textit{CASE 2: Exactly one of $Z_1$ and $Z_2$ is a box of $P$, and the other is a box of $N[P] \setminus P$.} For definiteness, assume that $Z_1 \in P$ and $Z_2 \in N[P] \setminus P$. Either $Z_1$ is an end-box of $P$ or it is not.

\null

\noindent \textit{CASE 2.1: $Z_1$ is an end-box of $P$.} Let $Z_2'$ be the $H$-neighbour of $Z_2$ in $P$. Then $Z_2'=Z_2+(0,-1)$, $Z_2'=Z_2+(-1,0)$, or $Z_2'=Z_2+(0,1)$. Since the second and third are symmetric, we only need to check the first two.

\null 

\noindent \textit{CASE 2.1(a): $Z_2'=Z_2+(0,-1)$.} The assumption of Case 2.1 implies that $Z_2'$ is not an end-box of $P$. By Lemma 4.6(b), $Z_2'+(0,1) \in P$ and $e(a;b,b+1) \in H$. But now, $e(a;b,b+1) \in P_0$, contradicts Lemma 4.3, and $e(a;b,b+1) \in E(P_1) \cup E(P_2)$ contradicts Lemma 4.10. See Figure 4.29. End of Case 2.1(a).

\endgroup 

\null 

\noindent \textit{CASE 2.1(b): $Z_2'=Z_2+(0,-1)$.} This is the same as Case 1.1(a). End of Case 2.1(b). End of Case 2.1.

\null 

\noindent \textit{CASE 2.2: $Z_1$ is not an end-box of $P$.} Let $Z_2'$ be the $H$-neighbour of $Z_2$ in $P$. Then $Z_2'=Z_2+(0,-1)$ or $Z_2'=Z_2+(-1,0)$ or $Z_2'=Z_2+(0,1)$. Since the second and third are symmetric, we only need to check the first two.

\null 

\begingroup 
\setlength{\intextsep}{0pt}
\setlength{\columnsep}{20pt}
\begin{wrapfigure}[]{r}{0cm}
\begin{adjustbox}{trim=0.5cm 0.5cm 0cm 0.5cm}

\end{adjustbox}
\end{wrapfigure}

\noindent \textit{CASE 2.2(a): $Z_2'=Z_2+(0,-1)$.} By Lemma 4.6 (b) and (c) we have that $e(a,a+1;b-2) \in H$. By Lemma 4.6(d), $e(a;b-1,b) \in H$ and $e(a+1;b-1,b) \in H$. Then $e(a;b,b+1) \notin H$ and $e(a+1;b,b+1) \notin H$. If $e(a,a+1;b+1) \in H$, then $Z_1$ is a switchable box of $P$ so we may assume that $e(a,a+1;b+1) \notin H$. Then $S_{\rightarrow}(a-1,b+1;a,b+2) \in H$ and $S_{\downarrow}(a+1,b+2;a+2,b+1) \in H$. Note that, by definition of $A_1$, $v(a+1,b)$ is not a corner of $A_1$. It follows that $e(a+1;b-1,b)\in E(P_0)$. Now, if $Z_2+(1,0) \in N[P] \setminus P$, then we're back to Case 1, so we may assume that $Z_2+(1,0) \in P$.

In the same way we find that $v(a,b)$ is not a corner of $A_1$, $e(a;b-1,b)\in E(P_0)$ and $Z_2+(-1,0) \in P$. See Figure 4.30. Now, either $Z_2'$ is an end-box of $P$ or it is not.

\begingroup 
\setlength{\intextsep}{0pt}
\setlength{\columnsep}{20pt}
\begin{wrapfigure}[]{l}{0cm}
\begin{adjustbox}{trim=0cm 0.5cm 0cm 0.5cm}

\end{adjustbox}
\end{wrapfigure}

\null

\noindent \textit{CASE 2.2($a_1$): $Z_2'$ is an end-box of $P$.} By Lemma 4.6(c), exactly one of $e(a;b-2,b-1)$ and $e(a+1;b-2,b-1)$ belong to $H$. By symmetry, we may assume WLOG that $e(a+1;b-2,b-1) \in H$. By ($(**)$), we have that $F$ is northern with $v(a+1,b-2)$ and $v(a+2,b-2)$ being the corners of $A_1$. But then $Z_2+(1,0)$ is a switchable box of $P$. See Figure 4.31 (a). End of Case 2.2($a_1$).

\null

\noindent \textit{CASE 2.2($a_2$): $Z_2'$ is not an end-box of $P$.} By Lemma 4.6(b), $Z_2'+(1,0) \in P$ and $e(a+1,a+2;b-1)\in H$. By definition of $A_1$, $v(a+1,b-1)$, cannot be a corner of $A_1$. This implies that $e(a+1,a+2;b-1)\in E(P_0)$. But this contradicts Lemma 4.3.  See Figure 4.31 (b). End of Case 2.2($a_2$). End of Case 2.2(a).

\null

\noindent \textit{CASE 2.2(b): $Z_2'=Z_2+(0,-1)$.} By Lemma 4.6(d), $e(a,a+1;b-1) \in H$. By Lemma 4.6 (b) and (c),  $e(a-1; b-1,b) \in H$. Now, either $Z_2'$ is an end-box of $P$ or it is not. 

\null 

\noindent \textit{CASE 2.2($b_1$): $Z_2'$ is not an end-box of $P$.} By Lemma 4.6(b), $Z_2'+(0,1) \in P$ and $e(a;b,b+1) \in H$. The assumption of Case 2.2 and the definition of $A_1$ imply that $v(a,b)$ is not a corner of $A_1$. But then $e(a;b,b+1) \in E(P_0)$, which contradicts Lemma 4.3.  See Figure 4.32 (a). End of Case 2.2($b_1$).

%pagemarker

\begingroup 
\setlength{\intextsep}{0pt}
\setlength{\columnsep}{20pt}
\begin{wrapfigure}[]{r}{0cm}
\begin{adjustbox}{trim=0cm 0.0cm 0cm 0.0cm}

\end{adjustbox}
\end{wrapfigure}

\noindent \textit{CASE 2.2($b_2$): $Z_2'$ is end-box of $P$.} By Lemma 4.6(c), exactly one of $e(a-1,a;b-1)$ and $e(a-1,a;b)$ belongs to $H$. Note that if the former, then we arrive at the same contradiction as we did in Case 2.2($b_1$), so we may assume the latter. Now, by ($(**)$),  $F$ is eastern. An argument similar to that in Case 2.2 ($a_2$) in Lemma 4.10 shows that the corners of $A_1$ must be $v(a-1,b)$ and $v(a-1,b+1)$. But then $Z_1$ is a middle-box of $A_1$ in $P$, contradicting Corollary 4.2 (b). See Figure 4.32 (b).  End of Case 2.2($b_2$). End of Case 2.2(b). End of Case 2.2. End of Case 2. $\square$

%pagemarker 

\noindent \textbf{Proposition 4.11.} Let $H$ be a Hamiltonian path of a polyomino $G$ and let $F$ be a \index{looping fat path}looping fat path in $G$. Then $F$ has no polyking junctions whenever:

(i) \ $G$ is a simply connected polyomino and both end-vertices of $H$ are on $B(G)$, or

(ii) $F$ is sturdy.

\begingroup 
\setlength{\intextsep}{0pt}
\setlength{\columnsep}{20pt}
\begin{wrapfigure}[]{l}{0cm}
\begin{adjustbox}{trim=0cm 0cm 0cm 0cm}

\end{adjustbox}
\end{wrapfigure}

\noindent \textit{Proof.} The proofs of (i) and (ii) overlap a lot so we combine them. We will prove the contrapositive. Assume that $F$ has a polyking junction. We will show that this assumption leads to a contradiction under (i) and do the same for (ii). We remark that we won't have to use (i) and (ii) until Case 2.

Let $v(a,b)$ be a polyking junction of $F$ with $Z=R(a-1,b-1) \in F$ and $Z+(1,1) \in F$. Either $v(a,b)$ is an end-vertex of $H$ or it is not. 

\noindent Suppose that $v(a,b)$ is an end-vertex of $H$. By symmetry, we may assume that $e(a-1,a;b) \in H$. Then $e(a;b,b+1)\notin H$ and $e(a,a+1;b)\notin H$. But now, $Z+(1,1) \in P$ conflicts with FPC-1 and $Z+(1,1) \in N[P] \setminus P$ contradicts Lemma 4.6 (d). See Figure 4.33 (a). It remains to check the case where $v(a,b)$ is not an end-vertex of $H$.

\noindent By symmetry, there are three possibilities to consider: $S_{\rightarrow}(a-1,b;a,b-1) \in H$, $S_{\rightarrow}(a-1,b;a,b+1) \in H$, and the case where $e(a-1,a;b)$ and $e(a,a+1;b)$ belong to $H$. We note that if $S_{\rightarrow}(a-1,b;a,b-1) \in H$ then $e(a;b,b+1)\notin H$ and $e(a,1+;b)\notin H$. See Figure 4.33 (b). Now, $Z+(1,1) \in N[P] \setminus P$ contradicts Lemma 4.6(d); and if $Z+(1,1) \in P$, then, by Lemma 4.6(a), at least one of $Z+(1,0)$ and $Z+(0,1)$ belongs to $P$, contradicting that $v(a,b)$ is a polyking junction. It remains to check the other two possibilities.

\endgroup 

\null

\begingroup 
\setlength{\intextsep}{0pt}
\setlength{\columnsep}{20pt}
\begin{wrapfigure}[]{r}{0cm}
\begin{adjustbox}{trim=0cm 0.25cm 0cm 0.5cm}
\begin{tikzpicture}[scale=1.5]

\begin{scope}[xshift=0cm]{
\draw[gray,very thin, step=0.5cm, opacity=0.5] (0,0) grid (1.5,1.5);

\fill[blue!40!white, opacity=0.5] (0,0) rectangle (0.5,0.5);
\fill[blue!40!white, opacity=0.5] (1,1) rectangle (1.5,1.5);

\fill[green!40!white, opacity=0.5] (1,0.5) rectangle (1.5,1);

\fill[green!40!white, opacity=0.5] (0.5,0) rectangle (1,0.5);

\draw[blue, line width=0.5mm] (0.5,0)--++(0.5,0);
\draw[blue, line width=0.5mm] (0.5,0.5)--++(0.5,0)--++(0,0.5);
\draw[blue, line width=0.5mm] (1.5,0.5)--++(0,0.5);

%\draw[fill=blue, opacity=1] (1,0.5) circle [radius=0.05];

%black lines
{
\draw[black, line width=0.15mm] (0.95,0.2)--++(0.1,0);
\draw[black, line width=0.15mm] (0.95,0.25)--++(0.1,0);
\draw[black, line width=0.15mm] (0.95,0.3)--++(0.1,0);

\draw[black, line width=0.15mm] (1.2,0.45)--++(0,0.1);
\draw[black, line width=0.15mm] (1.25,0.45)--++(0,0.1);
\draw[black, line width=0.15mm] (1.3,0.45)--++(0,0.1);

}

% nodes
{
\node[left] at (0,0) [scale=1]
{\tiny{-1}};
\node[left] at (0,0.5) [scale=1]
{\tiny{$b$}};
\node[left] at (0,1) [scale=1]
{\tiny{+1}};

\node[below] at (0.5,0) [scale=1]
{\tiny{-1}};
\node[below] at (1,0) [scale=1]
{\tiny{$a$}};
\node[below] at (1.5,0) [scale=1]
{\tiny{+1}};

\node at (0.75,0.25) [scale=0.8] {\small{$Z$}};
}

\node[below, align=center, text width=4.5cm] at (0.75, -0.25) { Fig. 4.34. Case 1.};

} \end{scope}

\end{tikzpicture}
\end{adjustbox}
\end{wrapfigure}

\noindent \textit{CASE 1: $S_{\rightarrow}(a-1,b;a,b+1) \in H$.} It follows that $e(a;b-1,b)\notin H$ and $e(a,a+1;b)\notin H$. Since $Z+(1,0)$ is not in $F$, and it is $H$-adjacent to $Z$, $Z$ must be in $N[P] \setminus P$. Similarly, $Z+(1,1) \in N[P] \setminus P$. Now Lemma 4.6 (d) implies that $e(a-1,a;b-1) \in H$ and $e(a+1;b,b+1) \in H$. Then we must have $Z+(-1,0) \in P$ and $Z+(1,2) \in P$. But then there is an $H$-cycle $Z+(-1,0), \ldots, Z+(1,2), Z+(1,1), Z+(1,0), Z,Z+(-1,0)$ which contradicts Proposition 1.2.1. See Figure 4.34. End of Case 1.

\null 

\endgroup

\noindent \textit{CASE 2: $e(a-1,a;b)$ and $e(a,a+1;b)$ belong to $H$.} Then $e(a;b-1,b)\notin H$ and $e(a;b,b+1)\notin H$. As in Case 1, we must have that $Z \in N[P] \setminus P$ and $Z+(1,1) \in N[P] \setminus P$. By Lemma 4.6 (d), $e(a-1,a;b-1) \in H$ and $e(a,a+1;b+1) \in H$. Then the $H$-neighbour of $Z$ in $P$ must be $Z+(-1,0)$ and the $H$-neighbour of $Z+(1,1)$ in $P$ must be $Z+(2,1)$.

\null

\begingroup 
\setlength{\intextsep}{0pt}
\setlength{\columnsep}{20pt}
\begin{wrapfigure}[]{l}{0cm}
\begin{adjustbox}{trim=0cm 0cm 0cm 0cm}

\end{adjustbox}
\end{wrapfigure}

\noindent Assume (i) holds. By Lemma 1.1.11, since $G$ is a simply connected polyomino, $G$ has no polyking junctions. Then at least one of $Z+(1,0)$ and $Z+(0,1)$ must belong to $G$. By symmetry we may assume without loss of generality that $Z+(1,0)$ belongs to $G$. Let $J_F$ be the $H$-component of $G$ that contains $F$. Then both $Z+(1,1)$ and $Z+(1,0)$ must belong to $J_F$. Then $J_F$ is self-adjacent. By Lemma 1.3.12, some edge of a main trail of $J_F$ has an end-vertex of $H$ incident on it in $G \setminus B(G)$. But this contradicts the assumption that both end-vertices of $H$ are on $B(G)$. See Figure 4.35. End of proof for (i).

\endgroup 

\null

\noindent Orient $H$. Let $v(a+1,b+1)=v_s$ and $v(a+1,b)=v_t$. Without loss of generality assume that $s<t$. Either $v(a,b)=v_{t+1}$ or $v(a,b)=v_{t-1}$.

\null 

\noindent \textit{CASE 2.1: $v(a,b)=v_{t+1}$.} Let $Q$ be the cycle in $G$ consisting of the subpath $P(v_s,v_t)$ of $H$ and the edge $\{v_s,v_t\} \in G \setminus H$, and let $U$ be the region bounded by $Q$. Now, either $Z+(2,1) \in U$, or $Z+(2,1) \notin U$.

\begingroup 
\setlength{\intextsep}{0pt}
\setlength{\columnsep}{20pt}
\begin{wrapfigure}[]{r}{0cm}
\begin{adjustbox}{trim=0cm 0cm 0cm 0.25cm}

\end{adjustbox}
\end{wrapfigure}

\null 

\noindent \textit{CASE 2.1(a): $Z+(2,1) \in U$.} Consider the $H$-subpath $P(Z_1, Z_r)$ of $P(X,Y)$, where $Z_1=Z+(2,1)$ and $Z_r=Z+(-1,0)$. For $j \in \{1, \ldots,r\}$ let $c_j$ be the center of the box $Z_j$. Note that the edges $(v_t, v_{t+1})$ and $(v_{t+1}, v_{t+2})$ of $H$ do not belong to $Q$. Then the segment $[c_1, c_r]$ (red in Figure 4.36) intersects $Q$ exactly once at the edge $\{v_s,v_t\}$. By Corollary 1.1.5, $c_1$ and $c_r$ are on different sides of $Q$. This implies that $Z_r$ is in $G \setminus U$. Consider the polygonal path $P(c_1,c_r)$ in the plane that has as edges the segments $[c_1,c_2], \ldots, [c_{r-1}, c_r]$. Note that each edge $[c_j, c_{j+1}]$ of $P(c_1,c_r)$ bisects the gluing edge between $Z_j$ and $Z_{j+1}$, and intersects no other edge of $G$. Since $c_1$ and $c_r$ are on different sides of $Q$, by JCT, $P(c_1,c_r)$ intersects $Q$. Then $P(c_1,c_r)$ either intersects $Q$ at an edge of $P(v_s,v_t)$ or at $\{v_s,v_t\}$. The former is not possible, because gluing edges of $P(Z_1,Z_r)$ are not in $H$. The latter implies that $Z+(1,1) \in P$, contradicting the assumption that $Z+(1,1) \in N[P] \setminus P$. End of Case 2.1(a). 

\endgroup 

\null 

\noindent \textit{CASE 2.1(b): $Z+(2,1) \notin U$.} The proof uses the same argument as above, so we omit it. End of Case 2.1(b). End of Case 2.1.

\null

\begingroup 
\setlength{\intextsep}{0pt}
\setlength{\columnsep}{20pt}
\begin{wrapfigure}[]{l}{0cm}
\begin{adjustbox}{trim=0cm 0cm 0cm 0.25cm}

\end{adjustbox}
\end{wrapfigure}

\noindent \textit{CASE 2.2: $v(a,b)=v_{t-1}$.} Assume that (ii) holds. Note that at least one of $Z+(-1,0)$ and $Z+(2,1)$ is not an end-box of $P$. By symmetry we may assume WLOG that $Z+(-1,0)$ is not an end-box of $P$. By Lemma 4.6 (b) and (c), $Z+(-1,1)\in P$, $Z+(-1,-1)\in P$, $e(a-2;b-1,b)\in H$, and $e(a+2;b,b+1)\in H$. Now, either $Z+(2,1)$ is an end-box of $P(X,Y)$, or it is not.

\null

\noindent \textit{CASE 2.2(a): $Z+(2,1)$ is not an end-box of $P(X,Y)$.} Then $Z+(2,2) \in P$ and $Z+(2,0) \in P$. If $Z+(1,1)$ is parallel, then, after $\text{Sw}(Z+(1,1))$, $Z+(2,1) \in P$ is switchable, violating FPC-4, so we may assume that $Z+(1,1)$ is anti-parallel. Note that this implies that $v(a,b+1)$ is not an end-vertex of $H$. Then exactly one of $e(a-1,a;b+1)$ and $e(a;b+1,b+2)$ belongs to $H$. If the former, then $Z+(1,1) \mapsto Z$ is a valid double-switch move after which $Z+(-1,0)\in P$ is switchable; and if the latter, then $Z+(1,1) \mapsto Z+(1,2)$, is a valid flip move after which $Z+(2,1) \in P$ is switchable. See Figure 4.37. End of Case 2.2(a).

\endgroup 

\null

\noindent \textit{CASE 2.2(b): $Z+(2,1)$ is an end-box of $P(X,Y)$.}  By Lemma 4.6(c), exactly one of $e(a+1,a+2;b)$ and $e(a+1,a+2;b+1)$ belongs to $H$.

\null

\begingroup 
\setlength{\intextsep}{0pt}
\setlength{\columnsep}{20pt}
\begin{wrapfigure}[]{r}{0cm}
\begin{adjustbox}{trim=0cm -0.5cm 0cm 0cm}

\end{adjustbox}
\end{wrapfigure}

\noindent \textit{CASE 2.2$(b_1)$: $e(a+1,a+2;b) \in H$ and $e(a+1,a+2;b+1) \notin H$.} Then $F$ can be northern or western. By $(**)$ we may assume that $F$ is western. Then the switchable middle-box of $F$ is $W=Z+(2,0)$. Let $Z+(2,1)$ be the end-box $Y$ of $P(X,Y)$. Note that $e(a+1,a+2;b+1) \notin H$, so $e(a+1;b+1,b+2) \in H$. If $Z+(1,1)$ is parallel then after $\text{Sw}(Z+(1,1))$, $W \mapsto Y$ is a valid flip move, violating FPC-4. So we may assume that $Z+(1,1)$ is anti-parallel. See Figure 4.38. As before, this implies that $v(a,b+1)$ is not an end-vertex of $H$. Then exactly one of $e(a-1,a;b+1)$ and $e(a;b+1,b+2)$ belongs to $H$. If the former, then $Z+(1,1) \mapsto Z$ is a valid move after which $Z+(-1,0) \in P$ is switchable; and if the latter, then after $Z+(1,1) \mapsto Z+(1,2)$, $W \mapsto Y$ is a valid flip move, violating FPC-4. End of Case 2.2$(b_1)$.

\endgroup

\null

\noindent \textit{CASE 2.2$(b_2)$: $e(a+1,a+2;b) \notin H$ and $e(a+1,a+2;b+1) \in H$.} Then $F$ can be northern or western. By $(**)$ we may assume that $F$ is western. Then the switchable middle-box of $F$ is $W=Z+(2,2)$. Let $Z+(2,1)$ be the end-box $X$ of $P(X,Y)$. If $Z+(1,1)$ is parallel then after $\text{Sw}(Z+(1,1))$, $W \mapsto X$ is a valid flip move, violating FPC-4.
\begingroup 
\setlength{\intextsep}{0pt}
\setlength{\columnsep}{15pt}
\begin{wrapfigure}[]{l}{0cm}
\begin{adjustbox}{trim=0cm 0cm 0cm 0cm}

\end{adjustbox}
\end{wrapfigure}
So we may assume that $Z+(1,1)$ is anti-parallel. Once again, this implies that $v(a,b+1)$ is not an end-vertex of $H$. Then exactly one of $e(a-1,a;b+1)$ and $e(a;b+1,b+2)$ belongs to $H$.

If the former, then $Z+(1,1) \mapsto Z$ is a valid double-switch move after which $Z+(-1,0) \in P$ is switchable (Figure 4.39 (a)), so we may assume the latter (Figure 4.39 (b)).  Now, if $e(a,a+1;b+3) \in H$, then $W \mapsto W+(-1,1)$ is a valid transpose move, violating FPC-4 so we may assume that $e(a,a+1;b+3) \notin H$. Then we must have that $S_{\downarrow}(a+1,b+4;a+2,b+3) \in H$ and that $W+(-1,1) \in P$. Then $S_{\rightarrow}(a-1,b+3;a,b+4)$ is in $H$. The $H$-neighbour of $W+(-1,1)$ in $P$ cannot be $W+(-1,2)$, since then $W+(-1,2)$ would have to be switchable, so the $H$-neighbour of $W+(-1,1)$ in $P$ must be $W+(-2,1)$. Observe that this implies that $W+(-2,0)$ belongs to $G$, since all four vertices of $W+(-2,0)$ are in $G$. Similarly,  $Z+(0,1)$ must belong to $G$ as well. But then $P$ must be the $H$-path $X, X+(-1,0), X+(-2,0), X+(-2,1), X+(-2,2),X+(-1,2), X+(0,2)$, contradicting that $Z+(-1,0) \in P$. See Figure 4.39(b). End of Case 2.2$(b_2)$. End of Case 2.2(b). End of Case 2.2. End of Case 2. $\square$

\endgroup

\null 

\noindent \textbf{Definitions.} Let $G$ be a polyomino, let $H$ be a Hamiltonian path or cycle of $G$ and let
$J$ be an $H$-subtree of an $H$-component of $G$. We define a \index{shadow edge|textbf}\textit{shadow edge} of $J$ to be a boundary edge of $J$ that is not in $H$. We denote the set of all shadow edges of $J$ by $E_{\text{sh}}(J)$. 

Suppose that $J$ is a looping fat path $G\langle N[P(X,Y)] \rangle$ that is also non-self adjacent, and has no polyking junctions (but is not necessarily sturdy). Then we say that $J$ is a \index{standard looping fat path|textbf}\textit{standard looping fat path}. We note that shadow edges of $J$ cannot be incident on boxes of $P$.

\null

\noindent \textbf{Corollary 4.12.}  Let $G$ be a polyomino, let $H$ be a Hamiltonian path or cycle of $G$ and let $F$ be a \index{sturdy looping fat path}looping fat path of $G$.

(a) If $F$ is a \index{sturdy looping fat path}sturdy looping fat path, then $F$ is a standard looping fat path.

(b) If $F$ is a standard looping fat path, then $F$ is a simply connected polyomino.

\null

\noindent \textit{Proof of (a).} Let $F$ be a sturdy looping fat path. By Propositions 4.7 and 4.11, $F$ is non-self-adjacent and has no polyking junctions, so it is a standard looping fat path.

%pagemarker

\noindent \textit{Proof of (b).} Let $F$ be a standard looping fat path. By definition, $F$ is an $H$-subtree in an $H$-component of $G$, that is non-self adjacent and has no polyking junctions. By Lemma 1.3.6, it follows that $F$ is a simply connected polyomino. 

\null 

\noindent \textbf{Lemma 4.13.} Let $G$ be a polyomino, let $H$ be a Hamiltonian path of $G$ and let $F=G\langle N[P(X,Y)] \rangle $ be a standard looping fat path in $G$. Then $E_{\text{sh}}(F) \cup (E(F) \cap E(H))$ is a Hamiltonian cycle of $F$ and the boundary of $F$.

\null 

\noindent \textit{Proof.} Let $E'=E_{\text{sh}}(F)\cup (E(F) \cap E(H))$. It follows from Lemma 4.6 (b), (c) and (d) that each vertex of $F$ must be incident on some edge of $E(F) \cap E(H)$. So, it is enough to show that $E'$ is a cycle.

By definition, $F$ is not self-adjacent and has no polyking junctions. Corollary 4.12 (b), $F$ is a simply connected polyomino. By Corollary 1.1.10, $B(F)=B_0(F)$. By Lemma 1.1.1, $B(F)$ is a cycle. Then the conclusion will follow if we can show that $B(F)=E'$.

Suppose that $e \in B(F)$ and assume that $e$ is incident on boxes $Z \in \text{Boxes}(F)$ and $Z' \in \text{Boxes}(G_{-1}) \ $ $\setminus \text{Boxes}(F)$. Then $e \in H$ or $e \notin H$. If $e \in H$, since $e$ is incident on the box $Z$ of $F$, we have that $e \in E(F) \cap E(H) \subset E'$; and if $e \notin H$, by definition, $e \in E_{\text{sh}}(F) \subset E'$.

Suppose that $e \in E'$, and assume that $e$ is incident on boxes $Z$ and $Z'$ of $G_{-1}$. For definiteness, assume that $Z \in \text{Boxes}(F)$. If $e \notin H$, then $e$ is a shadow edge of $F$, so by definition, $e\in B(F)$. If $e \in H$, then, by Proposition 4.7, $F$ is non-self-adjacent, so we must have $Z' \in \text{Boxes}(G_{-1}) \setminus \text{Boxes}(F)$. Thus $e \in B(F)$. $\square$

\null 

\noindent \textbf{Proposition 4.14.} Let $G$ be a polyomino, let $H$ be a Hamiltonian path of $G$ and let $F=G\langle N[P(X,Y)] \rangle $ be a \index{standard looping fat path}standard looping fat path in $G$ following an edge $e_F$. Then $B(F)$ does not have consecutive colinear edges other than the left and right collinear edges in the $A_1$-type of $F$ following $e_F$.

\null

\noindent \textit{Proof.}  We need to check that $B(F)$ does not have consecutive colinear edges in the case where one of those colinear edges is one of the  \index{left colinear edges of an A1 type@left colinear edges of an $A_1$-type}left or \index{right colinear edges of an A1 type@right colinear edges of an $A_1$-type}right colinear edges of the $A_1$-type of $F$ and in the case where neither of those colinear edges is one of the left or right colinear edges of the $A_1$-type of $F$. We divide the proof into Lemmas 4.15 and 4.16.

\null

\noindent \textbf{Lemma 4.15.} The shadow of $F$ does not have a pair of consecutive colinear edges in the case where one edge of the pair is one of the left or right colinear edges of the $A_1$-type that follows $e_F$.

\setlength{\intextsep}{0pt}
\setlength{\columnsep}{20pt}
\begin{center}
\begin{adjustbox}{trim=0cm 0.5cm 0cm 0.5cm}
\begin{tikzpicture}[scale=1.5]
\begin{scope}[xshift=0cm]{
\draw[gray,very thin, step=0.5cm, opacity=0.5] (0,0) grid (1.5,1);

\draw[green!50!black, dotted, line width=0.5mm] (1,0.5)--++(0,0.5);

%\draw[orange, dotted, line width=0.5mm] (0.5,1)--++(0.5,0);
%\draw[orange, dotted, line width=0.5mm] (1,1.5)--++(0.5,0);

\draw[blue, line width=0.5mm] (0,0.5)--++(0.5,0)--++(0,-0.5);
\draw[blue, line width=0.5mm] (1.5,0.5)--++(-0.5,0)--++(0,-0.5);
\draw[blue, line width=0.5mm] (0.5,1)--++(0.5,0);

\node[left] at (0,0.5) [scale=1]
{\tiny{$\ell{+}2$}};

\node[below] at (1, 0) [scale=1]
{\tiny{$k$}};

\node at (1.25,0.25) [scale=0.8] {\small{X}};
\node at (0.25,0.25) [scale=0.8] {\small{Y}};

 \node[right, align=left, text width=9cm] at (1.75, 0.5) { Fig. 4.40. The pair of  consecutive colinear edges is $e(k;l+1,l+2)$ and $e(k;l+2,l+3)$.};

} \end{scope}

\end{tikzpicture}
\end{adjustbox}
\end{center}

\noindent \textit{Proof.} For definiteness, assume that $F$ is a standard southern looping fat path following the edge $e(k-1,k;l+3)=e_F$ southward, and let $X=R(k,l+1)$. BWOC assume that $B(F)$ does have consecutive colinear edges and one of those edges is one of the right or left colinear edges of the $A_1$-type that follows $e_F$. For definiteness, assume that one of those consecutive colinear edges of $B(F)$ is one of the right colinear edges of the $A_1$-type that follows $e_F$. Then either $e(k;l+1,l+2)$, $e(k;l+2,l+3)$ is a pair of consecutive colinear edges or $e(k;l-1,l)$, $e(k;l,l+1)$ is a pair of consecutive colinear edges. If the former, we note that $e(k;l+2,l+3) \in B(F)$ contradicts Lemma 4.13, so we only need to check the latter. See Figure 4.40. 

Suppose then that $e(k;l-1,l)$, $e(k;l,l+1)$ is a pair of consecutive colinear edges of $B(F)$. Note that $X$ and $X+(0,-1)$ belong to $F$. If the edge $e(k;l-1,l)$ is in $B(F)$ then it is either a shadow edge or it belongs to $H$. We will show that both cases lead to contradictions.

\begingroup 
\setlength{\intextsep}{0pt}
\setlength{\columnsep}{20pt}
\begin{wrapfigure}[]{l}{0cm}
\begin{adjustbox}{trim=0cm 0.5cm 0cm 0.25cm}

\end{adjustbox}
\end{wrapfigure}
 
\null 

\noindent \textit{CASE 1: $e(k;l-1,l) \in H$.} Then exactly one of $X+(0,-1)$ and $X+(1,0)$ must belong to $P$.

\null

\noindent \textit{CASE 1.1: $X+(0,-1)\in P$.} Note that if $e(k+1;l,l+1) \in H$, then $X+(0,-1)$ is switchable, so we only need to check the case where $e(k+1;l,l+1) \notin H$. Then $S_{\uparrow}(k+1,l-1;k+2,l) \in H$ and $S_{\downarrow}(k+1,l+2;k+2,l+1) \in H$. Now exactly one of $X+(0,-2)$ and $X+(1,-1)$ belong to $P$. But, since neither can be the other end-box of $P$, either one would have to be a switchable box. See Figure 4.41 (a). End of Case 1.1

{
%Suppose that $X+(0,-2) \in P$. Then $e(k,k+1;l-1) \in H$ or $e(k,k+1;l-1) \notin H$. If the former, then $X+(-2,0)$ must be an end-box of $P$, but this contradicts the fact that the other end-box of $P$ is $Y=X+(-2,0)$; and if the latter then $X+(0,-2)$ is switchable. The case where $X+(1,-1)$ belongs to $P$ is similar so we omit the proof. End of Case 1.1
}

\endgroup

\null 

\noindent \textit{CASE 1.2: $X+(1,0)\in P$.} Then $X+(0,-1) \in N[P] \setminus P$ and, by Lemma 4.6(a), $X+(0,-2) \notin F$. Then $e(k,k+1;l)$ is a shadow edge of $F$. But then $\deg_{B(F)}(v(k,l))=3$, contradicting Lemma 4.13. See Figure 4.41 (b). End of Case 1.2.

\null

\noindent \textit{CASE 2: $e(k;l-1,l) \in E_{\text{sh}}(F)$.}  Then exactly one of $X-(1,2)$ and $X-(0,2)$ belongs to $N[P] \setminus P$.

\null 

\begingroup 
\setlength{\intextsep}{0pt}
\setlength{\columnsep}{20pt}
\begin{wrapfigure}[]{r}{0cm}
\begin{adjustbox}{trim=0cm 0cm 0cm 0cm}

\end{adjustbox}

\end{wrapfigure}

%pagemarker

\noindent If $X-(1,2)\in N[P] \setminus P$ then Lemma 4.13 implies that $e(k-1,k;l) \notin B(F)$, but this is a contradiction to Lemma 4.6(d); and if $X-(0,2) \in N[P] \setminus P$, by Lemma 4.6(d), $e(k,k+1;l) \in H$. But then, since $X-(0,1)$ must belong to $F$, this contradicts the assumption that $F$ is non-self-adjacent. See Figure 4.42. End of Case 2. $\square$

\endgroup 

\null

\noindent \textbf{Lemma 4.16.} The shadow of $F$ does not have a pair of consecutive colinear edges in the case where neither edge of the pair is one of the left or right colinear edges of the $A_1$-type that follows $e_F$.

\null

\noindent  \textit{Proof.} For a contradiction, assume that there is a pair of consecutive collinear edges in $B(F)$ where neither edge of the pair is a left or right colinear edge of the $A_1$-type that follows $e_F$. For definiteness, we may assume that these edges are the horizontal edges $e(a,a+1;b)$ and $e(a+1,a+2;b)$, and that $Z=R(a,b)$ belongs to $F$. By Lemma 4.13, $e(a+1;b,b+1) \notin B(F)$, otherwise $\deg_{B(F)}(v(a+1,b))=3$. Now, either $e(a,a+1;b) \in E_{\text{sh}}(F)$, or $e(a,a+1;b) \in H$.

\begingroup 
\setlength{\intextsep}{0pt}
\setlength{\columnsep}{20pt}
\begin{wrapfigure}[]{l}{0cm}
\begin{adjustbox}{trim=0cm 0.5cm 0cm 0cm}

\end{adjustbox}
\end{wrapfigure}

\null 

\noindent \textit{CASE 1: $e(a,a+1;b)$ belongs to $E_{\text{sh}}(F)$.} Since $Z$ is incident on a shadow edge, 
$Z$ belongs to $N[P] \setminus P$, and $Z+(0,-1) \notin F$. 
%Since $e(a+1,a+2;b) \in B(F)$, at least one of $Z+(1,0)$ and $Z+(1,-1)$ belongs to $F$. The assumption that $F$ has no polyking junctions implies that $Z+(1,0) \in F$. 
Since $e(a+1; b, b+1) \notin B(F)$, $e(a+1; b, b+1)$ is not in $H$ either. But then $Z$ can be neither switchable nor a leaf, which contradicts Lemma 4.6 (d). See Figure 4.43 (a). End of Case 1.

\null 

\noindent  \textit{CASE 2. $e(a,a+1;b)$ belongs to $H$.} Since $F$ is non-self-adjacent, $Z+(0,-1) \notin F$. Since $F$ is non-self-adjacent, exactly one of $Z+(1,0)$ and $Z+(1,-1)$ belongs to $G$. 
%(Note that this holds regardless of e(a+1,a+2;b) being in H or in E_{\text{sh}}(F)).
Since $F$ has no polyking junctions, it must be the case that $Z+(1,0)\in F$ and $Z+(1,-1) \notin F$. Note that if $e(a+1,a+2;b) \in E_{\text{sh}}(F)$, then $Z+(1,0)$ must belong to $N[P] \setminus P$. But then, as in Case 1, $Z+(1,0)$ can be neither a leaf nor a switchable box, contradicting Lemma 4.6(d). See Figure 4.43 (b). Therefore we may assume that $e(a+1,a+2;b)$ also belongs to $H$.

\endgroup

Now, if $Z$ and $Z+(1,0)$ both belong to $N[P] \setminus P$, then there is an $H$-cycle on $Z,Z'\ldots,Z'', Z+(1,0)$, where $Z'$ and $Z''$ are the $H$-neighbours in $P$ of $Z$ and $Z{+}(1,0)$, respectively, which contradicts Proposition 1.2.1. It follows, either $Z$ and $Z{+}(1,0)$ both belong to $P$, or exactly one of them belongs to $P$ and the other belongs to $N[P] \setminus P$.

\null 

\begingroup 
\setlength{\intextsep}{0pt}
\setlength{\columnsep}{20pt}
\begin{wrapfigure}[]{r}{0cm}
\begin{adjustbox}{trim=0cm -0.25cm 0cm 0.5cm}
\begin{tikzpicture}[scale=1.5]
\begin{scope}[xshift=0cm]{
\draw[gray,very thin, step=0.5cm, opacity=0.5] (0,0) grid (1,0.5);

\fill[blue!50!white, opacity=0.5](0,0) rectangle (1,0.5);

\draw[blue, line width=0.5mm] (0,0)--++(1,0);
\draw[blue, line width=0.5mm] (0,0.5)--++(0.5,0);

%black lines
{
\draw[black, line width=0.15mm] (-0.05,0.2)--++(0.1,0);
\draw[black, line width=0.15mm] (-0.05,0.25)--++(0.1,0);
\draw[black, line width=0.15mm] (-0.05,0.3)--++(0.1,0);

}

\node[left] at (0,0) [scale=1]
{\tiny{$b$}};

\node[left] at (0,0.5) [scale=1]
{\tiny{$+1$}};

\node[below] at (0,0) [scale=1]{\tiny{$a$}};
\node[below] at (0.5,0) [scale=1]{\tiny{$+1$}};
\node[below] at (1,0) [scale=1]{\tiny{$+2$}};

\node at (0.25,0.25) [scale=0.8]{\small{Z}};

\node[above] at (0.5,0.5) [scale=1.25]{(a)};

 \node[below, align=left, text width=5cm] at (1.5, -0.25) { Fig. 4.44. Case 2.1. (a) $e(a;b,b{+}1) \notin H$. (b) $e(a;b,b{+}1) \in H$};

} \end{scope}

\begin{scope}[xshift=2cm]{
\draw[gray,very thin, step=0.5cm, opacity=0.5] (0,0) grid (1,0.5);

\fill[blue!50!white, opacity=0.5](0,0) rectangle (1,0.5);

\draw[blue, line width=0.5mm] (0,0)--++(1,0);
\draw[blue, line width=0.5mm] (0,0)--++(0,0.5)--++(0.5,0);

\draw[blue, line width=0.5mm] (0.5,0)--++(0.5,0);

%black lines
{
\draw[black, line width=0.15mm] (0.7,0.45)--++(0,0.1);
\draw[black, line width=0.15mm] (0.75,0.45)--++(0,0.1);
\draw[black, line width=0.15mm] (0.8,0.45)--++(0,0.1);
}

\node[left] at (0,0) [scale=1]
{\tiny{$b$}};
\node[left] at (0,0.5) [scale=1]
{\tiny{$+1$}};

\node[below] at (0,0) [scale=1]{\tiny{$a$}};
\node[below] at (0.5,0) [scale=1]{\tiny{$+1$}};
\node[below] at (1,0) [scale=1]{\tiny{$+2$}};

\node at (0.25,0.25) [scale=0.8]{\small{Z}};

\node[above] at (0.5,0.5) [scale=1.25]{(b)};

} \end{scope}

\end{tikzpicture}
\end{adjustbox}
\end{wrapfigure}

\noindent \textit{CASE 2.1. $Z$ and $Z+(1,0)$ both belong to $P$.} We have that at least one of $e(a,a+1;b+1)$ and $e(a+1,a+2;b+1)$ belongs to $H$. For definiteness, assume that $e(a, a+1; b+1) \in H$. Now, either $e(a;b,b+1) \notin H$ or $e(a;b,b+1) \in H$. If $e(a;b,b+1) \notin H$ then $Z$ is switchable (Figure 4.44 (a)), so we only need to check the case where $e(a;b,b+1) \in H$. In that case (Figure 4.44 (b)), $Z$ is an end-box of $P$, and $F$ must be eastern. Note that if $e(a+1,a+2;b+1) \in H$, then $Z+(1,0) \in P$ is switchable, so we may assume that $e(a+1,a+2;b+1) \notin H$. But then, $e(a,a+1;b)$ and $e(a+1,a+2;b)$ must be the right colinear edges of the $A_1$-type that follows $e_F$, which contradicts the initial assumption of the Lemma.  End of Case 2.1

\null 

\begingroup 
\setlength{\intextsep}{0pt}
\setlength{\columnsep}{20pt}
\begin{wrapfigure}[]{l}{0cm}
\begin{adjustbox}{trim=0cm 0cm 0cm 0cm}
\begin{tikzpicture}[scale=1.5]
\begin{scope}[xshift=0cm]{
\draw[gray,very thin, step=0.5cm, opacity=0.5] (0,0) grid (1,0.5);

\fill[green!50!white, opacity=0.5](0,0) rectangle (0.5,0.5);
\fill[blue!50!white, opacity=0.5](0.5,0) rectangle (1,0.5);

\draw[blue, line width=0.5mm] (0,0)--++(1,0)--++(0,0.5);
\draw[blue, line width=0.5mm] (0,0.5)--++(0.5,0);

%black lines
{
\draw[black, line width=0.15mm] (0.7,0.45)--++(0,0.1);
\draw[black, line width=0.15mm] (0.75,0.45)--++(0,0.1);
\draw[black, line width=0.15mm] (0.8,0.45)--++(0,0.1);
}

\node[left] at (0,0) [scale=1]
{\tiny{$b$}};
\node[left] at (0,0.5) [scale=1]
{\tiny{$+1$}};

\node at (0.25,0.25) [scale=0.8] {\small{Z}};

\node[below] at (0,0) [scale=1]{\tiny{$a$}};
\node[below] at (0.5,0) [scale=1]{\tiny{$+1$}};
\node[below] at (1,0) [scale=1]{\tiny{$+2$}};

\node at (0.25,0.25) [scale=0.8]{\small{Z}};

 \node[below, align=center, text width=2.5cm] at (0.5, -0.25) { Fig. 4.45. Case 2.2.};

} \end{scope}

\end{tikzpicture}
\end{adjustbox}
\end{wrapfigure}

\noindent \textit{CASE 2.2.  Exactly one of $Z$ and $Z+(1,0)$ belongs to $N[P] \setminus P$}. For definiteness, assume that $Z \in N[P] \setminus P$. By Lemma 4.6 (d), $e(a,a+1;b+1) \in H$. Lemma 4.6 (b) implies that $Z+(1,0)$ is an end-box of $P$. By Lemma 4.6 (c), $e(a+2;b,b+1) \in H$, and $e(a+1,a+2;b+1) \notin H$. Then $F$ must be western.
But then $e(a,a+1;b)$ and $e(a+1,a+2;b)$ must be the left colinear edges of the $A_1$-type that follows $e_F$, contradicting the initial assumption of the lemma. See Figure 4.45. $\square$

\endgroup 

\null 

\noindent This completes the proof of Proposition 4.14. An immediate consequence of it is that the $A_1$-type following $e_F$ is the only $A_1$-type in $F$, so we can refer to it as \textit{the} $A_1$-type of $F$.

\subsection{Turns} 

\noindent In this section we show that every standard looping fat path $F$ must have a turn \index{turn}. We do this by showing that the boundary $B(F)$ of a looping fat path $F$ must have a (necessarily closed) turn and note that this would immediately imply that $F$ must have a turn. We will often use the definition of the shadow of standard southern looping fat path, Proposition 4.14, and the fact that the $A_1$-type of $F$ is unique in $F$, and write (DsFP) whenever we appeal to them.

\null

\endgroup

\noindent \textbf{Lemma 4.17.} Let $H$ be a Hamiltonian path or cycle of a polyomino $G$ and let $B(F)$ be the boundary of a \index{standard looping fat path}standard looping fat path $F$ of $G$. Then $B(F)$ has at least one turn $T_1$ such that both leaves of $T_1$ belong to $F$.

\null

\begingroup
\setlength{\intextsep}{0pt}
\setlength{\columnsep}{20pt}
\begin{wrapfigure}[]{l}{0cm}
\begin{adjustbox}{trim=0cm 0cm 0cm 0.25cm}

% northeastern turn
\begin{tikzpicture}[scale=1.5]

\begin{scope}[xshift=0cm, yshift=0cm]
{
\draw[gray,very thin, step=0.5cm, opacity=0.5] (0,0) grid (1.5,1);

\fill[blue!50!white, opacity=0.5] (0.0,0.0) rectangle (0.5,1.0);
\fill[blue!50!white, opacity=0.5] (1.0,0.0) rectangle (1.5,1.0);

\draw[blue, line width=0.5mm] (0.0,1.0)--++(0.5,0)--++(0,-1);
\draw[blue, line width=0.5mm] (1.0,1.0)--++(0,-1);

\draw[->, blue,line width=0.5mm] (1.0,1.0)--++(0.5,0);
%\draw[green!50!black, line width=0.5mm] (0.0,0.5)--++(0.5,0);

\node[below] at (0.0,0) [scale=1]{\tiny{$-2$}};
\node[below] at (0.5,0) [scale=1]{\tiny{$-1$}};
\node[below] at (1.0,0) [scale=1]{\tiny{$k'$}};
\node[below] at (1.5,0) [scale=1]{\tiny{$+1$}};

\node[left] at (0,1.0) [scale=1]{\tiny{$\ell'$}};
\node[left] at (0,0.5) [scale=1]{\tiny{$-1$}};
\node[left] at (0,0.0) [scale=1]{\tiny{$-2$}};

\node[above] at (1.25,1) [scale=1]{\small{$e_s$}};

\node at (0.25,0.75) [scale=0.8]{Y};
\node at (1.25,0.75) [scale=0.8]{X};

 \node[below, align=center, text width=4cm] at (0.75, -0.25) { Fig. 4.46. $\text{Boxes}(F)$ shaded blue.};

}
\end{scope}

\end{tikzpicture}
\end{adjustbox}
\end{wrapfigure}

\noindent \textit{Proof.} For definiteness, assume that $F=G\langle N[P(X,Y)] \rangle $ is a standard southern looping fat path following $e(k'-1,l';l'+1)$ southward, with $X=R(k',l'-1)$ and $Y=R(k'-2,l'-1)$. By Lemma 4.13, $B(F)$ is a cycle. Orient  $B(F)$ into a directed trail $\overrightarrow{K}$ so that the first edge of $\overrightarrow{K}$ 
is  $(v(k',l'), v(k'+1,l'))$. With this orientation we can give a direction  - $N,S,E$ or $W$ - to edges of $\overrightarrow{K}$, defined as the position of the head of an edge relative to its tail. Our choice of direction for the first edge and Lemma 1.3.3 (b) imply that $\text{Boxes}(\Phi(\overrightarrow{K}, \text{right})) \subset \text{Boxes}(F)$ so the boxes of $F$ are on the right side of the oriented edges of $K$. We call this fact (RSK) for reference. See Figure 4.46. We sweep the edges of $K$ in the direction of the orientation starting at $v(k',l')$. We observe that we must encounter at least one west edge $e_W$, since $Y$ is west of $X$. Let $e_W=e_0=e(k-1, k;l)$ be the first west edge encountered, let $e_1$ be the edge preceding $e_0$ in the sweep, let $e_j$ be the edge preceding the edge $e_{j-1}$ in the sweep and let $(v(k',l'), v(k'+1,l'))=e_s$.

\null 

\noindent In this proof we will use the fact that $e_W$ is the first west edge encountered (1stW) often. 

\null 

\begingroup
\setlength{\intextsep}{0pt}
\setlength{\columnsep}{20pt}
\begin{wrapfigure}[]{r}{0cm}
\begin{adjustbox}{trim=0cm 0.75cm 0cm 0cm}

% northeastern turn

\end{adjustbox}
\end{wrapfigure}

\noindent By (DsFP),  $e_0$ was immediately preceded by a south edge or a north edge.

\null 

\noindent \textit{CASE 1: $e_1$ is southern.} We shall find a northeastern turn. By (DsFP) and (1stW), the preceding edge $e_2$ must be eastern; By (DsFP) and (1stW), $e_3$  has to be southern. By (1stW) $e_4$ cannot be western. Then $e_4$ is southern or $e_4$ is eastern.

\null 

\noindent \textit{CASE 1.1: $e_4$ is southern.} Then, by (DsFP) and (RSK), we have that $\Phi(e_4, \text{right})=X$ or $\Phi(e_4, \text{right})=Y$. But the former contradicts Proposition 4.14, and the latter implies that $\text{deg}_{B(F)}(v(k,l+1)=3$, contradicting Lemma 4.13. See Figure 4.47. Thus, $e_4$ must be eastern. End of Case 1.1. 

\null

\endgroup 

\noindent \textit{CASE 1.2: $e_4$ is eastern.} By (DsFP), $e_5$ is not eastern. Then $e_5$ is northern or $e_5$ is southern. 

\null

\noindent \textit{CASE 1.2(a): $e_5$ is northern.} Then there is a northeastern turn $T_1$ on the edges $e_0, \ldots, e_5$ with both leaves contained in $F$. See Figure 4.48 (a). End of Case 1.2(a).

\setlength{\intextsep}{0pt}
\setlength{\columnsep}{20pt}
\begin{center}
\begin{adjustbox}{trim=0cm 0cm 0cm 0cm}

\end{adjustbox}
\end{center}

\noindent \textit{CASE 1.2 (b): $e_5$ is southern.} Let $Q(j)$ be the statement: ``$e_j$ is southern and $e_{j+1}$ is eastern''. Now, either $Q(j)$ is true for each $j \in \{1,3, \ldots, s-1\}$ (Case 1.2($b_1$)), or there is some $j_0 \in \{ 5,9, \ldots, s-1\}$  such that $Q(j)$ for each odd $j<j_0$, but $Q(j_0)$ is not true (Case 1.2($b_2$)).

\null 

\noindent \textit{CASE 1.2($b_1$).} Then we have a northeastern turn $T_1$ on the edges $e_0, \ldots, e_s, e(k'; l'-1, l')$ with both leaves contained in $F$. Figure 4.48 (b). End of Case 1.2($b_1$).

\null

\noindent \textit{CASE 1.2($b_2$).} By (DsFP), $e_{j_0}$ is not eastern. If $e_{j_0}$ is southern, then we run into the same contradiction as in Case 1.1; and if $e_{j_0}$ is northern then we have a northeastern turn $T_1$ on the edges $e_0, \ldots, e_{j_0}$ with both leaves contained in $F$. See Figure 4.48 (c). End of Case  1.2($b_2$). End of Case 1.2(b). End of Case 1.2.

\null

\begingroup
\setlength{\intextsep}{0pt}
\setlength{\columnsep}{20pt}
\begin{wrapfigure}[]{r}{0cm}
\begin{adjustbox}{trim=0cm 0cm 0cm 0cm}

% northeastern turn

\end{adjustbox}
\end{wrapfigure}

\noindent \textit{CASE 2: $e_1$ is northern.} We shall find a southeastern turn with both leaves contained in $F$. By (1stW) and (DsFP), $e_2$ is eastern. By (DsFP), $e_3$ is northern. Note that (RSK) and Lemma 4.6 (d) imply that $\Phi(e_1, \text{right}) \in P$ and $e_1 \in H$. By Lemma 4.13, $e(k,k+1;l-1) \notin H$ and $e(k,k+1;l) \notin H$. If $e(k+2;l-1,l) \in H$, then $\Phi(e_1, \text{right})$ is switchable and in $P$, contradicting FPC-3, so we may assume that $e(k+2;l-1,l) \notin H$. Then we must have that $S_{\downarrow}(k+1,l+1;k+2,l) \in H$,  $S_{\uparrow}(k+1,l-2;k+2,l-1) \in H$ and that $e(k+2;l-1,l) \in B(F)$.

By (1stW), $e_4$ is not western. If $e_4$ is northern then (DsFP) and (RSK) imply that $\Phi(e_3, \text{right}) =X$. But then $L=R(k-1,l)$ and then $\deg_H(v(k-1,l))=3$, contradicting that $H$ is Hamiltonian. Then $e_4$ must be eastern.  By (DsFP), $e_5$ is not eastern. Then $e_5$ is southern or northern. See Figure 4.49.

\endgroup

\begingroup
\setlength{\intextsep}{0pt}
\setlength{\columnsep}{20pt}
\begin{wrapfigure}[]{l}{0cm}
\begin{adjustbox}{trim=0cm 0cm 0cm 0.25cm}

% northeastern turn

\end{adjustbox}
\end{wrapfigure}

\noindent \textit{CASE 2.1: $e_5$ is southern.} By (RSK) and Lemma 4.6 (d), we have that $\Phi(e_4, \text{right}) \in P$, and that $e_4 \in H$.  By Lemma 4.13, $e(k-2;l-3,l-2) \notin H$ and $e(k-1;l-3,l-2) \notin H$. If $e(k-2,k-1;l-3) \in H$, then $\Phi(e_4, \text{right})$ is switchable and in $P$, contradicting FPC-3, so we may assume that $e(k-2,k-1;l-3) \notin H$. It follows that  $S_{\rightarrow}(k-3,l-3;k-2,l-4) \in H$ and that $S_{\uparrow}(k-1,l-4;k+2,l-1) \in H$. Then there is a southeastern turn on $e(k-2;l-4,l-3)$, $S_{\uparrow}(k-1,l-4;k+2,l-1)$, $e(k+1,k+2;l)$ with both leaves in $F$. See Figure 4.50 (b). End of Case 2.1

\endgroup
 
\null

\begingroup
\setlength{\intextsep}{0pt}
\setlength{\columnsep}{20pt}
\begin{wrapfigure}[]{r}{0cm}
\begin{adjustbox}{trim=0cm 0.5cm 0cm 1cm}

\end{adjustbox}
\end{wrapfigure}

\noindent \textit{CASE 2.2: $e_5$ is northern.} Let $Q(j)$ be the statement: ``$e_j$ is northern and $e_{j+1}$ is eastern''. Now, either $Q(j)$ is true for each $j \in \{1,3, \ldots, s-1\}$ (Case 2.2(a)), or there is some $j_0 \in \{ 5,7, \ldots, s-1\}$  such that $Q(j)$ for each odd $j<j_0$, but $Q(j_0)$ is not true (Case 2.2(b)).

\null

\noindent \textit{CASE 2.2(a).} Then (DsFP) and (RSK) imply that $\Phi(e_s, \text{right})=X$. This means that $e(k';l'-2,l'-1) \in H$, $e(k';l'-1,l') \in H$, and that $\Phi(e_{s-1}, \text{right}) \in F$. By Proposition 4.14, $e(k';l'-3,l'-2) \notin H$. Note that if $e(k'-1,k';l'-2) \in H$, then $P(X,Y)$ is the $H$-path $X, X+(0,-1), X+(0,-2), X+(-1,-2), X+(-2,-2), X+(-2,-1), Y$. This contradicts our finding that $\Phi(e_{s-1}, \text{right}) \in F$. Then it must be the case that $e(k'-1,k';l'-2) \notin H$. Then we must have $e(k', k'+1;l'-2) \in H$. It follows that $S_{\uparrow}(k'+1,l'-2;k+2,l-1) \in H$. See Figure 4.51. Then there is a southeastern turn $T_1$ on $e(k';l'-2,l'-1)$, $S_{\rightarrow}(k'+1,l'-2;k+2,l-1)$, $e(k+1,k+2;l)$ with both leaves contained in $F$ (see Figure 4.72). End of Case 2.2(a). 

\endgroup

\null 

\begingroup
\setlength{\intextsep}{0pt}
\setlength{\columnsep}{25pt}
\begin{wrapfigure}[]{l}{0cm}
\begin{adjustbox}{trim=0cm 1.5cm 0cm 0.5cm}

% northeastern turn

\end{adjustbox}
\end{wrapfigure}

\noindent \textit{CASE 2.2(b).} By (DsFP), $e_{j_0}$ is not eastern. Suppose that $e_{j_0}$ is northern (in orange in Figure 4.25). The assumption that $Q(j_0)$ is false implies that $e_{j_0+1}$ is not eastern and (1stW) implies that $e_{j_0+1}$ is not western. It must be the case that $e_{j_0+1}$ is also northern. By (DsFP) and (RSK), we have that $\Phi(e_{j_0}, \text{right})=X$. But then $j_0-1=s$, contradicting the assumption that $j_0 \in \{5, \ldots, s-1\}$ (in orange in Figure 4.52). Thus $e_{j_0}$ cannot be northern. It follows that $e_{j_0}$ is southern. Let $e_{j_0}=e(k'';l'',l''+1)$.

Using the same arguments as in Case 2.1, we find that $e(k'';l''-1,l'') \notin H$, $e(k''+1;l''-1,l'') \notin H$, $e(k'',k''+1;l''-1) \notin H$, and that $\Phi({j_0-1}, \text{right}) \in P$. It follows that $S_{\rightarrow}(k''-1,l''-1;k'',l''-2) \in H$ and that $S_{\uparrow}(k''+1,l''-2;k+2,l-1) \in H$. Then there is a southeastern turn $T_1$ on $e(k'';l''-2,l''-1)$, $S_{\uparrow}(k''+1,l''-2;k+2,l-1)$, $e(k+1,k+2;l)$ with both leaves contained in $F$ (in blue in Figure 4.52). End of Case 2.2(b). End of Case 2.2. End of Case 2. $\square$ 

\null

\endgroup 

\noindent Recall once more the setup from the sketch of the proof of Lemma 3.13 at the start of this chapter: there is an $H$-path $P(X,Y)$ following a leaf $L$ contained in a looping fat path $F$. We then promised to show that $F$ must have a turn $T$, and that we can find a cascade (called a \textit{weakening}) that collects one of its leaves. As we shall see in Sections 4.3-4.6 we cannot control which leaf of the turn is collected. In some cases this leads to unwanted behavior. Recall that one of the steps toward our goal is to make $P(X,Y)$ gain a switchable box. In particular, if one of the turn's leaves happens to be an end-box of $P(X,Y)$, then the weakening does not give a useful outcome: instead of allowing us to extend the weakening and collect $L$ (the original goal), it can cement $L$ into a stair\index{stairs} subgraph. To avoid such situations, we restrict attention to well-behaved turns, called \textit{admissible turns}, defined below.

\null 

\noindent\textbf{Definition.} Let $G$ be an  $m\times n$ grid graph, let $H$ be a Hamiltonian cycle of $G$, and let $F$ be a standard looping fat path in $G$. We say that a turn $T$ of $B(F)$ is \index{admissible turn|textbf}\textit{admissible} if:

(i) \  no leaf of $T$ is an end-box of $F$, and 

(ii) both leaves of $T$ belong to $F$.  

\null 

\endgroup

\noindent \textbf{Lemma 4.18.} Let $H$ be a Hamiltonian path or cycle of a polyomino $G$ and let $B(F)$ be the boundary of a \index{standard looping fat path}standard looping fat path of $G$. Then $B(F)$ has an admissible turn. 

\null 

\noindent \textit{Proof.} Let $F$, $X$, $Y$, $\overrightarrow{K}$ and $e_W, e_1, \ldots e_s$ be as in Lemma 4.17, including the assumption that $F$ is southern and (RSK).

\null 

\noindent \textit{CASE 1: $e_1$ is southern.} By Case 1 in Lemma 4.17, there is a northeastern turn $T_1$. 
%with eastern leaf $L_E(T_1)= \Phi(e_1, \text{right})$. 
We continue sweeping $\overrightarrow{K}$, beginning from $e_W$, until we find the first northern edge $e_N$ in the subtrail $\overrightarrow{K}(e_W,e_N)$ of $\overrightarrow{K}$, where $e_N=(\widehat{k}; \widehat{l}, \widehat{l}+1)= \widehat{e_0}$. 
We write (1stN) to refer to the fact that $e_N$ is the first northern edge encountered after $e_W$, whenever we appeal to it. Let $\widehat{e_1}$ be the edge preceding $\widehat{e_0}$ in the sweep, let $\widehat{e_j}$ be the edge preceding the edge $\widehat{e_{j-1}}$ in the sweep and let $\widehat{e_{t+1}}=e_W$. Then $\widehat{e_1}$ is western or $\widehat{e_1}$ is eastern. 

%[pagemarker] 

\begingroup 
\setlength{\intextsep}{0pt}
\setlength{\columnsep}{20pt}
\begin{wrapfigure}[]{r}{0cm}
\begin{adjustbox}{trim=0cm 0cm 0cm 0cm}
\begin{tikzpicture}[scale=1.5]
\begin{scope}[xshift=0cm]{
\draw[gray,very thin, step=0.5cm, opacity=0.5] (0,0) grid (3,2.5);

%\fill[green!50!white, opacity=0.5] (4,0)--++(-1,1)--++(0,1)--++(1,1);
\fill[blue!50!white, opacity=0.5]  (0,0)--++(1.5,1.5)--++(1.5,-1.5);
\fill[orange!50!white, opacity=0.5](2.5,0)--++(-1,1)--++(1.5,1.5)--++(0,-2.5);

%\draw[green!50!black, line width=0.5mm, opacity=0.5] (4,0)--++(-1,1)--++(0,1)--++(1,1);

\draw[black, line width=0.25mm] (0,2)--++(0.5,0)--++(0,-1);
\draw[black, line width=0.25mm] (1.5,2)--++(-0.5,0)--++(0,-1);
\draw[orange, line width=0.5mm, opacity=1] (2.5,0)--++(-1,1)--++(1.5,1.5);
\draw[blue, line width=0.5mm, opacity=1] (0,0)--++(1.5,1.5)--++(1.5,-1.5);
\draw[red, dashed, line width=0.5mm, opacity=1] (0,0)--++(2.5,2.5);

%nodes
{
\node[left] at (0,2) [scale=1]
{\tiny{0}};
\node[left] at (0,1.5) [scale=1]
{\tiny{-1}};
\node[left] at (0,1) [scale=1]
{\tiny{-2}};

\node[above] at (1,2.5) [scale=1]
{\tiny{0}};
\node[above] at (1.5, 2.5) [scale=1]
{\tiny{1}};
\node[above] at (2,2.5) [scale=1]
{\tiny{2}};

\node at (0.25,1.75) [scale=0.8] {\small{$Y$}};
\node at (1.25,1.75) [scale=0.8] {\small{$X$}};

\draw[fill=blue] (1.5,1) circle [radius=0.05];

\node at (1.5,0.25) [scale=0.8] {$U_2$};

\node at (2.5,1) [scale=0.8] {{$U_{1, \text{end}}$}};

\node[left] at (1.5,0.9) [scale=0.8] {\small{(a,b)}};
}

\node[below, align=left, text width=6cm] at (1.5, -0.1) { Fig. 4.53. Case 1. The line \\ $y-x+2=0$ in  red; $U_{1, \text{end}}$\\  shaded orange, $U_2$ shaded blue.};

} \end{scope}

\end{tikzpicture}
\end{adjustbox}
\end{wrapfigure}

\noindent Before we consider each case, we will check that the subtrail $\overrightarrow{K}(\widehat{e_0}, \widehat{e_t})$ of $\overrightarrow{K}$ does not contain the right or left colinear edges of the $A_1$-type of $F$. To this end, we will translate $H$ by $(-k',-l')$ to simplify calculations. (DsFP) and (1stW) imply that for every eastern edge in the subtrail $\overrightarrow{K}(e_s, e_1)$ of $\overrightarrow{K}$ there is at most one northern or southern edge. Denote by ${v_\text{end}}$ the head of the edge $e_W$.
The assumption that $e_1$ is southern and the fact that a shortest turn has length two imply that ${v_\text{end}}$ is contained in the region $U_{1, \text{end}}$, determined by $x \geq 1 $ and $|y+2| \leq x-1$ (Eq.1). Let ${v_\text{end}}=v(a,b)$. It follows that $\overrightarrow{K}(\widehat{e_t}, \widehat{e_0})$ is contained in the region $U_2$ bounded by
$y \leq b+1$ and $|x-a| \leq b+1-y$. (Eq.2) %\footnote{Note that this accounts for the possibility that $\overrightarrow{K}(\widehat{e_0}, \widehat{e_t})$ contains a pair of consecutive colinear edges of the $A_1$-type.}

We will check that $U_2$ and the colinear edges of the $A_1$-type of $F$ lie on two different sides of the line $y=x-2$. See Figure 4.53. By (Eq.1) we have that $b \leq a-3$ and by (Eq.2) we have that $y \leq x-a+b+1$. Let $(x,y) \in U_2$.  Then $y-x+2 \leq -a+b+3 \leq 0$, so $U_2$ lies below the line $y-x+2$. Plugging in the values of the coordinates of the vertices of the $A_1$-type, we see that they lie above the line $y=x-2$. This shows that $\overrightarrow{K}(\widehat{e_t}, \widehat{e_0})$ does not contain colinear edges. We will write (NCE) whenever we appeal to this fact. Note that (NCE) implies (i).

\endgroup

\null

\noindent \textit{CASE 1.1: $\widehat{e_1}$ is western.} Note that $\widehat{e_1} \neq e_W$, otherwise we get a cycle on $e_1, e_W, e_1, e_2$. By (DsFP), $\widehat{e_2}$ is not western and by (1stN), $\widehat{e_2}$ is not northern, so $\widehat{e_2}$ must be southern. By (NCE) $\widehat{e_3}$ cannot be southern. Then $\widehat{e_3}$ must be western. If $\widehat{e_3}=e_W$, then we have a southeastern turn $T_2$ on $\widehat{e_0}, \ldots, \widehat{e_3}, e_1, e_2$, satisfying (i) and, by (RSK), (ii) (in orange in Figure 4.54), so we may assume that $\widehat{e_3} \neq e_W$. By (DsFP), $\widehat{e_4}$ is not western and by (1stN), $\widehat{e_4}$ is not northern. Then $\widehat{e_4}$ is southern.

\begingroup
\setlength{\intextsep}{10pt}
\setlength{\columnsep}{20pt}
\begin{wrapfigure}[]{r}{0cm}
\begin{adjustbox}{trim=0cm 0cm 0cm 0cm}

% northeastern turn

\end{adjustbox}
\end{wrapfigure}

Let $Q(j)$ be the statement: ``$\widehat{e_j}$ is western and $\widehat{e_{j+1}}$ is southern''. Now, either $Q(j)$ is true for each $j \in \{1,3, \ldots, t+1\}$, or there is some $j_0 \in \{ 5,7, \ldots, t+1\}$  such that 
$Q(j)$ for each odd $j<j_0$, but $Q(j_0)$ is not true. If the former, then we have a southeastern turn $T_2$ on $\widehat{e_0}, \widehat{e_1}, \ldots, \widehat{e_t}, e_W, e_1,e_2$ satisfying (i) and (ii) (in blue in Figure 4.54), so assume the latter. By (NCE), $\widehat{e_{j_0}}$ is not southern. If $\widehat{e_{j_0}}$ is eastern then we have a southeastern turn on $\widehat{e_0}, \widehat{e_1}, \ldots, \widehat{e_{j_0}}$  satisfying (i) and (ii) (blue in Figure 4.75). Suppose then that $\widehat{e_{j_0}}$ is western (in green in Figure 4.54). This is impossible: by (1stN). $\widehat{e_{j_0+1}}$ is not northern; since $Q(j_0)$ is false, $\widehat{e_{j_0+1}}$ is not southern; and by (DsFP), $\widehat{e_{j_0+1}}$ is not western. End of Case 1.1

\endgroup 

\null

\begingroup
\setlength{\intextsep}{0pt}
\setlength{\columnsep}{20pt}
\begin{wrapfigure}[]{l}{0cm}
\begin{adjustbox}{trim=0cm 0cm 0cm 0cm}

% northeastern turn

\end{adjustbox}
\end{wrapfigure}

\noindent \textit{CASE 1.2: $\widehat{e_1}$ is eastern.} By (DsFP), $\widehat{e_2}$ is not eastern and by (1stN), $\widehat{e_2}$ is not northern, so $\widehat{e_2}$ must be southern. By (NCE) $\widehat{e_3}$  is not southern and by (DsFP), $\widehat{e_3}$ is not western. Then $\widehat{e_3}$ must be eastern. (DsFP) and (1stN) imply that $\widehat{e_4}$ must be southern.

Let $Q(j)$ be the statement: ``$\widehat{e_j}$ is eastern and $\widehat{e_{j+1}}$ is southern''. Now, either $Q(j)$ is true for each $j \in \{1, 3, \ldots, t-1\}$, or there is some $j_0 \in \{ 5,7, \ldots, t-1\}$  such that $Q(j)$ is true for each odd $j<j_0$, but $Q(j_0)$ is not true. 
 
\null 

\noindent \textit{CASE 1.2(a): $Q(j)$ is true for each $j \in \{1, 3, \ldots, t-1\}$.}  Then we have a southwestern turn on $\widehat{e_0}, \widehat{e_1}, \ldots, \widehat{e_{t-1}}, \widehat{e_t}, e_W$. Recall that $e_W=e(k-1, k;l)$. Observe that $R(k-2; l-1) \in P$, so $e(k-2;l-1,l) \notin H$.

\endgroup

\noindent Similarly, $R(\widehat{k}-1,\widehat{l}-1) \in P$ and $e(\widehat{k}-1, \widehat{k};\widehat{l}-1) \notin H$. It follows that $R(k-3, l-1) \in F$, $R(\widehat{k}-1,\widehat{l}) \in F$ and that there is a southeastern

\noindent turn $T_2$ on $e(k-3,k-2;l)$, $S_{\rightarrow}(k-3,l-1;\widehat{k}-1, \widehat{l}-2)$, $e(\widehat{k}; \widehat{l}-2;\widehat{l}-1)$ satisfying (i) and (ii) (in blue in Figure 4.55). End of Case 1.2(a).

\null 

\noindent \textit{CASE 1.2(b): There is some $j_0 \in \{ 5,7, \ldots, t-1\}$  such that $Q(j)$ is true for each odd $j<j_0$, but $Q(j_0)$ is not true.} If $\widehat{e_{j_0}}$ is western then we have a southeastern turn on $\widehat{e_0}, \widehat{e_1}, \ldots, \widehat{e_{j_0}}$. As in Case 1.2 (a), this turn satisfies (i) and (ii) (in blue in Figure 4.55, with $e_{j_0}$ dotted green). 

By (DsFP), $\widehat{e_{j_0}}$ is not southern. Suppose then $\widehat{e_{j_0}}$ is eastern. This is impossible: by (1stN), $\widehat{e_{j_0+1}}$ is not northern; since $Q(j_0)$ is false, $\widehat{e_{j_0+1}}$ is not southern; and by (DsFP), $\widehat{e_{j_0+1}}$ is not eastern (dotted orange in Figure 4.55). End of Case 1.2(b). End of Case 1.2.

\null 

\noindent \textit{CASE 2: $e_1$ is northern.}   By Case 2 in Lemma 4.17, $\overrightarrow{K}$ has a southeastern turn $T_1$. We continue sweeping $K$, beginning from $e_W$, until we find the first southern edge $e_S$ in the subtrail $\overrightarrow{K}(e_W,e_S)$ of $\overrightarrow{K}$, where $e_S=(\widehat{k}; \widehat{l}, \widehat{l}+1)= \widehat{e_0}$. 
We write (1stS) to refer to the fact that $e_S$ is the first southern edge encountered after $e_W$, whenever we appeal to it. Let $\widehat{e_1}$ be the edge preceding $\widehat{e_0}$ in the sweep, let $\widehat{e_j}$ be the edge preceding the edge $\widehat{e_{j-1}}$ in the sweep and let $\widehat{e_{t+1}}=e_W$. Then $\widehat{e_1}$ is western or $\widehat{e_1}$ is eastern.

\null

\noindent \textit{CASE 2.1: $\widehat{e_1}$ is western.} Note that the assumption that $e_1$ is northern implies that $\widehat{e_1} \neq e_W$, otherwise there is a cycle $e_2, e_1, e_W, \widehat{e_0}$. By (DsFP), $\widehat{e_2}$ is not western and by 1stS, $\widehat{e_2}$ is not southern, so $\widehat{e_2}$ must be northern. By (DsFP) %and RSK,
$\widehat{e_3}$ is not northern or eastern. Then $\widehat{e_3}$ must be western. By 1stS $\widehat{e_4}$ is not southern, and by (DsFP), $\widehat{e_4}$ is not western. Then $\widehat{e_4}$ must be northern. Let $Q(j)$ be the statement: ``$\widehat{e_j}$ is western and $\widehat{e_{j+1}}$ is northern. Then either $Q(j)$ is true for each $j \in \{1,3, \ldots, t-1\}$ or there is some $j_0 \in \{ 5,7, \ldots, t-1\}$  such that $Q(j)$ for each odd $j<j_0$, but $Q(j_0)$ is not true.

\begingroup
\setlength{\intextsep}{0pt}
\setlength{\columnsep}{20pt}
\begin{wrapfigure}[]{l}{0cm}
\begin{adjustbox}{trim=0cm 0cm 0cm 0cm}

\end{adjustbox}
\end{wrapfigure}

\null

\noindent \textit{CASE 2.1(a): $Q(j)$ is true for each $j \in \{1,3, \ldots, t-1\}$.} Then there is a northeastern turn on $\widehat{e_0}, \widehat{e_1}, \ldots, \widehat{e_{t-1}},\widehat{e_t}, e_W, e_1, e_2$. As in Case 1.2 (a), we have that $R(\widehat{k},\widehat{l}+1) \in P$, $R(k,l-1) \in P$, $e(\widehat{k}, \widehat{k}+1; \widehat{l}+2) \notin H$ and $e(k+1;l-1,l) \notin H$. Then there is a northeastern turn $T_2$ on $e(\widehat{k};  \widehat{l}+2,\widehat{l}+3)$, $S_{\downarrow}(\widehat{k}+1, \widehat{l}+3; k+2,l)$, $e(k+1,k+2;l-1)$ satisfying (i) and (ii) (in blue in Figure 4.56). End of Case 2.1(a).

\null 

\noindent \textit{CASE 2.1(b): There is some $j_0 \in \{ 5,7, \ldots, t-1\}$  such that $Q(j)$ for each odd $j<j_0$, but $Q(j_0)$ is not true.} If $\widehat{e_{j_0}}$ is eastern then we have a northeastern turn on $\widehat{e_0}, \widehat{e_1}, \ldots, \widehat{e_{j_0}}$. Then, as in Case 2.1(a), there is a northeastern turn $T_2$ satisfying (i) and (ii). 

By (DsFP), $\widehat{e_{j_0}}$ is not southern. Suppose then that $\widehat{e_{j_0}}$ is western. This is impossible: by 1stS. $\widehat{e_{j_0+1}}$ is not southern; since $Q(j_0)$ is false, $\widehat{e_{j_0+1}}$ is not northern; and by (DsFP), $\widehat{e_{j_0+1}}$ is not western (in orange in Figure 4.56). End of Case 2.1. End of Case 2.1(b). End of Case 2.1.

\endgroup

\null

\noindent \textit{CASE 2.2: $\widehat{e_1}$ is eastern.} By 1stS and (DsFP), $\widehat{e_2}$ is northern.  By (DsFP), $\widehat{e_3}$ is not western. An argument analogous to (NCE-1) in Case 1 can be used to show that $T_2$ and the $A_1$-type lie on two different sides of the line $y=2-x$. In this case, we have that the region $U_{1,\text{end}}$ containing $v_{\text{end}}$ is determined by $x \geq 1$ and $|y-2| \leq x-1$, and the region $U_2$ containing $T_2$, as defined in the next paragraph, is determined by $y\geq b-1$ and $|x-a| \leq y-b+1$. We will refer to this argument as (NCE-2). Note that by (NCE-2), $\widehat{e_3}$ is not northern. Then $\widehat{e_3}$ must be eastern. By (DsFP) and (1stS) $\widehat{e_4}$ is not southern or eastern. Then $\widehat{e_4}$ must be northern.

\begingroup
\setlength{\intextsep}{0pt}
\setlength{\columnsep}{20pt}
\begin{wrapfigure}[]{r}{0cm}
\begin{adjustbox}{trim=0cm 0cm 0cm 0cm}

% northeastern turn


\end{adjustbox}
\end{wrapfigure}

%pagemarker
\noindent Let $Q(j)$ be the statement: ``$\widehat{e_j}$ is eastern and $\widehat{e_{j+1}}$ is northern''. Now, either $Q(j)$ is true for each $j \in \{1,3, \ldots, t-1\}$, or there is some $j_0 \in \{ 5,7, \ldots, t-1\}$  such that $Q(j)$ for each odd $j<j_0$, but $Q(j_0)$ is not true. If the former then we have a northwestern turn $T_2$ on $\widehat{e_0}, \widehat{e_1}, \ldots, \widehat{e_{t-1}},\widehat{e_t}, e_W$ satisfying (i) and (ii), (in blue in Figure 4.57 with $\widehat{e_{t-1}},\widehat{e_t}, e_W$ dotted orange) so assume the latter. If $\widehat{e_{j_0}}$ is western then again we have a northwestern turn $T_2$ on $\widehat{e_0}, \widehat{e_1}, \ldots, \widehat{e_{j_0}}$ satisfying (i) and (ii) (in blue in Figure 4.57).  By (NCE-2), $\widehat{e_{j_0}}$ is not southern. Then, suppose that $\widehat{e_{j_0}}$ is western (green in Figure 4.57). This is impossible: by 1stS. $\widehat{e_{j_0+1}}$ is not southern; since $Q(j_0)$ is false, $\widehat{e_{j_0+1}}$ is not northern; and by (DsFP), $\widehat{e_{j_0+1}}$ is not eastern. End of Case 2.2. $\square$

\endgroup 

\null 

\noindent \textbf{Corollary 4.19.} Let $G$ be a polyomino, let $H$ be a Hamiltonian cycle or path of $G$, let $F=G\langle N[P(X,Y)] \rangle $ be a standard looping fat path in $G$, let $T$ be an admissible turn of $F$, let $L$ be a leaf of $T$, and let $L'$ be the $H$-neighbour of $L$ in $F$. Then:

(a) $d(T) \geq 3$, and 

(b)  $L \in N[P]\setminus P$ and $L' \in P$. 

\null 

\noindent  \textit{Proof of (a).} We prove the contrapositive. For definiteness, assume that $T$ is northeastern with northern leaf $L_N=R(k,l-1)$. Suppose that $d(T) < 3$. Then $d(T)=2$ and the eastern leaf of $T$ must be $L_E=R(k+1,l-2)$. Note that  $L_N+(0,-1) \in F$, otherwise $L_E, \ldots, L_N, L_N+(0,-1)$ is an $H$-cycle of boxes in $G$, contradicting Proposition 1.2.1.

Now, by Proposition 4.14, $e(k;l-2,l-1)$ and $e(k,k+1;l-2)$ cannot both belong to $H$. Then, either exactly one of $e(k;l-2,l-1)$ and $e(k,k+1;l-2)$ belongs to $H$, or neither does.

\begingroup 
\setlength{\intextsep}{0pt}
\setlength{\columnsep}{20pt}
\begin{wrapfigure}[]{l}{0cm}
\begin{adjustbox}{trim=0cm 0cm 0cm 0cm}

\end{adjustbox}
\end{wrapfigure}

\noindent \textit{CASE 1: exactly one of $e(k;l-2,l-1)$ and $e(k,k+1;l-2)$ belongs to $H$.} By symmetry, we may assume WLOG that  $e(k;l-2,l-1) \in H$ and $e(k,k+1;l-2) \notin H$. Note that the assumption that  $e(k,k+1;l-2) \notin H$ implies that $F$ is northern. It follows that $L_N$ is an end-box of $P(X,Y)$. Therefore $T$ is not admissible. See Figure 4.58 (a). End of Case 1.

\null 

\noindent \textit{CASE 2: neither $e(k;l-2,l-1)$ nor $e(k,k+1;l-2)$ belongs to $H$.} Then $e(k-1,k;l-1)$, $e(k+1;l-3,l-2)$ and \ $S_{\rightarrow}(k-1,l-2;k,l-3)$ belong to $H$. By Lemma 4.6(d), $L_N+(0,-1)$ must belong to $P(X,Y)$. It follows that at least one of $L_N$, $L_N+(0,-2)$, $L_E$ and $L_E+(-2,0)$ belongs to $P(X,Y)$ and is switchable, contradicting the assumption that $F$ is a looping fat path. See Figure 4.58 (b). End of Case 2. End of proof for (a).

\endgroup

\null 

\begingroup
\setlength{\intextsep}{0pt}
\setlength{\columnsep}{20pt}
\begin{wrapfigure}[]{r}{0cm}
\begin{adjustbox}{trim=0cm 0cm 0 0cm}

\end{adjustbox}
\end{wrapfigure}

\noindent \textit{Proof of (b).} Let $T$ be an admissible turn. For definiteness, assume that $T$ is northeastern with northern leaf $L_N=R(a,b)$. By Corollary 4.19 (a), $d(T) \geq 3$. Then we have that $e(a,a+1;b-1) \in H$, and that $e(a+1;b-1,b) \notin H$. By (RSK), $L_N+(0,-1)$ and $L_N+(1,-1)$ belong to $F$. This means that $L_N+(0,-1)=L_N'$ is the $H$-neighbour of $L_N$ in $F$. Now, either $e(a,a+1;b+1) \in H$ or $e(a,a+1;b+1) \notin H$. If $e(a,a+1;b+1) \in H$ (Figure 4.59 (a)), then, since $L_N$ is not an end-box of $P$, $L_N \in N[P]\setminus P$. Then, by Lemma 4.6 (a), $L_N' \in P$. And if $e(a,a+1;b+1) \notin H$ (Figure 4.59 (b)), then $L_N$ is switchable, so $L_N \in N[P] \setminus P$. Since $L_N' \in F$, by Lemma 4.6 (a), $L_N' \in P$. Either way we have that $L_N' \in P$ and  and $L_N \in N[P]\setminus P$. End of proof for (b). $\square$

\null

%{\noindent \textbf{Definition.} Let $H$ be a Hamiltonian path or cycle of a polyomino $G$, and let $F=G\langle N[P(X,Y)] \rangle$ be a looping fat path of $G$ following an edge $e_F$. Assume that $e_F$ is contained in the $j^{\textrm{th}}$ $A_0$ of a $j$-stack of $A_0$'s, whenever $j>0$. Let $T$ be a turn of $F$ such that:

%(i) \hspace{0.06cm}  neither leaf of $T$ is an end-box of $P(X,Y)$,

%(ii) \hspace{0.00cm} both leaves of $T$ are contained in $F$, and

%(iii) $\text{Sector}(T)$ and the $j$-stack of $A_0$'s are disjoint.

%Then we call $T$ an \textit{admissible} turn.

\endgroup

\noindent In this section we have shown that every standard looping fat path contains at least one admissible turn.

\subsection{Turn weakenings.}

\textbf{Definitions.} Let $H$ be a Hamiltonian path or cycle of an $m \times n$ grid graph $G$. Let $T$ be a turn\index{turn} of $H$ on $\{e(k;l-1,l), S_{\downarrow}(k+1,l;k',l'+1), e(k'-1,k';l')\}$ and let $L_N$ and $L_E$ be the northern and eastern leaves of $T$, respectively. Define the \index{sector of T@sector of $T$|textbf}\textit{sector of T} to be the induced subgraph of $G$ bounded by $e(k,k+1;l)$, $S_{\downarrow}(k+1,l;k',l'), e(k'-1,k';l')$, and the segments $[(k',l'), (m-1,l')]$, $[(m-1,l'), (m-1,n-1)]$, $[(m-1,n-1), (k,n-1)]$, $[(k,n-1), (k,l)]$, and denote it by Sector($T$). See Figure 4.60.  Analogous definitions apply to sectors of southeastern, southwestern and northwestern turns. Define the \index{end-vertex weakening terminal of a turn|textbf}\textit{end-vertex weakening terminal of $T$} to be the set of vertices $V(S_{\rightarrow}(k+1,l-1;k'-1,l'+1)) \cup \{ v(k,l-1), v(k'-1,l') \}$ and denote it by ew-Terminal($T$). We will call the edges $e(k, k+1; l-1)$ and $e(k'-1; l', l'+1)$ of $G \setminus H$ the \textit{northern} and \index{terminal edge of a turn|textbf}\textit{eastern terminal edge of $T$}, respectively.

Define a \index{weakening|textbf}\textit{weakening} of $T$ to be a cascade $\mu_1, \ldots, \mu_s$ such that $\mu_s$ is the first move after which one of the following holds:

1. One of the terminal edges of $T$ is in the resulting Hamiltonian path of $G$, in which 

\hspace{0.4cm} case we say that $T$ has an \index{edge weakening|textbf}\textit{edge} weakening,  or

2. An end-vertex of the resulting Hamiltonian path of $G$ is incident on 

\hspace{0.4cm} $\text{ew-Terminal}(T)$, in which case we say that $T$ has an \index{end-vertex weakening|textbf}\textit{end-vertex} weakening.

\null 

\noindent We call a weakening of $T$ consisting of three or less moves a \index{short weakening|textbf}\textit{short weakening of $T$}. We remark that if $H$ is a Hamiltonian path of $G$, then cookies are not defined and we need not be concerned with preserving their count after applying $\mu_1, \ldots, \mu_s$. We call the subgraph $S_{\downarrow}(k+1,l;k',l'+1)$ \textit{the stairs-part of T} and denote it by \index{stairs}$\text{stairs}(T)$. We say that $T$ has a \index{lengthening|textbf}\textit{lengthening} $T'$ if $T'$ is a turn of $H$ such that: 
 
 a) $d(T') \geq d(T)$ and
 
 b) $\text{stairs}(T') \supseteq \text{stairs}(T)+(1,1)$. 

\null

\noindent Define a set of lengthenings $\mathcal{T}(T_0)$ as follows:

1. The turn $T_0 \in \mathcal{T}(T_0)$.

2. The turn $T_j \in \mathcal{T}(T_0)$ if and only if $T_j$ is a lengthening of the turn $T_{j-1}$.

\noindent Note that by definition, $\mathcal{T}(T_0)$ is maximal.

\null

\noindent Define the \index{flank of a turn|textbf}\textit{flank of $T$} to be the set of vertices $V(\text{stairs}(T)+(1,1)) \cup V(L_N+(0,1)) \cup V(L_E+(1,0))$. We call $V(L_N+(0,1)) \cup V(L_E+(1,0))$ the \index{leaf flank of a turn|textbf}\textit{leaf flank of $T$} and denote it $\ell$-Flank$(T)$, and we call $V(\text{stairs}(T)+(1,1))$ the \index{stairs flank of a turn|textbf}\textit{stairs flank of $T$} and denote it $s$-Flank$(T)$. See Figure 4.60.

\setlength{\intextsep}{0pt}
\setlength{\columnsep}{20pt}
\begin{center}
\begin{adjustbox}{trim=0cm 0cm 0cm 0.5cm} 

\end{adjustbox}
\end{center}

\noindent \textbf{Observation 4.20.} Let $F=G\langle N[P(X,Y)] \rangle$ be a standard looping fat path, let $T$ be an admissible turn of $F$, and let $T'$ be a lengthening of $T$. Then:

(a) Every vertex of ew-Terminal($T$) is incident on a box of $P(X,Y)$. 

(b) ew-Terminal($T'$) $\subset$ Flank($T$).

(c) $\text{Flank}(T) \subset \text{Sector}(T)$.

(d) $\mathcal{T}(T) \subset \text{Sector}(T)$.

%pagemarker 

\noindent \textbf{Remark.} Recall from the definition of looping fat paths at the beginning of this chapter that our goal, roughly, is to find a cascade after which the edge $e'$ is in the resulting Hamiltonian path (see Figure 4.3). What happens if there is an end-vertex weakening? By Observation 4.20(a), every vertex of ew-Terminal($T$) is incident on a box of $P(X,Y)$. So, after applying the end-vertex weakening, FPC-1 fails in the resulting Hamiltonian path. We show in Chapter 5 (Claims 5.18-5.21) that when FPC-1 fails, we can find a cascade after which $e'$ is in the resulting Hamiltonian path. This explains why attaining an end-vertex weakening is sufficient for our purposes. This type of weakening is an effective way to deal with the Section 4.6 case, that considers Hamiltonian paths with $\text{Sector}(T)$ containing at most one end-vertex $u$ with $u \notin R_0$.

\null

\noindent \textbf{Lemma 4.21.} Let $G$ be an $m \times n$ grid graph, and let $H$ be a Hamiltonian cycle of $G$. Let $F$ be a looping fat path of $G$, anchored at some \index{outermost small cookie}outermost small cookie $C$. Then $F$ has an admissible turn $T$ such that Sector$(T)$ and the $j$-stack of $A_0$'s following $C$ are disjoint.

\begingroup
\setlength{\intextsep}{0pt}
\setlength{\columnsep}{20pt}
\begin{wrapfigure}[]{l}{0cm}
\begin{adjustbox}{trim=0cm 0cm 0cm 0cm} 

\end{adjustbox}
\end{wrapfigure}

\noindent  \textit{Proof.} For definiteness, assume that $C$ is a small southern cookie followed by a $j$-stack of $A_0$'s, which is then followed by the southern looping fat path $F=G\langle N[P(X,Y)] \rangle$, and let $X=R(k',l'-1)$. Let $\overrightarrow{K}$ and $e_W, e_1, \ldots e_s$ be as in Lemma 4.17. By Lemma 4.17, $F$ has a turn $T_1$. Either $T_1$ is a northeastern turn with $X$ as its northern leaf, or it is not. 

\null 

\noindent  \textit{CASE 1: $T_1$ is a northeastern turn with $X$ as its northern leaf.} Then $e_1$ is southern. By the proof of Lemma 4.18, $T_2$ is either southeastern or southwestern. In either case, we note that $\text{Sector}(T_2)$ is south of the stack of $A_0$'s. By Lemma 4.18, $T_2$ is admissible. See Figure 4.61. End of Case 1. 

\endgroup 

\null
\null 

\noindent  \textit{CASE 2: $T_1$ is not a northeastern turn with $X$ as its northern leaf.} By the proof of Lemma 4.17, $T_1$ is a northeastern turn or $T_1$ is a southeastern turn.

\begingroup
\setlength{\intextsep}{0pt}
\setlength{\columnsep}{20pt}
\begin{wrapfigure}[]{r}{0cm}
\begin{adjustbox}{trim=0cm 0.5cm 0cm 1cm} 

\end{adjustbox}
\end{wrapfigure}

\noindent  \textit{CASE 2.1: $T_1$ is a northeastern turn.} It follows from the proof of Lemma 4.17 that $\text{Sector}(T_1)$ is east of the stack of $A_0$'s, that both leaves of $T_1$ are in $F$, and that neither leaf of $T_1$ is an end-box of $P(X,Y)$. End of Case 2.1. See Figure 4.62 (a).

\null

\noindent  \textit{CASE 2.2: $T_1$ is a southeastern turn.}  It follows from the proof of Lemma 4.17 that $\text{Sector}(T_1)$ is southeast of the stack of $A_0$'s. More precisely, we can check that $\text{Sector}(T_1)$ is below the line $y = x + (l' - k')$, and the j-stack of $A_0$'s is above the line $y = x + (l' - k')$. By the proof of Lemma 4.17, we have that both leaves of $T_1$ are in $F$, and that no leaf of $T_1$ is an end-box of $P(X,Y)$. See Figure 4.62 (b). End of Case 2.2. End of Case 2. $\square$

\null

\noindent \textbf{Corollary 4.22.} Let $H$ be a Hamiltonian path of an $m \times n$ grid graph $G$, and let $e_F = e(k,k+1;l)$. Assume that $e_F$ is followed by a southern $j$-stack of $A_0$'s, which is then followed by a standard southern looping fat path $F = G\langle N[P(X,Y)] \rangle$. Then $F$ has an admissible turn $T$ such that $\text{Sector}(T)$ avoids the $j$-stack of $A_0$'s as well as all boxes incident on $e_F$. Furthermore, $T$ may be chosen to lie below the line $y = x - k + l -2j - 2$, or below the line $y = -x + k + l -2j - 1$.

\null 

\noindent\textit{Proof.} One such turn $T$ lying below the line $y = x - k + l -2j - 2$ can be found by inspecting the turns identified in Lemma~4.21. By the same argument and symmetry, we may find another turn that lies below the line $y = -x + k + l -2j - 1$. $\square$

\null 

\noindent In the next three sections we show that every looping fat path contains an admissible turn $T$ that has a weakening. We will need three separate versions of this fact addressing three distinct cases:

1.   $H$ is a Hamiltonian cycle (Section~4.4);

2.  $H$ is a Hamiltonian path and $\text{Sector}(T)$ contains at most one end-vertex $u$ with 

\hspace{0.45cm} $u \in R_0$ (Section~4.5); and

3.   $H$ is a Hamiltonian path and $\text{Sector}(T)$ contains at most one end-vertex $u$ with 

\hspace{0.45cm} $u \notin R_0$ (Section~4.6).

\noindent The first will be needed for the proof of Lemma 3.13 while thee second and third are used in Chapter 5.

%pagemarker

\subsection{Weakenings when $H$ is a Hamiltonian cycle of $G$}

\noindent \textbf{Lemma 4.23.} Let $H$ be a Hamiltonian cycle of an $m\times n$ grid graph $G$, and let $T$ be a turn in $H$ with $d(T)\geq 3$. Then:

I.  \ $T$ has a \index{short weakening}short weakening or $T$ has a \index{lengthening}lengthening. 

II. If $T'$ is a lengthening of $T$ and $T'$ has a weakening of length at most $s$, then $T$ has 

\hspace{0.45cm} a \index{weakening}weakening of length at most $s+1$, with $s+1 \leq \min(m,n)$.

\null

\noindent We prove Lemma 4.23 after we use it to prove Proposition 4.24. 

\null 

\noindent \textbf{Proposition 4.24.} Let $H$ be a Hamiltonian cycle of an $m\times n$ grid graph $G$, and let $T$ be a \index{turn}turn in $H$ with $d(T)\geq 3$. Then $T$ has a weakening of length at most $\min(m,n)$.

\null 

\noindent \textit{Proof.} Let $T=T_0$ be a turn of $H$ with $d(T) \geq 3$. 
If $T_0$ has a short weakening, then we're done, so we assume $T_0$ has no short weakening. By I in Lemma 4.23, $T_0$ has a lengthening $T_1$. So, $T_1\in \mathcal{T}$, where $\mathcal{T}=\mathcal{T}(T_0)$. Since $m,n < \infty$, we have that $|\mathcal{T}| < \infty$. Let $\mathcal{T}=\{T_0, T_1, \ldots, T_j\}$. %Then $|\mathcal{T}|=j+1$, $T_j \in \mathcal{T}$ and $T_{j+1} \notin \mathcal{T}$. 
Then $T_j$ has no lengthening; thus, by I of Lemma 4.23, it must have a short weakening. Then, by induction and II on Lemma 4.23, $T_0$ has a weakening. The bound follows immediately. $\square$

\null 

\noindent \textit{Proof of Lemma 4.23.} We first remark that none of the moves we use throughout this proof fit the description of the moves in Observation 3.4 (i) and (ii) in Section 3. We will use this fact repeatedly and implicitly.

Let $H$ be a Hamiltonian cycle of $G$ and let $T$ be a turn of $H$ with $d(T) \geq 3$. For definiteness, assume that $T$ is northeastern and that $T$ is on $\{e(k;l-1,l), S_{\downarrow}(k+1,l;k',l'+1), e(k'-1,k';l')\}$. Let $L_{\textrm{N}}$ be the northern leaf of $T$ and let $L_{\textrm{E}}$ be the eastern leaf of $T$. Since $d(T)\geq 3$, $m-1 \geq k+3$ and $0\leq l-3$. $L_{\textrm{N}}$ can be open or closed, so there are two cases to check.

\begingroup
\setlength{\intextsep}{0pt}
\setlength{\columnsep}{20pt}
\begin{wrapfigure}[]{l}{0cm}
\begin{adjustbox}{trim=0cm 0.5cm 0cm 0.5cm}

\end{adjustbox}
\end{wrapfigure}

\noindent \textit{CASE 1: $L_{\textrm{N}}$ is closed. Proof of I.}  First we note $n-1 \neq l$, otherwise $H$ misses $v(k+2,l)$. Then we must have $S_{\downarrow}(k+2,l+1;k+3,l) \in H$. See Figure 4.63 (a). Now, $n-1 = l+1$, $n-1 = l+2$, or $n-1 \geq l+3$.

\null 

\noindent \textit{CASE 1.1: n-1=l+1.} By Corollary 1.3.15(b), $L_N+(0,1) \in \text{int}(H)$. This implies that $L_N+(2,1)$ is a small cookie of $H$, so $e(k+3;l,l+1) \in H$. Then $e(k+3;l-2,l-1) \in H$. It follows that $L_N+(2,-2)=L_E$. But then $L_E \mapsto L_E+(0,1)$ is a short weakening of $T$. See Figure 4.63 (b).  End of Case 1.1.

\null 

\endgroup

\begingroup
\setlength{\intextsep}{0pt}
\setlength{\columnsep}{20pt}
\begin{wrapfigure}[]{r}{0cm}
\begin{adjustbox}{trim=0cm 0.5cm 0cm 1.5cm}

\end{adjustbox}
\end{wrapfigure}

\noindent \textit{CASE 1.2: n-1= l+2.} Either $e(k,k+1;l+1) \in H$ or $e(k,k+1;l+1) \notin H$. 

\null 

\noindent \textit{CASE 1.2(a): $e(k,k+1;l+1) \in H$.} Either $L_N+(0,2) \in \text{int}(H)$ or $L_N+(0,2) \in \text{ext}(H)$. If $L_N+(0,2) \in \text{int}(H)$, then $L_N+(0,1) \mapsto L_N$ is a short weakening of $T$. See Figure 4.64 (a). Suppose then that $L_N+(0,2) \in \text{ext}(H)$. This implies that $L_N+(0,2)$ is a small cookie. Then $L_N+(0,1) \mapsto L_N+(0,2)$, $L_N \mapsto L_N+(1,1)$ is a short weakening. See Figure 4.64 (b).  End of Case 1.2(a).

\null 

\begingroup
\setlength{\intextsep}{0pt}
\setlength{\columnsep}{20pt}
\begin{wrapfigure}[]{l}{0cm}
\begin{adjustbox}{trim=0cm 0cm 0cm 0cm}

\end{adjustbox}
\end{wrapfigure}

\noindent \textit{CASE 1.2(b): $e(k,k+1;l+1) \notin H$.} Then $S_{\rightarrow}(k-1,l+1;k,l+2) \in H$ and $S_{\downarrow}(k+1,l+2;k+2,l+1) \in H$. Note that if $e(k,k+1;l+2) \in H$, then $H$ misses $v(k+2,l+2)$ (Figure 4.65 (a)), so we may assume that $e(k,k+1;l+2) \notin H$ (Figure 4.65 (b)). It follows that $e(k+1,k+2;l+2) \in H$. Then $L_N+(0,2) \mapsto L_N+(1,2)$\footnotemark[1]\footnotetext[1]{This move will be referenced in Chapter 5, in the proof of Proposition 5.4.}, $L_N+(0,1) \mapsto L_N$ is a short weakening. End of Case 1.2(b). End of Case 1.2.

\endgroup 

\null 

\noindent \textit{CASE 1.3: $n-1 \geq l+3$.} By Case 1.2, we may assume that $e(k,k+1;l+1) \notin H$, $S_{\rightarrow}(k-1,l+1;k,l+2) \in H$ and that $S_{\downarrow}(k+1,l+2;k+2,l+1) \in H$. Now $L_{\textrm{E}}$ is either open or closed.

\begingroup
\setlength{\intextsep}{0pt}
\setlength{\columnsep}{15pt}
\begin{wrapfigure}[]{r}{0cm}
\begin{adjustbox}{trim=0cm 0.75cm 0cm 0.25cm}


\end{adjustbox}
\end{wrapfigure}

\null 

\noindent \textit{CASE 1.3(a): $L_{\textrm{E}}$ is closed.} By previous cases and symmetry we may assume that $m-1 \geq k'+3$. Using symmetry once more, we may assume that $e(k'+1;l',l'+1) \notin H$. Then the turn $\hat{T}$  on $\{e(k; l+1,l+2), S_{\downarrow}(k+1,l+2;k'+2,l'+1), e(k'+1,k'+2;l')\}$ is in $H$ and it is a lengthening of $T$. See Figure 4.66. End of proof of I for Case 1.3(a).

\null 

\noindent \textit{Proof of II for Case 1.3 (a).} WLOG assume that the last move $\mu_s$ of a weakening $\mu_1,\ldots, \mu_s$ of $\hat{T}$ is $Z \mapsto \hat{L}_{\textrm{N}}$, where $\hat{L}_{\textrm{N}}$ is the northern leaf of $\hat{T}$.  Then $\mu_1,\ldots, \mu_s, L_N+(0,1) \mapsto L_N$, is a weakening of $T$. 

\noindent It remains to check that $s+1 \leq \min(m,n)$. Since the $j^{\text{th}}$ lengthening $T_j$ in $\mathcal{T}(T)$ is $j$ units north and east of $T$, and $d(T) \geq 3$, there can be at most $\min(m,n)-3$ such lengthenings. Since a turn with no lengthening has a short weakening, $s+1 \leq 3+ \min(m,n)-3 =\min(m,n)$. End of proof of II for Case 1.3(a). End of Case 1.3(a).

\endgroup 

\null 

\noindent \textit{CASE 1.3(b): $L_{\textrm{E}}$ is open.} Either $m-1 < k'+2$ or $m-1 \geq k'+2$. It will follow from Case 2 that if a turn has an open leaf adjacent to the boundary or at distance one away from the boundary, then we can find a weakening outright. Therefore, we may assume that $m-1 \geq k'+2$. 

If $e(k';l'+1,l'+2) \in H$, then there is a weakening $L_E \mapsto L_E+(0,1)$, so we may assume that $e(k';l'+1,l'+2) \notin H$. Then the turn $\hat{T}$ on $\{e(k; l+1,l+2), S_{\downarrow}(k+1,l+2;k'+1,l'+2), e(k',k'+1;l'+1)\}$ is in $H$ and it is a lengthening of $T$. See Figure 4.67. End of proof of I for Case 1.3.

\begingroup
\setlength{\intextsep}{0pt}
\setlength{\columnsep}{20pt}
\begin{wrapfigure}[]{l}{0cm}
\begin{adjustbox}{trim=0cm 0cm 0cm 0cm}


\end{adjustbox}
\end{wrapfigure}

\null 

\noindent \textit{Proof of II for Case 1.3(b).} Let $\hat{L_N}$ and $\hat{L_E}$ be the northern and eastern leaves of $\hat{T}$ respectively, and let $\mu_1,\ldots, \mu_s$ be a weakening of $\hat{T}$. If $\mu_s$ is the move $X \mapsto \hat{L}_N$, then, as in Case 1.3(a), $\mu_1,\ldots, \mu_s, L_N+(0,1) \mapsto L_N$, is a weakening of $T$. Suppose then that $\mu_s$ is the move $Z' \mapsto \hat{L}_E$. Then $\mu_1,\ldots, \mu_s, L_E \mapsto L_E+(0,1)$, is a weakening of $T$. The argument that $s+1 \leq \min(m,n)$ is the same as the one in Case 1.3(a), so we omit it. End of proof of II for Case 1.3(b) End of Case 1.3(b). End of Case 1.3. End of Case 1.

\endgroup

\null 

\begingroup
\setlength{\intextsep}{0pt}
\setlength{\columnsep}{20pt}
\begin{wrapfigure}[]{r}{0cm}
\begin{adjustbox}{trim=0cm 0cm 0cm 0cm}
\begin{tikzpicture}[scale=1.75]

\begin{scope}[xshift=0cm, yshift=0cm]
\draw[gray,very thin, step=0.5cm, opacity=0.5] (0,0) grid (1,1);
\draw[yellow, line width =0.2mm] (0,1) -- (1,1);

%%%%%%%%%% Integer Coordinates %%%%%%%%

\node[right] at (1,1) [scale=1]{\tiny{$\ell$}};
\node[above] at (0, 1) [scale=1] {\tiny{k}};
\node[above] at (0.5, 1) [scale=1] {\tiny{+1}};

\fill[blue!40!white, opacity=0.5] (0,0)  rectangle (1,0.5);
\fill[blue!40!white, opacity=0.5] (0,0.5)  rectangle (0.5,1);

\draw[blue, line width=0.5mm] (0,0.5)--++(0,0.5);
\draw[blue, line width=0.5mm] (1,1)--++(-0.5,0)--++(0,-0.5)--++(0.5,0);

\draw[black, line width=0.15mm] (0.25,0.95)--++(0,0.1);
\draw[black, line width=0.15mm] (0.2,0.95)--++(0,0.1);
\draw[black, line width=0.15mm] (0.3,0.95)--++(0,0.1);

\node at  (0.25,0.75 ) [scale=0.8]{\small{$L_{\textrm{N}}$}};

\node[right, align=left, text width=2.5cm] at (1.25, 0.5) { Fig. 4.68. Case  2, \\ $n-1=l$.};

%   n =l+1 
\end{scope}

\end{tikzpicture}
\end{adjustbox}
\end{wrapfigure}

\noindent \textit{CASE 2: $L_{\textrm{N}}$ is open. Proof of I.} If $n-1=l$ then we must have $e(k+1, k+2;l) \in H$. Then $L_{\textrm{N}} \mapsto L_{\textrm{N}}+(1,0)$\footnotemark[1] is a weakening. See Figure 4.68.  Therefore, we may assume that $n-1 >l$.

\null 

\noindent \textit{CASE 2.1: $n-1=l+1$.} Either $e(k+1, k+2;l) \in H$ or  $e(k+1, k+2;l) \notin H$.

\endgroup 

\null 

\noindent \textit{CASE 2.1(a): $e(k+1, k+2;l) \in H$.} Then $e(k,k+1;l+1) \in H$ and $e(k+1,k+2;l+1) \in H$. Then by Corollary 1.3.15 (b), $L_N+(0,1) \in \textrm{int}H$, and so $L_N \in \textrm{int}H$ and $L_N+(1,0) \in \textrm{ext}H$. Either $k=0$, or $k>0$.

\begingroup
\setlength{\intextsep}{0pt}
\setlength{\columnsep}{20pt}
\begin{wrapfigure}[]{l}{0cm}
\begin{adjustbox}{trim=0cm 0cm 0cm 1.25cm}

\end{adjustbox}
\end{wrapfigure}

\noindent \textit{CASE 2.1($a_1$): $k>0$.} Then $L_N+(-1,0) \in \textrm{ext}H$, and $L_N+(-1,0)$ is either a small cookie of $H$ or it is not. If the former, then $L_N \mapsto L_N+(-1,0)$ is a short weakening; and if the latter then $L_N \mapsto L_N+(1,0)$ is a short weakening. See Figure 4.69. End of Case 2.1($a_1$).

\endgroup 

\null 
\null

\begingroup
\setlength{\intextsep}{0pt}
\setlength{\columnsep}{20pt}
\begin{wrapfigure}[]{r}{0cm}
\begin{adjustbox}{trim=0cm 0cm 0cm 1cm}
%
\end{adjustbox}
\end{wrapfigure}

\noindent \textit{CASE 2.1($a_2$): $k=0$.} Then $e(0;l,l+1) \in H$, $e(0;l-2,l-1) \in H$, and $S_{\rightarrow}(0,l-2;1,l-3) \in H$. This implies that $0\leq l-4$ and that $S_{\leftarrow}(1,l-3;0,l-4) \in H$. Then we must have $S_{\uparrow}(2,l-4;3,l-3) \in H$ as well, that $L_E=L_N+(2,-2)$ and that $L_E+(0,-1) \in \textrm{ext}H$. See Figure 4.70. Note that if $L_E+(0,-1)$ is a small cookie, then $L_E \mapsto L_E+(0,-1)$ is a short weakening, so we may assume that $L_E+(0,-1)$ is not a small cookie. Note that this implies that $0\leq l-5$. 

\noindent If $e(3;l-4,l-3) \in H$, then again $L_E \mapsto L_E+(0,-1)$ is a short weakening. Similarly, if $e(3;l-2,l-1) \in H$, then $L_E \mapsto L_E+(0,1)$ is a short weakening. Therefore we only need to check the case where $e(3;l-4,l-3) \notin H$ and $e(3;l-2,l-1) \notin H$. Then $L_E+(1,-1) \in \textrm{ext}H$, and by the assumption that $H$ is Hamiltonian, $m-1 \geq 4$. Then we have $S_{\downarrow}(3,l;4,l-1) \in H$. 

\begingroup
\setlength{\intextsep}{0pt}
\setlength{\columnsep}{20pt}
\begin{wrapfigure}[]{l}{0cm}
\begin{adjustbox}{trim=0cm 0.5cm 0cm 0cm}

\end{adjustbox}
\end{wrapfigure}

\noindent Note that either $e(2,3;l+1) \in H$ and $e(2,3;l) \in H$, or $e(2;l,l+1) \in H$ and $e(3;l,l+1) \in H$. Either way, we must have $e(3,4;l) \notin H$ and $e(3,4;l+1) \in H$. Now, either $e(3;l-3,l-2) \in H$ or $e(3;l-3,l-2) \notin H$.

\null

\noindent \textit{CASE 2.1($a_2$).(i): $e(3;l-3,l-2) \in H$.} Then $L_E+(1,0) \in \textrm{ext}H$. By Corollary 1.3.15(b), this implies that $m-1 \geq 5$. If $e(4;l-3,l-2) \in H$, then $L_E+(1,0) \mapsto L_E$ is a short weakening, so we may assume that $e(4;l-3,l-2) \notin H$. Then $S_{\downarrow}(4,l-1;5,l-2) \in H$ and $S_{\uparrow}(4,l-4;5,l-3) \in H$. We must also have that $S_{\downarrow}(4,l+1;5,l)\in H$ and $e(5;l,l+1)\in H$. Then $e(5;l-2,l-1)\in H$ as well. See Figure 4.71. Then $L_E+(2,0) \mapsto L_E+(2,1)$, $L_E+(1,0) \mapsto L_E$ is a short weakening of $T$. End of Case 2.1($a_2$).(i).

\endgroup 

\null

\noindent \textit{CASE 2.1($a_2$).(ii): $e(3;l-3,l-2) \notin H$.} Then $e(3,4;l-3) \in H$ and $e(3,4;l-2) \in H$. Now, either $e(2,3;l+1) \in H$ and $e(2,3;l) \in H$, or $e(2;l,l+1) \in H$ and $e(3;l,l+1) \in H$.

\begingroup
\setlength{\intextsep}{00pt}
\setlength{\columnsep}{20pt}
\begin{wrapfigure}[]{r}{0cm}
\begin{adjustbox}{trim=0cm 0cm 0cm 0.25cm}

\end{adjustbox}
\end{wrapfigure}

\noindent \textit{CASE 2.1($a_2$).$(ii)_1$:  $e(2,3;l+1) \in H$ and $e(2,3;l) \in H$.} Then $e(4;l-2,l-1) \notin H$. Note that if $e(4;l-1,l)\in H$, then $L_E+(1,1) \mapsto L_E+(1,2)$, $L_E \mapsto L_E+(0,1)$ is a short weakening of $T$, so we may assume that $e(4;l-1,l)\notin H$. Then $e(4,5;l-1) \in H$, $S_{\downarrow}(4,l+1;5,l) \in H$, and $e(5;l,l+1) \in H$. Then $L_E+(2,2) \mapsto L_E+(2,3)$, $L_E+(1,1) \mapsto L_E+(1,2)$, $L_E \mapsto L_E+(0,1)$ is a short weakening of $T$. See Figure 4.72 (a). End of Case 2.1($a_2$).$(ii)_1$.

\null 

\noindent \textit{CASE 2.1($a_2$).$(ii)_2$:  $e(2;l,l+1) \in H$ and $e(3;l,l+1) \in H$.} Now, if $e(4;l-2,l-1) \notin H$, then we can use the same argument and find the same cascades as in Case 2.1($a_2$).$(ii)_1$ (Figure 4.72 (a)); and if $e(4;l-2,l-1) \in H$, then we must have that $e(5;l-2,l-1) \in H$ as well. Then $L_E+(2,1) \mapsto L_E+(1,1)$, $L_E \mapsto L_E+(0,1)$ is a short weakening of $T$. See Figure 4.72 (b). End of Case 2.1($a_2$).$(ii)_2$. End of Case 2.1($a_2$).$(ii)$.
End of Case 2.1($a_2$). End of Case 2.1(a).

\endgroup

\null 

\noindent \textit{CASE 2.1(b): $e(k+1, k+2;l) \notin H$.} Then we must have $e(k+1;l,l+1) \in H$ and $S_{\downarrow}(k+2,l+1;k+3,l) \in H$. Now, either $e(k+1,k+2;l+1) \notin H$ or $e(k+1,k+2;l+1) \in H$.

\begingroup
\setlength{\intextsep}{0pt}
\setlength{\columnsep}{20pt}
\begin{wrapfigure}[]{r}{0cm}
\begin{adjustbox}{trim=0cm 1cm 0cm 0cm}

\end{adjustbox}
\end{wrapfigure}

\null 

\noindent \textit{CASE 2.1($b_1$): $e(k+1, k+2;l+1) \notin H$.} Then we must have $e(k+2,k+3;l+1) \in H$ and that $L_N+(1,1) \in \textrm{ext}H$ is the neck of the large cookie. Then $k>0$, or $k=0$.

If $k>0$, then after $L_N+(1,1) \mapsto L_N+(2,1)$ \footnotemark[1]\footnotetext[1]{This move will be referenced in Chapter 5, in the proof of Proposition 5.4.}, we are back to Case 2.1($a_1$). And if $k=0$ then we are effectively in the same scenario as in Case 2.1($a_2$), except that in this case, in the analogues of Cases 2.1$(a_2).(i)$ and $(ii)$, we get lengthenings instead of short weakenings. See Figure 4.73 (a). End of Case 2.1($b_1$). 

\null 

\noindent \textit{CASE 2.1($b_2$): $e(k+1, k+2;l+1) \in H$.} Then $e(k,k+1;l+1) \notin H$ and $S_{\downarrow}(k+2,l+1;k+3,l) \in H$. Then we must have that $e(k+3;l,l+1)\in H$ as well. This implies that $e(k+3;l-2,l-1)\in H$. It follows that $L_N+(2,-2)=L_E$, and that $L_E$ is open. Then $L_E \mapsto L_E+(0,1)$ is a short weakening. See Figure 4.73 (b). End of Case 2.1($b_2$). End of Case 2.1(b). End of Case 2.1.

\endgroup 

\null 

\begingroup
\setlength{\intextsep}{0pt}
\setlength{\columnsep}{20pt}
\begin{wrapfigure}[]{l}{0cm}
\begin{adjustbox}{trim=0cm -0.25cm 0cm 0cm}
\begin{tikzpicture}[scale=1.5]

\draw[gray,very thin, step=0.5cm, opacity=0.5] (0,0) grid (2,2);

%%%%%%%%%% Integer Coordinates %%%%%%%%
{
\foreach \x in {1,...,1}
\node[left] at (0,0.5*\x+1.5) [scale=1]{\tiny{\x}};
\node[left] at (0,1.5) [scale=1]{\tiny{$\ell$}};

\node[left] at (0,0.5) [scale=1]{\tiny{+1}};
\node[left] at (0,0) [scale=1]{\tiny{$\ell'$}};

\node[above] at (0.0, 2) [scale=1]{\tiny{k}};
\node[above] at (0.5, 2) [scale=1]{\tiny{+1}};

\foreach \x in {1,...,1}
\node[above] at (0.5*\x+1.5, 2) [scale=1]{\tiny{\x$'$}};
\node[above] at (1.5, 2) [scale=1]{\tiny{k$'$}};
}

%%%%%%%%%%%%     DFP-0   %%%%%%%%%%%%  
{
\draw[blue, line width=0.5mm] (0.0,1)--++(0,0.5);
\draw[blue, line width=0.5mm] (0.5,1)--++(0,0.5);

\draw[blue, line width=0.5mm] (1.0,0.5)--++(0.5,0);
\draw[blue, line width=0.5mm] (1.0,0)--++(0.5,0);
}

%green stuff and 
{
\draw[orange, line width=0.5mm](0.5,1.5)--++(0,0.5);
\draw[orange, line width=0.5mm](1.0,1.5)--++(0,0.5);

\draw[orange, line width=0.5mm](1.5,1)--++(0.5,0);
\draw[orange, line width=0.5mm](1.5,0.5)--++(0.5,0);
}

%%%%%%%%%% Circles %%%%%%%%
{
\tikzset{myCircle/.style={orange}}

\foreach \x in {0,...,2}
\fill[myCircle, orange] (1.125+0.125*\x,1.375-0.125*\x) circle (0.05);

\foreach \x in {0,...,2}
\fill[myCircle, blue] (0.625+0.125*\x,0.875-0.125*\x) circle (0.05);
}

% no edge signs
{
\draw[black, line width=0.15mm] (0.7,1.45)--++(0,0.1);
\draw[black, line width=0.15mm] (0.75,1.45)--++(0,0.1);
\draw[black, line width=0.15mm] (0.8,1.45)--++(0,0.1);

\draw[black, line width=0.15mm] (1.45,0.7)--++(0.1,0);
\draw[black, line width=0.15mm] (1.45,0.75)--++(0.1,0);
\draw[black, line width=0.15mm] (1.45,0.8)--++(0.1,0);
}

% labeling
{
\node at (0.25,1.25) [scale=0.8]{\small{$L_{\textrm{N}}$}};
\node at (0.75,1.75) [scale=0.8]{\small{$\hat{L}_{\textrm{N}}$}};
\node at (1.25,0.25) [scale=0.8]{\small{$L_{\textrm{E}}$}};

\node[below, align=center, text width=5cm] at (1, -0.1) { Fig. 4.74. Case 2.2.};

}

\end{tikzpicture}

\end{adjustbox}
\end{wrapfigure}

\noindent \textit{CASE 2.2: $n-1 \geq l+2$.} By previous cases we may assume that $m-1 \geq k'+2$, $e(k+1, k+2;l) \notin H$, $e(k+1;l,l+1) \in H$ and $S_{\downarrow}(k+2,l+1;k+3,l) \in H$. If $L_{\textrm{E}}$ is closed, then we're done by Case 1, so we may assume that $L_{\textrm{E}}$ is open. By Case 2.1, we may assume that $e(k';l'+1,l'+2) \notin H$. Then the turn $\hat{T}$ on $\{e(k+1; l,l+1), S_{\downarrow}(k+2,l+1; k'+1,l'+2), e(k',k'+1;l'+1)\}$ is in $H$ and it is a lengthening of $T$. See Figure 4.74. End of proof of I for Case 2.2.

The proof of II for Case 2.2 is the same as the proof of II for Case 1.3(a). End of proof for Case 2. $\square$

\endgroup

\null

\noindent \textbf{Observation 4.25.} All turn weakenings found in Lemma 4.23 are contained in Sector$(T)$.

\subsection{Weakenings when $H$ is a Hamiltonian path of $G$ with $u \in R_0$}

\noindent \textbf{Lemma 4.26.} Let $G$ be an $m\times n$ grid graph, let $H$ be a Hamiltonian path with end-vertices $u$ and $v$, and let $T$ be a \index{turn}turn in $H$ with $d(T)\geq 3$. Assume that $u$ and $v$ are not adjacent, that $u \in R_0$, and that $v$ is not in Sector$(T)$. Then:

I.  \ $T$ has a  \index{short weakening}short double-switch weakening\footnote{That is, a weakening consisting entirely of double-switch moves.} (SDSW) or $T$ has a lengthening.

II. If $T'$ is a \index{lengthening}lengthening of $T$ and $T'$ has a double-switch weakening of length at most 

\hspace{0.4cm}  $s$, then $T$ has a double-switch \index{weakening}weakening of length at most $s+1$, with 

\hspace{0.4cm} $s+1 \leq \min(m,n)$.

\null

\noindent \textbf{Proposition 4.27.} Let $G$ be an $m\times n$ grid graph, let $H$ be a Hamiltonian path with end-vertices $u$ and $v$, and let $T$ be a turn in $H$ with $d(T)\geq 3$. Assume that $u$ and $v$ are not adjacent, that $u \in R_0$, and that $v$ is not in Sector$(T)$. Then $T$ has a double-switch weakening of length at most $\min(m,n)$ moves.

\null 

\noindent \textit{Proof.} Let $T=T_0$ be a turn of $H$ with $d(T) \geq 3$. Assume that $u \in R_0$, and that $v \notin \textrm{Sector}(T)$.
If $T_0$ has a short double-switch weakening, we are done, so we may assume that $T_0$ does not have a short double-switch weakening. By I in Lemma 4.26, $T_0$ has a lengthening $T_1$. So, $T_1\in \mathcal{T}$, where $\mathcal{T}=\mathcal{T}(T_0)$. Since $m,n < \infty$, we have that $|\mathcal{T}| < \infty$. Let $\mathcal{T}=\{T_0, T_1, \ldots, T_j\}$. 
Then $T_j$ has no lengthening. So, by I in Lemma 4.26, $T_j$ must have a short double-switch weakening. Then, by induction and II on Lemma 4.26, $T_0$ has a double-switch weakening. The bound follows immediately. $\square$ 

\null

\noindent \textit{Proof of Lemma 4.26.} For definiteness, assume that $T$ is the northeastern turn on $\{e(k;l-1,l), S_{\downarrow}(k+1,l;k',l'+1), e(k'-1,k',l')\}$. Let $L_N$ and $L_E$ be the northern and eastern leaves of $T$, respectively. We will consider three cases: $L_N$ is closed and $n-1<l+2$, $L_N$ is open and $n-1<l+2$, and $n-1 \geq l+2$.

\begin{center}
\begin{adjustbox}{trim=0cm 0cm 0cm 0cm} 

\end{adjustbox}
\end{center}

\noindent \textit{CASE 1. $L_N$ is closed and $n-1<l+2$}. Suppose that $n-1=l+1$. Observe that (NAA) implies that $e(k,k+1;l+1) \in H$. Then $L_N+(0,1) \mapsto L_N$ is a SDSW. See Figure 4.75 (a). It remains to check the case where $n-1=l$.

Suppose that $n-1=l$. Then $v(k+2,l)=u$, and $e(k+2,k+3;l) \in H$. Note that if $e(k+3;l-2,l-1) \in H$, then we must have that $L_E=L_N+(2,-2)$, and that $L_E$ is open. Then $L_E \mapsto L_E+(0,1)$ is a SDSW. See Figure 4.75 (b). Therefore we may assume that $e(k+3;l-2,l-1) \notin H$. Since $v \notin \text{Sector}(T)$, this implies that $m-1 >k+3$. Then we have that $S_{\downarrow}(k+3,l;k+4,l-1) \in H$, and $e(k+4;l-1,l) \in H$. Now, $e(k+3,k+4;l-2)\in H$, or $e(k+3,k+4;l-2)\notin H$.

\endgroup 

\null 

\noindent \textit{CASE 1.1: $e(k+3,k+4;l-2)\in H$.} This also implies that  $L_E=L_N+(2,-2)$, and that $L_E$ is open. Then $L_E+(1,1) \mapsto L_E+(1,2)$, $L_E \mapsto L_E+(0,1)$ is a SDSW. See Figure 4.75 (c). End of Case 1.1.

\begingroup
\setlength{\intextsep}{0pt}
\setlength{\columnsep}{20pt}
\begin{center}
\begin{adjustbox}{trim=0cm 0cm 0cm 0cm} 

\end{adjustbox}
\end{center}

\noindent \textit{CASE 1.2: $e(k+3,k+4;l-2) \notin H$.} Then $e(k+3;l-3,l-2) \in H$ and $e(k+4;l-3,l-2) \in H$. Now, either $e(k+2,k+3;l-3)\in H$, or $e(k+3,k+4;l-3)\in H$. If the former, then $L_E=L_N+(2,-2)$, $L_E$ is closed, and $L_E+(1,0) \mapsto L_E$ is a SDSW. See Figure 4.76 (a). And if the latter, $L_E=L_N+(3,-3)$, $L_E$ is open, and $L_E \mapsto L_E+(0,1)$ is a SDSW. Figure 4.76 (b). End of Case 1.2. End of Case 1.

\null

\endgroup

\begingroup
\setlength{\intextsep}{0pt}
\setlength{\columnsep}{20pt}
\begin{wrapfigure}[]{r}{0cm}
\begin{adjustbox}{trim=0cm 0cm 0cm 0cm} 
\begin{tikzpicture}[scale=1.75]

\begin{scope}[xshift=0cm] 
{
\draw[gray,very thin, step=0.5cm, opacity=0.5] (0,0) grid (1,0.5);

% turn

\draw[blue, line width=0.5mm] (0,0)--++(0,0.5); 
\draw[blue, line width=0.5mm] (0.5,0.5)--++(0,-0.5)--++(0.5,0); 

\draw[fill=blue] (0.5,0.5) circle (0.05);
\draw[fill=blue] (1,0.5) circle (0.05);

%black lines
{
\draw[black, line width=0.15mm] (0.2,0.45)--++(0,0.1);
\draw[black, line width=0.15mm] (0.25,0.45)--++(0,0.1);
\draw[black, line width=0.15mm] (0.3,0.45)--++(0,0.1);

\draw[black, line width=0.15mm] (0.7,0.45)--++(0,0.1);
\draw[black, line width=0.15mm] (0.75,0.45)--++(0,0.1);
\draw[black, line width=0.15mm] (0.8,0.45)--++(0,0.1);

}

% labellings 
{

\node[left] at (0,0.5) [scale=1]{\tiny{$\ell$}};

\node[below] at (0,0) [scale=1]
{\tiny{k}};
\node[below] at (0.5,0) [scale=1]
{\tiny{+1}};
\node[below] at (1,0) [scale=1]
{\tiny{+2}};

\node at (0.25, 0.25) [scale=0.8]{$L_N$};

\node[above] at (0.5,0.5) [scale=1.25]{(a)};

\node[below, align=left, text width=5cm] at (1.5, -0.25) {Fig. 4.77. (a) $e(k{+}1,k{+}2;l) {\in} H$. (b) $e(k{+}1,k{+}2;l) {\notin}) H$.};

}

}

\end{scope}

\begin{scope}[xshift=2cm] 
{
\draw[gray,very thin, step=0.5cm, opacity=0.5] (0,0) grid (1,0.5);

% turn

\draw[blue, line width=0.5mm] (0,0)--++(0,0.5); 
\draw[blue, line width=0.5mm] (0.5,0.5)--++(0,-0.5)--++(0.5,0); 

\draw[blue, line width=0.5mm] (0.5,0.5)--++(0.5,0);

%black lines
{
\draw[black, line width=0.15mm] (0.2,0.45)--++(0,0.1);
\draw[black, line width=0.15mm] (0.25,0.45)--++(0,0.1);
\draw[black, line width=0.15mm] (0.3,0.45)--++(0,0.1);

}

% labellings 
{

\node[left] at (0,0.5) [scale=1]{\tiny{$\ell$}};

\node[below] at (0,0) [scale=1]
{\tiny{k}};
\node[below] at (0.5,0) [scale=1]
{\tiny{+1}};
\node[below] at (1,0) [scale=1]
{\tiny{+2}};

\node at (0.25, 0.25) [scale=0.8]{$L_N$};

\node[above] at (0.5,0.5) [scale=1.25]{(b)};

}

}

\end{scope}

\end{tikzpicture}
\end{adjustbox}
\end{wrapfigure}

\noindent \textit{CASE 2. $L_N$ is open and $n-1<l+2$}. Suppose that $n-1=l$. Either $e(k+1,k+2;l) \in H$, or $e(k+1,k+2;l) \notin H$.

\null 

\noindent Note that the latter implies that the end-vertices of $H$ are on $v(k+1,l)$ and $v(k+2,l)$, contradicting the assumption that $v\notin \text{Sector}(T)$. Then, it must be the case that $e(k+1,k+2;l) \in H$.  Then $L_N \mapsto L_N+(1,0)$ is a SDSW. See Figure 4.77. It remains to check the case where $n-1=l+1$.

\endgroup 

\null 

\begingroup
\setlength{\intextsep}{0pt}
\setlength{\columnsep}{20pt}
\begin{wrapfigure}[]{l}{0cm}
\begin{adjustbox}{trim=0cm 0cm 0cm 0cm} 

\end{adjustbox}
\end{wrapfigure}

\noindent Suppose that $n-1=l+1$. Note that if $e(k+1,k+2;l) \in H$ then $L_N \mapsto L_N+(1,0)$ is a SDSW, so we may assume that $e(k+1,k+2;l) \notin H$. Then $e(k+1;l,l+1) \in H$. It follows that $S_{\downarrow}(k+2,l+1;k+3,l) \in H$. There are three possibilities: $e(k+2,k+3;l+1) \in H$ (Figure 4.78 (a)), $v(k+2,l+1)=u$ (Figure 4.78 (b)), and $e(k+1,k+2;l+1) \in H$ (Figure 4.79 (a)).

\noindent If the first, then $L_N+(1,1) \mapsto L_N+(2,1)$, $L_N \mapsto L_N+(1,0)$ is a SDSW. If the second, then we must have $e(k+3;l,l+1)\in H$. It follows that $e(k+3;l-2,l-1)\in H$ as well. It follows that  $L_E=L_N+(2,-2)$, and that $L_E$ is open. Then $L_E \mapsto L_E+(0,1)$ is a SDSW. It remains to check the third.

Suppose that $e(k+1,k+2;l+1) \in H$. As above, if $e(k+3;l-2,l-1)\in H$, then $L_E=L_N+(2,-2)$, $L_E$ is open, and $L_E \mapsto L_E+(0,1)$ is a SDSW, so we may assume that $e(k+3;l-2,l-1)\notin H$. It follows that $m-1 \geq k+4$, that $S_{\downarrow}(k+3,l;k+4;l-1) \in H$, that $u=v(k+3,l+1)$, that $e(k+3,k+4;l+1)\in H$. Now, $e(k+4;l-1,l) \in H$, or $e(k+4;l-1,l) \notin H$.

\endgroup

\begingroup
\setlength{\intextsep}{0pt}
\setlength{\columnsep}{20pt}
\begin{wrapfigure}[]{r}{0cm}
\begin{adjustbox}{trim=0cm 0cm 0cm 0.5cm} 

\end{adjustbox}
\end{wrapfigure}

\null 

\noindent \textit{CASE 2.1. $e(k+4;l-1,l) \in H$}. Now, either $e(k+3,k+4;l-2) \in H$ or $e(k+3,k+4;l-2) \notin H$. 

\null 

\noindent \textit{CASE 2.1(a). $e(k+3,k+4;l-2) \in H$}. Then $L_E=L_N+(2,-2)$, $L_E$ is open, and $L_E+(1,1) \mapsto L_E+(1,2)$, $L_E \mapsto L_E+(0,1)$ is a SDSW. See Figure 4.79 (b). End of Case 2.1(a).

\endgroup

\begingroup
\setlength{\intextsep}{0pt}
\setlength{\columnsep}{20pt}
\begin{wrapfigure}[]{l}{0cm}
\begin{adjustbox}{trim=0cm 0.5cm 0cm 0cm} 
%
\end{adjustbox}
\end{wrapfigure}

\noindent \textit{CASE 2.1(b). $e(k+3,k+4;l-2) \notin H$}. Then $e(k+3;l-3,l-2)\in H$, and $e(k+4;l-3,l-2)\in H$. Now either $e(k+2,k+3;l-3) \in H$, or $e(k+2,k+3;l-3) \notin H$.

\null 

\noindent \textit{CASE 2.1($b_1$). $e(k+2,k+3;l-3) \in H$}. Then $L_E=L_N+(2,-2)$, $L_E$ is closed, and $L_E+(1,0) \mapsto L_E$ is a SDSW. See Figure 4.80 (a). End of Case 2.1($b_1$).

\null 

\noindent \textit{CASE 2.1($b_2$). $e(k+2,k+3;l-3) \notin H$}. Then we must have $e(k+3,k+4;l-3) \in H$, that $L_E=L_N+(3,-3)$, that $L_E$ is open, and that  $L_E \mapsto L_E+(0,1)$ is a SDSW. See Figure 4.80 (b). End of Case 2.1($b_2$). End of Case 2.1(b). End of Case 2.1.

\endgroup 

\null

\noindent \textit{CASE 2.2: $e(k+4;l-1,l) \notin H$.} It follows that $m-1 \geq k+5$, that $S_{\downarrow}(k+4,l+1;k+5,l) \in H$, and that $e(k+5;l,l+1) \in H$. Note that if $e(k+4,k+5;l-1) \in H$, then, after $L_N+(4,0) \mapsto L_N+(4,1)$, we are back to Case 2.1. Note that the weakenings in Case 2.1 required at most two moves. So, including $L_N+(4,0) \mapsto L_N+(4,1)$, we would still need only three moves for a weakening.  See Figure 4.81 (a).

Assume that $e(k+4,k+5;l-1) \notin H$. It follows that $e(k+4;l-2,l-1) \in H$, $e(k+5;l-2,l-1) \in H$, $e(k+3,k+4;l-2) \notin H$, and $e(k+3,l-3,l-2) \in H$. Now, either $e(k+4;l-3,l-2) \in H$ or $e(k+4;l-3,l-2) \notin H$.

\begin{center}
\begin{adjustbox}{trim=0cm 0cm 0cm 0cm} 

\end{adjustbox}
\end{center}

\noindent \textit{CASE 2.2(a): $e(k+4;l-3,l-2) \in H$.} Then $e(k+3,k+4;l-3) \notin H$. It follows that $L_E=L_N+(2,-2)$, and that $L_E$ is closed. Then $L_E+(1,0) \mapsto L_E$ is a SDSW. See Figure 4.81 (b). End of Case 2.2(a)

\begin{center}
\begin{adjustbox}{trim=0cm 0cm 0cm 0cm} 
%
\end{adjustbox}
\end{center}

\noindent \textit{CASE 2.2(b): $e(k+4;l-3,l-2) \notin H$.} Then $e(k+4,k+5;l-2) \in H$. Note that if $e(k+2,k+3;l-3)\in H$, then $L_E=L_N+(2,-2)$, $L_E$ is closed, and $e(k+4,k+5;l-3) \in H$. Then $L_E+(2,0) \mapsto L_E+(2,1)$, $L_E+(1,0) \mapsto L_E$ is a SDSW. See Figure 4.82 (a). So, we may assume that $e(k+2,k+3;l-3)\notin H$ (Figure 4.82 (b)). Then $e(k+3,k+4;l-3) \in H$. Now, either $e(k+4,k+5; l-3) \in H$, or  $e(k+4,k+5; l-3) \notin H$.

\endgroup 

\null

\noindent \textit{CASE 2.2($b_1$): $e(k+4,k+5; l-3) \in H$.} Then $L_E=L_N+(3,-3)$, and $L_E+(1,1) \mapsto L_E+(1,2)$, $L_E \mapsto L_E+(0,1)$ is a SDSW. See Figure 4.83 (a). End of Case 2.2($b_1$).

\null

\noindent \textit{CASE 2.2($b_2$): $e(k+4,k+5; l-3) \notin H$.} It follows that $e(k+4;l-4,l-3) \in H$, and $e(k+5;l-4,l-3) \in H$. Now either $e(k+3,k+4;l-4) \in H$, or $e(k+4,k+5;l-4) \in H$. If the former, then $L_E=L_N+(3,-3)$, $L_E$ is closed, and $L_E+(1,0) \mapsto L_E$ is a SDSW. And if the latter, then $L_E=L_N+(4,-4)$, $L_E$ is open, and $L_E \mapsto L_E+(0,1)$ is a SDSW. See Figure 4.83 (b) and (c). End of Case 2.2($b_2$). End of Case 2.2(b). End of Case 2.2. End of Case 2.

\begingroup
\setlength{\intextsep}{0pt}
\setlength{\columnsep}{20pt}

\begin{center}
\begin{adjustbox}{trim=0cm 0cm 0cm 0.5cm} 

\end{adjustbox}
\end{center}

\noindent \textit{CASE 3: $n-1 \geq l+2$.} This is very similar to Cases 1.3 and 2.2 in Lemma 4.23, so we omit it. End of Case 3. $\square$

\null

%\textcolor{red}{[Odd Observation like 4.25 that all turn weakenings in lemma are contained in Sector(T).]}

\subsection{Weakenings when $H$ is a Hamiltonian path of $G$ with $u \notin R_0$}

Turn weakenings are useful for the reconfiguration of Hamiltonian paths, which are discussed in Chapter 5. So far, we have found turn weakenings using only double-switch moves. When an end-vertex of $H$ lies within the sector of a turn, we will find a weakening that also allows single-switches and backbites. While a double-switch weakening may still exist, it seems that it would require several additional pages of configuration checking, and we believe it would be no shorter than the one we find.

\null

\noindent \textbf{Proposition 4.28.} Let $G$ be an $m\times n$ grid graph, let $H$ be a Hamiltonian path with end-vertices $u$ and $v$, and let $T$ be a turn in $H$ with $d(T)\geq 3$. Assume that $u$ and $v$ are not adjacent, that $u \notin R_0$, and that $v$ is not in Sector$(T)$. Then $T$ has a weakening of length at most $\min(m,n)$ moves.

\null 

\noindent The proof of Proposition 4.28 requires the following lemma.

\null 

\noindent \textbf{Lemma 4.29.} Let $G$ be an $m\times n$ grid graph, let $H$ be a Hamiltonian path with end-vertices $u$ and $v$, and let $T$ be a \index{turn}turn in $H$ with $d(T)\geq 3$. Assume that $u$ and $v$ are not adjacent, that $u \notin R_0$, and that $v$ is not in Sector$(T)$. Then:

I.  \ $T$ has a \index{short weakening}short weakening or $T$ has a \index{lengthening}lengthening.

II. If $T'$ is a lengthening of $T$ and $T'$ has a weakening of length at most $s$, then $T$ has  

\hspace{0.4cm} a \index{weakening}weakening of length at most $s+1$, with  $s+1 \leq \min(m,n)$.

\null 

\noindent \textit{Proof of Proposition 4.28.} This is the same as the proof of Proposition 4.24, once Lemma 4.29 has been proved, so we omit it. $\square$. 

\null

\noindent The proof of Lemma 4.29 requires Lemma 4.30 below. The definitions of edge weakening, end-vertex weakening, end-vertex weakening terminal, and Flank, introduced at the start of Section 4.3, will be used throughout the remainder of this section.

\null 

\noindent \textbf{Lemma 4.30.}  Let $H$ be a Hamiltonian path of $G$ with non-adjacent end-vertices, and let $T$ be a turn of $H$ with $d(T)\geq 3$. Then:

(a) If $\textrm{Flank}(T)$ contains an end-vertex of $H$, then $T$ has a short weakening.

(b) If $T'$ is a lengthening of $T$, and $T'$ has an \index{end-vertex weakening}end-vertex weakening of length at most

\hspace{0.45cm}  $s$, then $T$ has a weakening of length at most $s+1$, with $s+1 \leq \min(m,n)$.

\null 

\noindent \textit{Proof of Lemma 4.29.} For definiteness, let $T$ be the northeastern turn of $H$ on $\{e(k;l-1,l), S_{\downarrow}(k+1,l;k',l'+1), e(k'-1,k',l')\}$. Since $v \notin \text{Sector}(T)$, by Observation 4.20(c), $v \notin \text{Flank}(T)$. If $u \in \text{Flank}(T)$, by Lemma 4.30 (a), $T$ has a short weakening, so we may assume that $u \notin \text{Flank}(T)$. We observe that the proof of Cases 1 and 2 in Lemma 4.26 show that if $n-1<l+2$, or if $m-1<k'+2$, then $T$ has a weakening, so we may assume that $n-1 \geq l+2$ and that $m-1 \geq k'+2$. Let $L_N=R(k,l-1)$ and $L_E=R(k'-1,l')$. Then $L_N$ is either open or closed.

\null 

\noindent \textit{CASE 1. $L_N$ is closed.}  As in Lemma 4.23, if $e(k,k+1;l+1) \in H$, then $T$ has a weakening, so we may assume that $e(k,k+1;l+1) \notin H$. Then $e(k;l+1,l+2) \in H$. Now, $L_E$ is open, or $L_E$ is closed.

\null

\begingroup
\setlength{\intextsep}{0pt}
\setlength{\columnsep}{20pt}
\begin{wrapfigure}[]{r}{0cm}
\begin{adjustbox}{trim=0cm 0.5cm 0cm 0.5cm}


\end{adjustbox}
\end{wrapfigure}

\noindent \textit{CASE 1.1: $L_E$ is open.} If $e(k';l'+1,l'+2) \in H$, then $L_E \mapsto L_E+(0,1)$ is a short weakening, so we may assume that $e(k';l'+1,l'+2) \notin H$. Then $e(k', k'+1;l'+1) \in H$ and $S_{\downarrow}(k+1,l+2;k'+1,l'+2) \in H$. Now the turn $T'$ on $e(k;l+1,l+2)$,  $S_{\downarrow}(k+1,l+2;k'+1,l'+2)$ and $e(k',k'+1;l'+1)$ is a lengthening of $T$. See Figure 4.84. End of proof of I for Case 1.1. 

%pagemarker 

\noindent Assume that $T'$ has a weakening $\mu_1, \ldots, \mu_s$. Either it has an edge weakening or an end-vertex weakening. If $\mu_1, \ldots, \mu_s$ is an edge weakening, after applying it, either $e(k,k+1;l+1)$, or $e(k';l'+1,l'+2)$ is in the resulting Hamiltonian path. If the former then $L_N+(0,1) \mapsto L_N$ is a weakening of $T$; and if the latter, then $L_E \mapsto L_E+(0,1)$ is a weakening of $T$. In either case, the argument that $s+1\leq \min(m,n)$ is the same as the one used in Case 1.3(a) of Lemma 4.23, so we omit it.

If $T'$ $\mu_1, \ldots, \mu_s$ is an end-vertex weakening, then, by Lemma 4.30 (b), $T$ has a weakening of length $s+1 \leq \min(m,n)$. End of Case 1.1.

Case 1.2 and Case 2 are similar, so we omit the proofs. $\square$

\null 

\noindent \textit{Proof of Lemma 4.30(a).} For definiteness, let $T$ be the northeastern turn of $H$ on $\{e(k;l-1,l), S_{\downarrow}(k+1,l;k',l'+1), e(k'-1,k',l')\}$ and assume that $v$ is an end-vertex of $H$ in Flank($T$). Orient $H$ so that $v=v_1$. If either leaf of $T$ is parallel, then switching that leaf is an edge weakening, so we may assume that both leaves of $T$ are anti-parallel. We will use Lemmas 1.4.1 and 1.4.4 repeatedly and implicitly. Now, $v \in$  $s$-Flank$(T)$ or $v \in \ell$-Flank$(T)$, so there are two cases to check.

\begingroup
\setlength{\intextsep}{0pt}
\setlength{\columnsep}{20pt}
\begin{wrapfigure}[]{r}{0cm}
\begin{adjustbox}{trim=0cm 0cm 0cm 0cm} 

\end{adjustbox}
\end{wrapfigure}

\null 

\noindent \textit{CASE 1: $v \in$  $\ell$-Flank$(T)$}.  For definiteness, assume that $v \in V(L_N+(0,1))$. Now, $L_N$ is open or closed.

\null 

\noindent \textit{CASE 1.1: $L_N$ is closed.} If $e(k,k+1;l+1) \in H$ then  $L_N +(0,1)\mapsto L_N$ is a short edge weakening, so we may assume that $e(k,k+1;l+1) \notin H$. Then $v=v(k,l+1)$ or $v=v(k+1,l+1)$. If the former, then $bb_{v}(east)$, $L_N +(0,1) \mapsto L_N$ is a short edge weakening (see Figure 4.85); and if the latter, then  $bb_{v}(west)$, $L_N+(0,1) \mapsto L_N$ is a short edge weakening. End of Case 1.1.

\endgroup 

\null

\begingroup
\setlength{\intextsep}{0pt}
\setlength{\columnsep}{20pt}
\begin{wrapfigure}[]{l}{0cm}
\begin{adjustbox}{trim=0cm 0cm 0cm 0cm} 

\end{adjustbox}
\end{wrapfigure}

\noindent \textit{CASE 1.2: $L_N$ is open.} Then $v=v(k+1,l)$, $v=v(k,l)$,  $v=v(k+1,l+1)$, or $v=v(k,l+1)$.

\null 

\noindent \textit{CASE 1.2(a): $v=v(k+1,l)$.} Then $bb_{v}(east)$, $L_N \mapsto L_N+(1,0)$ is a short edge weakening. See Figure 4.86. End of Case 1.2+(a).

\endgroup

\null 

\begingroup
\setlength{\intextsep}{0pt}
\setlength{\columnsep}{20pt}
\begin{wrapfigure}[]{r}{0cm}
\begin{adjustbox}{trim=0cm 0cm 0cm 0cm} 
%
\end{adjustbox}
\end{wrapfigure}

\noindent \textit{CASE 1.2 (b): $v=v(k,l)$.} Then, after  $bb_{v}(east)$, the end-vertex is on $v(k+1,l-1) \in \text{ew-Terminal}(T)$, so $bb_{v}(east)$ is a short end-vertex weakening. See Figure 4.87. End of Case 1.2(b).

\endgroup

\begingroup
\setlength{\intextsep}{0pt}
\setlength{\columnsep}{20pt}
\begin{wrapfigure}[]{l}{0cm}
\begin{adjustbox}{trim=0cm 0cm 0cm 0cm} 
%
\end{adjustbox}
\end{wrapfigure}

\null 

\noindent \textit{CASE 1.2(c): $v=v(k+1,l+1)$.} By (NAA), $v(k+2,l+1)$ and $v(k+1,l)$ are not end-vertices of $H$. If $e(k+1,k+2;l) \in H$, then $L_N \mapsto L_N+(1,0)$ is a short edge weakening, so we may assume that $e(k+1,k+2;l) \notin H$. Then  $e(k+1;l,l+1) \in H$. Note that if $e(k+2;l-1,l)\in H$, then $L_N+(1,-1)=L_E$, which contradicts the assumption that $d(T) \geq 3$, so we must have that $e(k+2;l-1,l) \notin H$. It follows that $e(k+2;l,l+1) \in H$. Note that if $L_N+(1,1)$ is parallel, then Sw($L_N+(1,1)$), $L_N \mapsto L_N +(1,0)$ is a short edge weakening, so we may assume that $L_N+(1,1)$ is anti-parallel.

\noindent Let $H'$ be the Hamiltonian path obtained after $bb_{v}(east)$, and let $v'$ be the end-vertex of $H'$ to which $v$ is relocated by the move. Note that $v'=v(k+2,l)$, and that $e(k+1;l-1,l)$ precedes $e(k+1;l,l+1)$ in $H'$. Now we can apply $bb_{v'}(west)$ to $H'$, after which, $v'$ is relocated to $v(k+1,l-1) \in \text{ew-Terminal}(T)$, so  $bb_{v}(east)$ by $bb_{v'}(west)$ is a short end-vertex weakening. See Figure. 4.88. End of Case 1.2(c).

\begingroup
\setlength{\intextsep}{0pt}
\setlength{\columnsep}{20pt}
\begin{wrapfigure}[]{l}{0cm}
\begin{adjustbox}{trim=0cm 0cm 0cm 0cm} 

\end{adjustbox}
\end{wrapfigure}

\null 

\noindent \textit{CASE 1.2(d): $v=v(k,l+1)$.} By case 1.2(a), we may assume that $v(k+1,l)$ is not an end-vertex of $H$, and that $e(k+1,k+2;l) \notin H$. Then $e(k+1;l,l+1) \in H$. Now, either $e(k;l,l+1) \in H$ or $e(k;l,l+1) \notin H$. 

\null 

\noindent \textit{CASE 1.2($d_1$): $e(k;l,l+1) \in H$.} The assumption that $L_N$ is anti-parallel implies that $bb_{v}(east)$ relocates $v$ to $v(k+1,l)$, and we are back to Case 1.2(a). See Figure 4.89. End of Case 1.2($d_1$).

\endgroup

\begingroup
\setlength{\intextsep}{0pt}
\setlength{\columnsep}{20pt}
\begin{center}
\begin{adjustbox}{trim=0cm 0cm 0cm 0cm} 

\end{adjustbox}
\end{center}

\noindent \textit{CASE 1.2 ($d_2$): $e(k;l,l+1) \notin H$.} By NAA, $v(k,l)$ and $v(k+1,l+1)$ are not end-vertices of $H$. Let $v_s=v(k,l)$ and $v_t=v(k+1,l+1)$. Now, $v(k,l-1)=v_{s-1}$ (Figure 4.90 (a)), or $v(k,l-1)=v_{s+1}$ (Figure 4.90 (b)). If the former, then after $bb_{v}(south)$, $v$ is relocated to $v(k,l-1) \in \text{ew-Terminal}(T)$, so $bb_{v}(south)$ is a short end-vertex weakening. And if the latter, than the assumption that $L_N$ is anti-parallel implies that $v_{t-1}=v(k+1,l)$. Then, after $bb_{v}(east)$, $v$ is relocated to $v(k,l+1)$, and we are back to Case 1.2(a). End of Case 1.2($d_2$). End of Case 1.2(d). End of Case 1.2. End of Case 1.

\null

\noindent \textit{CASE 2: $v \in$  $s$-Flank$(T)$.} As we have seen before, if $e(k+1,k+2;l) \in H$, then $L_N \mapsto L_N+(1,0)$ is a short edge weakening, so we may assume that $e(k+1,k+2;l) \notin H$. Let:

$V_1= \{v(k+2,l), v(k+2,l+1), v(k+3,l), v(k',l'+2), v(k'+1,l'+2), v(k',l'+3) \}$,

$V_2= \{v(k+3,l-1), v(k+4,l-2), \ldots, v(k'-1,l'+3) \}$, and

$V_3= \{v(k+4,l-1), v(k+5,l-2), \ldots, v(k'-1,l'+4) \}$,

\noindent Note that $\{V_1, V_2, V_3\}$ is a partition of $s$-Flank$(T)$, that $V_3=\emptyset$ if $d(T)=4$, and that $V_2$ and $V_3$ are both empty if $d(T)=3$. Then there are three cases to check: $v \in V_1$, $v \in V_2$ and $v \in V_3$.

\null 

\noindent \textit{CASE 2.1: $v \in V_1$.} By symmetry we only need to check the cases where $v=v(k+2,l)$, $v=v(k+2,l+1)$ and $v=v(k+3,l)$.

\null 

\begingroup
\setlength{\intextsep}{0pt}
\setlength{\columnsep}{20pt}
\begin{wrapfigure}[]{l}{0cm}
\begin{adjustbox}{trim=0cm 0cm 0cm 0cm} 

\end{adjustbox}
\end{wrapfigure}

\noindent \textit{CASE 2.1 (a): $v=v(k+2,l)$.} If $e(k+2;l-1,l) \in H$ then $L_N+(1,-1)=L_E$, which contradicts the assumption that $d(T)\geq 3$, so we may assume that $e(k+2;l-1,l) \notin H$. Let $v(k+1,l)=v_s$. Then $v(k+1,l-1)=v_{s-1}$ or $v(k+1,l-1)=v_{s+1}$. 

\endgroup 

\null

\noindent \textit{CASE 2.1 ($a_1$): $v(k+1,l-1)=v_{s-1}$.} Then, after $bb_{v}(west)$, $v$ is relocated to $v(k+1,l-1) \in \text{ew-Terminal}(T)$, so $bb_{v}(west)$ is a short end-vertex weakening. See Figure 4.91. End of Case 2.1 ($a_1$).

\begingroup
\setlength{\intextsep}{0pt}
\setlength{\columnsep}{20pt}

\begin{center}
\begin{adjustbox}{trim=0cm 0cm 0cm 0cm} 

\end{adjustbox}
\end{center}

\noindent \textit{CASE 2.1 ($a_2$): $v(k+1,l-1)=v_{s+1}$.} If $L_N$ is open (Figure 4.92 (a)), then $bb_{v}(west)$ is the move $e(k+1;l,l+1) \mapsto e(k+1,k+2;l)$. Then we can apply $L_N \mapsto L_N+(1,0)$, and see that $bb_{v}(west)$, $L_N \mapsto L_N+(1,0)$ is a short edge weakening. And if $L_N$ is closed (Figure 4.92 (b)), then $bb_{v}(west)$ is the move $e(k,k+1;l) \mapsto e(k+1,k+2;l)$. Then we can apply $L_N \mapsto L_N+(1,0)$, and see that $bb_{v}(west)$, $L_N \mapsto L_N+(1,0)$ is a short weakening. Either way $bb_{v}(west)$, $L_N \mapsto L_N+(1,0)$ is a short edge weakening. End of Case 2.1 ($a_2$). End of Case 2.1(a).

\endgroup 

\null

\noindent \textit{CASE 2.1 (b): $v=v(k+2,l+1)$.} By NAA $v(k+2,l)$ and $v(k+1,l+1)$ are not end-vertices of $H$. Since $d(T)\geq 3$, $e(k+2;l-1,l) \notin H$. Then we must have that $e(k+2;l,l+1) \in H$.

\begingroup
\setlength{\intextsep}{0pt}
\setlength{\columnsep}{20pt}
\begin{center}
\begin{adjustbox}{trim=0cm 0cm 0cm 0cm} 

\end{adjustbox}
\end{center}

\noindent Note that if $L_N$ is closed (Figure 4.93 (a)), since $v(k+1,l+1)$ is not an end-vertex of $H$, then we must have that $e(k,k+1;l+1)\in H$. But then $L_N+(0,1) \mapsto L_N$ is a short edge weakening, so we may assume that $L_N$ is open. Then $e(k+1;l,l+1) \in H$. 

%pagemarker
\noindent Let $v(k+1,l+1)=v_s$. Then $v(k+1,l)=v_{s-1}$ (Figure 4.93 (b)) or $v(k+1,l)=v_{s+1}$ (Figure 4.93 (c)). If $v(k+1,l)=v_{s+1}$, then $L_N+(1,1)$ is parallel, and then Sw($L_N+(1,1)$), $L_N \mapsto L_N+(1,0)$ is a short edge weakening. Assume $v(k+1,l)=v_{s-1}$. Then after $bb_{v}(west)$, $v$ relocates to $v(k+1,l)$, and we are back to Case 1.2(a). End of Case 2.1(b).

\endgroup

\null

\begingroup
\setlength{\intextsep}{0pt}
\setlength{\columnsep}{20pt}
\begin{wrapfigure}[]{r}{0cm}
\begin{adjustbox}{trim=0cm 0cm 0cm 0.5cm} 

\end{adjustbox}
\end{wrapfigure}

\noindent \textit{CASE 2.1 (c): $v=v(k+3,l)$.} By NAA $v(k+2,l)$ and $v(k+3,l-1)$ are not end-vertices of $H$. Then we must have that $e(k+2,k+3;l) \in H$, and $e(k+3;l-2,l-1) \in H$. It follows that $L_E=L_N+(2,-2)$, and that $L_E$ is open. Then $L_E \mapsto L_E+(0,1)$ is a short edge weakening. See Figure 4.94. End of Case 2.1(c). End of Case 2.1.

\endgroup

\begingroup
\setlength{\intextsep}{0pt}
\setlength{\columnsep}{20pt}

\begin{center}
\begin{adjustbox}{trim=0cm 0cm 0cm 0cm} 
%
\end{adjustbox}
\end{center}

\noindent \textit{CASE 2.2: $v \in V_2$}. We may assume that $V_2 \neq \emptyset$. Note that this implies that $d(T) \geq 4$. Let $v=v(a,b)$. Now, after $bb_{v}(west)$, the end-vertex in the resulting Hamiltonian path is at $v(a-2,b)$ or on $v(a-1,b-1)$. See Figure 4.95 (a). Since $\{v(a-2,b), v(a-1,b-1)\} \subset \text{ew-Terminal}(T)$, $bb_{v}(west)$ is a short end-vertex weakening. End of Case 2.2.

\null

\noindent \textit{CASE 2.3: $v \in V_3$}. We may assume that $V_3\neq \emptyset$. Note that this implies that $d(T) \geq 5$. Let $v=v(a,b)$. By NAA $v(a-1,b)$ and $v(a,b-1)$ are not end-vertices of $H$. But then we must have that $e(a;b-1,b)\in H$ and $e(a-1,a;b)\in H$, which conflicts with our assumption that $v(a,b)$ is an end-vertex of $H$. See Figure 4.95 (b). End of Case 2.3. End of Case 2. End of proof for part (a).

\null

\noindent \textit{Proof of Lemma 4.30(b).} Let $T'$ be a lengthening of $T$ and assume that
$\mu_1, \ldots, \mu_s$ is an end-vertex weakening of $T'$. By Observation 4.20(b) ew-Terminal $(T')\subset \text{Flank}(T)$, so after $\mu_s$, the resulting Hamiltonian path has an end-vertex $u'$ in $\text{Flank}(T)$. Then, by part (a), $T$ has a weakening. The argument that $s+1\leq \min(m,n)$ is the same as the one used in Case 1.3(a) of Lemma 4.23, so we omit it. End of proof for part (b). $\square$

\null 

\noindent \textbf{Corollary 4.31.} Let $G$ be an $m \times n$ grid graph and let $H$ be a Hamiltonian path of $G$ with end-vertices $u$ and $v$ not adjacent. Let $T$ be a turn of $H$ with $d(T) \geq 3$. Then:

(i) \ If $v \notin \text{Sector}(T)$ then $T$ has a weakening contained in $\text{Sector}(T)$ that fixes $v$.

(ii) The last move in every edge weakening of $T$ is a non-backbite move involving one 

\hspace{0.6cm} of its leaves.

\null

\noindent \textit{Proof.} Both (i) and (ii) follow by inspecting all the weakenings used in Lemmas 4.30, 4.26, and 4.29, and Propositions 4.27 and 4.28. Note that part (i) makes implicit use of Observation 4.20 (d). $\square$

\null 

\noindent \textbf{Lemma 4.32.} Let $G$ be an $m \times n$ grid graph, and let $H$ be a Hamiltonian path or cycle of $G$. Let $e_F$ be an edge in $H$ followed by a $j$-stack of $A_0$'s, which is then followed by an $A_1$-type with switchable middle-box $W$. Let $X$ and $Y$ be the boxes adjacent to $W$ that are not its $H$-neighbours. Assume that $F = G\langle N[P(X,Y)] \rangle$ is a standard looping fat path. Let $T$ be an admissible turn of $F$, and let $\mu_1, \ldots, \mu_s$ be an edge-weakening of $T$ that is contained in Sector$(T)$, avoids the $j$-stack of $A_0$'s and the boxes incident on $e_F$, and is such that after the application of $\mu_{s-1}$, $F$ remains a standard looping fat path in the resulting Hamiltonian path or cycle. Then we can extend $\mu_1, \ldots, \mu_s$ by a non-backbite cascade of length at most two, after which $W$ is switched.

\null 

\noindent \textit{Proof.} For definiteness, assume that $T$ is northeastern with northern leaf $L_N=R(a,b)$ and (by Corollary 4.31 (ii)) that $\mu_s$ is the move $L_N \mapsto L_N'$, $L_N' \mapsto L_N$, or $\text{Sw}(L_N)$. By Corollary 4.19(b), $L_N \in N[P] \setminus P$ and $L_N+(0,-1) \in P$. Note that we must have $S_{\rightarrow}(a,b-1;a+1,b-2) \in H$, $e(a,a+1;b) \notin H$ and $e(a+1;b-1,b)\notin H$. Now, either $L_N+(0,-1)$ is an end-box of $P$ (Figure 4.96) or it is not (Figure 4.97).

\begingroup
\setlength{\intextsep}{10pt}
\setlength{\columnsep}{20pt}

\begin{adjustbox}{trim=0cm 0cm 0cm 0cm}

\end{adjustbox}

\noindent \textit{CASE 1: $L_N+(0,-1)$ is not an end-box of $P$.} Then $e(a;b-1,b) \notin H$. Then, after $\mu_s$, $L_N+(0,-1) \in P(X,Y)$ is switchable. Then, by Proposition 4.5 (a), there is a non-backbite cascade of length at most two, after which $W$ is switched. End of Case 1.

\null 

\noindent \textit{CASE 2: $L_N+(0,-1)$ is an end-box of $P$.} Assume that $L_N+(0,-1)=Y$. Then we have that $e(a;b-1,b)\in H$, that $F$ is northern, that $e_F=e(a-1,a;b-2)$, and that $W=Y+(-1,0)$. Then $\mu_s$ can be followed by $W \mapsto Y$. $\square$

\endgroup 

\null

\noindent Now we are ready to give a proof of Lemma 3.13.

\null 

\noindent \textbf{Lemma 3.13.} Let $G$ be an $m \times n$ grid graph, and let $H$ be a Hamiltonian cycle of $G$. Let $C$ be a small cookie of $G$.  Assume that $G$ has only one large cookie, and that there is a $j$-stack of $A_0$ starting at the $A_0$-type containing $C$. Let $L$ be the leaf in the top ($j^{\text{th}}$) $A_0$ of the stack, and assume that $L$ is followed by an $A_1$-type. Let $X$ and $Y$ be the boxes adjacent to the middle-box of the $A_1$-type that are not its $H$-neighbours. If $P(X,Y)$ has no switchable boxes, then either:

(i) there is a cascade of length at most $\min(m,n)$, which avoids the stack of $A_0$'s, and 

\hspace{0.4 cm} after  which $P(X,Y)$ gains a switchable box, or 

(ii) there is a cascade of length at most $\min(m,n)+1$, that collects $L$ and avoids the 

\hspace{0.4 cm} stack of $A_0$'s.

\null

\noindent \textit{Proof.} Suppose that $P(X,Y)$ has no switchable boxes. Note that $X$ and $Y$ belong to the same $H$-component. Since $H$ is a cycle, it follows that $P(X,Y)$ is contained in a \index{standard looping fat path}standard looping fat path $F=G\langle N[P(X,Y)] \rangle$. By Lemma 4.21, $F$ has an admissible turn $T$ such that $\text{Sector}(T)$ and the $j$-stack of $A_0$'s are disjoint. Then, by Corollary 4.19(a), $d(T) \geq 3$. By Proposition 4.24, $T$ has a weakening $\mu_1, \ldots, \mu_s$ of length at most $\min(m,n)$. By Observation 4.25, $\mu_1, \ldots, \mu_s$ is contained in $\text{Sector}(T)$, and thus it avoids the j-stack of $A_0$'s. By Lemma 4.32, and the fact that $H$ is a cycle, after $\mu_s$, either  $P(X,Y)$ gains a switchable box $Z$, or $W \mapsto W'$ is a valid move, where $W'=X$ or $W'=Y$. If the former, then (i) holds; and if the latter, apply $W \mapsto W'$, and then (ii) holds. $\square$.

\null 

%\noindent \textbf{Observation 4.33}. Let $G$ be an $m \times n$ grid graph, and let $H$ be a Hamiltonian cycle of $G$ with more than one cookie. Then it takes at most $\frac{1}{2}\max(m,n)+\min(m,n)+2$ moves to reduce the number of cookies of $G$ by one.

%\null 

%\noindent \textit{Proof.} By Proposition 3.8, if there are at least two large cookies, then we can reduce the number of large cookies by one with at most two moves. Assume then that there is only one large cookie. Let $C$ be an outermost small cookie, followed by a $j$-stack of $A_0$-types, which is followed by an $A_1$-type with looping $H$-path $P(X,Y)$\footnote{Recall that if the $j$-stack of $A_0$s is either not followed by an $A_1$-type or is a full stack, then these cases are handled by Lemmas 3.9 and 3.11 in Chapter 3. Moreover, the assumption that there is only one large cookie implies $X$ and $Y$ must belong to the same $H$-component.}.

%Let $L$ be the leaf contained in the top $A_0$-type of the stack, and assume that the $P(X,Y)$ is contained in a standard looping fat path $F$ that has an admissible turn $T$ with set of lengthenings $\mathcal{T}(T)$. The worst case scenario occurs when $j$ and $|\mathcal{T}(T)|$ are as large as possible, and the best we can do is to make $P(X,Y)$ gain a switchable box. Assume we are in this scenario. Note that $j$ can be no larger than $\frac{1}{2}\max(m,n)$, and that by Lemma 3.12, we need at most $\min(m,n)+3$ moves to collect $L$. Then can use $(j{-1})$-flips to collect $C$. Thus, we need at most $\frac{1}{2}\max(m,n)+\min(m,n)+2$ moves to collect $C$. $\square$ 

\noindent \textbf{Remark.} As per the proof of Lemma 3.13, the first valid move that $F$ will admit, either results in a switchable box for $P(X,Y)$, or makes available one of the moves $W \mapsto X$ and $W \mapsto Y$. Let $\min(m,n)=s$. It can be shown that there is no $j_0 \in \{1, \ldots, s\}$ such that after applying the cascade $\mu_1, \ldots, \mu_{j_0}$, $P(X,Y)$ gains a switchable box, or one of $W \mapsto X$ and $W \mapsto Y$ becomes available. However, it is easier to note that if such a $j_0$ were to exist, then we can just as well use this shorter cascade.

\null

\subsection{Summary}

In this chapter we examined the structure of fat paths and used it to prove Lemma 3.13, which guarantees the existence of cascades required for the 1LC algorithm.
Key definitions include fat paths, turns, weakenings, and the fat path conditions (FPC-1 through FPC-4). 

In Section 4.1 we proved some structural properties of looping fat paths. Let $F$ be such a path. Then $F$ has no colinear edges other than those in its $A_1$-type (Proposition 4.14), and its boundary $B(F)$ is a Hamiltonian cycle of $F$ (Lemma 4.13). Proposition 4.5 generalized the characterization of double-switch moves for Hamiltonian cycles of $m \times n$ grid graphs (Proposition 3.3) to Hamiltonian paths in polyominoes.

The two main structural results are that every fat path contains an admissible turn (Proposition 4.18, Section 4.2) and that every admissible turn has a weakening (Propositions 4.24, 4.27, and 4.28 in Sections 4.4-4.6). We identified two types of weakenings: edge weakenings and end-vertex weakenings. Both are useful, with end-vertex weakenings required specifically when $H$ is a Hamiltonian path and certain end-vertex configurations arise Section 4.6.

\null 
\null 
\null

\null 
\null

\null 
\null 
\null 
\null

\newpage

\section{Reconfiguration algorithm  for Hamiltonian paths}

In this chapter, we prove that any two Hamiltonian paths on an $m \times n$ grid graph can be reconfigured into one another (Theorem 5.9). The proof relies on Theorem 2.1, which handles reconfiguration of Hamiltonian cycles. Our strategy is to reconfigure an arbitrary path into one whose end-vertices lie in corners of $R_0$. Once both end-vertices are positioned in corners, Proposition 5.4 shows that we can effectively treat the resulting Hamiltonian path as a Hamiltonian cycle, and so Theorem 2.1 applies.

We move the end-vertices to corners in two stages. First, the End-vertex-to-Boundary (EtB) algorithm moves an end-vertex to a side of $R_0$ (Proposition 5.6). Second, the End-vertex-to-Corner (EtC) algorithm moves the end-vertex from a side to a corner (Proposition 5.7). Both algorithms draw on Sections 1.3--1.5 of Chapter 1 and the fat path-turn-weakening pipeline from Chapter 4.

\null 

\noindent \textbf{Definitions.} Let $G$ be an $m \times n$ grid graph. Throughout this paper, we position grid graphs in the first quadrant of the plane, so that $v(0,0)$ is the bottom-left corner and $v(m{-}1,n{-}1)$ is the top-right corner. We call the vertices $v(0,0)$, $v(0,n{-}1)$, $v(m{-}1,n{-}1)$, and $v(m{-}1,0)$ the \index{corner of R0@corner of $R_0$|textbf}\textit{corners} of $R_0$. A vertex $v(a,b)$ is called \index{even vertex|textbf}\textit{even} if $a + b$ is even, and \index{odd vertex|textbf}\textit{odd} if $a + b$ is odd. We write $p(v)$ to denote the parity of a vertex $v$.

Let $H$ be a Hamiltonian path of $G$ with end-vertices $u$ and $v$, with parities $p(u)$ and $p(v)$, respectively. We write $e(v)$ to denote the edge of $H$ incident on $v$, and $e(u)$ to denote the edge of $H$ incident on $u$. 

Let $s$ be a side of $G$ with corners $a$ and $b$. If $p(u) \in \{p(a), p(b)\}$ and $p(v) \in \{p(a), p(b)\}$, we call that side a \index{parity compatible side of R0@parity compatible side of $R_0$|textbf}\textit{parity-compatible} side of $G$.

%If $u$ and $v$ are any two vertices of $G$, we say that the \textit{Hamiltonian path problem $(G, u, v)$} is \textit{solvable} if there exists a Hamiltonian path from $u$ to $v$ in $G$.

\null

\noindent The existence of Hamiltonian paths on grid graphs was studied by  Itai et al. in \cite{itai1982hamilton}. The following is a consequence of the main result (Theorem 3.2) in that paper. 

\null

\noindent \textbf{Corollary 5.1.} Let $H$ be a Hamiltonian path in an $m \times n$ grid graph $G$ with end-vertices $u$ and $v$. Then: 

(i) \ \ If $m$ and $n$ are both odd, then both $u$ and $v$ have even parity.

(ii) \ If at least one of $m$ and $n$ is even, then $u$ and $v$ have different parities.  $\square$

\null

\noindent We record below Observation 5.2, which follows immediately from Corollary 5.1.

\null

\noindent \textbf{Observation 5.2.} Let $H$ be a Hamiltonian path in an $m \times n$ grid graph $G$, such that $m,n \geq 3$. Then $G$ has at least one pair of opposite sides that are parity-compatible. $\square$

\null

\noindent  \textbf{Definitions.} Let $G^+_{\mathrm{north}}$ be the $m \times (n{+}1)$ grid graph with vertex set $\{v(i,j) : 0 \le i \le m{-}1,\, 0 \le j \le n\}$ and edge set consisting of all edges between pairs of vertices at Euclidean distance $1$. We refer to $G^+_{\mathrm{north}}$ as the \index{northern extension of an m x n grid graph@northern extension of an $m \times n$ grid graph|textbf}\textit{northern extension} of $G$. Define eastern extensions of $G$ analogously.

\begingroup
\setlength{\intextsep}{0pt}
\setlength{\columnsep}{20pt}
\begin{center}
\begin{adjustbox}{trim=0cm 0cm 0cm 0cm} 
\begin{tikzpicture}[scale=1.5]

\begin{scope}[xshift=0cm]
  % grid unchanged
\draw[gray,very thin, step=0.5cm, opacity=0.5] (0,0) grid (2,1.5);

\draw[blue, line width=0.5mm] (0,1)--++(0,-1)--++(0.5,0)--++(0,1)--++(0.5,0)--++(0,-1)--++(1,0)--++(0,0.5)--++(-0.5,0)--++(0,0.5)--++(0.5,0); 

\draw[orange, line width=0.5mm] (0,1)--++(0,0.5)--++(2,0)--++(0,-0.5);

\draw[fill=blue, opacity=1] (0,1) circle [radius=0.05];
\draw[fill=blue, opacity=1] (2,1) circle [radius=0.05];

\node[above] at (1,1.5) [scale=1.25]{(a)};

\end{scope}

\begin{scope}[xshift=2.75cm]
  % grid unchanged
\draw[gray,very thin, step=0.5cm, opacity=0.5] (0,0) grid (2.5,1.5);

\draw[blue, line width=0.5mm] (0,1)--++(0,-1)--++(0.5,0)--++(0,1)--++(0.5,0)--++(0,-1)--++(0.5,0)--++(0,1)--++(0.5,0)--++(0,-1); 

\draw[orange, line width=0.5mm] (0,1)--++(0,0.5)--++(2.5,0)--++(0,-1.5)--++(-0.5,0);

\draw[fill=blue, opacity=1] (0,1) circle [radius=0.05];
\draw[fill=blue, opacity=1] (2,0) circle [radius=0.05];

\node[above] at (1.25,1.5) [scale=1.25]{(b)};

\node[right, align=left, text width=7cm] at (3, 0.75) {Fig. 5.1. (a) A northern path with a northern handle. (b) A northwest-southeast path with a northeastern handle.};

\end{scope}

\end{tikzpicture}
\end{adjustbox}
\end{center}

\noindent A \index{northern Hamiltonian path|textbf}\textit{northern Hamiltonian path} in $G$ is a Hamiltonian path $H_{\mathrm{north}}$ with end-vertices $v(0,n{-}1)$ and $v(m{-}1,n{-}1)$, the top-left and top-right corners of $G$ (Figure 5.1 (a)). Define eastern, southern, and western Hamiltonian paths of $G$ analogously. A \index{northwest-southeast Hamiltonian path|textbf}\textit{northwest-southeast Hamiltonian path} in $G$ is a Hamiltonian path $H_{\mathrm{NW-SE}}$ with end-vertices $v(0,n{-}1)$ and $v(m{-}1,0)$, the top-left and bottom-right corners of $G$ (Figure 5.1 (b)). Define a {southwest-northeast Hamiltonian path} analogously.

We call the set of edges $\{e(0;n{-}1,n),\, e(0,1;n),\, \ldots,\, e(m{-}2,m{-}1;n), e(m{-}1;n{-}1,n)\}$ the \index{northern handle|textbf}\textit{northern handle} of $H_{\mathrm{north}}$ and denote it by $\mathrm{Handle}_{\mathrm{north}}(H)$. Define the eastern, southern, and western handle of a Hamiltonian path analogously. We call the set of edges $\{e(0;n{-}1,n), e(0,1;n), \ldots,\, e(m{-}2,m{-}1;n),$ $ e(m-1,m;n), e(m;n-1,n), \ldots,$ $e(m;0,1),$ $ e(m-1,m;0)\}$ the \index{northeastern handle|textbf}\textit{northeastern handle} of $H_{\mathrm{NW-SE}}$ and denote it by $\mathrm{Handle}_{\mathrm{NE}}(H)$. Define the \textit{southwestern handle} of $H_{\mathrm{SW-NE}}$ analogously

We call the Hamiltonian cycle of $G^+_{\mathrm{north}}$ obtained by adding $\mathrm{Handle}_{\mathrm{north}}(H)$ to a northern Hamiltonian path $H$ the \index{handle-path cycle of a Hamiltonian path|textbf}\textit{handle-path cycle} of $H$, and we denote it by $\overline{H}$. Note that $\overline{H}$ lies in $G^+_{\mathrm{north}}$ and is defined only when $H$ is a northern Hamiltonian path. Analogous statements hold for eastern, southern, and western handle-path cycles.

\null

\noindent \textbf{Lemma 5.3.} Let $G$ be an $m \times n$ grid graph with $m,n \geq 3$, let $H$ be a northern Hamiltonian path in $G$ and let $X$ and $Y$ be boxes of $G$. Assume that $X \mapsto Y$ is a valid move for $\overline{H}$, after which we obtain $\overline{H}'$. Then $X \mapsto Y$ is valid for $H$, yielding a northern Hamiltonian path $H'$, where $H'=\overline{H}' \setminus \mathrm{Handle}_{\mathrm{north}}(H)$. An analogous statement holds for eastern, southern, and western Hamiltonian paths of $G$.

\null 

\noindent \textit{Proof.} Let $H=v_1, \ldots,v_r$ be a northern Hamiltonian path of $G$, and let $v_{r+1} \ldots, v_1$ be its northern handle. Let $\overline{H}=v_1, .., v_r, v_{r+1} \ldots, v_1$, and let $\overline{H}'$ be the Hamiltonian cycle obtained after applying $X \mapsto Y$ to $\overline{H}$. Let $(v_s, v_{s+1})$ and $(v_t, v_{t+1})$ be the edges of $X$ in $\overline{H}$. By Corollary 1.5.2, $X$ is anti-parallel. Then $v_s$ is adjacent to $v_{t+1}$ and $v_{s+1}$ is adjacent to $v_t$. After switching $X$ we obtain two cycles $H_1=$ $P(v_{s+1},v_t), \{v_{s+1},v_t\}$ and $H_2=  P(v_1,v_s), \{v_s,v_{t+1}\},P(v_{t+1},v_r),$ $P(v_r,v_1)$. Since the $(H_1, H_2)$-port $Y$ is contained in $G$, neither of its edges may belong to $P(v_r,v_1)$. Then $Y$ must be a port between the cycle $H_1$ and the sub-path $P(v_1,v_s), \{v_s,v_{t+1}\},P(v_{t+1},v_r)$ of $H_2$. Since $X \mapsto Y$ does not interact with $\mathrm{Handle}_{\mathrm{north}}(H)$, it follows that $X \mapsto Y$ is a valid move for $H$ and that $H'=\overline{H}' \setminus \mathrm{Handle}_{\mathrm{north}}(H)$, where $H'$ is the northern Hamiltonian path obtained after applying $X \mapsto Y$ to $H$. $\square$

\null

\noindent \textbf{Proposition 5.4.} Let $H$ and $K$ be northern Hamiltonian paths in an $m \times n$ grid graph $G$ with $n \geq m \geq 3$. Then there is a sequence of at most $n^2m$ double-switch moves that reconfigures $H$ into $K$. An analogous statement holds for eastern, southern, and western Hamiltonian paths of $G$.

\null 

\noindent \textit{Proof.} By Theorem 2.1, there is a sequence $\mu_1, \ldots, \mu_s$ of double-switch moves that reconfigures $\overline{H}$ into $\overline{K}$. Recall that no new cookies are created by the RtCF algorithm, and that by definition, cascades create no new cookies. Because there are no cookie necks on $\text{Boxes}(G_{\text{north}}^+) \setminus \text{Boxes}(G)$, and because even when the neck of a large cookie is relocated, it remains on the same side of the boundary,\footnote{These moves appear in the proof of Lemma 4.23 in Case 1.2(b), the first paragraph of Case 2, and Case 2.1$(b_1)$.} the entire sequence $\mu_1, \ldots, \mu_s$ is contained within $G$. Then, by Lemma 5.3 and induction, after applying $\mu_1, \ldots, \mu_s$ to $H$, we obtain $K$. The proof for the case where the Hamiltonian paths are eastern, southern, or western is similar, so we omit it. $\square$

\null

\noindent \textbf{Corollary 5.4.1.} Let $H$ and $K$ be either both northwest–southeast Hamiltonian paths or both southwest–northeast Hamiltonian paths of an $m \times n$ grid graph $G$, with $n \geq m \geq 3$. Then there is a sequence of at most $n^2m$ double-switch moves that reconfigures $H$ into $K$.

\null 

\noindent \textit{Proof.} For definiteness, assume that $H$ is northwest-southeast. We can use a northeastern handle to obtain a handle-path cycle of $H$. Then the result follows by using essentially the same arguments as in Lemma 5.3 and Proposition 5.4. $\square$

\null

\noindent \textbf{Proposition 5.5.} Let $G$ be an $m \times 2$ grid graph, and let $H$ and $K$ be any two Hamiltonian paths of $G$. Then there is a sequence of at most $m$ backbite moves that reconfigures $H$ into $K$.

\null

\noindent \textit{Proof.} Let $H$ be a Hamiltonian path of the $m \times n$ grid graph $G$ with end-vertices $u$ and $v$. By Corollary 5.1, since $n=2$, $u$ and $v$ have distinct parities. This means that $G$ does have Hamiltonian e-cycles. We define the canonical forms of $G$ to be the set of e-cycles of $G$. We show below (Lemma 5.10) that any two e-cycles can be reconfigured into one another by at most $m$ backbite moves.  It remains to show that a Hamiltonian path of $G$ can be reconfigured into an e-cycle by a sequence of at most $m$ backbite moves.

If $H$ is already an e-cycle we are done, so assume that $u$ and $v$ are not adjacent. For definiteness, assume that $u$ is east of $v$, and that $v=v(k,1)$ lies on the northern side of $R_0$. It is enough to check that there is a backbite move that relocates $v$ to a vertex with x-coordinate one unit greater than $v$. Let $e(v)$ be the edge of $H$ incident on $v$. Either $k=0$, or $k>0$.

\begingroup
\setlength{\intextsep}{0pt}
\setlength{\columnsep}{20pt}
\begin{wrapfigure}[]{l}{0cm}
\begin{adjustbox}{trim=0cm 0cm 0cm 0cm} 
\begin{tikzpicture}[scale=1.75]

\begin{scope}[xshift=0cm] 
{
\draw[gray,very thin, step=0.5cm, opacity=0.5] (0,0) grid (1,0.5);

\draw [->,black, very thick] (1.2,0.25)--(1.8,0.25);
\node[above] at  (1.5,0.25) [scale=0.8]{\small{$bb_v(east)$}};

\draw [blue, line width=0.5mm] (0,0.5)--++(0,-0.5)--++(0.5,0)--++(0,0.5)--++(0.5,0);

\draw[fill=blue] (0,0.5) circle (0.05);
\node[right] at (0,0.5) [scale=0.8]{\small{$v$}};

% labellings 
{

\node[left] at (0.05,-0.1) [scale=1]{\tiny{$0$}};

\node[left] at (0,0.5) [scale=1]
{\tiny{$1$}};
\node[below] at (0.5,0) [scale=1]
{\tiny{$1$}};
\node[below] at (1,0) [scale=1]
{\tiny{$2$}};

}

}

\node[below, align=center, text width=5cm] at (1.5, -0.25) { Fig. 5.2. Case 1.};

\end{scope}

\begin{scope}[xshift=2.1cm] 
{
\draw[gray,very thin, step=0.5cm, opacity=0.5] (0,0) grid (1,0.5);

\draw [blue, line width=0.5mm] (0.5,0)--++(-0.5,0)--++(0,0.5)--++(1,0);

\draw[fill=blue] (0.5,0) circle (0.05);

% labellings 
{

\node[left] at (0.05,-0.1) [scale=1]{\tiny{$0$}};

\node[left] at (0.025,0.5) [scale=1]
{\tiny{$1$}};
\node[below] at (0.5,0) [scale=1]
{\tiny{$1$}};
\node[below] at (1,0) [scale=1]
{\tiny{$2$}};

}

}

\end{scope}

\end{tikzpicture}
\end{adjustbox}
\end{wrapfigure}

\noindent \textit{CASE 1: $k=0$.} Note that if $e(v)$ is eastern, then $u=v(0,0)$, and $H$ is an e-cycle, so we may assume that  $e(v)$ is not eastern. Then $e(v)$ must be southern. Then $e(0,1;0) \in H$ and $S_{\uparrow}(1,0;2,1) \in H$. Then $bb_v(east)$ relocates $v$ east, to $v(1,0)$. See Figure 5.2. End of Case 1.

\endgroup 

\null

\noindent \textit{CASE 2: $k>0$.} Then $e(v)$ is eastern, southern, or western. We will see that the first two possibilities lead to invalid configurations.

\null

\noindent \textit{CASE 2.1: $e(v)$ is eastern.} Then we must have that $e(k-1,k;0) \in H$ and $e(k,k+1;0) \in H$. Let $v_s=v(k,0)$.  Either $e(k+1;0,1) \in H$, or $e(k+1;0,1) \notin H$.

\null 

\begingroup
\setlength{\intextsep}{0pt}
\setlength{\columnsep}{10pt}
\begin{wrapfigure}[]{r}{0cm}
\begin{adjustbox}{trim=0cm 0cm 0cm 0cm} 
\begin{tikzpicture}[scale=1.75]

\begin{scope}[xshift=0cm] 
{
\draw[gray,very thin, step=0.5cm, opacity=0.5] (0,0) grid (1,0.5);

\begin{scope}[very thick,decoration={
    markings,
    mark=at position 0.6 with {\arrow{>}}}
    ] 
    \draw[postaction={decorate}, blue, line width=0.5mm] (0.5,0.5)--++(0.5,0);
    \draw[postaction={decorate}, blue, line width=0.5mm] (1,0.5)--++(0,-0.5);
    \draw[postaction={decorate}, blue, line width=0.5mm] (1,0)--++(-0.5,0);
    \draw[postaction={decorate}, blue, line width=0.5mm] (0.5,0)--++(-0.5,0);

\end{scope}

\draw[fill=blue] (0.5,0.5) circle (0.05);
\node[left] at (0.5,0.5) [scale=0.8]{\small{$v$}};

\draw[fill=blue] (0.5,0) circle (0.035);
\node[above] at (0.5,0) [scale=0.8]{\small{$v_s$}};

\draw[fill=blue] (0,0) circle (0.035);
\node[above] at (0,0) [scale=0.8]{\small{$v_{s+1}$}};

% labellings 
{

\node[left] at (0.05,-0.1) [scale=1]{\tiny{$0$}};

\node[left] at (0,0.5) [scale=1]
{\tiny{$1$}};
\node[below] at (0.5,0) [scale=1]
{\tiny{$1$}};
\node[below] at (1,0) [scale=1]
{\tiny{$2$}};

}

}

\node[below, align=center, text width=4cm] at (0.5, -0.25) { Fig. 5.3. Case 2.1(a).};

\end{scope}

\end{tikzpicture}
\end{adjustbox}
\end{wrapfigure}

\noindent \textit{CASE 2.1(a): $e(k+1;0,1) \in H$.} Then $v_{s+1}=v(k-1,0)$. Since the x-coordinate of $u$ is greater than $k$, and the x-coordinate of $v_{s+1}$ is $k-1$, the subpath $P(v_{s+1}, u)$ must contain a vertex of $H$ with x-coordinate $k$. But the only two vertices of $G$ with x-coordinate $k$ are $v_1$ and $v_s$, and neither belongs to $P(v_{s+1}, u)$. See Figure 5.3. End of Case 2.1(a).

\endgroup 

\null 

\noindent \textit{CASE 2.1(b): $e(k+1;0,1) \notin H$.} Then, by Corollary 1.5.2, $R(k,0)$ is anti-parallel. Then $v_{s+1}=v(k-1,0)$, and now, we can use the same argument as in Case 2.1(a) to check that this configuration is invalid. End of Case 2.1(b). End of Case 2.1.

\null

\noindent \textit{CASE 2.2: $e(v)$ is southern.} The arguments for this case are similar to the one used in Case 2.1(a), so we omit the proof. End of Case 2.2.

\null 

\begingroup
\setlength{\intextsep}{0pt}
\setlength{\columnsep}{20pt}
\begin{wrapfigure}[]{l}{0cm}
\begin{adjustbox}{trim=0cm 0cm 0cm 0cm} 
\begin{tikzpicture}[scale=1.75]

\begin{scope}[xshift=0cm] 
{
\draw[gray,very thin, step=0.5cm, opacity=0.5] (0,0) grid (1.5,0.5);

\draw [->,black, very thick] (1.7,0.25)--(2.3,0.25);
\node[above] at  (2,0.25) [scale=0.8]{\small{$bb_v(east)$}};

\begin{scope}[very thick,decoration={
    markings,
    mark=at position 0.6 with {\arrow{>}}}
    ] 
    \draw[postaction={decorate}, blue, line width=0.5mm] (0.5,0.5)--++(-0.5,0);
    \draw[postaction={decorate}, blue, line width=0.5mm] (0,0)--++(0.5,0);
    \draw[postaction={decorate}, blue, line width=0.5mm] (0.5,0)--++(0.5,0);
    \draw[postaction={decorate}, blue, line width=0.5mm] (1,0)--++(0,0.5);
    \draw[postaction={decorate}, blue, line width=0.5mm] (1,0.5)--++(0.5,0);

\end{scope}

\draw[fill=blue] (0.5,0.5) circle (0.05);
\node[right] at (0.5,0.5) [scale=0.8]{\small{$v$}};

\draw[fill=blue] (0.5,0) circle (0.035);
\node[above] at (0.5,0) [scale=0.8]{\small{$v_s$}};

\draw[fill=blue] (1,0) circle (0.035);
\node[right] at (1,0) [scale=0.8]{\small{$v_{s+1}$}};

% labellings 
{

\node[left] at (0.05,-0.1) [scale=1]{\tiny{$0$}};

\node[left] at (0,0.5) [scale=1]
{\tiny{$1$}};
\node[below] at (0.5,0) [scale=1]
{\tiny{$1$}};
\node[below] at (1,0) [scale=1]
{\tiny{$2$}};

}

}

\node[below, align=center, text width=4cm] at (2, -0.25) { Fig. 5.4. Case 2.3.};

\end{scope}

\begin{scope}[xshift=2.5cm] 
{
\draw[gray,very thin, step=0.5cm, opacity=0.5] (0,0) grid (1.5,0.5);

\draw[ blue, line width=0.5mm] (0,0.5)--++(1.5,0);
\draw[ blue, line width=0.5mm] (0,0)--++(1,0);

\draw[fill=blue] (1,0) circle (0.05);

% labellings 
{

\node[left] at (0.05,-0.1) [scale=1]{\tiny{$0$}};

\node[left] at (0.025,0.5) [scale=1]
{\tiny{$1$}};
\node[below] at (0.5,0) [scale=1]
{\tiny{$1$}};
\node[below] at (1,0) [scale=1]
{\tiny{$2$}};

}

}

\end{scope}

\end{tikzpicture}
\end{adjustbox}
\end{wrapfigure}

\noindent \textit{CASE 2.3: $e(v)$ is western.} Then we must have $S_{\uparrow}(k+1,0;k+2,2) \in H$, $e(k-1,k;0) \in H$, and $e(k,k+1;0) \in H$. Let $v_s=v(k,0)$. Note that, using an argument similar to the one in Case 2.1, we can show that, whether $e(k-1;0,1) \in H$ or $e(k-1;0,1) \notin H$, $v_{s+1}=v(k+1,0)$. Then, $bb_v(east)$ relocates $v$ east, to $v(k+1,0)$. See Figure 5.4. End of Case 2.3. End of Case 2. $\square$

\endgroup

% A MORE FOCUSED VERSION OF PROPOSITIONS 5.6 AND 5.7.
{
%\noindent \textbf{Proposition 5.6.} Let $H$ be a Hamiltonian path from $u$ to $v$ in an $m \times n$ grid graph $G$ with $m,n \geq 3$. Suppose $s$ is a side of $R_0$ with corners $c_u$ and $c_v$ such that the parity of $c_u$ matches that of $u$, and the parity of $c_v$ matches that of $v$. Let $s_u\neq s$ be the side of $R_0$ incident on $c_u$ and let $s_v\neq s$ be the side of $R_0$ incident on $c_v$. Then:

%(a) There is a cascade that reconfigures $H$ into a Hamiltonian path with end-vertices $u'$ and $v'$ matching the parities of $u$ and $v$, respectively, such that $u'$ lies on $s_u$ and $v'$ lies on $s_v$. We refer to this procedure as the \textit{End-vertex-to-Boundary (EtB) algorithm}.

%(b) If $u$ lies on $c_u$ and $v$ lies on $c_v$, then there is a cascade that reconfigures $H$ into a Hamiltonian path whose end-vertices lie at $c_u$ and $c_v$. We refer to this procedure as the \textit{End-vertex-to-Corner (EtC) algorithm}.
}

\null

\noindent \textbf{Proposition 5.6.} \emph{(End-vertex-to-Boundary (EtB) Algorithm)} Let $H$ be a Hamiltonian path of an $m \times n$ grid graph $G$, with end-vertices $u$ and $v$, and let $s$ and $s'$ be a pair of opposite parity-compatible sides of $R_0$\footnote{By Observation 5.2, there is such a pair.}. Then:

(i) \ If neither of $u$ and $v$ is on $s$ or $s'$, then there is a cascade after which one of $u$ and 

\hspace{0.6cm} $v$ is on $s$ or $s'$.

(ii) If $u$ is on $s$ and $v$ is not on $s'$, then there is a cascade that fixes $u$, after which $v$ 

\hspace{0.6cm} is on $s'$.

\noindent Either cascade has length at most
$n^2+mn+4n$.

\null

\noindent \textbf{Proposition 5.7.} \emph{(End-vertex-to-Corner (EtC) Algorithm)} Let $G$ be an $m \times n$ grid graph, and let $H$ be a Hamiltonian path of $G$ with end-vertices $u$ and $v$. Let $s_u$ and $s_v$ be a pair of opposite sides of $R_0$, and assume that $u$ is on $s_u$, and $v$ is on $s_v$. Assume further that $c_v$ is a corner of $s_v$ that has the same parity as $v$, and that $c_u$ is a corner of $s_u$ that has the same parity as $u$. Then:

(i) \ If $u$ is not at $c_u$ and $v$ is not at $c_v$, then there is a cascade after which $v$ is at $c_v$ 

\hspace{0.6cm} and $u$ is on $s_u$.

(ii) If $u$ is at $c_u$ and $v$ is not at $c_v$, then there is a cascade after which $v$ is at $c_v$ and $u$ 

\hspace{0.6cm} is at $c_u$.

\noindent Either cascade has length at most
$\frac{n^2}{2}+mn+\frac{5n}{2}$.

\null

\noindent The proofs of Propositions 5.6 and 5.7 are given in Sections 5.1 and 5.2, respectively.

\null 

\noindent \textbf{Proposition 5.8.} Let $H$ be a Hamiltonian path in an $m \times n$ grid graph $G$ with $n \geq m \geq 3$. Then there is a cascade of length at most
$\frac{7n^2}{2}+5mn+\frac{31n}{2}$ that reconfigures $H$ into a northern or eastern Hamiltonian path. 

\null

\noindent \textit{Proof.} Since $m$ and $n$
must be both odd, both even, or differ in parity, we consider three cases.

\null 

\noindent \textit{CASE 1: $m$ and $n$ are odd.} Then all four corners of $R_0$ have even parity, and by Corollary 5.1(i), so do $u$ and $v$. Using EtB, we send $u$ and $v$ to the western and eastern sides of $R_0$, respectively. Using EtC we send $u$ and $v$ to the northwestern and northeastern corners of $R_0$, respectively. Then the resulting Hamiltonian path is northern. End of Case 1.

\null 

\noindent \textit{CASE 2: $m$ and $n$ are even.} By Corollary 5.1(ii), $u$ and $v$ have distinct parities. WLOG assume that $u$ is odd and $v$ is even. Note that the corners $v(0,0)$ and $v(m-1,n-1)$ are even, and the corners $v(0,n-1)$ and $v(m-1,0)$ are odd. This means that each side of $R_0$ is incident on an even corner and an odd one. Using  EtB, we send $u$ and $v$ to opposite sides of $R_0$. WLOG assume that we send $u$ to the western side and $v$ to the eastern side. Using EtC, we send $u$ to $v(0,n-1)$ and $v$ to $v(m-1,n-1)$.  Then the resulting Hamiltonian path is northern. End of Case 2.

\null 

\noindent \textit{CASE 3: $m$ and $n$ have different parities.} By Corollary 5.1(ii), $u$ and $v$ have distinct parities. WLOG assume that $u$ is even and $v$ is odd, and that $m$ is even and $n$ is odd. Then $v(0,0)$ and $v(0,n-1)$ are even, and $v(m-1,0)$ and $v(m-1,n-1)$ are odd. Using  EtB, we send one of $u$ and $v$ to the northern side of $R_0$, and the other to the southern side. Either $u$ has been sent to the northern side or to the southern side.

\null 

\noindent \textit{CASE 3.1: $u$ has been sent to the northern side.} Using EtC, we send $u$ to $v(0,n-1)$ and $v$ to $v(m-1,0)$. Using EtC again, we send $v$ to $v(m-1,n-1)$. Then the resulting Hamiltonian path is northern. End of Case 3.1.

\null 

\noindent \textit{CASE 3.2: $u$ has been sent to the southern side.} Using EtC, we send $u$ to $v(0,0)$ and $v$ to $v(m-1,n-1)$. Using EtC again, we send $u$ to $v(0,n-1)$. Then the resulting Hamiltonian path is northern. End of Case 3.2. End of Case 3.

\null

\noindent Note that it takes at most two applications of EtB and three applications of EtC to reconfigure a Hamiltonian path into an eastern or a northern one. So we need at most $2(n^2+mn+4n) +3 (\frac{n^2}{2}+mn+\frac{5n}{2})=\frac{7n^2}{2}+5mn+\frac{31n}{2}$ moves to complete the reconfiguration. $\square$

\null

\noindent \textbf{Theorem 5.9.} Let $G$ be an $m \times n$ grid graph with $n\geq m$, and let $H$ and $K$ be any two Hamiltonian paths of $G$. Then there is a cascade with at most $n^2m+7n^2+10mn+31n$ moves that reconfigures $H$ into $K$.

\null

\noindent \textit{Proof.} By Proposition 5.5, we may assume that $m,n \geq 3$. By Proposition 5.8, we may assume WLOG that $H$ and $K$ can be reconfigured into the northern Hamiltonian paths $H'$ and $K'$ by the cascades $\mu_1, \ldots, \mu_s$ and $\nu_1', \ldots, \nu_t'$. By Proposition 5.4,  $H'$ can be reconfigured into $K'$ by a cascade $\eta_1, \ldots, \eta_q$. For $j \in \{1, \ldots,t \}$ let $\nu_j$ be the inverse of the move $\nu_j'$. Then $\mu_1, \ldots, \mu_s$, $\eta_1, \ldots, \eta_q$, $\nu_t, .., \nu_1$ is a cascade that reconfigures $H$ into $K$. 

The bound is the sum of the bound of Theorem 2.1 with twice that of Proposition 5.8. $\square$

\null 

\noindent The rest of the chapter is structured as follows. Section 5.1  handles certain configurations that appear frequently in the proofs of both Propositions 5.6 and 5.7. We isolate these scenarios as separate lemmas, which we then invoke as needed in the main proofs. The proof of Proposition 5.6 is given in Section 5.2, and the proof of Proposition 5.7 is given in Section 5.3. Both proofs draw on Sections 1.3--1.5 of Chapter 1, the fat path-turn-weakening pipeline from Chapter 4, and the lemmas from Section 5.1.

\null

\subsection{Recurring Scenarios}

\noindent \textbf{Lemma 5.10.} Let $G$ be an $m \times n$ grid graph, and let $H$ be Hamiltonian path with adjacent end-vertices. Then there is a cascade of backbite moves that reconfigures $H$ into an e-cycle.

\begingroup
\setlength{\intextsep}{0pt}
\setlength{\columnsep}{20pt}
\begin{center}
\begin{adjustbox}{trim=-0.25cm 0cm 0cm 0cm}

\end{adjustbox}
\end{center}

\noindent \textit{Proof.}
Orient $H=H_1$ as $v_1, \ldots,v_r$.  For each $j \in  \{1, \ldots, r-1 \}$, let $H_j = v_j,v_{j+1},$ $ \ldots ,v_r,v_1, \ldots,v_{j-1}$, which is a Hamiltonian path. For each $j$, let $(v_j,v_{j+1}) \mapsto \{ v_{j-1},v_j\}$ be the backbite move $\beta_j$. Then that $\beta_j$ reconfigures $H_j$ into $H_{i+1}$. Let $j_0$ be the smallest index in $\{1, \ldots, r\}$ such that the edge $\{v_{j_0}, v_{j_0+1}\}$ of $G \setminus H_{j_0+1}$  is in $R_0$. Note that $H_{j_0+1}$ is an e-cycle. So, using $j_0$ backbite moves, we have reconfigured $H_1$ into the e-cycle $H_{j_0+1}$. See Figure 5.5.

\endgroup 

\null 

\noindent \textbf{Lemma 5.11.} Let $H$ be a Hamiltonian path of an $m \times n$ grid graph $G$. Let $e(k,k+1;l)$ be an edge of $H$ followed southward by a $j$-stack of $A_0$'s, which is not followed by an $A_1$-type. Let $R(k,b)$ be the leaf in the $j^{\text{th}}$ $A_0$ of the stack, and assume that the rectangle $R(k-1,k+2;b-2,b-1)$ does not contain an end-vertex of $H$. Then there is a cascade of flips, after which $e(k,k+1;l-1)$ is in the resulting Hamiltonian path of $G$. Analogous results apply for the cases where the $j$-stack follows an edge in another direction.

\begingroup
\setlength{\intextsep}{0pt}
\setlength{\columnsep}{20pt}
\begin{wrapfigure}[]{r}{0cm}
\begin{adjustbox}{trim=0cm 0.75cm 0cm 0.25cm}


\end{adjustbox}
\end{wrapfigure}

\null

\noindent \textit{Proof.} Let $R(k,b)$ be the southern leaf of the last ($j^{\text{th}}$) $A_0$-type in the $j$-stack of $A_0$'s that follows $e(k,k+1;l)$. Now either $e(k,k+1;b-1) \in H$, or $e(k,k+1;b-1) \notin H$. 

\null 

\noindent \textit{CASE 1: $e(k,k+1;b-1) \in H$.} Then $R(k,b-1) \mapsto R(k,b)$, $R(k,b+1) \mapsto R(k,b+2)$, $\ldots,$ $R(k,l-3) \mapsto R(k,l-2)$, is a cascade, after which $e(k,k+1;l-1)$ is in the resulting Hamiltonian path. See Figure 5.6 (a). End of Case 1.

\null

\noindent \textit{CASE 2: $e(k,k+1;b-1) \notin H$.} Since the rectangle $R(k-1,k+2;b-2,b-1)$ does not contain an end-vertex of $H$, we must have $S_{\rightarrow}(k-1,b-1; k,b-2) \in H$ and $S_{\uparrow}(k+1,b-2; k+2,b-1) \in H$. Since $R(k,b)$ is not followed by an $A_1$-type, at least one of $e(k-1,k;b-2)$ and $e(k+1,k+2;b-2)$ is in $H$. By symmetry, we may assume WLOG that $e(k-1,k;b-2) \in H$. Then, after $R(k,b-2) \mapsto R(k-1,b-2)$, $e(k,k+1;b-1)$ is in the resulting Hamiltonian cycle, and we are back to Case 1.  See Figure 5.6 (b). End of Case 2. $\square$

\endgroup

\null

\noindent \textbf{Lemma 5.12.}  Let $H$ be a Hamiltonian path of an $m \times n$ grid graph $G$, and let $n-1$ be even. Fix $k \in \{0, .., m-1\}$, and assume that the end-vertices of $H$ have x-coordinate at most $k$. A configuration in which for each $i \in \{0,2, \ldots, n-1\}$ we have that $e(k,k+1;i) \in H$ and that $e(k,k+1;i+1) \notin H$ is not a possible configuration of $H$. An analogous statement applies for the case where these edges are vertical.

\null

\noindent \textit{Proof.} Let $EE(q)$ be the statement ``For every $i \in \{0,2, \ldots, n-1\}$, $e(q,q+1;i) \in H$ and $e(q,q+1;i+1) \notin H$''. Let ($*_1$) be the statement ``$EE(q)$ implies $EE(q+2)$, whenever $q \geq k$'', and assume for now that $(*_1)$ is true. For a contradiction, assume that $EE(k)$ is true. Then $m-1-k$ is odd or $m-1-k$ is even.

\begingroup
\setlength{\intextsep}{0pt} %this shifts tikz up-down
\setlength{\columnsep}{25pt}
\begin{wrapfigure}[]{r}[0cm]{0cm}
\begin{adjustbox}{trim=0cm 0cm 0cm 0cm}
\begin{tikzpicture}[scale=1.5]

% CASE I
\begin{scope}[xshift=0cm]
{
\draw[gray,very thin, step=0.5cm, opacity=0.5] (0,0) grid (1.5,2);

{
\foreach \x in {0,...,2}
\draw[blue, line width=0.5mm] (1,\x)--++(0.5,0);

\foreach \x in {0,...,1}
\draw[black, line width=0.15mm] (1.2,\x+0.45)--(1.2,\x+0.55);

\foreach \x in {0,...,1}
\draw[black, line width=0.15mm] (1.25,\x+0.45)--(1.25,\x+0.55);

\foreach \x in {0,...,1}
\draw[black, line width=0.15mm] (1.3,\x+0.45)--(1.3,\x+0.55);

\draw[green!50!black, line width=0.5mm] (1.5,0)--++(0,1);

\draw[opacity=1] (1.5,1.5) circle [radius=0.075];

\node[left] at (0,0) [scale=1]{\tiny{0}};
\node[left] at (0,0.5) [scale=1]{\tiny{1}};
\node[left] at (0,1) [scale=1]{\tiny{2}};
\node[left] at (0,1.5) [scale=1]{\tiny{3}};
\node at (1.5,-0.15) [scale=1]{\tiny{m-1}};
}

\node[above] at (0.75,2) [scale=1.25]{(a)};

\node[below, align=center, text width=7cm] at (2, -0.25) { Fig. 5.7. (a) Case I. (b) Case II.};

}
\end{scope}

% CASE II
% CASE 1
\begin{scope}[xshift=2.5cm]
{
\draw[gray,very thin, step=0.5cm, opacity=0.5] (0,0) grid (1.5,2);

{
\foreach \x in {0,...,2}
\draw[blue, line width=0.5mm] (0.5,\x)--++(0.5,0);

\foreach \x in {0,...,1}
\draw[black, line width=0.15mm] (0.7,\x+0.45)--(0.7,\x+0.55);

\foreach \x in {0,...,1}
\draw[black, line width=0.15mm] (0.75,\x+0.45)--(0.75,\x+0.55);

\foreach \x in {0,...,1}
\draw[black, line width=0.15mm] (0.8,\x+0.45)--(0.8,\x+0.55);

\draw[green!50!black, line width=0.5mm] (1,0)--++(0.5,0)--++(0,0.5)--++(-0.5,0)--++(0,0.5);

\draw[opacity=1] (1.5,1) circle [radius=0.075];

\node[left] at (0,0) [scale=1]{\tiny{0}};
\node[left] at (0,0.5) [scale=1]{\tiny{1}};
\node[left] at (0,1) [scale=1]{\tiny{2}};
\node at (1.5,-0.15) [scale=1]{\tiny{m-1}};
}

\node[above] at (0.75,2) [scale=1.25]{(b)};

}
\end{scope}

\end{tikzpicture}
\end{adjustbox}
\end{wrapfigure}

\null 

\noindent \textit{CASE I. $m-1-k$ is odd.} By ($*_1$) and induction, $EE(m-2)$ is true. Then we must have $e(m-1;0,1) \in H$ and $e(m-1;1,2) \in H$. But now, since $v_1$ and $v_r$ are west of $x=k$, $H$ misses $v(m-1,3)$. See Figure 5.7(a). End of proof for Case I.

\null 

\noindent \textit{CASE II. $m-1-k$ is even.} By ($*_1$) and induction, $EE(m-3)$ is true. Then we must have $S_{\rightarrow}(m-2,0;m-1,1) \in H$. It follows that $S_{\downarrow}(m-2,2;m-1,1) \in H$. But now, since $v_1$ and $v_r$ are west of $x=k$, $H$ misses $v(m-1,2)$.  See Figure 5.7 (b). End of proof for Case II. It remains to prove ($*_1$).

\endgroup 

\null

\begingroup
\setlength{\intextsep}{0pt} %this shifts tikz up-down
\setlength{\columnsep}{20pt}
\begin{wrapfigure}[]{l}[0cm]{0cm}
\begin{adjustbox}{trim=0cm 0cm 0cm 0cm}

\end{adjustbox}
\end{wrapfigure}

\noindent \textit{Proof of ($*_1$).} Assume that $EE(q)$ is true. Then, either $e(q+1,q+2;0) \in H$ or $e(q+1;0,1) \in H$. See Figure 5.8. Initial configurations in blue.

\null 

\noindent \textit{CASE 1: $e(q+1,q+2;0) \in H$.} Then we must have that $S_{\downarrow}(q+1,2;q+2,1) \in H$, which implies $S_{\downarrow}(q+1,4;q+2,3) \in H$, which implies $\ldots,$ which implies $S_{\downarrow}(q+1,n-1;q+2,n-2) \in H$. It follows that $S_{\uparrow}(q+2,n-2;q+3,n-1) \in H$, which implies $\ldots, $ which implies $S_{\uparrow}(q+2,1;q+3,2) \in H$. Then $e(q+2,q+3; 0)$ must be in $H$ as well, and so, $EE(q+2)$ is true. End of proof for Case 1.

\null 

\noindent \textit{CASE 2: $e(q+1;0,1) \in H$.} Either $e(q+1;1,2) \in H$ or $e(q+1,q+2;1) \in H$.

\null 

\noindent \textit{CASE 2.1: $e(q+1;1,2) \in H$.} As in Case 1, we must have $S_{\downarrow}(q+1,4;q+2,3) \in H$, $\ldots, $ $S_{\downarrow}(q+1,n-1;q+2,n-2) \in H$. It follows that $S_{\uparrow}(q+2,n-2;q+3,n-1) \in H$, $\ldots,$ $S_{\uparrow}(q+2,1;q+3,2) \in H$. Then we must have $S_{\downarrow}(q+2,1;q+3,0)$ in $H$ as well, and so, $EE(q+2)$ is true. End of proof for Case 2.1.

\begingroup
\setlength{\intextsep}{0pt} %this shifts tikz up-down
\setlength{\columnsep}{20pt}
\begin{wrapfigure}[]{r}[0cm]{0cm}
\begin{adjustbox}{trim=0cm 0cm 0cm 0cm}
\begin{tikzpicture}[scale=1.5]

% CASE 2.2(a)
\begin{scope}[xshift=0cm]
{
\draw[gray,very thin, step=0.5cm, opacity=0.5] (0,0) grid (1.5,3);

{
\foreach \x in {0,...,3}
\draw[blue, line width=0.5mm] (0,\x)--++(0.5,0); % was (0.5,\x)

\foreach \x in {0,...,2}
\draw[black, line width=0.15mm] (0.2,\x+0.45)--(0.2,\x+0.55); % was 0.7

\foreach \x in {0,...,2}
\draw[black, line width=0.15mm] (0.25,\x+0.45)--(0.25,\x+0.55); % was 0.75

\foreach \x in {0,...,2}
\draw[black, line width=0.15mm] (0.3,\x+0.45)--(0.3,\x+0.55); % was 0.8

\draw[blue, line width=0.5mm] (0.5,0)--++(0,0.5)--++(0.5,0); % was (1,0)

\foreach \x in {1,...,2}
\draw[blue, line width=0.5mm] (0.5,\x)--++(0,0.5)--++(0.5,0); % was (1,\x)

\draw[orange!75!black, line width=0.5mm] (0.5,3)--++(1,0); % was (1,3)

\foreach \x in {1,...,3}
\draw[green!50!black, line width=0.5mm] (1.0,\x-0.5)--++(0,-0.5)--++(0.5,0); % was 1.5

\node[left] at (0,0.5) [scale=1]{\tiny{1}};
\node[left] at (0,1) [scale=1]{\tiny{2}};

\node[left] at (0,2.5) [scale=1]{\tiny{$n{-}2$}};
\node[left] at (0,3) [scale=1]{\tiny{$n{-}1$}};

\node at (0,-0.15) [scale=1]{\tiny{$q$}};
\node at (0.5,-0.15) [scale=1]{\tiny{+1}};
\node at (1,-0.15) [scale=1]{\tiny{+2}};
\node at (1.5,-0.15) [scale=1]{\tiny{+3}};
}

\node[below, align=center, text width=3cm] at (0.75, -0.25) { Fig. 5.9. Case 2.2(a).};

}
\end{scope}

\end{tikzpicture}
\end{adjustbox}
\end{wrapfigure}

\null

\noindent \textit{CASE 2.2: $e(q+1,q+2;1) \in H$.} Let $Q(j)$ be the statement ``$e(q+1;j,j+1) \in H$ and $e(q+1, q+2;j+1) \in H$''. Note that $Q(0)$ is true. Then either $Q(j)$ is true for each $j \in \{0,2, \ldots, n-3  \}$ or there is some $j_0 \in \{0,2, \ldots, n-3  \}$ such that $Q(j)$ is true for all $j \leq j_0$ but $Q(j_0+2)$ is not true.

\null 

\noindent \textit{CASE 2.2(a): $Q(j)$ is true for each $j \in \{0,2, \ldots, n-3  \}$.} This implies that $S_{\downarrow}(q+2,1;q+3,0), \dots,$ $S_{\downarrow}(q+1,n-2;q+2,n-3)$ belong to $H$. Then we must have $e(q,q+1;n-1) \in H$ and $e(q+1,q+2;n-1) \in H$ as well, and so, $EE(q+2)$ is true. See Figure 5.9. End of proof for Case 2.2(a).

\endgroup

\begingroup
\setlength{\intextsep}{0pt} %this shifts tikz up-down
\setlength{\columnsep}{20pt}
\begin{wrapfigure}[]{l}[0cm]{0cm}
\begin{adjustbox}{trim=0cm 0.25cm 0cm 0cm}

\end{adjustbox}
\end{wrapfigure}

\null 

\noindent \textit{CASE 2.2(b): There is some $j_0 \in \{0,2, \ldots, n-3  \}$ such that $Q(j)$ is true for all $j \leq j_0$ but $Q(j_0+2)$ is not true.} Then $S_{\downarrow}(q+2,1;q+3,0), \ldots, $ $S_{\downarrow}(q+2,j_0+1;q+3,j_0)$ belong to $H$. Note that if $e(q+1;j_0+2,j_0+3)\notin H$ and $e(q+1,q+2;j_0+3)\notin H$ then $H$ misses $v(q+1,j_0+3)$. Therefore we only need to check two cases:  Case $b_1$: $e(q+1;j_0+2,j_0+3)\in H$ and $e(q+1, q+2;j_0+3)\notin H$ and $b_2$: $e(q+1;j_0+2,j_0+3)\notin H$ and $e(q+1, q+2;j_0+3)\in H$.

\null

\noindent \textit{CASE 2.2($b_1$): $e(q+1;j_0+2,j_0+3)\in H$ and $e(q+1, q+2;j_0+3)\notin H$.} Then $e(q+1; j_0+3,j_0+4) \in H$,  $S_{\downarrow}(q+2,j_0+3;q+3,j_0+2) \in H$ and $S_{\downarrow}(q+1,j_0+6;q+2,j_0+5), \ldots,$ $S_{\downarrow}(q+1,n-1;q+2,n-2)$ all belong to $H$. It follows that $S_{\uparrow}(q+2,n-2;q+3,n-1), \ldots, $ $S_{\uparrow}(q+2,j_0+3;q+3,j_0+4)$ all belong to $H$ as well, and so $EE(q+2)$ is true. See Figure 5.10 (a). End of proof for Case 2.2($b_1$).

\null 

\noindent \textit{CASE 2.2($b_2$): $e(q+1;j_0+2,j_0+3)\notin H$ and $e(q+1, q+2;j_0+3)\in H$.}  We must have $e(q+1,q+2;j_0+2)\in H$ and $e(q+1;j_0+3, j_0+4) \in H$. It follows that $S_{\downarrow}(q+1,j_0+6;q+2,j_0+5), \ldots,$ $S_{\downarrow}(q+1,n-1;q+2,n-2)$ all belong to $H$. Then we must have $S_{\uparrow}(q+2,n-2;q+3,n-1), \ldots$, $S_{\uparrow}(q+2,j_0+3;q+3,j_0+4)$ in $H$ as well. Then we must have $e(q+2,q+3;j_0+2)\in H$, and so $EE(q+2)$ is satisfied. See Figure 5.10 (b). End of proof for Case 2.2($b_2$). End of Case 2.2(b). End of Case 2.2. End of Case 2. End of proof for ($*_1$). $\square$

\null

\noindent \textbf{Lemma 5.13.}  Let $H$ be a Hamiltonian path of an $m \times n$ grid graph  $G$. Fix $k \in \{0, .., m-1\}$, and assume that the end-vertices of $H$ have x-coordinate at most $k$. Assume that $n-1$ is odd, and let $R(k+1,0)$ be a northern leaf on the southern boundary. The configuration in which there is a full stack of $A_0$'s starting at $R(k+1,0)$ is not a possible configuration of $H$. 

Analogous statements apply for the other compass directions.

\null

\noindent \textit{Proof.} Note that if $k+3=m-1$, then $e(k+3;0,1)\in H$ and $e(k+3;1,2)\in H$. Then $e(k+3;3,4) \in H$,  or $e(k+3;4,5) \in H$. If the former, then $v(k+3,3)$ is an end-vertex of $H$, contradicting the assumption that the end-vertices have x-coordinate at most $k$; and if the latter, then  $H$ misses $v(k+3,3)$. Then we must have that, $m-1 \geq k+4$. See Figure 5.11(a).

If $k+4=m-1$, then we must have that $S_{\rightarrow}(k+3,0;k+4,1) \in H$. This implies that $S_{\downarrow}(k+3,2;k+4,1) \in H$. But then again, as in the previous paragraph, either $H$ misses $v(k+4,2)$ or $v(k+4,2)$ is an end-vertex. So, we may assume that $m-1 \geq k+5$. See Figure 5.11(b).

\begingroup
\setlength{\intextsep}{0pt} %this shifts tikz up-down
\setlength{\columnsep}{15pt}
\begin{wrapfigure}[]{r}[0cm]{0cm}
\begin{adjustbox}{trim=0cm 0cm 0cm 0cm}
\begin{tikzpicture}[scale=1.5]

% CASE 1
\begin{scope}[xshift=0cm]
{
\draw[gray,very thin, step=0.5cm, opacity=0.5] (0,0) grid (1.5,2.5);

\draw[blue, line width=0.5mm] (0,0)--++(0.5,0)--++(0,0.5)--++(0.5,0)--++(0,-0.5)--++(0.5,0);

\draw[blue, line width=0.5mm] (0,1)--++(0.5,0)--++(0,0.5)--++(0.5,0)--++(0,-0.5)--++(0.5,0);

\draw[blue, line width=0.5mm] (0,2)--++(0.5,0)--++(0,0.5)--++(0.5,0)--++(0,-0.5)--++(0.5,0);

\draw[blue, line width=0.5mm] (1.5,0)--++(0,1);

\draw[opacity=1] (1.5,1.5) circle [radius=0.05];

\node[left] at (0,0) [scale=1]{\tiny{0}};
\node[left] at (0,2.5) [scale=1]{\tiny{$n-1$}};
\node[left] at (0,2) [scale=1]{\tiny{$n-2$}};
\node at (0,-0.15) [scale=1]{\tiny{$k$}};
\node at (0.5,-0.15) [scale=1]{\tiny{+1}};
\node at (1,-0.15) [scale=1]{\tiny{+2}};
\node at (1.5,-0.15) [scale=1]{\tiny{+3}};

\node[above] at (0.75,2.5) [scale=1.25]{(a)};

\node[below, align=center, text width=7cm] at (2.25, -0.25) { Fig. 5.11.};

}
\end{scope}

\begin{scope}[xshift=2.5cm]
{
\draw[gray,very thin, step=0.5cm, opacity=0.5] (0,0) grid (2,2.5);

\draw[blue, line width=0.5mm] (0,0)--++(0.5,0)--++(0,0.5)--++(0.5,0)--++(0,-0.5)--++(0.5,0);

\draw[blue, line width=0.5mm] (0,1)--++(0.5,0)--++(0,0.5)--++(0.5,0)--++(0,-0.5)--++(0.5,0);

\draw[blue, line width=0.5mm] (0,2)--++(0.5,0)--++(0,0.5)--++(0.5,0)--++(0,-0.5)--++(0.5,0);

\draw[blue, line width=0.5mm] (1.5,0)--++(0.5,0)--++(0,0.5)--++(-0.5,0)--++(0,0.5);

\draw[opacity=1] (2,1) circle [radius=0.05];

\node[left] at (0,0) [scale=1]{\tiny{0}};
\node[left] at (0,2.5) [scale=1]{\tiny{$n-1$}};
\node[left] at (0,2) [scale=1]{\tiny{$n-2$}};
\node at (0,-0.15) [scale=1]{\tiny{$k$}};
\node at (0.5,-0.15) [scale=1]{\tiny{+1}};
\node at (1,-0.15) [scale=1]{\tiny{+2}};
\node at (1.5,-0.15) [scale=1]{\tiny{+3}};

\node[above] at (1,2.5) [scale=1.25]{(b)};
}
\end{scope}

\end{tikzpicture}
\end{adjustbox}
\end{wrapfigure}

Let $SS(q)$ be the statement ``There is a full stack of $A_0s$ starting at $R(q,0)$''. Let $(*_2)$ be the statement ``$SS(q)$ implies $SS(q+2)$, whenever $q>k$", and assume for now that $(*_2)$ is true. For a contradiction, assume that $SS(q+1)$ is true. Then, $m-1-k$ is odd or it is even. In each case, the arguments are similar to those in Cases I and II of Lemma 5.12, so we omit them.

It remains to prove $(*_2)$.

\null

\noindent \textit{Proof of $(*_2)$.} Assume $SS(q)$. Then we must have that $S_{\uparrow}(q+2,n-2;q+3,n-1) \in H$, which implies $S_{\uparrow}(q+2,n-4;q+3,n-3) \in H$, which implies $\ldots,$ which implies $S_{\uparrow}(q+2,0;q+3,1) \in H$. It follows that $S_{\downarrow}(q+3,1;q+4,0) \in H$, which implies $\ldots,$ which implies $S_{\uparrow}(q+3,n-1;q+4,n-2) \in H$. So, $SS(q+2)$ is true. $\square$

\endgroup

\null

\noindent \textbf{Observation 5.14.} Let $G$ be an $m \times n$ grid graph, and let $H$ be a Hamiltonian path of $G$ with end-vertices $u$ and $v$ not adjacent (NAA). Let $J$ be an $H$-component with more than one box, and let $N_J$ be the neck of $J$. If $N_J$ is not switchable, then:

(a) exactly one of $u$ and $v$ (say, $u$) is incident on the neck-edge of $N_J$,

(b) $u$ is not at a corner of $R_0$, and 

(c) $e(u) \in R_0$.

\null

\noindent \textbf{Lemma 5.15.} Let $G$ be an $m \times n$ grid graph, and let $H$ be a Hamiltonian path of $G$ with end-vertices $u$ and $v$ not adjacent. Let $W$ be a switchable box in $H$, and let $X$ and $Y$ be the boxes adjacent to $W$ that are not its $H$-neighbours. Assume that $X$ and $Y$ belong to distinct $H$-components $J_X$ and $J_Y$. Then:

\begin{itemize}
    \item[$(a)$] If at least one of $u$ and $v$ does not lie on $R_0$, then there exists a cascade $\mu_1, \dots, \mu_s$, such that $\mu_s$ switches $W$ and fixes both $u$ and $v$, and $s \in \{1,2\}$.

    \item[$(b)$] If at least one of $u$ and $v$ lies on a corner of $R_0$, then there exists a cascade $\mu_1, \dots, \mu_s$ such that $\mu_s$ switches $W$ and fixes both $u$ and $v$, and $s \in \{1,2\}$.
    
    \item[$(c)$] If both end-vertices $u$ and $v$ lie on $R_0$, then there exists a cascade that keeps $v$ on the same side of $R_0$ and that switches $W$ and fixes $u$, and there exists another cascade that keeps $u$ on the same side of $R_0$ and that switches $W$ and fixes $v$. Each cascade has length at most three.

\end{itemize}

\noindent \textit{Proof.} Orient $H$ as $v=v_1, v_2, \ldots, v_r=u$. Note that if $W$ is parallel then we are done after $\text{Sw}(W)$, so we may assume that $W$ is anti-parallel. Let $e_X$ and $e_Y$ be the edges of $W$ in $H$, let $X=\Phi(e_X, \text{right})$, let $\Phi(e_Y, \text{right})=Y$, and WLOG assume that $e_X$ precedes $e_Y$ in $H$. Let $N_{J_X}$ and $N_{J_Y}$ be the necks of $J_X$ and $J_Y$, respectively, and let $\{v_s,v_t\}$ be the neck-edge of $J_X$. As in Corollary 1.3.8, the neck-edge of $J_X$ splits $H$ into three subpaths $P_1=P(v_1,v_s)$, $P_2=P(v_s,v_t)$ and $P_3=P(v_t,v_r)$. Now, either at least one of $N_{J_X}$ and $N_{J_Y}$ is switchable, or neither is. 

\null

\noindent\textit{CASE 1. At least one of $N_{J_X}$ and $N_{J_Y}$ is switchable.} It follows from Observation 5.14 that if at least one of $u$ and $v$ does not lie on $R_0$ or at least one of $u$ and $v$ lies on a corner of $R_0$, then at least one of $N_{J_X}$ and $N_{J_Y}$ is switchable. Therefore, parts (a) and (b) will follow from this case alone.

WLOG, assume that $N_{J_X}$ is switchable. Then $N_{J_X}$ is parallel or $N_{J_X}$ is anti-parallel.

\null 

\noindent\textit{CASE 1.1: $N_{J_X}$  is anti-parallel.} Then, $J_X$ contains neither end-vertex, or it contains exactly one end-vertex, or it contains both.

\begingroup
\setlength{\intextsep}{0pt}
\setlength{\columnsep}{20pt}
\begin{wrapfigure}[]{l}{0cm}
\begin{adjustbox}{trim=0cm 0cm 0cm 0cm}

\end{adjustbox}
\end{wrapfigure}

\null

\noindent \textit{CASE 1.1(a): $J_X$ contains neither end-vertex.}  
We will first check that $e_X \in P_2$ and $e_Y \notin P_2$. Since $X=\Phi(e_X, \text{right})$, by Corollary 1.3.10, $e_X \in \overrightarrow{K}_{J_X}$. By Corollary 1.3.11, since $J_X$ contains no end-vertices of $H$, $\overrightarrow{K}_{J_X}=\overrightarrow{K}_{P_2}$. By construction, the directed subpath $P_2$ of $H$ is identical to $\overrightarrow{K}_{P_2}$. Thus, $e_X \in P_2$. By Corollary 1.3.11,  $\Phi(P_2, \text{right}) = J_X$.  For contradiction, assume that $e_Y \in P_2$. %,and WLOG assume that $e_X$ precedes $e_Y$. 
By Observation 1.2.4 (a), $\Phi(P_2, \text{right})\supset \Phi(\overrightarrow{K}(e_X,e_Y),\text{right})$. But then there is an $H$-path $P(X,Y) \subseteq \Phi(\overrightarrow{K}(e_X,e_Y),\text{right}) \subset \Phi(P_2, \text{right})=J_X$, contradicting that $X$ and $Y$ belong to distinct components. Thus, $e_Y \notin P_2$.

Next we check that the edges of $N_{J_X}$ in $H$ are $(v_s, v_{s+1})$ and $(v_{t-1}, v_t)$. Since $\Phi(P_2, \text{right}) = J_X$, $N_{J_X}$ must be added by an edge in $P_2$. Since the neck-edge of $N_{J_X}$ is $\{v_s,v_t\}$, the edge of $H$ incident on $v_s$ adding $N_{J_X}$ is either $(v_{s-1}, v_s)$ 
or $(v_s, v_{s+1})$. But $v_{s-1} \notin P_2$. Thus $N_{J_X}$ must be added by $(v_s, v_{s+1})$. Using the same argument, we see that the other edge of  $N_{J_X}$ in $H$ must be  $(v_{t-1}, v_t)$. See Figure 5.12.

Now, after switching $N_{J_X}$, $e_X \in H_c=P(v_{s+1},v_{t-1}), \{v_{t-1},v_{s+1}\}$ and $e_Y \in H_p=P(v_1,v_s)$, $\{v_s,v_t\}, P(v_t,u)$. Then $W$ is an $(H_c, H_p)$-port and so  $N_{J_X} \mapsto W$ is a valid move that switches $W$ and fixes both end-vertices of $H$. End of Case 1.1(a). 

\endgroup

\null

\noindent \textit{CASE 1.1(b): $J_X$ contains exactly one end-vertex.} WLOG assume that $J_X$ contains $v$. Corollary 1.5.5 implies that $v \notin J_X \setminus B(J_X)$. Then $v \in B(J_X)$. By Corollary 1.3.11(d),  $B(J_X)=P_2, \{v_s,v_t\}$. It follows that $v=v_s$, so $P_1=v$ is trivial. Then, by Corollary 1.3.11(b), $\overrightarrow{K}_{J_X}=\overrightarrow{K}_{P_1},\overrightarrow{K}_{P_2}=\overrightarrow{K}_{P_2}$. The rest of the argument is the same as in Case 1.1(a), so we omit it. End of Case 1.1(b)

\null 

\noindent\textit{CASE 1.1(c): $J_X$ contains both end-vertices.} By Corollary 1.5.5, either $u$ and $v$ are both in $J_X \setminus B(J_X)$, or neither is.

\null

\noindent\textit{CASE 1.1$(c_1)$: Neither $u$ nor $v$ is in $J_X \setminus B(J_X)$.} Then both $u$ and $v$ must belong to $B(J_X)$. Using Corollary 1.1.13(d) again, we see that we must have $v=v_s$ and $u=v_t$. As in Case 1.1(b), this means that $\overrightarrow{K}_{J_X}=\overrightarrow{K}_{P_2}$. The rest of the argument is the same as in Case 1.1(a), so we omit it. End of Case 1.1$(c_1)$.

\null

\noindent\textit{CASE 1.1$(c_2)$: both $u$ and $v$ belong to $J_X \setminus B(J_X)$.} Then $J_Y$ contains neither end-vertex. By Observation 5.14, $N_{J_Y}$ is switchable. Corollary 1.5.3 implies that $N_{J_Y}$ is anti-parallel. Then, by Case 1.1(a), $N_{J_Y} \mapsto W$ is a valid move, that switches $W$ and fixes both end-vertices of $H$. End of Case 1.1$(c_2)$. End of Case 1.1(c). End of Case 1.1.

\null 

\noindent \textit{CASE 1.2: $N_{J_X}$  is parallel.} Then there are two subcases to consider:

(a) $e_X\in P_i$, $e_Y\in P_j $, with $(i,j) \in \{ (1,1), (2,2), (3,3), (1,3), (3,1) \}$. 

(b) $e_X\in P_i$, $e_Y\in P_j $, with $(i,j) \in \{ (1,2), (2,1), (2,3), (3,2)  \}$, and 

\null

\noindent \textit{CASE 1.2(a).} We will check that each subcase in (a) is impossible. 

\null 

\noindent \textit{CASE 1.2($a_1$): $(i,j)=(1,1)$ or $(i,j)=(3,3)$}. By symmetry, assume WLOG that $(i,j)=(1,1)$. By Corollary 1.3.8 (b), $\Phi(\overrightarrow{K}_{P_1}, \text{right}) \subset J_X$ or $\Phi(\overrightarrow{K}_{P_1},\text{right}) \subset G_{-1} \setminus J_X$. If the former, then there is an $H$-path $P(X,Y)$ contained in $\Phi(\overrightarrow{K}(e_X,e_Y),\text{right}) \subseteq \Phi(\overrightarrow{K}_{P_1}\text{right})$, contradicting that $X$ and $Y$ belong to distinct components. If the latter, then $X=\Phi(e_X, \text{right}) \subset \Phi(\overrightarrow{K}_{P_1},\text{right})$ implies that $X \in G_{-1} \setminus J_X$, contradicting $X \in J_X$. End of Case 1.2($a_1$).

\null 

\noindent \textit{CASE 1.2($a_2$): $(i,j)=(2,2)$.} Since $X=\Phi(e_X,\text{right})$, Corollary 1.3.8 (a) $\Phi(P_2, \text{right}) \subseteq J_X$. For contradiction, assume the $e_Y \in P_2$ as well. Then there is an $H$-path $P(X,Y)$ contained in $\Phi(\overrightarrow{K}(e_X,e_Y),\text{right}) \subset \Phi(P_2, \text{right})$, contradicting that $X$ and $Y$ belong to distinct components. End of Case 1.2($a_2$).

\null

\noindent \textit{CASE 1.2($a_3$): $(i,j)=(1,3)$ or $(i,j)=(3,1)$.}  By symmetry, assume WLOG that $(i,j)=(1,3)$. Let $e_X=(v_x,v_{x+1})$ and $e_Y=(v_{y-1},v_y)$. By definition of $\overrightarrow{K}_{P_1}$ (recall definition in Section 1.3), the edge $(v_{x+1},v_x)$ of $H^*$ (recall definition in Section 1.2) belongs to $\overrightarrow{K}_{P_1}$. Since $X=\Phi((v_x,v_{x+1}), \text{right})$, we have that $W=\Phi((v_{x+1},v_x), \text{right})$. Similarly, $(v_y,v_{y-1})$ belongs to $\overrightarrow{K}_{P_3}$, and $W=\Phi((v_y,v_{y-1}), \text{right})$. So we have that $W \in \Phi(\overrightarrow{K}_{P_1}, \text{right}) \cap \Phi(\overrightarrow{K}_{P_3}, \text{right})$. Since $X \in J_X \cap \Phi(\overrightarrow{K}_{P_1}, \text{right})$, by Corollary 1.3.8 (b), $\Phi(\overrightarrow{K}_{P_1}, \text{right}) \subset J_X$. Since $W \in \Phi(\overrightarrow{K}_{P_1}, \text{right})$, $W \in J_X$. Since $W \in J_X \cap \Phi(\overrightarrow{K}_{P_3},\text{right})$, by Corollary 1.3.8 (b), $\Phi(\overrightarrow{K}_{P_3},\text{right}) \subset J_X$. But then $Y=\Phi(e_Y, \text{right})$ belongs to $\Phi(\overrightarrow{K}_{P_3},\text{right})\subset J_X$, contradicting that $X$ and $Y$ belong to distinct $H$-components. End of Case 1.2($a_3$) End of Case 1.2 (a).

\null

\noindent \textit{CASE 1.2(b).} By Corollary 1.4.2 (ii), if we apply Sw($N_{J_X}$), $W$ becomes parallel. Then, Sw($N_{J_X}$), Sw($W$), is a cascade that switches $W$ and fixes both end-vertices of $H$. End of Case 1.2(b). End of Case 1.2. End of Case 1. End of proof for parts (a) and (b), as well as for the scenario in part (c) where at least one of $N_{J_X}$ and $N_{J_Y}$ is switchable.

\null

\begingroup
\setlength{\intextsep}{0pt}
\setlength{\columnsep}{20pt}
\begin{wrapfigure}[]{l}{0cm}
\begin{adjustbox}{trim=0cm 0cm 0cm 0cm}

\end{adjustbox}
\end{wrapfigure}

\noindent \textit{CASE 2. Neither $N_{J_X}$ nor $N_{J_Y}$ is switchable.} By Observation 5.14, the neck-edges of $J_Y$ and $J_X$ are incident on distinct end-vertices of $H$. In keeping with the notation of Corollary 1.3.8, we note that in this case, $v_1=v_s$, so $P_1$ is the trivial path consisting of just the vertex $v_1$. Thus the neck-edge of $J_X$ splits $H$ into the subpaths $P_2=P(v_1,v_t)$ and $P_3=P(v_t,v_r)$, where $v_r$ is incident on the neck edge of $N_{J_Y}$. By Observation 5.14, $(v_1,v_2)$ must be in $R_0$ and neither  $u$ nor $v$ may be at a corner of $R_0$. See Figure 5.13. There are two cases to consider:

1) $e_X\in P_i, e_Y\in P_j $, with $i=j$ 

2) $e_X\in P_i, e_Y\in P_j $, with $i \neq j$ 

\null 

\noindent \textit{CASE 2.1: $e_X\in P_i, e_Y\in P_j $, with $i=j$.} Since $P_1$ is the trivial path, there are two possibilities: $e_X\in P_2, e_Y\in P_2$, and $e_X\in P_3, e_Y\in P_3$. The case where $e_X\in P_2, e_Y\in P_2$ was shown to be impossible in Case 1.2($a_2$), so we omit it. We will check that the case where $e_X\in P_3, e_Y\in P_3$ is impossible as well. 

\endgroup

%deleted
{
%\textit{CASE 2.1(a): $e_X\in P_2, e_Y\in P_2$. } By Observation 1.2.4 (a), $\Phi(P_2, \text{right})\supset \Phi(\overrightarrow{K}(e_X,e_Y),\text{left})$. Since $G$ is simply connected, by Corollary 1.3.8 (a),  $\Phi(P_2, \text{right}) \subset J_X$, or $\Phi(P_2, \text{right}) \subset G_{-1} \setminus J_X$. If the former, then there is an $H$-path $P(X,Y) \subseteq \Phi(\overrightarrow{K}(e_X,e_Y), \text{right}) \subseteq  \Phi(P_2, \text{right}) \subseteq J_X$, contradicting the fact that $X$ and $Y$ belong to distinct components. And if the latter, then $X=\Phi(e_X, \text{right}) \in \Phi(\overrightarrow{K}(e_X,e_Y),\text{right}) \subseteq  \Phi(P_2, \text{right}) \subset G_{-1} \setminus J_X$, contradicting $X \in J(X)$. End of Case 2.1(a).
}

Assume that $e_X\in P_3$ and $e_Y\in P_3$. By Observation 1.2.4 (a), $\Phi(\overrightarrow{K}_{P_3}, \text{right})\supset \Phi(P_3, \text{right})$ $\supset \Phi(\overrightarrow{K}(e_X,e_Y),\text{right})$. By Corollary 1.3.8 (b), $\Phi(\overrightarrow{K}_{P_3}, \text{right}) \subset J(X)$ or $\Phi(\overrightarrow{K}_{P_3}, \text{right}) \subset G_{-1} \setminus J_X$. If the former, there is an $H$-path $P(X,Y) \subseteq \Phi(\overrightarrow{K}(e_X,e_Y),\text{right}) \subseteq \Phi(P_3, \text{right}) \subseteq 
\Phi(\overrightarrow{K}_{P_3}, \text{right}) \subset J_X$, contradicting that $X$ and $Y$ belong to distinct components. And if the latter then $X=\Phi(e_X, \text{right}) \in$
$ \Phi(P_3, \text{right}) \subseteq 
\Phi(\overrightarrow{K}_{P_3}, \text{right}) \subset G_{-1} \setminus J_X$, contradicting that $X \in J_X$. End of Case 2.1.

\null

\noindent \textit{CASE 2.2: $e_X\in P_i, e_Y\in P_j $ with $i \neq j$.} Since both end-vertices of $H$ are in $R_0$, $J_X$ contains neither end-vertex. By NAA, $v_t \neq u$. Then, by Case 1.1(a), we have that $e_X \in P_2$, $e_Y \notin P_2$, and that $(v_{t-1},v_t)$ is the edge of $N_{J_X}$ in $H$ that is incident on $v_t$. It follows that $e_Y \in P_3$, and that $(v_t,v_{t+1}) \in R_0$.

\begingroup
\setlength{\intextsep}{0pt}
\setlength{\columnsep}{20pt}
\begin{center}
\begin{adjustbox}{trim=0cm 0cm 0cm 0cm}

\end{adjustbox}
\end{center}

\noindent For definiteness, assume that $v=v_1=v(k,0)$. It follows that $v_2=v(k-1,0)$ and $v_t=v(k+1,0)$. Note that $W \neq N_{J_X}$ and $W \neq N_{J_X}+(1,0)$. Now, by Corollary 1.4.5 (ii), after $bb_{v_1}(east)$, $W$ is parallel in the resulting Hamiltonian path $H'$, the end-vertex $v$ is relocated to $v(k+1,1)$, and the other end-vertex $u$ is fixed. Since $W$ is neither $N_{J_X}$, nor adjacent to $N_{J_X}$, applying Sw($W$) to $H'$ to obtain $H''$ does not remove the edges $e(k,k+1;0)$ or $e(k+1,k+2;0)$. Then we can apply $bb_{v(k+1,1)}(south)$ to $H''$, after which the end-vertex on $v(k+1,1)$ is relocated to $R_0$, and $u$ is fixed. It follows that, starting at $H$, $bb_{v_1}(east)$, Sw($W$), $bb_{v(k+1,1)}(south)$, is a cascade that switches $W$, keeps $v$ on the same side of $R_0$, and fixes $u$. See Figure 5.14.

There is a similar cascade that starts with a backbite move on the end-vertex $u$ incident on the neck-edge of $N_{J_Y}$ that switches $W$,  keeps $u$ on the same side of $R_0$, and fixes $v$.  End of Case 2.2 $\square$

\null

\noindent \textbf{Lemma 5.16.}  Let $G$ be an $m \times n$ grid graph and let $H$ be a Hamiltonian path of $G$ with an end-vertex $v$ incident on $R_0$, and assume that the end-vertices of $H$ are not adjacent. Let the edge $e$ of $H$ be followed by an $A_1$-type with switchable middle-box $W$, and let $X$ and $Y$ be the boxes adjacent to $W$ that are not its $H$-neighbours. Assume that $X$ and $Y$ belong to the same $H$-component, and let $P(X,Y)$ be the looping $H$-path of $W$. If $v$ is incident on $P(X,Y)$, then $P(X,Y)$ has a switchable box $Z$ incident on $R_0$. 

%\textcolor{red}{and such that $Z$ is either incident on $e(v)$, or incident the edge $(v_1,v_2)$.}

\null

\begingroup
\setlength{\intextsep}{0pt}
\setlength{\columnsep}{20pt}
\begin{wrapfigure}[]{r}{0cm}
\begin{adjustbox}{trim=0cm 0cm 0cm 0cm} 
\begin{tikzpicture}[scale=1.75]

\begin{scope}[xshift=0cm] 
{
\draw[gray,very thin, step=0.5cm, opacity=0.5] (0,0) grid (1,0.5);

\fill[blue!50!white, opacity=0.5] (0,0)--++(0.5,0)--++(0,0.5)--++(-0.5,0);

\draw[blue, line width=0.5mm] (0,0)--++(0,0.5);
\draw [blue, line width=0.5mm] (0.5,0.5)--++(0,-0.5)--++(0.5,0);

\draw[fill=blue] (0,0) circle (0.05);
\node[right] at (0,0) [scale=0.8]{\small{$v$}};

% labellings 
{

\node at (0.25,0.25) [scale=0.8]
{\small{$Z$}};

\node[above] at (0,0.5) [scale=1]
{\tiny{$k$}};
\node[above] at (0.5,0.5) [scale=1]
{\tiny{$+1$}};

\node[right] at (1,0) [scale=1]
{\tiny{$0$}};
\node[right] at (1,0.5) [scale=1]
{\tiny{$1$}};

}

}

\node[below, align=center, text width=3cm] at (0.5, -0.15) { Fig. 5.15. $e(v)$ is northern.};

\end{scope}

\end{tikzpicture}
\end{adjustbox}
\end{wrapfigure}

\noindent \textit{Proof.} We will use NAA repeatedly and implicitly. First note that since the $A_1$-type follows $e$, $W$ cannot be incident on $R_0$. Let $v$ be an end-vertex of $H$ incident on a box $Z \in P(X,Y)=P$. For definiteness, assume that $Z=R(k,0)$ and $v=v(k,0)$. Note that if $e(v)$ is northern, then we must have that $S_{\downarrow}(k+1,1;k+2,0)\in H$. This implies that $Z$ is an end-box of $P$, which is not possible since $Z$ does not fit the description of an end-box of $P$. See Figure 5.15. The observation that $W$ is not incident on $R_0$ implies that $W\neq Z+(\pm 1,0)$; and $W\neq Z+(0,1)$, because that means that $P(X,Y)=Z$. See Figure 5.17. Then $e(v)$ must be either eastern, or western.

\endgroup 

\null

\noindent \textit{CASE 1: $e(v)$ is eastern.} By NAA, it is easy to observe that $k \geq 2$. Then $e(k-1;0,1) \in H$. Either $Z+(-1,0) \in P$, or  $Z+(-1,0) \notin P$.

\null 

\noindent \textit{CASE 1.1: $Z+(-1,0) \notin P$.} Note that the observation that $W$ is not incident on $R_0$, together with the assumption that $W$ is the switchable middle-box of an $A_1$-type, imply that $Z$ cannot be an end-box of $P$. Then the $H$-neighbours of $Z$ in $P$ must be $Z+(0,1)$ and $Z+(1,0)$. Then $S_{\rightarrow}(k-1,1;k,2) \in H$, and $e(k+1,k+2;0) \in H$. There are three possibilities: $v(k+1,1)$ is the other end-vertex $u$ of $H$ and $e(u)=e(k+1,k+2;1)$, $v(k+1,1)$ is the other end-vertex $u$ of $H$ and $e(u)=e(k+1;1,2)$, and $v(k+1,1)$ is not the other end-vertex of $H$.

\null 

\noindent \textit{CASE 1.1(a): $u=v(k+1,1)$, $e(u)=e(k+1,k+2;1)$.} Note that $Z+(1,0)$ does not fit the description of an end-box of $P$. Then $Z+(1,0) \in P$ must be switchable. See Figure 5.16 (a). End of Case 1.1(a).

\null 

\noindent \textit{CASE 1.1(b): $u=v(k+1,1)$, $e(u)=e(k+1;1,2)$.} Note that $Z+(0,1)$ does not fit the description of an end-box of $P$. Then $Z+(0,1) \in P$ must be switchable. See Figure 5.16 (b). End of Case 1.1(b).

\null 

\noindent \textit{CASE 1.1(c): $v(k+1,1) \neq u$.} Then $S_{\downarrow}(k+1,2;k+2,1) \in H$. Now, $W \neq Z+(2,0)$, and since  $Z+(1,1)$ is not switchable, $W \neq Z+(1,1)$. It follows that $Z+(1,0)$ is not an end-box of $P$. Then $Z+(1,0)$ must be switchable. See Figure 5.16 (c). End of Case 1.1(c). End of Case 1.1.

\begingroup

\begin{center}
\setlength{\intextsep}{0pt}
\setlength{\columnsep}{20pt}
\begin{adjustbox}{trim=0cm 0cm 0cm 0cm} 

\end{adjustbox}
\end{center}

\endgroup

\noindent \textit{CASE 1.2: $Z+(-1,0) \in P$.} Now, either $Z+(-1,0)$ is an end-box of $P$ or it is not. 

\null 

\noindent \textit{CASE 1.2(a): $Z+(-1,0)$ is an end-box of $P$.} This implies that $W=Z+(-1,1)$, and that $(k-1,k;1) \in H$, $(k,k+1;1) \in H$, and that $Z$ is not an end-box of $P$. Then $Z$ must be switchable. See Figure 5.16 (d). End of Case 1.2(a).

\null 

\noindent \textit{CASE 1.2(b): $Z+(-1,0)$ is not an end-box of $P$.} Then we must have that $e(k-1,k;1) \notin H$, and that $S_{\downarrow}(k,2;k+1,1) \in H$. Now, $W \neq Z+(1,0)$, and since $Z+(0,1)$ is not switchable, $W \neq Z+(0,1)$, so $Z$ is not an end-box of $P$. It follows that $Z$ must be switchable. See Figure 5.16 (e). End of Case 1.2(b). End of Case 1.2. End of Case 1.

\null

\noindent \textit{CASE 2: $e(v)$ is western.} Then $e(k+1;0,1) \in H$. Either $Z$ is an end-box of $P$, or it is not.

\null

\begingroup
\setlength{\intextsep}{0pt}
\setlength{\columnsep}{20pt}
\begin{wrapfigure}[]{r}{0cm}
\begin{adjustbox}{trim=0cm 0cm 0cm 0cm} 
\begin{tikzpicture}[scale=1.75]

\begin{scope}[xshift=0cm] 
{
\draw[gray,very thin, step=0.5cm, opacity=0.5] (0,0) grid (1,1);

\fill[blue!50!white, opacity=0.5] (0,0)--++(1,0)--++(0,0.5)--++(-1,0);

\draw[blue, line width=0.5mm] (0,0)--++(0.5,0);
\draw [blue, line width=0.5mm] (0,0.5)--++(1,0)--++(0,-0.5);

\draw[fill=blue] (0.5,0) circle (0.05);
\node[above] at (0.5,0) [scale=0.8]{\small{$v$}};

% labellings 
{

\node at (0.75,0.25) [scale=0.8]
{\small{$Z$}};
\node at (0.75,0.75) [scale=0.8]
{\small{$W$}};

\node[above] at (1,1) [scale=1]
{\tiny{$+1$}};
\node[above] at (0.5,1) [scale=1]
{\tiny{$k$}};
\node[above] at (0,1) [scale=1]
{\tiny{$-1$}};

\node[right] at (1,0) [scale=1]
{\tiny{$0$}};
\node[right] at (1,0.5) [scale=1]
{\tiny{$1$}};
\node[right] at (1,1) [scale=1]
{\tiny{$2$}};
}

}

\node[above] at (0.5,1.15) [scale=1.25]{(a)};

\node[below, align=left, text width=5cm] at (1.5, -0.25) { Fig. 5.17. (a) Case 2.1. (b) Case 2.2(b).};

\end{scope}

\begin{scope}[xshift=2cm] 
{
\draw[gray,very thin, step=0.5cm, opacity=0.5] (0,0) grid (1,1);

\fill[blue!50!white, opacity=0.5] (0,0)--++(1,0)--++(0,0.5)--++(-1,0);
\fill[blue!50!white, opacity=0.5] (0.5,0.5)--++(0.5,0)--++(0,0.5)--++(-0.5,0);

\draw[blue, line width=0.5mm] (0,0)--++(0.5,0);
\draw [blue, line width=0.5mm] (0,0.5)--++(0.5,0)--++(0,0.5);
\draw [blue, line width=0.5mm] (1,0.5)--++(0,-0.5);

\draw[fill=blue] (0.5,0) circle (0.05);
\node[above] at (0.5,0) [scale=0.8]{\small{$v$}};

% labellings 
{

\node at (0.75,0.25) [scale=0.8]
{\small{$Z$}};

\node[above] at (1,1) [scale=1]
{\tiny{$+1$}};
\node[above] at (0.5,1) [scale=1]
{\tiny{$k$}};
\node[above] at (0,1) [scale=1]
{\tiny{$-1$}};

\node[right] at (1,0) [scale=1]
{\tiny{$0$}};
\node[right] at (1,0.5) [scale=1]
{\tiny{$1$}};
\node[right] at (1,1) [scale=1]
{\tiny{$2$}};
}

}

\node[above] at (0.5,1.15) [scale=1.25]{(b)};

\end{scope}

\end{tikzpicture}
\end{adjustbox}
\end{wrapfigure}

\noindent \textit{CASE 2.1: $Z$ is an end-box of $P$.} Then we must have that $W=Z+(0,1)$, so $e(k-1,k;1)\in H$,  and $e(k,k+1;1)\in H$. Then $Z+(-1,0) \in P$ is not an end-box of $P$ and is switchable. See Figure 5.17 (a). End of Case 2.1.

\null

\noindent \textit{CASE 2.2: $Z$ is not an end-box of $P$.} Then the $H$-neighbours of $Z$ in $P$ must be $Z+(-1,0)$ and $Z+(0,1)$, and so we must have that $S_{\rightarrow}(k-1,1;k,2) \in H$. Once again, $W \neq Z+(-2,0)$, and since $Z+(-1,1)$ is not switchable, $W\neq Z+(-1,1)$. It follows that $Z+(-1,0) \in P$ is not an end-box and is switchable. See Figure 5.17 (b). End of Case 2.2. End of Case 2. $\square$

\null

\subsection{Algorithm for relocating an endpoint to the boundary (EtB)}

\noindent In this section, we prove Proposition 5.6. There is a lot of overlap between the proofs of (i) and (ii), so we will combine the proofs. Several times during the proof, we encounter subcases that require detailed analysis. We handle these separately as Claims 5.17--5.21, which can be read after the main proof if preferred. A flowchart (Flowchart 2) illustrating the proof structure appears at the end of this section.

\null 

\noindent \textit{Proof of Proposition 5.6.} Orient $H$ as $v=v_1, v_2, \ldots, v_r=u$, and recall that $e(v)$ is the edge of $H$ incident on $v$. We may assume WLOG that the eastern and western sides of $R_0$ are parity-compatible, and that $v_1=v(k,l)$ is an easternmost end-vertex of $H$.

Note that if one of $v$ or $u$ is incident on a side that is not parity-compatible, we can apply one backbite move that relocates that end-vertex to $G \setminus R_0$. So, for the proof of (i), we may assume that $v$ and $u$ are not incident on $R_0$, and neither is incident on the eastern or western side of $R_0$. In order to prove (i), it is enough to show that if neither of $u$ and $v$ is incident on $R_0$, then there is a cascade after which one of them lies closer to the eastern side of $R_0$ than either did prior to the cascade. We call this statement $(\dag_1)$ for reference. 

Similarly, for the proof of (ii), we may assume that $u$ is on the western side of $R_0$, and that $v$ is not on $R_0$. To prove (ii), it is enough to show that there is a cascade that fixes $u$, after which $v$ has moved closer to the eastern side of $R_0$. We call this statement $(\dag_2)$ for reference.

First we check that the case where the end-vertices are adjacent reduces to the case where (i) holds and end-vertices are not adjacent. Suppose the end-vertices of $H$ are adjacent. By Lemma 5.10, we can reconfigure $H$ into an e-cycle with one end-vertex $u$ on the western side of $R_0$, so (i) holds. Now, it is easy to check that we can apply a backbite move that keeps $u$ fixed and moves $v$ to a vertex not adjacent to $u$, lying strictly closer than $u$ to the eastern side of $R_0$. Therefore, from here on we assume that the end-vertices of $H$ are not adjacent. Furthermore, if $u$ and $v$ are equally easternmost, we assume WLOG that $v$ is north of $u$. Now, either $e(v)$ is eastern, or it is not.

\null 

\noindent \textit{CASE 1: $e(v)$ is not eastern.} We check that $(\dag_1)$  and $(\dag_2)$ hold. Since $e(v)$ is not eastern, $e(k,k+1;l) \notin H$. Let $v+(1,0)=v_s$. Then $v_{s-1}=v_1+(1,1)$, $v_{s-1}=v_1+(2,0)$ or $v_{s-1}=v_1+(1,1)$. See Figure 5.21. Then, after $bb_{v}(east)$, by Lemma 1.4.4, the end-vertex $v'$ of the resulting Hamiltonian path is located on the vertex of $G$ where $v_{s-1}$ was located in $H$. Thus, $v'$ is at least one unit closer to the eastern boundary, and so $(\dag_1)$  and $(\dag_2)$ are satisfied. See Figure 5.18 (a). End of Case 1. 

\endgroup

\begingroup
\setlength{\intextsep}{0pt}
\setlength{\columnsep}{20pt}
\begin{center}
\begin{adjustbox}{trim=0cm 0cm 0cm 0cm} 

\end{adjustbox}
\end{center}

\noindent \textit{CASE 2: $e(v)$ is eastern.} Then, either at least one of $e(v)+(0,1)$ and $e(v)+(0,-1)$ belongs to $H$, or neither does. 

\null

\noindent \textit{CASE 2.1: At least one of $e(v)+(0,1)$ and $e(v)+(0,-1)$ belongs to $H$.} We check that $(\dag_1)$ and $(\dag_2)$ hold. By symmetry, we may assume WLOG that $e(v)+(0,1)\in H$. Orient $H$ so that $v=v_1$ and let $v_s=v_1+(0,1)$. Then $v(k+1,l+1)=v_{s-1}$ (Figure 5.18 (b)) or $v(k+1,l+1)=v_{s+1}$  (Figure 5.18 (c)) . If $v(k+1,l+1)=v_{s+1}$, then $e(k+1;l,l+1)\notin H$. It follows that $R(k,l)$ is parallel. Apply Sw($R(k,l)$). Now the edge incident on $v$ in the resulting Hamiltonian path is northern, and we are back to Case 1. 

Assume that $v(k+1,l+1)=v_{s-1}$. Then, after $bb_{v}(north)$, $v$ is relocated to $v(k+1,l+1)$, and so  $(\dag_1)$ and $(\dag_2)$ are satisfied. End of Case 2.1

\endgroup 

\null

\noindent \textit{CASE 2.2: Neither $e(v)+(0,1)$ nor $e(v)+(0,-1)$ belongs to $H$.} Note that if $l+1=n-1$, then the other end-vertex of $H$ must be on $v(k,l+1)$. Since this conflicts with NAA, we may assume that $n-1 \geq l+2$, and by symmetry, that $l-2 \geq 0$. See Figure 5.19 (a). Also, if $m-1=k+1$, then the other end-vertex of $H$ must be on $v(k+1,l+1)$ or on $v(k+1,l-1)$. See Figure 5.19 (b) for an illustration. Since both possibilities conflict with the assumption that $v$ is an easternmost end-vertex, we must have that $m-1 \geq k+2$. Since $v \notin R_0$, we may assume that $ k-1 \geq 0$. See Figure 5.19 (c). Then, by NAA, $S_{\rightarrow}(k-1,l-1; k,l-2) \in H$ and $S_{\rightarrow}(k-1,l+1; k,l+2) \in H$. By the assumption that $v$ is an easternmost end-vertex of $H$, $v(k+1,l+1)$ and $v(k+1,l-1)$ are not end-vertices of $H$. Now, either $e(k+1;l,l+1) \in H$, or $e(k+1;l,l+1) \notin H$. If the former, then we must have that $S_{\uparrow}(k+1,l-2; k+2,l-1) \in H$; and if the latter, then we must have that $S_{\downarrow}(k+1,l+2; k+2,l+1) \in H$. Either way, $e(v)$ is followed by an $A$-type southward or northward. See Figures 5.23 (d) and (e).

\begingroup
\setlength{\intextsep}{0pt}
\setlength{\columnsep}{20pt}
\begin{center}
\begin{adjustbox}{trim=0cm 0cm 0cm 0cm} 

\end{adjustbox}
\end{center}

\endgroup

\noindent WLOG we may assume that $e(v)$ is followed by an $A$-type southward (as in Figure 5.19 (d)). Note that if $e(k-1,k;l-2) \in H$, then, after $R(k,l-2) \mapsto R(k-1,l-2)$, we are back to Case 2.1, so we may assume that $e(k-1,k;l-2) \notin H$. By symmetry, we may also assume that $e(k+1,k+2;l-2) \notin H$. Note that if $v(k,l-2)=u$ is the other end-vertex of $H$, then we are back to Case 1. Therefore we may assume that $v(k,l-2)$ is not the other end-vertex of $H$. Thus we have that $e(v)$ is followed southward by an $A_0$-type or by an $A_1$-type.

\null 

\noindent By Claim 5.17, we may assume that $m-1 \geq k+3$ and $0\leq k-2$. Now, there are three possibilities \footnote{We added the abbreviations as this trio comes up several times in later cases.}: 

(a) $e(v)$ is followed southward by a $j$-stack of $A_0$s, not followed by an $A_1$-type 

\hspace{0.6cm} $(e(v), js\text{-}A_0, \neg A_1, \text{South})$,

(b) $e(v)$ is followed southward by a $j$-stack of $A_0$s followed by an $A_1$-type 

\hspace{0.6cm} $(e(v), js\text{-}A_0, A_1,\text{South})$,

(c) $e(v)$ is followed southward by a full stack of $A_0$s $(e(v), fs\text{-}A_0, \text{South})$.

\endgroup

\null

\noindent \textit{CASE 2.2(a): $(e(v), js\text{-}A_0, \neg A_1, \text{South})$.} Let $R(k,b)$ be the southern leaf of the last ($j^{\text{th}}$) $A_0$-type in the $j$-stack of $A_0$'s that follows $e(v)$. 

If the other end-vertex $u$ of $H$ is not incident on the rectangle $R(k-1,k+2;b-2,b-1)$, then, by Lemma 5.11, there is a cascade of flips, after which $e(k,k+1;l-1)$ is in the resulting Hamiltonian path, and we are back to Case 2.1. Therefore we may assume that $u$ is incident on the rectangle $R(k-1,k+2;b-2,b-1)$.

\begingroup
\setlength{\intextsep}{0pt} %this shifts tikz up-down
\setlength{\columnsep}{20pt}
\begin{center}
\begin{adjustbox}{trim=0cm 0cm 0cm 0cm}

\end{adjustbox}
\end{center}

\noindent If $e(k,k+1;b-1) \in H$, then we can use Case 1 of Lemma 5.11 to return to Case 2.1 again, so we may assume that $e(k,k+1;b-1) \notin H$. Now, the assumption that $v$ is an eastmost end-vertex implies that $S_{\uparrow}(k+1,b-2;k+2,b-1) \in H$. 

\noindent If $v(k,b-1)$ is the other end-vertex $u$ of $H$, then we can apply $bb_{u}(east)$, after which $(\dag_1)$ and $(\dag_2)$ are satisfied (Figure 5.20 (a)). So, we may assume that $v(k,b-1)$ is not the other end-vertex of $H$. Then  $S_{\rightarrow}(k-1,b-1;k,b-2) \in H$. 

Now there are three possibilities: $v(k,b-2)=u$, $e(k,k+1;b-2) \in H$, or  $e(k-1,k;b-2) \in H$. If the first, then we can apply $bb_{u}(east)$, after which $(\dag_1)$ and $(\dag_2)$ are satisfied.  (Figure 5.20 (b)). If the second, then we can apply  then $R(k,b-2) \mapsto R(k-1,b-2)$, and then we can use Lemma 5.11 again to return to Case 2.1 (Figure 5.20 (b)). The third is not possible, as it would mean that the length of the $j$-stack of $A_0$'s is $j+1$, contradicting the assumption that the $j$-stack of $A_0$'s has length $j$ (marked off in orange in Figure 5.20 (c)). End of Case 2.2(a).

\endgroup 

\null

\noindent \textit{CASE 2.2(b):} $(e(v), js\text{-}A_0, A_1,\text{South})$. Let $e_c$ be the edge in $G \setminus H$ between the two corners of the $A_1$-type, and let $W=R(k,l-2j-2)$. Note that $W$ is the switchable box of the $A_1$-type following the $j$-stack of $A_0$'s. If $W$ is parallel, then, after $\text{Sw}(W)$, we are back to Case 2.2(a), so we may assume that $W$ is anti-parallel. Let $X$ and $Y$ be the boxes adjacent to $W$ that are not its $H$-neighbours. Note that Claim 5.17, the fact that the $A_1$-type has edges south of $W$, and the fact that the $A_1$-type is south of $e(v)$, imply that $X$ and $Y$ are not incident on $R_0$. Now, either $X$ and $Y$ belong to different $H$-components, or they belong to the same $H$-component.

\null

\noindent If $X$ and $Y$ belong to distinct $H$-components, then, by Lemma 5.15 (a), there is a non-backbite cascade, after which $W$ is switched, and we are back to Case 2.2(a).
%Then, by Case 2.2(a), either there is a non-backbite cascade after which $e(k,k+1;l-1)$ is in the resulting Hamiltonian path, we are back to Case 2.1, or, $(\dag_1)$ and $(\dag_2)$ are satisfied. 
Therefore we may assume that $X$ and $Y$ belong to the same $H$-component. In this case, $P(X,Y)$ is unique. Let $F=G\langle N[P(X,Y)] \rangle$. Then either $F$ is a sturdy looping fat path or it is not. 

If $P(X,Y)$ has a switchable box $Z$, then, by Proposition 4.5, either $\text{Sw}(Z)$, $\text{Sw}(W)$ is a cascade or $Z \mapsto W$ is a valid move, and we are back to Case 2.2(a).  If $F$ fails FPC-4, apply the available switch or flip or transpose move. If now, one of $W \mapsto X$, or $W \mapsto Y$ is valid, after applying that move, we are back to Case 2.2(a); and if $P(X,Y)$ has gained a switchable box $Z$, we can use  Proposition 4.5 again to return to Case 2.2(a). Therefore we may assume that $F$ satisfies FPC-2,3 and 4. Now, either $F$ satisfies FPC-1 or it does not.

\null

\noindent \textit{CASE 2.2$(b_1)$: $F$ does not satisfy FPC-1.} Recall that $u \neq v$ is the other end-vertex of $H$. Then, at least one of $u$ and $v$ is incident on $P(X,Y)$. Then, either $u$ is incident on $P(X,Y)$, or it is not.

\null

\noindent \textit{CASE 2.2$(b_1).(i)$: $u$ is incident on $P(X,Y)$.}  Either $u$ is on $R_0$ or it is not. If $u$ is in $R_0$, then by Lemma 5.16, $P(X,Y)$ has a switchable box, contradicting the assumption that $F$ satisfies FPC-3.

Assume that $u$ is not in $R_0$. Then, by Claim 5.18, either there is a cascade after which $(\dag_1)$ and $(\dag_2)$ hold, or there is a cascade that fixes $v$, after which we are back to Case 2.2(a).  End of Case 2.2$(b_1).(i)$.

\null

\noindent \textit{CASE 2.2$(b_1).(ii)$: $u$ is not incident on $P(X,Y)$} Then $v$ must be incident on $P(X,Y)$. Then, by Claim 5.20, $(\dag_1)$ and $(\dag_2)$ hold. End of Case 2.2$(b_1).(ii)$.  End of Case 2.2$(b_1)$

\null

\noindent \textit{CASE 2.2$(b_2)$: $F$ satisfies FPC-1.} This means that $F$ is a \index{sturdy looping fat path}sturdy looping fat path. By Corollary 4.12(a), $F$ is a standard looping fat path. 
%Note that it is sufficient to find a cascade after which, either $W$ is switched, or $e_c$ is in the resulting Hamiltonian path.
By Corollary 4.22, $F$ has an admissible \index{turn}turn $T$ such that $\text{Sector}(T)$ avoids the $j$-stack of $A_0$'s and the boxes incident on $e(v)$, and is below the lines $y=x+l-2j-k-2$. Note that this implies that $\text{Sector}(T)$ avoids $v$ as well. By Corollary 4.19(a), $d(T) \geq 3$. Now, either $u \in R_0$ or $u \notin R_0$.

\null

\noindent \textit{CASE 2.2$(b_2).(i)$: $u \in R_0$.} By Proposition 4.27, $T$ has a double-switch weakening $\mu_1, \ldots, \mu_s$. By Corollary 4.31(i), the weakening $\mu_1, \ldots, \mu_s$ fixes $v$, and is contained in Sector($T$). By Lemma 4.32, $\mu_s$ can be followed by a non-backbite cascade of length at most two, after which $W$ is switched. End of Case 2.2$(b_2). (i)$. 

\null

\noindent \textit{CASE 2.2$(b_2).(ii)$: $u \notin R_0$.}  By Proposition 4.28, $T$ has a weakening $\mu_1, \ldots, \mu_s$.  By Corollary 4.31(i), $\mu_1, \ldots, \mu_s$ is contained in Sector($T$).

If there is a $q$ in $\{1, \ldots, s-1\}$ such that after $\mu_{q}$, $X$ and $Y$ belong to distinct $H$-components, then we can use Lemma 5.15(a) to extend $\mu_1, \ldots,\mu_{q}$ into a cascade of length at most $q+2$ that switches $W$ and fixes $v$. So, we may assume that no such $q$ exists. If there is a $p$ in $\{1, \ldots, s-1\}$ such that after $\mu_{p}$, $u$ is incident on $P(X,Y)$, then we can use Claim 5.18 to extend $\mu_1, \ldots,\mu_{p}$ into a cascade of length at most $p+2$, after which we are back to Case 2.2(a), or into a cascade of length at most $p+3$, after which $(\dag_1)$ and $(\dag_2)$ hold. Thus, we may assume that no such $p$ exists either. 

Let $H'$ be the resulting Hamiltonian path after $\mu_s$. Since $v$ and $\text{Sector}(T)$ are disjoint, the cascade $\mu_1, \ldots, \mu_s$ fixes $v$. Let $u'$ be the other end-vertex of $H'$. Since $\mu_1, \ldots, \mu_s$ is a weakening of $T$, after applying $\mu_s$, either one of the terminal edges of $T$ is in $H'$, or $u'$ is in ew-Terminal$(T)$. If the former, then by Lemma 4.32, $\mu_s$ can be followed by a non-backbite cascade of length at most two, after which $W$ is switched; and if the latter, by Observation 4.20(a), $u'$ is incident on a box of $P(X,Y)$. Then, we can use Claim 5.18 again to extend $\mu_1, \ldots, \mu_s$ into a cascade of length at most $s+2$ after which we are back to Case 2.2(a), or a cascade of length at most $s+3$ after which $(\dag_1)$ and $(\dag_2)$ hold. End of Case 2.2$(b_2).(ii)$. End of Case 2.2$(b_2)$. End of Case 2.2(b). End of Case 2.2.

\null 

\noindent \textit{CASE 2.2(c): $(e(v), fs\text{-}A_0, \text{South})$.} Note that if $l$ is odd, by NAA, $e(k,k+1;0) \in H$, and we are back to Case 2.2(a), so we may assume that $l$ is even. We have that $S_{\downarrow}(k+2,1;k+3,0) \in H$. Note that  $l=2j$, and $j \geq 1$. If $j \geq 2$, then for each $i \in \{0,2, \ldots,l-4\}$, $S_{\downarrow}(k+2,i+1;k+3,i) \in H$ implies $S_{\downarrow}(k+2,i+3;k+3,i+2) \in H$. Thus we have that for each $i \in \{0,2, \ldots,l-2\}$, $S_{\downarrow}(k+2,i+1;k+3,i) \in H$. See Figure 5.26 (a). Now, $e(k+1,k+2;l)\in H$ or $e(k+1;l,l+1)\in H$.

\begingroup
\setlength{\intextsep}{0pt} 
\setlength{\columnsep}{20pt}
\begin{center}
\begin{adjustbox}{trim=0cm 0cm 0cm 0cm}

\end{adjustbox}
\end{center}

\noindent {CASE 2.2($c_1$): $e(k+1,k+2;l)\in H$.} Then $S_{\downarrow}(k+1,l+2;k+2,l+1) \in H$. There are three possibilities: 

1. $(e(v), js\text{-}A_0, \neg A_1, \text{North})$,

2. $(e(v), js\text{-}A_0, A_1, \text{North})$,

3. $(e(v), fs\text{-}A_0, \text{North})$.

\noindent If the first, then, by symmetry, this is the same as Case 2.2(a). If the second, then, by symmetry, this is the same as Case 2.2(b). Assume the third. If $n-1$ is odd, we have that $e(k,k+1;n-1) \in H$. Then, after the cascade of flips, $R(k,n-2) \mapsto R(k,n-3), \ldots,R(k,l+2) \mapsto R(k,l+1)$, we are back to Case 2.1. And if $n-1$ is even, this is impossible by Lemma 5.12. See Figure 5.26 (b). End of Case 2.2($c_1$). 

\endgroup 

\null

\noindent {CASE 2.2$(c_2)$: $e(k+1;l,l+1)\in H$.} 
Then we have that $S_{\downarrow}(k+2,l+1;k+3,l) \in H$. Again, there are three possibilities:

1. $(e(k+1,k+2;l-1), js\text{-}A_0, \neg A_1, \text{North})$,

2. $(e(k+1,k+2;l-1), js\text{-}A_0, A_1, \text{North})$,

3. $(e(k+1,k+2;l-1), fs\text{-}A_0, \text{North})$.

\noindent Assume the first. By the assumption that $v$ is north of $u$, the rectangle $R(k,k+3;l,n-1)$ does not contain $u$. Then, by Lemma 5.11 and symmetry, there is a cascade of flips after which $e(k+1,k+2;l) \in H$, and then we are back to Case 2.2($c_1$). If the third, then $n-1$ is even or odd. If $n-1$ is even, then we must have that $e(k+1,k+2;n-1) \in H$. Then, after the cascade of flips $R(k+1,n-2) \mapsto R(k+1,n-3), \ldots,R(k+1,l+1) \mapsto R(k+1,l)$, we are back to Case 2.2($c_1$). By Lemma 5.13, the case where $n-1$ is odd is not a possible configuration. 

\null 

\noindent Assume the second. Let $W=R(k+1,l+2j)$. Note that $W$ is the switchable box of the $A_1$-type following the $j$-stack of $A_0$'s. Let $X$ and $Y$ be the boxes adjacent to $W$ that are not its $H$-neighbours. 

%pagemarker

\noindent The fact that the $A_1$-type has edges located south of $W$, and the fact that the $A_1$-type is north of $e(v)$, imply that $X$ and $Y$ are not incident on the southern side or the northern side $R_0$. Since $k>0$, they are not incident on the eastern side of $R_0$. An argument like the one in Cases I and II of Lemma 5.12 can be used to show that $m-1 \geq k+4$, and thus $X$ and $Y$ cannot be incident on the western side of $R_0$.

\null 

\noindent Now, either $X$ and $Y$ belong to distinct $H$-components, or they belong to the same $H$-component. If $X$ and $Y$ belong to distinct $H$-components, then, by Lemma 5.15 (a) and (b), there is a non-backbite cascade, after which $W$ is switched. Then, by Lemma 5.11, there is a cascade of flips after which $e(k+1,k+2;l)$ is in the resulting Hamiltonian path, and we are back to Case 2.2($c_1$). Therefore we may assume that $X$ and $Y$ belong to the same $H$-component. In this case, $P(X,Y)$ is unique. Let $F=G\langle N[P(X,Y)] \rangle$.

If $W$ is parallel (that is, if $F$ fails FPC-2), then, after $\text{Sw}(W)$, by Case 2.2(a) and symmetry, there is a cascade of flips after which we are back to Case 2.2($c_1$), so we may assume that FPC-2 holds. If $F$ fails FPC-3, then $P(X,Y)$ has a switchable box $Z$. By Proposition 4.5, either $\text{Sw}(Z)$, $\text{Sw}(W)$ is a cascade or $Z \mapsto W$ is a valid move. Then, by Lemma 5.11 and symmetry, there is a cascade of flips, after which we are back to Case 2.2$(c_1)$. If $F$ fails FPC-4, apply the available switch or flip or transpose move. If now one of $W \mapsto X$ or $W \mapsto Y$ is valid, we apply it. Then, by Lemma 5.11 and symmetry, there is a cascade of flips, after which we are back to Case 2.2$(c_1)$. If otherwise $P(X,Y)$ has gained a switchable box $Z$, we can use  Proposition 4.5 and Lemma 5.11 again to return to Case 2.2$(c_1)$. Therefore we may assume that $F$ satisfies FPC-2, 3 and 4. Now, either $F$ satisfies FPC-1 or it does not.

\null 

\noindent \textit{CASE 2.2$(c_2).(i)$: $F$ does not satisfy FPC-1.} Then there is an end-vertex of $H$ incident on $P(X,Y)$. Recall that $u \neq v$ is the other end-vertex of $H$. Then, at least one of $u$ and $v$ is incident on $P(X,Y)$. If $u$ is incident  on $P(X,Y)$, then we can use Claim 5.19, to show that either there is a cascade of length three after which $(\dag_1)$ and $(\dag_1)$ hold, or there is a cascade of length two, after which, $W$ is switched. If the latter, then we can use Lemma 5.11 and symmetry to return to Case 2.2$(c_1)$. Therefore we may assume that $v$ is incident on $P(X,Y)$ and that $u$ is not. Then, by Claim 5.21, $(\dag_1)$ and $(\dag_2)$ hold. End of 2.2$(c_2).(i)$

\null 

\noindent \textit{CASE 2.2$(c_2).(ii)$: $F$ satisfies FPC-1.} This means that $F$ is a \index{sturdy looping fat path}sturdy looping fat path. The case where $u \in R_0$ uses the same argument as Case 2.2$(b_2).(i)$. The case where $u \notin R_0$ uses the same argument as Case 2.2$(b_2).(ii)$. In each case, we either find a cascade after which we are back to Case 2.2$(c_1)$, or we find a cascade after which $(\dag_1)$ and $(\dag_2)$ hold. End of Case 2.2$(c_2).(ii)$. End of Case 2.2$(c_2)$. End of Case 2.2(c). End of Case 2.2. End of Case 2. $\square$

\null 

\noindent It remains to prove Claims 5.17--5.21. Claim 5.17 deals with boundary proximity checking. Claims 5.18--5.21 all handle similar scenarios: an end-vertex is incident on a looping $H$-path that otherwise satisfies FPC-2, FPC-3, and FPC-4. We postpone proving the bound on the length of the cascade required by the EtB algorithm until after we prove these claims.

%Claim  5.17.1 - that maybe will never be written.
{
%We check that there is a cascade after which $u$ is fixed and $v$ is not adjacent to $u$. Either  $e(1,2;n-1) \in H$ or $e(1,2;n-1) \notin H$.

%\textit{CASE 1: $e(1,2;n-1) \notin H$.} Let $v_s=v_{2,n-1}$. Note that Corollary 1.5.2 implies that $v_{s-1}=v(2,n-2)$. Then after $bb_v(east)$, the end-vertex $v$ is relocated to $v(2,n-2)$
}

% Extra detail
{%If $u \notin R_0$ ($u \in R_0$) then, using the same argument as in Case 2.2$(b_1).(i)$ (Case 2.2$(b_1).(ii))$, we find a cascade that fixes $v$, and avoids the $j$-stack of $A_0s$ and the boxes incident on $e(k+1,k+2;l+1)$, after which $W$ is switched \textcolor{red}{[Show more here? Essentially, by Corollary 4.22, the weakenings are contained in the region above the line $y=-x+4+2+l+3j$, avoiding anything relevant.]}. Then, by Lemma 5.11, there is a cascade, after which $e(k+1,k+2;l)$ is in the resulting Hamiltonian path, and we are back to Case 2.2$(c_1)$. End of Case 2.2$(c_2).(ii)$. End of Case 2.2$(c_2)$. End of Case 2.2(c). End of Case 2.2. End of Case 2. $\square$
}

\begingroup
\setlength{\intextsep}{0pt}
\setlength{\columnsep}{10pt}
\begin{center}
\begin{adjustbox}{trim=0cm 0cm 0cm 0cm} 

\end{adjustbox}
\end{center}

\noindent \textbf{Claim 5.17}. Recall that this is a continuation of Case 2.2 of the proof or Proposition 5.6, with $e(v)$ followed southward by an $A_0$ or an $A_1$. We will check that if $k+2=m-1$, then $(\dag_1)$ and $(\dag_2)$ are satisfied. Assume that $k+2=m-1$. Note that if $e(k+2;l-1,l) \in H$ (Figure 5.22 (a)), then $v(k+2;l-2)$ must be an end-vertex, which contradicts the assumption that $v$ is an easternmost end-vertex. So, we may assume that $e(k+2;l-1,l) \notin H$ (Figure 5.22 (b)). It follows that $e(k+2;l-2,l-1) \in H$. 

If $l-2=0$, then again, $v(k+2,l-2)$ must be an end-vertex, which contradicts the assumption that $v$ is an easternmost end-vertex, so we may assume that $l-3 \geq 0$ (Figure 5.22 (c)). Then $e(k+2;l-3,l-2) \in H$.

Note that if $l-3=0$, then we must have $e(k+1,k+2;l-3)\in H$, $e(k,k+1;l-3)\in H$, and $e(k,k+1;l-2)\in H$. But then, after $R(k,l-3) \mapsto R(k,l-2)$, we are back to Case 2.1. So we may assume that $l-4 \geq 0$.

\endgroup

\null 

\begingroup
\setlength{\intextsep}{0pt}
\setlength{\columnsep}{20pt}
\begin{center}
\begin{adjustbox}{trim=0cm 0cm 0cm 0cm} 

\end{adjustbox}
\end{center}

\noindent If $e(k+1;l-3,l-2) \in H$, then, after $R(k+1,l-3) \mapsto R(k,l-2)$, we are back to Case 2.1, so we may assume that $e(k+1;l-3,l-2) \notin H$. It follows that $e(k,k+1;l-2) \in H$. As before, if $e(k,k+1;l-3) \in H$, then, after $R(k,l-3) \mapsto R(k,l-2)$, we are back to Case 2.1, so we may assume that $e(k,k+1;l-3) \notin H$. Then $S_{\uparrow}(k+1,l-4;k+2,l-3) \in H$. Since $v$ is easternmost, there are no end-vertices east of the line $x=k$. Now, either $e(k;l-4,l-3) \in H$, or $e(k;l-4,l-3) \notin H$. If the former (Figure 5.23(a)), then we must have $e(k+1,k+2;l-4) \in H$ as well. But then there is $R(k,l-4) \mapsto R(k+1,l-4)$, followed by $R(k,l-3) \mapsto R(k,l-2)$, we are back to Case 2.1 once more. If the latter (Figure 5.23(a)) then $u=v(k,l-3)$, and we can apply $bb_u(east)$, after which $(\dag_1)$ and $(\dag_2)$ hold.

Therefore, we may assume that $m-1 \geq k+3$, and by symmetry, that $0 \leq k-2$. $\square$ 

\endgroup 

\null

\noindent \textbf{Claim 5.18.} This is a continuation of Case 2.2$(b_1).(i)$. It will also be used by Case 2.2$(b_2).(ii)$.

We restate the setup here. $e(v)$ is followed by a $j$-stack of $A_0$s followed by an $A_1$-type southward. $W=R(k,l-2j-2)$ is the switchable middle-box of the $A_1$-type, $X$ and $Y$ are the boxes adjacent to $W$ that are not its $H$-neighbours, and $P(X,Y)$ is contained in an $H$-component of $G$. $F=G\langle N[P(X,Y)] \rangle$ satisfies FPC-2,3, and 4, and fails FPC-1. We assume that $u$ is incident on $P(X,Y)$, and that both $u$ and $v$ have x-coordinate at most $k$. We will show that there is a cascade of length at most three after which $(\dag_1)$ and $(\dag_2)$ hold, or that there is a cascade of length at most two, after which we return to Case 2.2(a).

\null 

\noindent \textit{Proof.} Let $l'=l-2j-2$, let $v_y=v(k,l'+1)$ and $v_x=v(k+1,l'+1)$. Either $u$ is incident on a box that contains $v_x$ or $v_y$, or it is not.

\begingroup
\setlength{\intextsep}{0pt}
\setlength{\columnsep}{20pt}
\begin{center}
\begin{adjustbox}{trim=0cm 0cm 0cm 0cm} 

\end{adjustbox}
\end{center}

\noindent \textit{CASE 1: $u$ is incident on a box that contains $v_x$ or $v_y$.} Since the x-coordinate of $u$ is at most $k$, there are four possibilities: $u=v(k,l'+2)$, $u=v(k-1,l'+2)$, $u=v(k-1,l'+1)$, or $u=v(k-1,l')$. See Figure 4.24 (a).

\null 

\noindent \textit{CASE 1.1: $u=v(k,l'+2)$.} If $j=0$, $u=v(k,l)=v$, which is impossible; and if $j>0$, then we must have $S_{\downarrow}(k,l'+3;k+1;l'+2) \in H$, which is again impossible. So  $u \neq v(k,l'+2)$. End of Case 1.1.

\null

\noindent \textit{CASE 1.2: $u=v(k-1,l'+2)$.} $j=0$ contradicts NAA, so assume $j>0$. Then after $bb_{u}(east)$, $u$ is either relocated to $v(k+1;l'+2)$ or to $v(k;l'+3)$. If the former, then $\dag_1$ and $\dag_2$ hold; and if the latter, then we can follow the first move with $bb_{v(k,l'+3)}(east)$, after which $(\dag_1)$ and $(\dag_2)$ hold.  See Figure 4.24 (b). End of Case 1.2.

\begin{center}
\begin{adjustbox}{trim=0.25cm 0cm 0cm 0cm} 

\end{adjustbox}
    
\end{center}

\noindent \textit{CASE 1.3: $u=v(k-1,l'+1)$.} Then we must have $S_{\rightarrow}(k-2,l';k-1,l'-1) \in H$. Let $v_q=v(k-1,l')$ and $v_{q'}=v(k-1,l'-1)$. Either $q'<q$ or $q'>q$.

\null 

\noindent \textit{CASE 1.3(a): $q'<q$.} It follows that $W+(-1,-1)$ is parallel. Then, after $\textrm{Sw}(W+(-1,-1))$, $W \mapsto W+(0-1,0)$, we are back to Case 2.2(a) of Proposition 5.6. See Figure 5.25 (a) End of Case 1.3(a).

\null 

\noindent \textit{CASE 1.3(b): $q'>q$.} Then, by Lemma 1.4.4, $bb_u(south)$ is the move $e(k-1;l'-1,l') \mapsto e(k-1;l',l'+1)$. By Proposition 4.4, we may partition $H$ into the subpaths $P_0=P(v_x,v_y)$, $P_1=P(v_1,v_x)$ and $P_2=P(v_y,v_r)$ such that every box of $P(X,Y)$ is incident on a vertex of $P_0$ and a vertex of $P_1$ or $P_2$. Now, either $v_{q'} \in P_0$, or $v_{q'}\notin P_0$.

\null 

\noindent \textit{CASE 1.3($b_1$): $v_{q'} \in P_0$.} This means that the order of vertices of $H$ is $v_1, \ldots v_x, \ldots, v_q, \ldots, v_y,$ $u$. Then, by Corollary 1.4.5, after $bb_u(south)$, $W$ becomes parallel, so we can apply $\textrm{Sw}(W)$, and return to Case 2.2(a) of Proposition 5.6. See Figure 5.25 (b). End of Case 1.3($b_1$).

\null 

\noindent \textit{CASE 1.3($b_2$): $v_{q'} \notin P_0$.} Then $v_{q'}$ must belong to $P_1$. Note that this implies that $e(k-1,k;l'-1) \notin H$, and, by Lemma 1.4.4, that $bb_u(south)$ fixes $P_0$. Then, after $bb_u(south)$, $bb_{v(k-1,l'-1)}(east)$ is the move $e(k;l'-1,l')\mapsto e(k-1,k;l'-1)$, and we apply it. Now, after $bb_{v(k,l')}(east)$, the end-vertex is relocated to $v(k+1,l'-1)$, so $(\dag_1)$ and $(\dag_2)$ hold. See Figure 5.26. Note that we have used Lemma 1.4.4 for the previous two moves. End of Case 1.3($b_2$). End of Case 1.3(b). End of Case 1.3.

\begingroup
\setlength{\intextsep}{0pt}
\setlength{\columnsep}{20pt}

\begin{center}
\begin{adjustbox}{trim=0.25cm 0cm 0cm 0cm} 

\end{adjustbox}
\end{center}
\endgroup

\noindent \textit{CASE 1.4: $u=v(k-1,l')$.} Let $v_q=v(k,l')$ and $v_{q'}=v(k,l'-1)$. Then $q<q'$ or $q>q'$.

\null 

\noindent \textit{CASE 1.4(a): $q<q'$.} By Lemma 1.4.4, $bb_u(east)$ is the backbite move $e(k;l'-1,l') \mapsto e(k-1,k;l')$, and we apply it. Then, after $W \mapsto W+(-1,0)$, we are back to Case 2.2(a) of Proposition 1.17. See Figure 5.27 (a). End of Case 1.4(a).

\begin{center}
\begin{adjustbox}{trim=0.25cm 0cm 0cm 0cm} 
%
\end{adjustbox}
\end{center}

\noindent \textit{CASE 1.4(b): $q>q'$.} By Lemma 1.4.4, $bb_u(east)$ is the backbite move $e(k;l',l'+1) \mapsto e(k-1,k;l')$, and we apply it. Now, after $bb_{v(k,l'+1)}(east)$, the end-vertex is relocated to $v(k+1,l')$, or to $v(k+2;l'+1)$. Either way, $(\dag_1)$ and $(\dag_2)$ hold.  See Figure 5.27 (b). End of Case 1.4(b). End of Case 1.4. End of Case 1.

\null

\noindent \textit{CASE 2: $u$ is not incident on a box that contains $v_x$ or $v_y$.} Let $Z$ be the box of $P(X,Y)$ on which $u$ is incident.  By Proposition 4.4, we may partition $H$ into the subpaths $P_0=P(v_x,v_y)$, $P_1=P(v_1,v_x)$ and $P_2=P(v_y,v_r)$ such that every box of $P(X,Y)$ is incident on a vertex of $P_0$ and a vertex of $P_1$ or $P_2$. Let $v_s$ be the vertex of $Z$ that belongs to $P_0$. WLOG, assume that $Z=R(a,b)$, and $u=v(a+1,b)$. Either $v_s$ is adjacent to $u$, or it is not.

\null 

\noindent \textit{CASE 2.1: $v_s$ is adjacent to $u$.} For definiteness, assume that $v_s=v(a,b)$. Note that the edge$\{u,v_s\}$ must belong to $G \setminus H$, and that the order of vertices of $H$ is $v_1, \ldots v_x, \ldots, v_s, \ldots , v_y,$ $u$. Then, by Corollary 1.4.5, after $bb_u(west)$, $W$ becomes parallel, so we can apply $\textrm{Sw}(W)$, and return to Case 2.2(a) of Proposition 5.6. End of Case 2.1.

\null

\begingroup
\setlength{\intextsep}{0pt}
\setlength{\columnsep}{20pt}
\begin{wrapfigure}[]{l}{0cm}
\begin{adjustbox}{trim=0cm 0cm 0cm 0cm}
\begin{tikzpicture}[scale=1.75]
\usetikzlibrary{decorations.markings}
\begin{scope}[xshift=0cm, yshift=0cm] 

\draw[gray,very thin, step=0.5cm, opacity=0.5] (0,0) grid (1.5,1.5);

\draw[blue, line width=0.5mm] (0,1)--++(0.5,0)--++(0,0.5);
\draw[orange!90!black, line width=0.5mm] (1.5,1)--++(-0.5,0)--++(0,0.5);

\draw[fill=orange, opacity=1] (0.5,0.5) circle [radius=0.035];

\draw[fill=orange, opacity=1] (1.5,1) circle [radius=0.035];

\draw[fill=orange, opacity=1] (1,0.5) circle [radius=0.05];
\node[right] at  (1,0.5) [scale=0.8]{\small{$u$}};

\draw[fill=orange, opacity=1] (1,1) circle [radius=0.035];

\draw[fill=blue, opacity=1] (0.5,1) circle [radius=0.035];
\node[above] at  (0.4,1) [scale=0.8]{\small{$v_s$}};

% black lines
{
\draw[black, line width=0.15mm] (0.45,0.70)--++(0.1,0);
\draw[black, line width=0.15mm] (0.45,0.75)--++(0.1,0);
\draw[black, line width=0.15mm] (0.45,0.80)--++(0.1,0);

\draw[black, line width=0.15mm] (0.95,0.70)--++(0.1,0);
\draw[black, line width=0.15mm] (0.95,0.75)--++(0.1,0);
\draw[black, line width=0.15mm] (0.95,0.80)--++(0.1,0);

\draw[black, line width=0.15mm] (0.75,0.95)--++(0,0.1);
\draw[black, line width=0.15mm] (0.7,0.95)--++(0,0.1);
\draw[black, line width=0.15mm] (0.8,0.95)--++(0,0.1);
}

\node at (0.75,0.75) [scale=1]{\small{$Z$}};

\node[left] at (0,0.5) [scale=1]{\tiny{$b$}};
\node[left] at (0,1) [scale=1]{\tiny{$+1$}};

\node[above] at (0.5,1.5) [scale=1]{\tiny{$a$}};
\node[above] at (1,1.5) [scale=1]{\tiny{$+1$}};
\node[above] at (1.5,1.5) [scale=1]{\tiny{$+2$}};

\node[below, align=center, text width=4cm] at (0.75, -0.1) { Fig. 5.28.  Case 2.2. };

\end{scope}
\end{tikzpicture}
\end{adjustbox}
\end{wrapfigure}

\noindent \textit{CASE 2.2: $v_s$ is not adjacent to $u$.} For definiteness, assume that $v_s=v(a,b+1)$. If $v(a,b) \in P_0$ or $v(a+1,b+1) \in P_0$ then we are back to Case 1, so we may assume that $v(a,b)$ and $v(a+1,b+1)$ belong to $P_1$ or $P_2$, and by NAA, that neither vertex is $v_1$. By the premise of Case 2, the edges $e(a;b,b+1) \notin H$, $e(a,a+1;b+1) \notin H$. This implies that $S_{\rightarrow}(a-1,b+1;a,b+2) \in H$. Since $u$ is an end-vertex, at least one of $e(a,a+1;b)$ and $e(a+1;b,b+1)$ does not belong to $H$. By symmetry, we may assume WLOG that $e(a+1;b,b+1) \notin H$. Then we must have $S_{\downarrow}(a+1,b+2;a+2,b+1) \in H$, so $v(a+2,b+1)$ is in $P_1$ or $P_2$. See Figure 5.28. Now, either $e(a,a+1;b) \in H$ or $e(a,a+1;b) \notin H$.

\null 

\endgroup

\noindent \textit{CASE 2.2(a): $e(a,a+1;b) \in H$.} Either  $Z+(1,0) \in P(X,Y)$ or $Z+(1,0) \notin P(X,Y)$. 

\null 

\begingroup
\setlength{\intextsep}{0pt}
\setlength{\columnsep}{20pt}
\begin{wrapfigure}[]{r}{0cm}
\begin{adjustbox}{trim=0cm 0cm 0cm 0.25cm}

\end{adjustbox}
\end{wrapfigure}

\noindent \textit{CASE 2.2($a_1$): $Z+(1,0) \in P(X,Y)$.} Now, either $v(a+2;b) \in P_0$ (Figure 5.29 (a)), or $v(a+2;b) \notin P_0$ (Figure 5.29 (b)). If $v(a+2,b) \in P_0$, then we are back to Case 2.1. Assume that $v(a+2,b) \notin P_0$. Then, by Lemma 4.1, we must have that $v(a+2,b+1)\in P_0$. But this implies that $v(a+2,b+1)=v_x$ or $v(a+2,b+1)=v_y$, contradicting the premise of Case 2. End of Case 2.2($a_1$).

\endgroup 

\null

\begingroup
\setlength{\intextsep}{0pt}
\setlength{\columnsep}{20pt}
\begin{wrapfigure}[]{l}{0cm}
\begin{adjustbox}{trim=0cm 0cm 0cm 0cm}

\end{adjustbox}
\end{wrapfigure}

\noindent \textit{CASE 2.2($a_2$): $Z+(1,0) \notin P(X,Y)$.} Then $Z$ must have exactly two $H$-neighbours in $P(X,Y)$. In particular, $Z+(0,1) \in P(X,Y)$. Note that $Z+(0,1)$ can not be $X$ or $Y$. Therefore, it must also have exactly two $H$-neighbours in $P(X,Y)$. This implies that $e(a,a+1;b+2)\notin H$. But then $Z+(0,1) \in P$ is switchable, which contradicts that $F$ satisfies FPC-3. See Figure 5.30 (a). End of Case 2.2($a_2$). End of Case 2.2(a)

\null

\noindent \textit{CASE 2.2(b): $e(a,a+1;b) \notin H$.} This implies that $S_{\rightarrow} (a-1,b;a,b-1) \in H$. Now, exactly one of $e(a+1;b-1,b)$ and $e(a+1,a+2;b)$ belongs to $H$. By symmetry, we may assume without loss of generality that $e(a+1;b-1,b) \in H$. If $Z+(1,0) \in P$, then we're back to Case 2.2($a_1$), so we may assume that $Z+(1,0) \notin P$. Note that $Z$ cannot be $X$ or $Y$, so it must have at least two $H$ neighbours in $P$. Then, at least one of them must be $Z+(-1,0)$ or $Z+(0,1)$, and either possibility brings us back to Case 2.2($a_2$). See Figure 5.30 (b). End of Case 2.2(b). End of Case 2.2. End of Case 2. $\square$ 

\null

\noindent \textbf{Claim 5.19.} This is the version of Claim 5.18 suitable for use in Case 2.2$(c_2)$ of Proposition 5.6. The setup and cases here are very similar to the setup in Claim 5.18. We begin with the scenario $((k+1,k+2;l-1), js\text{-}A_0, A_1, \text{North})$, and the rest of the starting assumptions are the same. Case 2 is the same as in Claim 5.18. 

Case 1 cannot occur: The assumption that the x-coordinates of the end-vertices are at most $k$, precludes Case 1.1. The assumption that $v$ is located north of $u$ precludes Cases 1.2, 1.3, and 1.4. $\square$

\null

\noindent  \textbf{Claim 5.20.} Recall that this is a continuation of Case 2.2$(b_1).(ii)$ of Proposition 5.6. We are in the scenario where $(e(v), js\text{-}A_0, A_1,\text{South})$, FPC-2,3, and 4 hold, $u$ is not incident on $P(X,Y)$, but $v$ is. We will show that in this case, $(\dag_1)$ and $(\dag_2)$ hold.

\null 

\noindent \textit{Proof.} For definiteness, assume that $X=W+(1,0)$, $Y=W+(-1,0)$ and let $v_x$ and $v_y$ be the corners of the $A_1$-type incident on $X$ and $Y$, respectively. Consider the $j$-stack of $A_0$'s following $e(v)$. Either $j=0$, or $j>0$.

\null

\begingroup
\setlength{\intextsep}{0pt}
\setlength{\columnsep}{20pt}
\begin{wrapfigure}[]{r}{0cm}
\begin{adjustbox}{trim=0cm 0cm 0cm 0cm}

\end{adjustbox}
\end{wrapfigure}

\noindent \textit{CASE 1: $j=0$.} Then $W=R(k,l-2)$. Let $v_x=v(k+1, l-1)$, and $v_y=v(k, l-1)$ be the corners of the $A_1$-type of $P(X,Y)$. Now, either $v_{y-1}=v(k,l-2)$ (Figure 5.31), or $v_{y-1}=v(k-1,l-1)$. If the former, then $bb_{v}(south)$ must be the move $(v_{y-1}, v_y) \mapsto \{v_y,v_1\}$, and the we are back to Case 1 of Proposition 5.6. We will show that the latter is not possible.

For contradiction, assume that  $v_{y-1}=v(k-1,l-1)$. This means that $v_{x+1}=v(k+2,l-1)$. Then $x<y$, or $x>y$.

\endgroup 

\null 

\begingroup
\setlength{\intextsep}{0pt}
\setlength{\columnsep}{20pt}
\begin{wrapfigure}[]{l}{0cm}
\begin{adjustbox}{trim=0cm 0cm 0cm 0cm}

\end{adjustbox}
\end{wrapfigure}

\noindent \textit{CASE 1.1: $x<y$.} Let $U$ be the region bounded by the cycle $Q$ of $G$ consisting of the subpath $P(v_2,v_x)$ of $H$ and the edge $\{v_x,v_2\}$ of $G \setminus H$. JCT and Corollary 1.1.5 imply that $v_y$ and $v_{x+1}$ are on distinct sides of $Q$. Either  $v_{x+1} \in U$ or $v_{x+1} \in G \setminus U$.

Assume $v_{x+1} \in U$. By JCT, $P(v_{x+1}, v_r)$ is contained in $U$ as well. But since $x<y$, $v_y \in P(v_{x+1}, v_r)$, contradicting $v_y$ and $v_{x+1}$ are on distinct sides of $U$. See Figure 5.32. The case where $v_{x+1} \in G \setminus U$ is similar so we omit it.

\endgroup 

\null

\begingroup
\setlength{\intextsep}{0pt}
\setlength{\columnsep}{20pt}
\begin{wrapfigure}[]{r}{0cm}
\begin{adjustbox}{trim=0cm 0cm 0cm 0cm}

\end{adjustbox}
\end{wrapfigure}

\noindent \textit{CASE 1.2: $x>y$.} Since $v$ is incident on $P(X,Y)$ at least one of $W+(0,2)$, $W+(-1,2)$,  $W+(-1,1)$, and  $W+(0,1)$ belongs to $P(X,Y)$. First we check that at least one of $W+(0,2)$, $W+(-1,2)$, and $W+(-1,1)$ belongs to $P(X,Y)$. If $W+(0,1)\in P(X,Y)$, since $W+(0,1)$ is not $X$ or $Y$, both of its $H$-neighbours belong to $P(X,Y)$. In particular, $W+(-1,1)$ would belong to $P(X,Y)$. And if $W+(0,1)\notin P(X,Y)$, then at least one of $W+(0,2)$, $W+(1,2)$, and $W+(-1,1)$ must belong to $P(X,Y)$. Either way, at least one of $W+(0,2)$, $W+(1,2)$,  and $W+(-1,1)$ belongs to $P(X,Y)$.

Let $U$ be the region bounded by the cycle $Q$ of $G$ consisting of the subpath $P(v_1,v_y)$ of $H$ and the edge $\{v_y,v_1\}$ of $G \setminus H$. See Figure 5.41. JCT and Corollary 1.1.5 imply that $Y$ and $W+(-1,1)$ belong to distinct sides of $Q$. Similarly, $X$ and $W+(-1,1)$ belong to distinct sides of $Q$. Let $X'$ be a box from the set $\{W+(0,2), W+(1,2), W+(-1,1)\}$ that belongs to $P(X,Y)$, and let $P(X,X')$ and $P(X',Y)$ be such that $P(X,X') \cap P(X',Y) = X'$. See Figure 5.33.

Let $c_1, \ldots, c_s$ be the centers of the boxes $X=X_1$, \ldots, $X'=X_s$ of $P(X,X')$. Note that for each $j \in \{1, \ldots, q-1\}$, the segment $[c_j,c_{j+1}]$ intersects the gluing edge of $X_j$ and $X_{j+1}$, and no other edge of $G$. Since $c_1$ and $c_s$ are on distinct sides of $Q$, the polygonal path $[c_1,c_2],\ldots, [c_{s-1},c_s]$ must intersect $Q$. Since each segment intersects gluing edges of $P(X,Y)$, the polygonal path $[c_1,c_2],\ldots, [c_{s-1},c_s]$ must intersect $Q$ at the edge $\{v_y,v_1\}$. It follows that $W+(0,1)$ and $W+(-1,1)$ both belong to $P(X,X')$. A similar argument show that $W+(0,1)$ and $W+(-1,1)$ must belong to $P(X',Y)$. But this contradicts our assumption that $P(X,X') \cap P(X',Y) = X'$. End of Case 1.2. End of Case 1.

\null

\begingroup
\setlength{\intextsep}{0pt}
\setlength{\columnsep}{20pt}
\begin{wrapfigure}[]{l}{0cm}
\begin{adjustbox}{trim=0cm 0.5cm 0cm 0cm} 

\end{adjustbox}
\end{wrapfigure}

\noindent \textit{CASE 2: The $j$-stack of $A_0$'s has length strictly greater than zero.} Let $P_0$, $P_1$, $P_2$ be the $A_1$ partitioning of $H$ given in Lemma 4.1. That is, $P_1=P(v_1,v_x)$, $P_0=P(v_x,v_y)$, and $P_2=P(v_y,v_r)$. Then $v(k,l-1) \in P_0$ or $v(k,l-1) \notin P_0$.

\null

\noindent \textit{CASE 2.1: $v(k,l-1) \in P_0$.} Let $v(k,l-1)=v_s$. Then $1<x<s<y<r$, and one of the edges of $W$ in $H$ is in $P(v_1,v_s)$, while the other is in $P(v_s,v_r)$. Now, $v_{s-1}=v(k-1,l-1)$, or $v_{s-1}=v(k,l-2)$. 

\null

\noindent \textit{CASE 2.1(a):  $v_{s-1}=v(k-1,l-1)$.} Then, by Lemma 1.4.5, after $bb_{v_1}(south)$, $W$ becomes parallel. For each $i \in  \{1, \ldots, j\}$, let $\mu_i$ be the move $W+(0,2j-1) \mapsto W+(0,2j)$. Let $\Xi$ be the sequence of moves $bb_{v_1}(south)$, Sw($W$), $\mu_1, \ldots, \mu_j$. Note that $\Xi$ is a cascade and we apply it. Let $H'$ be the resulting Hamiltonian path. Then $v_1'=v(k-1,l-1)$ and $v_s'=v(k,l-1)$ are the first and $s^{\text{th}}$ vertices in $H'$, respectively. Then $v_{s-1}'=v(k+1,l-1)$ or $v_{s-1}'=v(k,l)$. If the former, then after $bb_{v_1'}(east)$, the end-vertex in the resulting Hamiltonian path is at $v(k+1,l-1)$, so $(\dag_1)$ and $(\dag_2)$ are satisfied. And if the latter, then after $bb_{v_1'}(east)$, we are back to Case 2.1 of Proposition 5.6. See Figure 5.34. End of Case 2.1(a).

\endgroup 

\null

\begingroup
\setlength{\intextsep}{0pt}
\setlength{\columnsep}{20pt}
\begin{wrapfigure}[]{r}{0cm}
\begin{adjustbox}{trim=0cm 0.5cm 0cm 0.5cm} 

\end{adjustbox}
\end{wrapfigure}

\noindent \textit{CASE 2.1(b):  $v_{s-1}=v(k,l-2)$.}  As in Case 2.1(a), after $bb_{v_1}(south)$, $W$ is switchable. If $j=1$, then after $bb_{v_1}(south)$, Sw($W$), we are back to Case 2.1 of Proposition 5.6, so assume that $j>1$. For each $i \in  \{1, \ldots, j-1\}$, let $\mu_i$ be the move $W+(0,2j-1) \mapsto W+(0,2j)$. Let $\Xi$ be the sequence of moves $bb_{v_1}(west)$, Sw($W$), $\mu_1, \ldots, \mu_{j-1}$. Note that $\Xi$ is a cascade and we apply it. Again, we are back to Case 2.1 of Proposition 5.6. See Figure 5.35. End of Case 2.1(b). End of Case 2.1.

\null 

\noindent \textit{CASE 2.2: $v(k,l-1) \notin P_0$.} Then the $A_0$-type following $e(v)$ southward must also belong to $P_1$ or $P_2$. Let $W'=R(k,l)$. We have that $V(W'+(0,-1)) \subset V(P_1) \cup V(P_2)$. Then, by Lemma 4.1, $W'+(0,-1) \notin P(X,Y)$. Now we check that $W' \in P(X,Y)$.

\endgroup

\begingroup
\setlength{\intextsep}{0pt}
\begin{center}
\begin{adjustbox}{trim=0cm 0cm 0cm 0cm} 

\end{adjustbox}
\end{center}

\noindent For a contradiction, assume that $W' \notin P(X,Y)$. Note that this implies that both $W'+(-1,0)$ and $W'+(-1,-1)$ belong to $P(X,Y)$. Since neither $W'+(-1,0)$ nor $W'+(-1,-1)$ is an end-box, their $H$-neighbours in $P(X,Y)$ must be $W'+(-2,0)$ and $W'+(-2,-1)$, respectively. But this implies that $v(k,l-1)$ is an end-vertex of $H$, contradicting the NAA. Thus we must have that  $W' \in P(X,Y)$. See Figure 5.36 (a). 

Next, we check that $e(k+1;l,l+1) \in H$. For a contradiction, assume that $e(k+1;l,l+1) \notin H$. Then we have that $S_{\downarrow}(k+1,l+2;k+2,l+1) \in H$, and that $e(k+1,k+2;l) \in H$. Then, at least one of $W'+(0,1)$ and $W'+(1,0)$ must be an $H$-neighbour of $W'$ in $P(X,Y)$. Suppose that $W'+(1,0)$ is an $H$-neighbour of $W'$ in $P(X,Y)$. Then $W'+(1,0)$ is either switchable, or an end-box of $P(X,Y)$, neither of which is possible. The same goes for $W'+(0,1)$. It follows that we must have $e(k+1;l,l+1) \in H$. See Figure 5.36 (b). Then the $H$-neighbours of $W'$ in $P(X,Y)$ must be $W'+(-1,0)$ and $W'+(0,1)$.

Next we check that $v(k,l+1)$ must belong to $P_0$. Note that the vertices $v_1$, $v_2$, and $v_3$ belong to $P_1$, while $W' \in P(X,Y)$. Then, by Lemma 4.1, $v(k,l+1) \in P_0$. See Figure 5.36 (c).

\noindent Let $v(k,l+1)=v_s$. Then $1<x<s<y<r$, and one of the edges of $W$ in $H$ is in $P(v_1,v_s)$, while the other is in $P(v_s,v_r)$.  Now, $v_{s-1}=v(k,l+2)$ or $v_{s-1}=v(k-1,l+1)$. 

\endgroup 

\begingroup
\begin{center}

\setlength{\intextsep}{0pt}
\setlength{\columnsep}{20pt}

\begin{adjustbox}{trim=0cm 0cm 0cm 0cm} 

\end{adjustbox}
\end{center}
\endgroup 

\noindent \textit{CASE 2.2(a): $v_{s-1}=v(k,l+2)$.} Then, by Corollary 1.4.5, after $bb_{v_1}(north)$, $W$ becomes parallel. For each $i \in  \{1, \ldots, j+1\}$, let $\mu_i$ be the move $W+(0,2j-1) \mapsto W+(0,2j)$. Let $\Xi$ be the sequence of moves $bb_{v_1}(north)$, Sw($W$), $\mu_1, \ldots, \mu_{j+1}$. Note that $\Xi$ is a cascade, and we apply it. Now, if $e(k,k+1;l+2) \notin H$, we are back to Case 1 of Proposition 5.6; and if $e(k,k+1;l+2) \in H$, we are back to Case 2.1 of Proposition 5.6. See Figure 5.37 (a). End of Case 2.2(a). End of Case 2.1.

\null 

\noindent \textit{CASE 2.2(b): $v_{s-1}=v(k-1,l+1)$.} By Corollary 1.4.5, after $bb_{v_1}(north)$, $W$ becomes parallel. For each $i \in  \{1, \ldots, j+1\}$, let $\mu_i$ be the move $W+(0,2j-1) \mapsto W+(0,2j)$. Let $\Xi$ be the sequence of moves $bb_{v_1}(north)$, Sw($W$), $\mu_1, \ldots, \mu_{j+1}$. Note that $\Xi$ is a cascade, and we apply it. Let $H'$ be the resulting Hamiltonian path. Then $v_1'=v(k-1,l+1)$ and $v_s'=v(k,l+1)$ are the first and $s^{\text{th}}$ vertices in $H'$, respectively. Then $v_{s-1}'=v(k+1,l+1)$ or $v_{s-1}'=v(k,l+2)$. See Figure 5.37 (b). If the former, then after $bb_{v_1'}(east)$, the end-vertex in the resulting Hamiltonian path is on $v(k+1,l+1)$, so $(\dag_1)$ and $(\dag_2)$ are satisfied. 

Suppose then that $v_{s-1}'=v(k,l+2)$. Apply $bb_{v_1}(east)$. Now, if $e(k,k+1;l+2) \notin H$, we are back to Case 1 of Proposition 5.6; and if $e(k,k+1;l+2) \in H$, we are back to Case 2.1 of Proposition 5.6. End of Case 2.2(b). End of Case 2.2. End of Case 2. $\square$

\null

\noindent \textbf{Claim 5.21.} Recall that this is a continuation of Case 2.2$(c_2).(i)$ of Proposition 5.6. We are in the scenario where $(e(k+1,k+2;l-1), js\text{-}A_0, A_1,\text{North})$, FPC-2,3, and 4 hold, $u$ is not incident on $P(X,Y)$, but $v$ is. We will show that in this case, $(\dag_1)$ and $(\dag_2)$ hold.

% OMITTED: the case where u is incident on P(X,Y), as it is too similar to the same scenario as in Case 2.2(b_1).(i).
{
%Suppose that FPC-2,3 and 4 hold but FPC-1 fails. First we will consider the case where $u$ is incident on $P(X,Y)$. By Lemma 5.16, if $u\in R_0$ and incident on $P(X,Y)$, then $P(X,Y)$ has a switchable box, contradicting the assumption that FPC-3 holds. And if  $u\notin R_0$ and incident on $P(X,Y)$, then By Proposition 4.5(b), there is a cascade that fixes $v$, after which $e_c \in H$. Then we may apply the $F_1$ sequence of flips, after which we are back to Case 2.2$(c_1)$. Therefore, we may assume that $u$ is not incident on $P(X,Y)$. Then it must be the case that $v$ is incident on $P(X,Y)$.
}

Let $v_x$ and $v_y$ be the corners of the $A_1$-type of $P(X,Y)$. Let $P_0$, $P_1$, $P_2$ be the $A_1$ partitioning of $H$ given in Lemma 4.1. That is, $P_1=P(v_1,v_x)$, $P_0=P(v_x,v_y)$, and $P_2=P(v_y,v_r)$. Then $v(k,l-1) \in P_0$ or $v(k,l-1) \notin P_0$.

\endgroup 

\null

\noindent \textit{CASE 1: $v(k,l-1) \in P_0$.} Let $v_s=v(k,l-1)$. Then $1<x<s<y<r$, and one of the edges of $W$ in $H$ is in $P(v_1,v_s)$, while the other is in $P(v_s,v_r)$. Now, $v_{s-1}=v(k,l-2)$ or $v_{s-1}=v(k-1,l-1)$.

\null

\noindent \textit{CASE 1.1: $v_{s-1}=v(k,l-2)$.} Then we apply $bb_v(south)$. By Corollary 1.4.5, $W$ is now parallel. Apply $\text{Sw}(W)$. Apply the sequence of flips $W+(0,-1) \mapsto W+(0,-2)$, $\ldots$, $R(k+1,l-1) \mapsto R(k+1,l-2)$\footnote{ \noindent Note that if the $j$-stack of $A_0$s is empty, this sequence consists of only the move $W+(0,-1) \mapsto W+(0,-2)$.}. See Figure 5.38. Now we are back to a translation by $(0,-2)$ of Case 2.2$(c_1)$ of Proposition 5.6. End of Case 1.1

\begin{center}
\begin{adjustbox}{trim=0cm 0cm 0cm 0cm}

\end{adjustbox}
\end{center}

\noindent \textit{CASE 1.2: $v_{s-1}=v(k-1,l-1)$.} We apply $bb_v(south)$. By Corollary 1.4.5, $W$ is now parallel. Apply $\text{Sw}(W)$. Let $\Xi$ be the sequence of flips $W+(0,-1) \mapsto W+(0,-2)$, $\ldots,$ $R(k+1,l+1) \mapsto R(k+1,l)$. Apply $\Xi$. See Figure 5.39.  In the resulting Hamiltonian

\begin{center}
\setlength{\intextsep}{0pt} 
\setlength{\columnsep}{20pt}
\begin{adjustbox}{trim=0cm 0cm 0cm 0cm}
%
\end{adjustbox}
\end{center}

\noindent  path $H'$, let $v(k-1,l-1)=v'$ and $v(k,l-1)=v_{s'}$. Then $v_{s'-1}=v(k,l)$ or $v_{s'-1}=v(k,l-2)$. If $v_{s'-1}=v(k,l)$ (Figure 5.40 (a)), then after applying $bb_{v(k-1,l-1)}(east)$, we are back to Case 2.2$(c_1)$ of Proposition 5.6; and if $v_{s'-1}=v(k,l-2)$ (Figure 5.40 (b)), then after applying $bb_{v(k-1,l-1)}(east)$, $R(k+1,l-1) \mapsto R(k+1,l-2)$, we are back to a translation by $(0,-2)$ of Case 2.2$(c_1)$ of Proposition 5.6.  End of Case 1.2. End of Case 1.

%%%%% BETWEEN TIKZPICS %%%%%%%

\begin{center}
\setlength{\intextsep}{0pt} 
\setlength{\columnsep}{20pt}
\begin{adjustbox}{trim=0.5cm 0cm 0cm 0cm}

\end{adjustbox}
\end{wrapfigure}

\noindent \textit{CASE 2: $v(k,l-1) \notin P_0$.} Then either $j=0$ or $j>0$. 

\null

\noindent \textit{CASE 2.1: $j=0$.} Note that if $X+(0,1) \in P(X,Y)$, then after $X+(0,1) \mapsto W$, we are back to Case 2.1 of Proposition 5.6, so we may assume that $X+(0,1) \notin P(X,Y)$. It follows that $X+(-1,0) \in P(X,Y)$. First we will show that $v(k,l+1)$ must belong to $P_0$. For a contradiction, assume that $v(k,l+1) \notin P_0$. Consider the box $X+(-1,0)$ of $P(X,Y)$.  By Lemma 4.1, at least one vertex of $X+(-1,0)$ must belong to $P_0$. We have that $v(k,l) \in P_1$, and since $v(k,l+1) \notin P_0$, $v(k-1,l+1)$ is not in $P_0$ either. Now, either $e(k-1;l,l+1)\in H$, or $e(k-1;l,l+1)\notin H$. If $e(k-1;l,l+1)\in H$, then $v(k-1,l)$ is not in $P_0$ either, contradicting Lemma 4.1. And if $e(k-1;l,l+1)\notin H$, then we must have $(k-1;l-1,l) \in H$, and so again, $v(k-1,l)$ cannot be in $P_0$, contradicting Lemma 4.1. It follows that we must have $v(k,l+1) \in P_0$. See Figure 5.41.

\null

\begingroup
\setlength{\intextsep}{0pt}
\setlength{\columnsep}{20pt}
\begin{wrapfigure}[]{r}{0cm}
\begin{adjustbox}{trim=0cm 0cm 0cm 0cm}

\end{adjustbox}
\end{wrapfigure}

\noindent Let $v_s=v(k,l+1)$. Then $1<x<s<y<r$, and one of the edges of $W$ in $H$ is in $P(v_1,v_s)$, while the other is in $P(v_s,v_r)$. Now, $v(k,l+2)=v_{s+1}$ or $v(k,l+2)=v_{s-1}$. The former can be shown to be impossible; however it is shorter to note that if it were true (Figure 5.42), after $\text{Sw}(X+(0,1))$, we are back to Case 2.1 of Proposition 5.6, so we may assume the latter. Now, either $e(k,k+1;l+2) \in H$, or $e(k,k+1;l+2) \notin H$.

\null

\noindent \textit{CASE 2.1(a): $e(k,k+1;l+2) \in H$.} Then, after $bb_v(north)$, $W+(0,1)$ is switchable; and after $\text{Sw}(W+(0,1))$, we are back to a translation by $(0,2)$ of Case 2.2$(c_1)$ of Proposition 5.6. See Figure 5.43 (a). End of Case 2.1(a).

\begingroup
\setlength{\intextsep}{0pt}
\setlength{\columnsep}{20pt}

\begin{center}
    
\begin{adjustbox}{trim=0cm 0.25cm 0cm 0.25cm}

\end{adjustbox}

\end{center}

\endgroup

\noindent \textit{CASE 2.1(b): $e(k,k+1;l+2) \notin H$.} Then, after $bb_v(north)$, we are back to a translation by $(0,2)$ of Case 1 of Proposition 5.6. See Figure 4.43(b). End of Case 2.1(b). End of Case 2.1.

\null

\noindent \textit{CASE 2.2: $j>0$.} Note that all vertices of $R(k,l)$, except perhaps for $v(k,l+1)$ belong to $P_1$. Then, by Lemma 4.1, $v(k,l+1) \in P_0$. Let $v(k,l+1)=v_s$. Then $v_{s-1}=v(k,l+2)$ or $v_{s-1}=v(k-1,l+1)$. 

\null

\noindent \textit{CASE 2.2(a): $v_{s-1}=v(k,l+2)$.}  Then we apply $bb_v(north)$. By Corollary 1.4.5, $W$ is now parallel. Apply $\text{Sw}(W)$. Apply the sequence of flips $W+(0,-1) \mapsto W+(0,-2)$, $\ldots,$ $R(k+1,l+3) \mapsto R(k+1,l+2)$\footnote{\label{fn:emptyseq}Note that if $j=1$, this sequence is empty.}. See Figure 5.44. Now we are back to a translation by $(0,2)$ of Case 2.2$(c_1)$ of Proposition 5.6. End of Case 2.2(a)

\null 

\begin{center}
\setlength{\intextsep}{0pt} %this shifts tikz up-down
\setlength{\columnsep}{0pt}
\begin{adjustbox}{trim=0cm 0cm 0cm 0cm}

\end{adjustbox}
\end{center}

\noindent \textit{CASE 2.2(b): $v_{s-1}=v(k-1,l+1)$.} Then we apply $bb_v(north)$. By Corollary 1.4.5, $W$ is now parallel. Apply $\text{Sw}(W)$. Let $\Xi$ be the sequence of flips $W+(0,-1) \mapsto W+(0,-2)$, $\ldots,$ $R(k+1,l+3) \mapsto R(k+1,l+2)$\footnotemark[\getrefnumber{fn:emptyseq}]. See Figure 5.45. Apply $\Xi$. In the resulting Hamiltonian path $H'$, let $v(k-1,l+1)=v'$ and $v(k,l+1)=v_{s'}$. Then $v_{s'-1}=v(k,l+2)$ or $v_{s'-1}=v(k,l)$. If $v_{s'-1}=v(k,l+2)$, then after applying $bb_{v'}(east)$, we are back to a translation by (0,2) of Case 2.2$(c_1)$ of Proposition 5.6 (Figure 5.46(a)); and if $v_{s'-1}=v(k,l)$, then after applying $bb_{v'}(east)$, $R(k+1,l+1) \mapsto R(k+1,l)$, we are back to Case 2.2$(c_1)$ of Proposition 5.6 (Figure 5.46(b)). End of Case 2.2(b). End of Case 2.2. End of Case 2. $\square$.

\begin{center}
\setlength{\intextsep}{0pt} 
\setlength{\columnsep}{20pt}
\begin{adjustbox}{trim=0cm 0.5cm 0cm 0cm}

\end{adjustbox}
\end{center}

\noindent \textbf{Bound for Proposition 5.6 (EtB).} The argument is similar to that of Observation 4.33. Let $n \geq m$. The longest cascade required to relocated an end-vertex closer to a desired side of $R_0$ arises when proceeding through Case 2.2$(c_2).(ii)$ of Proposition 5.6. That is, $(e(k+1,k+2;l-1), js\text{-}A_0, A_1)$, where $P(X,Y)$ is contained in a standard looping fat path that has a turn $T$ with a set of lengthenings $\mathcal{T}(T)$, with $|\mathcal{T}(T)|$ and $j$ as large as possible, and where $T$ has an edge weakening. Note that if $T$ has an end-vertex weakening instead, by Claims 5.18, 5.19, 5.20, and 5.21, we can find shorter cascades after which $(\dag_1)$ and $(\dag_2)$ hold. As in Observation 4.33, $j$ can be at most $\frac{n}{2}$. By Proposition 4.28, it takes at most $m$ moves for an edge weakening of $T$. After that, by Lemma 4.32, it takes at most two moves to switch $W$, and then at most $j-1 \leq \frac{n}{2}-1$ moves to return to Case 2.2$(c_1)$. Thus, it takes at most $\frac{n}{2}-1+m+2=\frac{n}{2}+m+1$ to return to Case 2.2$(c_1)$. By the same argument, we need at most another $\frac{n}{2}+m+1$ moves to return to Case 2.1, and then at most another two moves for $(\dag_1)$ and $(\dag_2)$ to hold. Therefore, in the worst case scenario, at most $2(\frac{n}{2}+m+1)+2=n+2m+4$ moves are required to relocate an end-vertex closer to the desired side of $R_0$. We wish to relocate the end-vertices to opposite sides of $R_0$, so EtB requires at most $n(n+2m+4)=n^2+2mn+4n$ moves.

\null 

\vspace*{-2cm}
\makebox[\textwidth][l]{%
  \hspace{-1.5cm}%
  \includegraphics[scale=0.85]{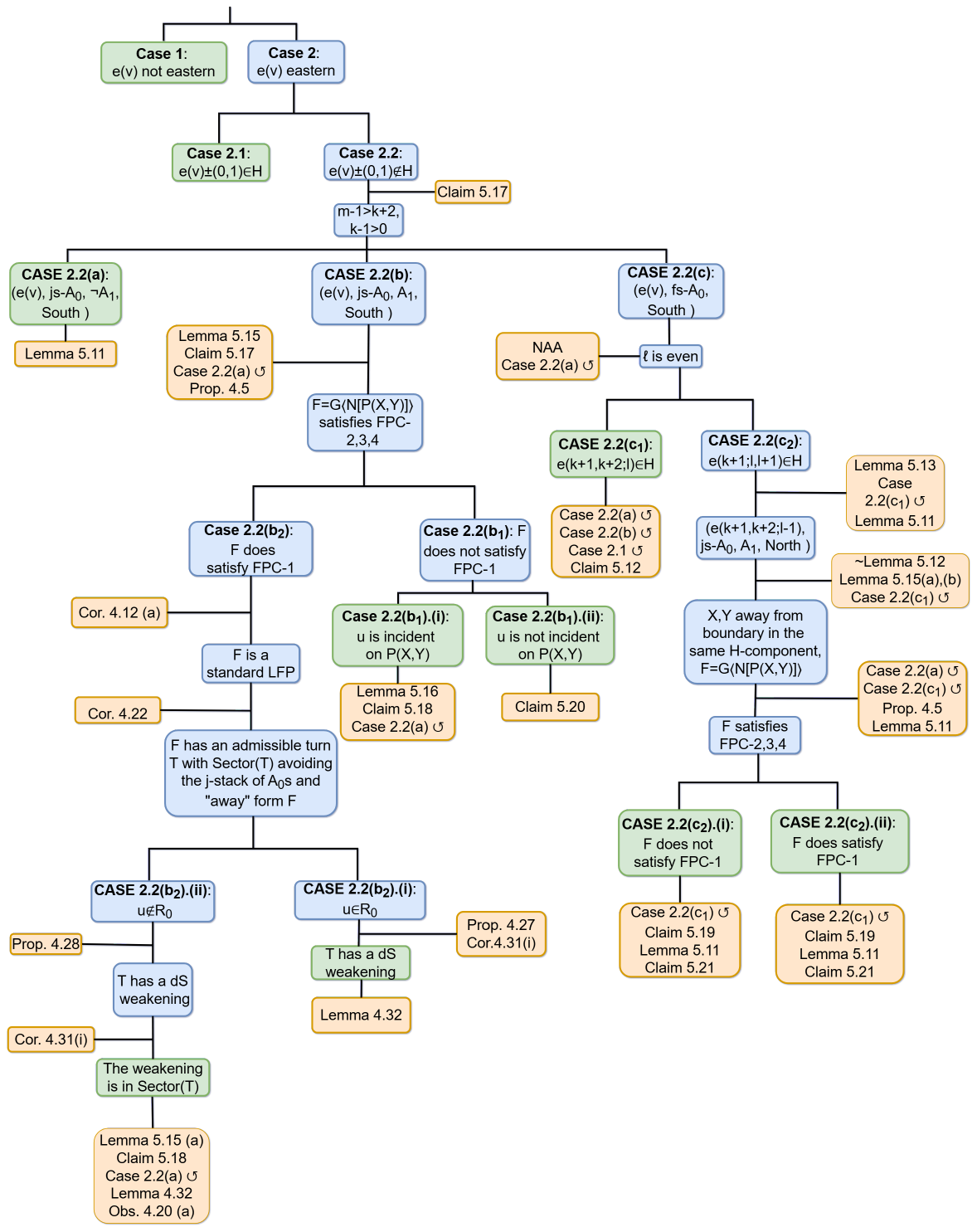}%
}

\vspace*{-2cm} 

\noindent \textbf{Flowchart 2.} Proof structure for Proposition 5.6. Green boxes indicate terminal cases, blue boxes indicate intermediate steps requiring further analysis. Orange boxes show the results applied at each step. When an orange box connects to a green box, that result completes the subcase. When an orange box appears on an edge between blue boxes, that result is used to transition between cases.

\subsection{Algorithm for relocating an endpoint to a corner (EtC)}

\noindent In this section, we prove Proposition 5.7. Four similar scenarios arise in the proof, each featuring an edge of $H$ followed by an $A$-type. We consider those scenarios together in Claim 5.25. As with Proposition 5.6, the claim can be read after the main proof. For Claim 5.25, we will also need Lemma 5.22.

\null 

\noindent \textbf{Lemma 5.22.} Let $G$ be an $m \times n$ grid graph, with $n \geq m \geq 3$. Let $H$ be a Hamiltonian path of $G$ with end-vertices $u$ and $v$. Assume that $v$ is on the eastern side of $R_0$, with y-coordinate at most $b$, and that $u$ is on the western side of $R_0$. Let $a \in \{0,1\}$ and let the edge $e=e(a;b,b+1)$ of $H$ be followed by a $j$-stack of $A_0$'s, which is then followed by an $A_1$-type with switchable middle-box $W$. Let $X$ and $Y$ be the boxes adjacent to $W$ that are not its $H$-neighbours, and assume that $F=G\langle N[P(X,Y)] \rangle$ is a \index{standard looping fat path}standard looping fat path. Then there is a double-switch cascade of length at most $m+2$ that switches $W$, and that avoids the $j$-stack of $A_0$'s and the boxes incident on $e$.

Analogous statements apply for the other compass directions.

\null 

\begingroup
\setlength{\intextsep}{0pt}
\setlength{\columnsep}{20pt}
\begin{wrapfigure}[]{r}{0cm}
\begin{adjustbox}{trim=0cm 0cm 0cm 0cm} 
\begin{tikzpicture}[scale=1.5]

\begin{scope}[xshift=0cm] 
{

\draw[gray,very thin, step=0.5cm, opacity=0.5] (0,0) grid (4,2.5);

\draw [blue, line width=0.5mm] (0,2)--++(0,-0.5)--++(0.5,0)--++(0,-0.5)--++(-0.5,0)--++(0,-0.5);

\draw [blue, line width=0.5mm] (1,2)--++(0,-0.5)--++(0.5,0)--++(0,-0.5)--++(-0.5,0)--++(0,-0.5);

\draw [blue, line width=0.5mm] (2,2)--++(0,-0.5)--++(0.5,0)--++(0,-0.5)--++(-0.5,0)--++(0,-0.5);

\draw [blue, line width=0.5mm] (3,2)--++(0,-0.5)--++(1,0);
\draw [blue, line width=0.5mm] (3,0.5)--++(0,0.5)--++(1,0);

\draw [orange, line width=0.5mm] (0.5,1)--++(0,0.5);

\draw[fill=blue] (0,0.5) circle (0.05);

\draw[red, dashed, line width=0.5mm, opacity=1] (4,0.5)--++(-2,2);

\node[left] at (0,0.5) [scale=1]{\tiny{$b$}};
\node[left] at (0,1) [scale=1]{\tiny{$+1$}};
\node[left] at (0,1.5) [scale=1]{\tiny{$+2$}};

\node[below] at (0,0) [scale=1] {\tiny{$0$}};
\node[below] at (0.5,0) [scale=1] {\tiny{$1$}};
\node[below] at (1,0) [scale=1] {\tiny{$2$}};
\node[below] at (1.5,0) [scale=1] {\tiny{$2$}};
\node[below] at (2,0) [scale=1] {\tiny{$3$}};
\node[below] at (2.5,0) [scale=1] {\tiny{$4$}};
\node[below] at (3,0) [scale=1] {\tiny{$5$}};

{
\node at (3.25,0.75) [scale=1]
{\small{$X$}};
\node at (3.25,1.25) [scale=1]
{\small{$W$}};
\node at (3.25,1.75) [scale=1]
{\small{$Y$}};

\node[right] at (0.45,1.25) [scale=1]
{\small{$e$}};

\node[below, align=left, text width=7cm] at (2, -0.25) { Fig. 5.47. An illustration with $a=1$ and $j=2$. The line $y=-x+a+b+2j+2$ in red.};

}

}

\end{scope}

\end{tikzpicture}
\end{adjustbox}
\end{wrapfigure}

\noindent \textit{Proof.} By Corollary 4.22, $F$ has an admissible \index{turn}turn $T$ such that $\text{Sector}(T)$ avoids the $j$-stack of $A_0$'s and lies above the line $y=-x+a+b+2j+2$. See Figure 5.47. By Corollary 4.19(a), $d(T)\geq 3$. By Proposition 4.27, $T$ has a double-switch weakening $\mu_1, \ldots, \mu_s$. By Corollary 4.31(i), $\mu_1, \ldots, \mu_s$ is contained in $\text{Sector}(T)$, and thus it avoids the $j$-stack of $A_0$'s. By Lemma 4.32 and the fact that there are no parallel switchable boxes, $\mu_1, \ldots, \mu_s$ can be extended by a double-switch move, after which $W$ is switched. The bound follows from Proposition 4.27 and Lemma 4.32. $\square$

\endgroup 

\null

\noindent Now we are ready to prove Proposition 5.7.

%pagemarker

\noindent \textit{Proof of Proposition 5.7.} Orient $H$ as $v=v_1, \ldots, v_r=u$. By Corollary 1.5.2, all switchable boxes are anti-parallel. We will use this fact implicitly and repeatedly throughout the proof. We will also use Lemmas and Corollaries 1.4.1-1.4.5 repeatedly, and often implicitly. For definiteness, assume that $v=v(0,l)$ is on the western side of $R_0$, that $u$ is on the eastern side, and that $c_v=v(0,n-1)$. To prove (i), it is enough to show that if $u$ is not at $c_u$ and $v$ is not at $c_v$, then there is a cascade after which $v$ is relocated to $v=v(0,l+2)$ or $v=v(0,l+4)$, and $u$ remains on $s_u$. We call this statement $(\ddag_1)$ for reference. To prove (ii), it is enough to show that if $u$ is at $c_u$ and $v$ is not at $c_v$, then there is a cascade after which $v$ is relocated to $v=v(0,l+2)$ or $v=v(0,l+4)$, and $u$ is at $c_u$. We call this statement $(\ddag_2)$ for reference.  

There is a lot of overlap between the proofs of (i) and (ii), so we will combine the proofs.

Either $e(0;l,l+1)\in H$ or $e(0;l,l+1)\notin H$.

\null 

\noindent \textit{CASE 1: $e(0;l,l+1)\notin H$.} Then $S_{\downarrow}(0,l+2;1,l+1) \in H$. Let $v(0,l+1)=v_q$. If $v(0,l+2)=v_{q-1}$, then $bb_v(north)$, $v$ is relocated to $v(0,l+2)$ and $(\ddag_1)$ and $(\ddag_2)$ hold. So we may assume that $v(0,l+2)=v_{q+1}$. Let $Z =R(0,l+2)$. Either $e(1;l+1,l+2) \in H$ or $e(1;l+1,l+2) \notin H$.

\null

\noindent \textit{CASE 1.1: $e(1;l+1,l+2) \notin H$.} Either $e(0,1;l+2)\in H$ or $e(0,1;l+2)\notin H$.

\begingroup
\setlength{\intextsep}{0pt}
\setlength{\columnsep}{20pt}
\begin{wrapfigure}[]{l}{0cm}
\begin{adjustbox}{trim=0cm 0cm 0cm 0cm} 

\end{adjustbox}
\end{wrapfigure}

\null

\noindent \textit{CASE 1.1(a): $e(0,1;l+2)\in H$.} Then, after $bb_v(north)$, $bb_{v(1,l+1)}(north)$, $(\ddag_1)$ and $(\ddag_2)$ hold. End of Case 1.1(a). See Figure 5.57.

\null

\noindent \textit{CASE 1.1(b): $e(0,1;l+2)\notin H$.} Since $u\neq v(0,n-1)$ and since $v$ has the same parity as the corner $v(0,n-1)$, $n-1 \geq l+4$. Then $e(0;l+2,l+3) \in H$ and $e(1;l+2,l+3) \in H$. Note that if $e(0,1;l+3) \in H$ (Figure 5.49 (a)), then we must have $e(0,1;l+4) \in H$. Then, after $Z+(0,1) \mapsto Z$, we are back to Case 1.1(a), so we may assume that $e(0,1;l+3) \notin H$. Then we must have that $n-1\geq l+4$, and that $e(0;l+3,l+4) \in H$. Let $v_s=v(1,l+1)$ and $v_t=v(1,l+2)$. Then, either $s<t$ or $s>t$.

\begin{center}
\begingroup
\setlength{\intextsep}{0pt}
\setlength{\columnsep}{20pt}
\begin{adjustbox}{trim=0cm 0cm 0cm 0cm} 

\end{adjustbox}
\endgroup 
\end{center}

\noindent \textit{CASE 1.1($b_1$): $s<t$.} Then, after $bb_{v}(north)$, $bb_{v(1,l+1)}(north)$, $bb_{v(1,l+3)}(west)$, $(\ddag_1)$ and $(\ddag_2)$ hold. See Figure 5.49 (b).End of Case 1.1 ($b_1$).

\begin{center}
\begingroup
\setlength{\intextsep}{0pt}
\setlength{\columnsep}{20pt}
\begin{adjustbox}{trim=0cm 0cm 0cm 0cm} 
%
\end{adjustbox}
\endgroup  
\end{center}

\noindent \textit{CASE 1.1($b_2$): $s>t$.} Then, after $bb_{v}(north)$, Sw$(Z)$, $bb_{v(1,l+1)}(west)$,  $(\ddag_1)$ and $(\ddag_2)$ hold. See Figure 5.50. End of Case 1.1($b_2$). End of Case 1.1(b). End of Case 1.1.

\endgroup

\null

\noindent \textit{CASE 1.2: $e(1;l+1,l+2) \in H$.} Then $e(0,1;l+2) \notin H$. Since $u\neq v(0,n-1)$ and since $v$ has the same parity as the corner $v(0,n-1)$, $n-1 \geq l+4$. Then $e(0;l+2,l+3) \in H$. See Figure 5.51 (a). Now, either $e(1;l+2,l+3) \in H$, or $e(1;l+2,l+3) \notin H$.

\begingroup
\setlength{\intextsep}{0pt}
\setlength{\columnsep}{20pt}
\begin{adjustbox}{trim=0cm 0cm 0cm 0cm} 

\end{adjustbox}

\null

\noindent \textit{CASE 1.2(a): $e(1;l+2,l+3) \in H$.} Then there is a cascade $bb_v(north), \text{Sw}(Z+(0,-1)), Z+(0,-1) \mapsto Z$, $bb_{v(1,l+1)}(west)$, after which $(\ddag_1)$ and $(\ddag_2)$ hold. See Figure 5.51 (b). End of Case 1.2(a).

\null 

\noindent \textit{CASE 1.2(b): $e(1;l+2,l+3) \notin H$.} Then $e(1,2;l+2) \in H$. See Figure 5.52 (a). Now, either $e(0;l+3,l+4) \in H$ or $e(0;l+3,l+4) \notin H$.

\null 

\noindent \textit{CASE 1.2($b_1$): $e(0;l+3,l+4) \in H$.} Then $S_{\downarrow}(1,l+4;2,l+3) \in H$. Note that if $e(2;l+3,l+4) \in H$ (Figure 5.52 (b)), then, after $Z+(1,0)\mapsto Z+(1,1)$, we are back to Case 1.2(a). And if $e(2;l+1,l+2) \in H$ (Figure 5.52 (c)), then, after $Z+(1,0)\mapsto Z+(1,-1)$, we are back to Case 1.1(b). So, we may assume that $e(2;l+3,l+4) \notin H$ and that $e(2;l+1,l+2) \notin H$. We call this configuration $I_1$ for reference (Figure 5.52 (d)). We will prove in Claim 5.23 below that in this case, $(\ddag_1)$ and $(\ddag_2)$ hold. End of Case 1.2$(b_1)$.

\begingroup
\setlength{\intextsep}{0pt}
\setlength{\columnsep}{20pt}
\begin{adjustbox}{trim=0cm 0cm 0cm 0cm} 

\end{adjustbox}

\noindent \textit{CASE 1.2($b_2$): $e(0;l+3,l+4) \notin H$.}  Then $e(0,1;l+3) \in H$. Since $u\neq v(0,n-1)$, $n-1>l+4$. Then $S_{\downarrow}(0,l+5;1,l+4) \in H$. Note that if $e(1,2;l+3)\in H$ (Figure 5.53(a)), then, after $Z+(1,0) \mapsto Z+(0,1)$, we are back to Case 1.2(a) so we may assume that $e(1,2;l+3)\notin H$. Then $e(1;l+3,l+4) \in H$.

If $e(2;l+3,l+4) \in H$ (Figure 5.53(b)), then, after $Z+(1,1) \mapsto Z+(0,1)$, we are back to Case 1.2($b_1$), so we may assume that $e(2;l+3,l+4) \notin H$. Note that this implies that $m-1 \geq 3$. Now, if $e(1;l+3,l+4)$ is not followed by an $A_0$-type or an $A_1$-type, then, by Lemma 5.11, there is a cascade of flips after which $e(2;l+3,l+4) \in H$. Then we can apply $Z+(1,1) \mapsto Z$, and once again we are back to Case 1.2$(b_1)$. Therefore we may assume that $e(1;l+3,l+4)$ is followed by an $A_0$-type or an $A_1$-type. We call this configuration $I_2$ for reference (Figure 5.53(c)). We will prove in Claim 5.23 below that in this case, $(\ddag_1)$ and $(\ddag_2)$ hold. End of Case 1.2($b_2$). End of Case 1.2(b). End of Case 1.2. End of Case 1. $\square$

\begin{center}
\begingroup
\setlength{\intextsep}{0pt}
\setlength{\columnsep}{20pt}
\begin{adjustbox}{trim=0cm 0cm 0cm 0cm} 

\end{adjustbox}
\end{center}

\noindent \textit{CASE 2.1(a): $e(0,1;l+1) \in H$.} Then we must have $e(0,1;l+2) \in H$. Then, after $Z+(0,1) \mapsto Z$, we are back to Case 1. See Figure 5.54 (a). End of Case 2.1(a). 

\null 

\noindent \textit{CASE 2.1(b): $e(0,1;l+1) \notin H$.} Then we must have $e(0;l+1,l+2) \in H$. Then we can apply $bb_v(east)$, $bb_{v(1,l+1)}(west)$, and then $(\ddag_1)$ and $(\ddag_2)$ hold. See Figure 5.54 (b). End of Case 2.1(b). End of Case 2.1. 

\endgroup

\null 

\noindent \textit{CASE 2.2: $e(1;l,l+1) \notin H$.} Since $m-1>1$, we must have $S_{\uparrow}(1,l-1;2,l) \in H$, and since $u$ is on the eastern side of $R_0$, we must have $S_{\uparrow}(0,l-2;1,l-1) \in H$. Now, either $e(0;l+1,l+2) \in H$, or $e(0;l+1,l+2) \notin H$.

%[pagemarker]

\noindent \textit{CASE 2.2(a): $e(0;l+1,l+2) \in H$.} Then $S_{\downarrow}(1,l+2;2,l+1) \in H$. We call this configuration $I_3$ for reference. See Figure 5.55 (a). We will prove in Claim 5.23 below, that in this case $(\ddag_1)$ and $(\ddag_2)$ hold. End of Case 2.2(a).

\begingroup
\setlength{\intextsep}{0pt}
\setlength{\columnsep}{20pt}

\begin{adjustbox}{trim=0cm 0cm 0cm 0cm} 

\end{adjustbox}

\noindent \textit{CASE 2.2(b): $e(0;l+1,l+2) \notin H$.} Then $S_{\downarrow}(0,l+3;1,l+2) \in H$. Let $W'=R(0,l+1)$, $X'=R(0,l)$ and $Y'=R(0,l+2)$. Now, either $e(1,2;l+1) \in H$ or $e(1,2;l+1) \notin H$.

{
%Note that if $e(1;l+2,l+3) \in H$, then we must have $e(1,2;l+1) \in H$. Then, after $W \mapsto Y$, we are back to Case 3.2(a). Therefore we may assume that $e(1;l+2,l+3) \notin H$. \textcolor{red}{[See Figure 4.42(b).]} .
}

\null

\noindent \textit{CASE 2.2($b_1$): $e(1,2;l+1) \in H$.} Then $e(1;l+1,l+2) \notin H$.  Now, $X'$ and $Y'$ either belong to the same $H$-component, or they do not.

\null 

\noindent \textit{CASE 2.2$(b_1).(i)$: $X'$ and $Y'$ belong to the same $H$-component.} This implies that $e(2;l,l+1) \notin H$. Let $J$ be the $H$ component containing $X'$ and $Y'$. Then $X'+(1,0)$ is a switchable box in the $H$-path $P(X',Y')$ contained in $J$. Then, by Proposition 4.5, $X'+(1,0) \mapsto W' $ is a valid move. Then we can apply $bb_v(east)$, $bb_{v(1,l+1)}(west)$, after which $(\ddag_1)$ and $(\ddag_2)$ hold. See Figure 5.55 (b). End of Case 2.2$(b_1).(i)$.

\null

\begingroup
\setlength{\intextsep}{0pt}
\setlength{\columnsep}{20pt}
\begin{wrapfigure}[]{r}{0cm}
\begin{adjustbox}{trim=0cm 0cm 0cm 0.75cm} 

\end{adjustbox}
\end{wrapfigure}

\noindent \textit{CASE 2.2$(b_1).(ii)$: $X'$ and $Y'$ belong to distinct $H$-components.} By Lemma 5.15 (c), there is a cascade that keeps $u$ on the eastern side of $R_0$, after which $W'$ is switched. Then, we are back to Case 2.2(a). See Figure 5.56.  End of Case 2.2$(b_1).(ii)$.  End of Case 2.2$(b_1)$.  

\null 

\endgroup

\begingroup
\setlength{\intextsep}{0pt}
\setlength{\columnsep}{20pt}
\begin{wrapfigure}[]{l}{0cm}
\begin{adjustbox}{trim=0cm 0cm 0cm 0cm} 
%
\end{adjustbox}
\end{wrapfigure}

\noindent \textit{CASE 2.2($b_2$): $e(1,2;l+1) \notin H$.} Then $e(1;l+1,l+2) \in H$. By Lemma 5.11, if $e(1;l+1,l+2)$ is not followed by an $A_0$-type or an $A_1$-type, then there is a cascade after which $e(2;l+1,l+2) \in H$. Then we can apply $W'+(1,0) \mapsto W'$, and then we are back to Case 2.2(a). Therefore we may assume that $e(1;l+1,l+2)$ is followed by an $A_0$-type or an $A_1$-type. We call this configuration $I_4$ for reference. See Figure 5.57 We will prove in Claim 5.23 below that in this case, $(\ddag_1)$ and $(\ddag_2)$ hold. End of Case 2.2($b_2$). End of Case 2.2(b). End of Case 2.2. End of Case 2. $\square$

\endgroup

\null 

\noindent It remains to prove Claim 5.23. We postpone proving the bound on the length of the cascade required by the EtC algorithm until after we prove Claim 5.23.

\null

\noindent \textbf{Claim 5.23.} Suppose that we are in one of the configurations $I_1, \ldots, I_4$, encountered in the proof of Proposition 5.7. Then there is a cascade after which $(\ddag_1)$ and $(\ddag_2)$ are satisfied.

\null 

\noindent \textit{Proof.} For $p \in \{1, \ldots, 4 \}$, each configuration $I_p$ (Figure 5.58) features an edge $e=e(a;b,b+1)$ with $a \in \{0,1\}$  followed by an $A_1$-type or an $A_0$-type, and such that the end-vertex $v$ has y-coordinate at most $b$. There are three possibilities:

1. $(e, js\text{-}A_0, \neg A_1, \text{East})$, \ \ \ 2. $(e, fs\text{-}A_0, \text{East})$, and \ \ \ 3. $(e, js\text{-}A_0, A_1, \text{East})$.

\begin{center}
    
\begingroup
\setlength{\intextsep}{0pt}
\setlength{\columnsep}{20pt}

\begin{adjustbox}{trim=0cm 0cm 0cm 0cm} 

\end{adjustbox}

\endgroup 

\end{center}

\noindent \textit{CASE 1: $(e, js\text{-}A_0, \neg A_1, \text{East})$.} Let $a'=2j$. Since the $j$-stack of $A_0$s is not a full stack, $m-1\geq a'+2$. We consider each $p \in \{ 1,2,3,4\}$ separately.

\null

\begingroup
\setlength{\intextsep}{0pt}
\setlength{\columnsep}{20pt}
\begin{wrapfigure}[]{r}{0cm}
\begin{adjustbox}{trim=0cm 0cm 0cm 0cm} 
%
\end{adjustbox}
\end{wrapfigure}

\noindent \textit{CASE 1.1: $p=1$.} Note that if $e(a'+1;l+2,l+3) \in H$, then after $R(a',l+2) \mapsto R(a'-1,l+2), \ldots, R(2,l+2) \mapsto R(1,l+2)$, $R(1,l+2)$ is switched, and we are back to Case 1.2(a) of Proposition 5.7. Therefore we may assume that $e(a'+1;l+2,l+3) \notin H$. See Figure 5.59. 

If the rectangle $R(a'+1,a'+2;l+1,l+4)$ does not contain $u$, then, by Lemma 5.11, there is a cascade of flips after which $R(1,l+2)$ is switched, and we are back to Case 1.2(a) of Proposition 5.7. Therefore we may assume that  $R(a'+1,a'+2;l+1,l+4)$ does contain $u$. This implies that $m-1=a'+2$. Note that if $e(m-1;l+2,l+3) \in H$, then stack of $A_0$s is a full one, contrary to the Case 1 assumption, so we may assume that  $e(m-1;l+2,l+3) \notin H$. Then, at least one of $e(m-1; l+1,l+2)$ and $e(m-1; l+3,l+4)$ belongs to $H$. WLOG assume that $e(m-1; l+1,l+2) \in H$. Then, after $R(a'+1;l+2) \mapsto R(a'+1;l+1)$, $e(a'+1;l+2,l+3) \in H$, and we can reuse the cascade in the first paragraph. End of Case 1.1.

\null 

\noindent \textit{CASE 1.2: $p=2$.} By an argument similar to the one in Case 1.1, we find a cascade of flips after which $R(2,l+3)$ is switched. Then we can apply $R(1,l+3) \mapsto R(0,l+3)$, and then we are back to Case 1.1. End of Case 1.2.

\null 

\noindent \textit{CASE 1.3: $p=3$.} By an argument similar to the one in Case 1.1, we find a cascade of flips after which $R(1,l)$ is switched, and then we are back to Case 2.1(b) of Proposition 5.7.  End of Case 1.3.

\null 

\noindent \textit{CASE 1.4: $p=4$.} By an argument similar to the one in Case 1.1, we find a cascade of flips after which $R(2,l+1)$ is switched. Then we can apply $R(1,l+1) \mapsto R(0,l+1)$, and then we are back to Case 1.3. End of Case 1.4. End of Case 1.

\null 

\noindent \textit{CASE 2: $(e, fs\text{-}A_0, \text{East})$.} Recall that one of the assumptions of Proposition 5.7 stated that $u$ and $c_u$, and $v$ and $c_v$, have matching parities. We will refer to this assumption as the ``end-vertex-corner parity'' assumption and write ``ECP'' whenever we refer to it. We will consider each $p$ separately.

\null

\noindent \textit{CASE 2.1. $p=1$.} If $m-1$ is odd (Figure 5.60 (a)), then we must have $e(m-1;l+2,l+3) \in H$. Then, after $R(m-2,l+2) \mapsto R(m-3,l+2), \ldots, R(2,l+2) \mapsto R(1,l+2)$, we are back to Case 1.2(a) of Proposition 5.7. Therefore we may assume that $m-1$ is even. 

If $m-1=2$ (Figure 5.60 (b)), then $u=v(2,l+1)$. Now, since $u=v+(2,1)$, $u$ and $v$ have different parities, while the corners $c_v$ and $c_u$ on $v(0,n-1)$ and $v(2,n-1)$ must have the same parity. Then either $c_u$ and $u$ have different parities, or $c_v$ and $v$ do. Either way, this contradicts the ECP assumption.  

If $m-1=4$ (Figure 5.60 (c)), then we must have $S_{\uparrow}(2,l;3,l+1) \in H$. It follows that $u=v(4,l+1)$, and this also contradicts the ECP assumption.

Suppose that $m-1\geq 6$. Then, for each $i \in \{2, 4, \ldots, m-5\}$, $S_{\uparrow}(i,l;i+1,l+1) \in H$ implies that $S_{\uparrow}(i+2,l;i+3,l+1) \in H$. It follows that for $i \in \{2, 4, \ldots, m-3\}$, $S_{\uparrow}(i,l;i+1,l+1) \in H$. Then we must have that $u=v(m-1,l+1)$. But then again, this contradicts the ECP assumption. End of Case 2.1.

\begingroup
\setlength{\intextsep}{0pt}
\setlength{\columnsep}{20pt}
\begin{adjustbox}{trim=0cm 0cm 0cm 0cm} 

\end{adjustbox}

\noindent \textit{CASE 2.2. $p=2$.} If $m-1$ is even, then we must have $e(m-1;l+3,l+4) \in H$.  Then, after $R(m-2,l+3) \mapsto R(m-3,l+3), \ldots, R(1,l+3) \mapsto R(0,l+3)$, we are back to Case 2.1.

Therefore we may assume that $m-1$ is odd. See Figure 5.60 (d). Then, using the same argument as in the last paragraph of Case 2.1, we can check that $u=v(m-1,l+2)$. See Figure 5.73(d). Now we see that both end-vertices have y-coordinate at most $l+2$, and there is a full stack of $A_0$'s starting at the small eastern leaf $R(0,l+3)$. Then, by Lemma 5.13, this is an impossible configuration. End of Case 2.2.

\null

\begingroup
\setlength{\intextsep}{0pt}
\setlength{\columnsep}{20pt}
\begin{adjustbox}{trim=0cm 0cm 0cm 0cm} 

\end{adjustbox}

\noindent \textit{CASE 2.3. $p=3$.} If $m-1$ is odd (Figure 5.61 (a)), then we must have $e(m-1;l,l+1) \in H$. Then, after $R(m-2,l) \mapsto R(m-3,l), \ldots, R(2,l) \mapsto R(1,l)$, we are back to Case 2.1(b) of Proposition 5.7.  And if $m-1$ is even (Figure 5.61 (b)), then, the same argument as in Case 2.1, shows that $u=v(m-1,l-1)$, which again, contradicts the ECP assumption. End of Case 2.3.

\null 

\noindent \textit{CASE 2.4. $p=4$.} If $m-1$ is odd (Figure 5.61 (c)), then we must have $e(m-1;l+1,l+2) \in H$. Then, after $R(m-2,l+1) \mapsto R(m-3,l+1), \ldots, R(2,l+1) \mapsto R(1,l+1)$, we are back to Case 2.3. Therefore, we may assume that  $m-1$ is even. See figure 5.61 (d). Then, using the same argument as in Case 2.1, we can check that $u=v(m-1,l)$. But then, as in Case 2.2, we can use Lemma 5.13 again to show that this is an impossible configuration. See Figure 5.75(b). End of Case 2.4. End of Case 2.

\null

\noindent \textit{CASE 3: $(e, js\text{-}A_0, A_1, \text{East})$.} Let $W$ be the switchable box of the $A_1$-type, and let $X$ and $Y$ be the boxes adjacent to $W$ that are not its $H$-neighbours. Either $X$ and $Y$ belong to distinct $H$-components, or they do not.

\null 

\noindent \textit{CASE 3.1: $X$ and $Y$ belong to distinct $H$-components.} By Lemma 5.15 (c), there is a cascade after which $W$ is switched, and which otherwise avoids $I_p$, for $p \in \{ 1,2,3,4\}$. Then we are back to Case 1. End of Case 3.1.

\null

\noindent \textit{CASE 3.2: $X$ and $Y$ belong to the same $H$-component.} Let $F=G\langle N[P(X,Y)] \rangle$. Either $F$ is a looping fat path or it is not.

\null

\noindent \textit{CASE 3.2(a): $F$ is a looping fat path.} Let $J(F)$ be the $H$-component of $G$ that contains $F$. Since both end-vertices of $H$ are in $R_0$, no edge of a main trail of $J(F)$ can have an end-vertex incident on it in $G \setminus R_0$. Then, by Lemma 1.3.12, $J(F)$ is non-self-adjacent. Since $F$ is an $H$-subtree of $J(F)$, $F$ is also non-self-adjacent.  By Proposition 4.11(i), $F$ has no polyking junctions. It follows that $F$ is a standard looping fat path. Then, by Lemma 5.22, there is a double-switch cascade $\mu_1, \ldots, \mu_q$ that switches $W$, that avoids the $j$-stack of $A_0$’s, and that is above the line $y=-x+a+b+2j+2$. See Figure 5.47. Then, after $\mu_q$, we are back to Case 1.

\null

\noindent \textit{CASE 3.2(b): $F$ is not looping fat path.} By Corollary 1.5.2, $W$ cannot be parallel, so FPC-2 must hold. If the end-vertex $u$ is incident on $P(X,Y)$, by Lemma 5.16, $P(X,Y)$ has a switchable $Z$ box incident on the eastern side of $R_0$. Then, after $Z \mapsto W$, we are back to Case 1. Now, either $v$ is incident on $P(X,Y)$, or it is not.

\null 

\noindent \textit{CASE 3.2$(b_1)$: $v$ is incident on $P(X,Y)$.}  First we check that for each $p \in \{ 1,2,3,4\}$, none of $R(a,b+1)$, $R(a,b)$, and $R(a+1,b)$ can belong to $P(X,Y)$.

\noindent Suppose that $p=1$. In this case, $a=0$ and $b=l+2$. Either $e(2;l+2,l+3) \in H$ or $e(2;l+2,l+3) \notin H$. 

Assume that $e(2;l+2,l+3) \in H$ (Figure 5.62 (a)). If $R(1,l+2) \in P$, then it must be an end-box of $P$, which it is not, so  $R(1,l+2) \notin P$. If $R(0,l+2) \in P$, it is not an end-box, so $R(0,l+1)$ must belong to $P$. But then $R(0,l+1)$  must be an end-box of $P$, which it is not, so $R(0,l+2)$ cannot be in $P$. A similar argument shows that  $R(0,l+3) \notin P$.

\begingroup
\setlength{\intextsep}{0pt}
\setlength{\columnsep}{20pt}
\begin{wrapfigure}[]{r}{0cm}
\begin{adjustbox}{trim=0cm 0cm 0cm 0cm} 

\end{adjustbox}
\end{wrapfigure}

Assume that $e(2;l+2,l+3) \notin H$ (Figure 5.62 (b)). Then $W=R(1,l+2)$. By Corollary 4.2(b), $R(1,l+2) \notin P$. The arguments showing that $R(0,l+2) \notin P$ and  $R(0,l+3) \notin P$ are similar to the ones used in the previous paragraph, so we omit them.

Hence, we have shown that if $p=1$, then none of $R(a,b+1)$, $R(a,b)$, and $R(a+1,b)$ can belong to $P(X,Y)$. The arguments for the cases where $p=2$, $p=3$, and $p=4$ are similar, so we omit them.

\noindent Now, the assumption on $v$ of Case 3.2$(b_1)$ implies that $p \neq 3$ and $p\neq 4$. Then $p=1$ or $p=2$. 

\endgroup 

\null

\begingroup
\setlength{\intextsep}{0pt}
\setlength{\columnsep}{20pt}
\begin{adjustbox}{trim=0cm 0cm 0cm 0cm} 

\end{adjustbox}

\noindent \textit{CASE 3.2$(b_1).(i)$: $p=1$.} First we check that if $e(v)$ is eastern, then $R(0,l)$ and $R(0,l-1)$ cannot belong to $P$.  Suppose that $e(v)$ is eastern. Then, if $R(0,l) \in P$, it would have to be an end-box, which it is not, so $R(0,l) \notin P$. Since $u$ is on the eastern side, we must have $S_{\uparrow}(0,l-2;1,l-1) \in H$. Then again, $R(0,l-1) \in P$, it would have to be an end-box, which it is not, so $R(0,l-1) \notin P$.

Thus we only need to check the case where $e(v)$ is southern. See Figure 5.63 (a). Then $R(0,l-1)$ must belong to $P$ and be switchable. Then, after $R(0,l-1) \mapsto W$, we apply the cascade of flips $W \mapsto W+(-1,0), \ldots, R(2,l+2) \mapsto R(1,l+2)$, followed by  $R(0,l) \mapsto R(2,l)$. See Figure 5.63 (a). Now we apply $bb_v(east)$, $bb_{v(1,l+1)}(west)$, (not shown in Figure 5.63 (a)) after which $(\ddag_1)$ and $(\ddag_2)$ hold. End of Case 3.2$(b_1).(i)$.

\null 

\noindent \textit{CASE 3.2$(b_1).(ii)$: $p=2$.} As in the previous case, we must have that $e(v)$ is southern, and that $R(0,l-1)$ must belong to $P(X,Y)$ and be switchable (Figure 5.63 (b). Then, after $R(0,l-1) \mapsto W$, we apply the cascade of flips $W \mapsto W+(-1,0), \ldots, R(1,l+3) \mapsto R(0,l+3)$, followed by $R(0,l) \mapsto R(0,l+3)$. See Figure 5.63 (b). Now we apply $bb_v(east)$, $bb_{v(1,l+1)}(west)$,  (not shown in Figure 5.63 (b)) after which $(\ddag_1)$ and $(\ddag_2)$ hold. End of Case 3.2$(b_1).(ii)$.  End of Case 3.2$(b_1)$.

\null 

\noindent \textit{CASE 3.2$(b_2)$: $v$ is not incident on $P(X,Y)$.} Then $P(X,Y)$ must have a switchable box $Z$. We consider each $p \in \{ 1,2,3,4\}$ separately. 

\null

\noindent \textit{CASE 3.2$(b_2).(i)$: $p=1$.} Since $v$ is not incident on $P(X,Y)$, $Z \neq R(0,l)$ and $Z \neq R(0,l-1)$. Then, after $Z \mapsto W$, we are back to Case 1.2(a) of Proposition 5.7. End of Case 3.2$(b_2).(i)$.

\null 

\noindent \textit{CASE 3.2$(b_2).(ii)$: $p=2$.} Since $v$ is not incident on $P(X,Y)$, $Z \neq R(0,l)$ and $Z \neq R(0,l-1)$. Then, after $Z \mapsto W$, we are back to Case 1.2. End of Case 3.2$(b_2).(ii)$.

\null 

\noindent \textit{CASE 3.2$(b_2).(iii)$: $p=3$.} Then, after $Z \mapsto W$, we are back to Case 1.3. End of Case 3.2$(b_2).(iii)$.

\null 

\noindent \textit{CASE 3.2$(b_2).(iv)$: $p=4$.} Then, after $Z \mapsto W$, we are back to Case 1.4. End of Case 3.2$(b_2).(iv)$. End of Case 3.2$(b_2)$. End of Case 3.2(b). End of Case 3.2. End of Case 3. $\square$

\null

\noindent \textbf{Bound for Proposition 5.7 (EtC).} Let $n \geq m$. The longest cascade required to relocate an end-vertex closer to a desired corner while keeping it on the same side of $R_0$ arises when proceeding through Case 3.2(a) of Claim 5.23, and considering configuration $I_2$. As before, $j \leq \frac{n}{2}$. By Lemma 5.22, switching $W$ requires at most $m+2$ moves. Then, after $m+2+\frac{n}{2}-1$ moves, we are back to configuration $I_1$. Using the same argument, we see that at most $m+1+\frac{n}{2}$ moves are required to return to Case 1.1$(b_1)$ of Proposition 5.7, and then three additional moves for $(\ddag_1)$ or $(\ddag_2)$ to be satisfied. Thus we need $n+2m+5$ moves to relocate an end-vertex closer to a desired corner while keeping it on the same side of $R_0$. Note that it takes at most $\frac{n}{2}$ such relocations to take an end-vertex from a vertex on a side of $R_0$ to a desired corner incident on that side. Therefore, EtC requires at most $\frac{n^2}{2}+mn+\frac{5n}{2}$ moves to complete (i) or (ii).

\null

\subsection{Summary}

\noindent In this chapter, we proved Theorem 5.9: any two Hamiltonian paths on an $m \times n$ grid graph can be reconfigured into one another. The proof uses Theorem 2.1, which handles reconfiguration of Hamiltonian cycles. We showed that to apply Theorem 2.1 to paths, it suffices to move the end-vertices to corners of $R_0$. The End-vertex-to-Boundary algorithm (Section 5.2) moves an end-vertex to a side of $R_0$, and the End-vertex-to-Corner algorithm (Section 5.3) moves it from a side to a corner. Together, these algorithms extend the reconfiguration result from cycles to paths.

\newpage

\section{Conclusion}

\noindent This chapter summarizes our main results (Section 6.1), discusses their computational complexity and implications for sampling algorithms (Sections 6.2 and 6.3), characterizes extremal reconfiguration cases (Section 6.4), explores the role of boundary structure (Section 6.5), and proposes several directions for future research (Section 6.6). Sections 6.4--6.6 present several conjectures for future investigation.

\subsection{Summary of Results}

In this dissertation we gave algorithms for reconfiguring Hamiltonian cycles and Hamiltonian paths on rectangular grid graphs. The main results are Theorem 2.1 and Theorem 5.9. Theorem 2.1 proves that any two Hamiltonian cycles on the same $m \times n$ rectangular grid graph with $m\leq n$ can be reconfigured into one another using only double-switch moves, with at most $O(n^2m)$ moves required. Theorem 5.9 proves the analogous result for Hamiltonian paths using  double-switch, single-switch, and backbite moves, again with an $O(n^2m)$ move bound. Both algorithms have $O(m^2n^2)$ time complexity with details given in Section 6.2.

\null 

\noindent We first analyzed the structure that a Hamiltonian path $H$ imposes on a polyomino, decomposing it into $H$-components. These $H$-components generalize the notion of cookies  defined by Nishat and Whitesides in \cite{nishat2017bend}, originally defined for Hamiltonian cycles, so that the concept applies to both Hamiltonian cycles and Hamiltonian paths. Then we established properties of $H$-components, which we then used to prove the existence of required moves for the reconfiguration algorithms.

\null

\noindent We introduce canonical forms for Hamiltonian cycles, designed to be simple yet sufficient as terminal configurations: once a cycle is reconfigured into a canonical form, the reconfiguration problem is essentially solved. 

\null 

\noindent The double-switch move used throughout this dissertation is a generalization of the flip and transpose moves used in \cite{nishat2017bend}. To extend the Hamiltonian cycle reconfiguration result to Hamiltonian paths, we included the single-switch and backbite moves.

\null 

\noindent The proof of Theorem 2.1 centers on the Reconfiguration to Canonical Form (RtCF) algorithm, supported by the OneLargeCookie (1LC) and ManyLargeCookies (MLC) algorithms. RtCF works iteratively from the outside in on a nested sequence of rectangles $R_0, R_1, \ldots, R_t$, reconfiguring the cycle so that all but two edges of each $R_i$ belong to the cycle before proceeding inward. The 1LC and MLC algorithms establish the existence of required moves at each step.

\null 

\noindent Theorem 5.9 reduces path reconfiguration to cycle reconfiguration through two new algorithms: End-to-Boundary (EtB) and End-to-Corner (EtC). These algorithms transform any Hamiltonian path into a form where Theorem 2.1 applies. The most involved cases in all four algorithms -- 1LC, MLC, EtB, and EtC -- are resolved through the fat-path-turn-weakening machinery developed in Chapter 4.

\subsection{Computational complexity}

\noindent In this section we give bounds for the time and space complexity of the reconfiguration algorithms presented in this dissertation. We assume that we are reconfiguring Hamiltonian paths or cycles on an $m \times n$ grid graph, with $n \geq m$.

\null

\noindent \textbf{Space complexity}: The grid graph can be represented by an adjacency list using $O(mn)$ space. An adjacency matrix representation would require $O(m^2n^2)$ space, but this is not necessary since the underlying graph is fixed. The current Hamiltonian cycle or path configuration requires $O(mn)$ space to store which edges are currently in use. Therefore, the total space complexity for the reconfiguration task is $O(mn)$. If one wishes to record the complete reconfiguration sequence, this requires additional space proportional to the number of moves. Since the reconfiguration requires at most $O(mn^2)$ moves, storing the complete sequence requires $O(mn^2)$ space. This can be done by storing an initial configuration and a sequence of length $mn^2$, where each term describes a move and costs $O(1)$.

\null

\noindent \textbf{Time complexity.} Consider the following scenario: an edge $e$ followed by a $j$-stack of $A_0$-types (with $j \leq n/2$), which is then followed by an $A_1$-type with switchable middle-box $W$. The middle-box $W$ determines an $H$-path $P(X,Y)$. A key computational bottleneck is determining which boxes belong to $P(X,Y)$. Using breadth-first search, this requires $O(mn)$ time in the worst case, as $P(X,Y)$ may traverse up to $O(mn)$ boxes. We now show that the remaining tasks -- checking FPC conditions, locating admissible turns, and finding turn weakenings -- each require at most $O(mn)$ time.

\null

\noindent\textit{Checking FPC conditions.} Once $P(X,Y)$ is identified, we check the FPC conditions as follows. FPC-2 (checking whether $W$ is parallel) requires checking which two edges of the Hamiltonian path are incident on $W$, which takes $O(mn)$ time, as well as one final check to see if the edges are parallel. So FPC-2 can be checked in $O(mn)$ time. FPC-1 (checking whether an end-vertex of the Hamiltonian path is incident on a box of $P(X,Y)$) requires checking each box of $P(X,Y)$ against the two end-vertices, taking $O(mn)$ time. FPC-3 (checking whether $P(X,Y)$ contains a switchable box) requires checking whether each box in $P(X,Y)$ is switchable, also $O(mn)$ time.

For FPC-4, we must check whether there exists a move (switch, flip, or transpose) after which $W$ is switched, $P(X,Y)$ gains a switchable box, or $W \to X$ or $W \to Y$ becomes valid. Note that only moves in the vicinity of boxes in $N[P(X,Y)] \setminus P(X,Y)$ can affect $P(X,Y)$ in the required manner. More precisely, it suffices to check moves near leaves or switchable boxes in $N[P(X,Y)] \setminus P(X,Y)$, which can be identified during the initial $O(mn)$ scan of $P(X,Y)$. For each such location, we check a $3 \times 3$ region in $O(1)$ time. Since there are at most $O(mn)$ such locations, the total time for FPC-4 is $O(mn)$.

\null 

\noindent\textit{Weakening of an admissible turn in a standard looping fat path.} Let$F = G\langle N[P(X,Y)] \rangle$ be a standard looping fat path and assume that we want to locate an admissible turn $T$ of the boundary $B(F)$. Recall that we store the Hamiltonian path as a sequence of vertices with each term identified by its coordinates, and that we begin our search for $T$ at an edge of the $A_1$-type of $F$. For definiteness, assume that $T$ is northeastern. This can be determined by checking the direction of the first two edges $e_1$ and $e_2$ of $T$ in $O(1)$ time. Next, we need to find the first western edge $e_W$ encountered after $e_1$, which can be no farther than $n$ units east of the first edge, requiring at most $2n$ direction checks. If the turn is admissible, we are done. If not, we repeat the procedure once more, but this time we need to find the first northern edge, which can be no farther than $m$ units south of the first western edge, requiring at most $2m$ direction checks. Thus, finding an admissible turn requires $O(m+n)$ time in the worst case. Given this admissible turn $T=T_1$, we can find $\text{Sector}(T)$ with cost $O(1)$ and check whether an end-vertex is incident on $\text{Sector}(T)$ with cost $O(1)$ as well. From this information, we can determine exactly which potential lengthening of $T$ may contain that end-vertex in its flank, avoiding the need to check every potential lengthening individually.

For each turn $T_i$ in the sequence of lengthenings of $T$ (where $i = 1, 2, \ldots, q$ and $q \leq m$), we check for a short weakening (a cascade of at most three moves) in $O(1)$ time by examining two $3 \times 3$ rectangles centered at the leaves of $T_i$. If no short weakening exists, there must be a lengthening $T_{i+1}$, whose exact structure (whether each leaf is open or closed) can be determined during the same check in $O(1)$ time. This process continues for at most $q \leq m$ iterations, requiring $O(m)$ time total. Once the final turn $T_q$ is reached (guaranteed to have a short weakening), we apply the weakenings in reverse order from $T_q$ back to $T_1$, requiring one move per turn for a total of $O(m)$ moves, each taking $O(1)$ time to identify. Thus the time complexity is $O(m+n)$. 

\null

\noindent\textit{The fat-path-turn-weakening procedure.} Consider once more the scenario described above: an $A_1$-type with switchable middle-box $W$ determining an $H$-path $P(X,Y)$, and the subgraph $F=G\langle N[P(X,Y)] \rangle$. We define the \textit{fat-path-turn-weakening procedure} as the complete process consisting of: (i) identifying which boxes belong to $P(X,Y)$, (ii) checking the FPC conditions for $F$, (iii) locating an admissible turn $T$ in $F$ whenever $F$ is a standard looping fat path, and (iv) finding a weakening of $T$. As shown above, this procedure requires $O(mn)$ time. We will reference this procedure in the complexity bounds that follow.

\null

\noindent \textit{Time complexity of the 1LC algorithm.} The 1LC algorithm handles Hamiltonian cycles with exactly one large cookie and at least one small cookie. 

Identifying all small cookies requires $O(m+n)$ time. For each outermost small cookie $C$, the algorithm first checks whether $C$ can be collected by a single move, which takes $O(1)$ time. If not, by Lemma 3.9, $C$ is followed by a $j$-stack of $A_0$-types (with $j \leq n/2$) and an $A_1$-type with switchable middle-box $W$. Checking for this costs $O(n)$ operations.

In this case, collecting $C$ requires the fat-path-turn-weakening procedure, which as shown above takes $O(mn)$ time. Following this, $j$ flips are executed to collect $C$, requiring $O(n)$ additional time. Thus, collecting an outermost small cookie requires at most $O(mn)$ time.

\null

\noindent\textit{Time complexity of Theorem 2.1.} The time complexity for MLC is $O(m+n)$ (see pseudocode at the end of Section 3.1), so the overall time complexity is dominated by the 1LC operations. 

The algorithm for Theorem 2.1 applies RtCF to both input Hamiltonian cycles to reduce them to canonical forms, then reconfigures between the canonical forms. Thus the total cost is twice the cost of RtCF plus the cost of canonical reconfiguration. The reconfiguration between canonical forms takes $O(m)$ time if we retain the cookie neck information from the RtCF executions, or $O(mn)$ time if computed from scratch. Either way, this is dominated by the RtCF cost.

Recall from the move bound for Theorem 2.1 at the end of Chapter 2 that at iteration $j$ of RtCF, there are at most $(n+m-4j)$ cookies on rectangle $R_j$ of size $(m-2j) \times (n-2j)$. Since the time complexity of 1LC on $R_0$ is $O(mn)$ per cookie, collecting a cookie via 1LC on $R_j$ requires $O((m-2j)(n-2j))$ time. Therefore, the total time complexity is:

$$T = 2 \cdot \sum_{j=0}^{\lfloor m/2 \rfloor} (n+m-4j) \cdot O((m-2j)(n-2j)) + O(m) = O(m^2n^2).$$

\noindent\textit{Time complexity of the EtB algorithm.} The EtB algorithm moves the two endpoints of a Hamiltonian path to opposite boundaries of the grid.
For each endpoint relocation step, either the endpoint can be moved directly in $O(1)$ time, or we are in the setting of Case 2.2 of the proof of Proposition 5.6. In Case 2.2, the costliest scenario requires at most two applications of the fat-path-turn-weakening procedure, along with at most $\frac{n}{2}$ flips. Since each application of the fat-path-turn-weakening procedure takes $O(mn)$ time and the flips take $O(n)$ time, the cost per relocation is $O(mn)$.
Moving both endpoints to opposite boundaries may require at most $O(n)$ such relocation steps. Therefore, the total time complexity of the EtB algorithm is $O(n) \cdot O(mn) = O(mn^2)$.

\null

\noindent\textit{Time complexity of the EtC algorithm.} The EtC algorithm moves the two endpoints of a Hamiltonian path to specified corners of the grid. Each relocation step either completes in $O(1)$ time or uses Claim 5.23. Similarly, in Claim 5.23, either we complete relocation in $O(1)$, or we are required to use Lemma 5.22, which invokes the fat-path-turn-weakening procedure. Its cost of $O(mn)$ dominates all other operations in the relocation.

To move a single endpoint from a boundary to a corner requires at most $\frac{n}{2}$ such relocation steps. Therefore, the time complexity to move each endpoint to its target corner is $O\left(\frac{n}{2}\right) \cdot O(mn) = O(mn^2)$.

\null

\noindent\textit{Time complexity of Theorem 5.9.} Theorem 5.9 reconfigures an arbitrary Hamiltonian path into another arbitrary Hamiltonian path on an $m \times n$ grid graph. The strategy, as described in the proof, is to apply Proposition 5.8 twice to transform both Hamiltonian paths into northern or eastern Hamiltonian paths, which can be viewed as Hamiltonian cycles. We then may apply Theorem 2.1 to reconfigure between them.

From the proof of Proposition 5.8, we know that transforming a path into northern form requires at most two applications of EtB and three applications of EtC. Since EtB has time complexity $O(mn^2)$ and EtC has time complexity $O(mn^2)$, one application of Proposition 5.8 requires $2 \cdot O(mn^2) + 3 \cdot O(mn^2) = O(mn^2)$. So the total time complexity is $2O(mn^2)+O(m^2n^2)=O(m^2n^2)$, accounting for one application of Proposition 5.8 per path, and one application of Theorem 2.1 to reconfigure the resulting northern or eastern paths.

\null 

\noindent\textbf{Summary.} We have analyzed the computational complexity of the reconfiguration algorithms presented in this dissertation. The fat-path-turn-weakening procedure, which is invoked in the most time-intensive scenarios, requires $O(mn)$ time. The endpoint relocation algorithms (EtB and EtC) each require $O(mn^2)$ time. The complete reconfiguration of Hamiltonian cycles (Theorem 2.1) and Hamiltonian paths (Theorem 5.9) both require $O(m^2n^2)$ time. Space complexity is $O(mn)$ for basic reconfiguration or $O(mn^2)$ if storing the complete move sequence. Notably, the time complexities for Theorems 2.1 and 5.9 are within a factor of $n$ of the move bounds for these theorems.

\null 

\begin{table}[h]
\centering
\renewcommand{\arraystretch}{1.3}
\begin{tabular}{lcc}
\toprule
\textbf{Algorithm} & \textbf{Move Bound} & \textbf{Time Complexity} \\
\midrule
Identifying $P(X,Y)$ & N/A & $O(mn)$ \\
MLC & $ 2 $ & $O(m+n)$ \\
1LC & $n^2/2+m+3$ & $O(mn)$ \\
EtB & $n^2+mn+4n$ & $O(mn^2)$ \\
EtC & $n^2/2+mn+5n/2$ & $O(mn^2)$ \\
Theorem 2.1 & $n^2m$ & $O(m^2n^2)$ \\
Theorem 5.9 & $n^2m+7n^2+10mn+31n$ & $O(m^2n^2)$ \\
\bottomrule
\end{tabular}
\caption{Move bounds and time complexities for the main algorithms.}
\label{tab:complexity-summary}
\end{table}

\subsection{Implications for sampling algorithms}  Consider the state space of all Hamiltonian cycles (or paths) on an $m \times n$ rectangular grid graph. Recall from the Introduction chapter that in Markov chain theory, a chain is called \textit{irreducible} if it is possible to reach any state from any other state through a sequence of transitions. For chains on finite state spaces with symmetric transition probabilities (as in our case), irreducibility -- together with the easily satisfied condition of aperiodicity -- implies \index{ergodicity|textbf}\textit{ergodicity}, the property that the chain eventually explores all possible configurations and converges to a unique equilibrium distribution. An \index{ergodicity class|textbf}\textit{ergodicity class} is defined as a maximal subset of the state space in which any two configurations can reach each other through a sequence of moves. A chain is irreducible if and only if the entire state space forms a single ergodicity class. Otherwise, the state space partitions into multiple ergodicity classes, and a Markov chain started in one class can never reach configurations in another class. This issue is particularly relevant when using Markov chain Monte Carlo to sample self-avoiding walks (including high-density polymer models, which correspond to Hamiltonian paths on lattices). Madras and Slade observe in \cite{madras2013self} that for certain Monte Carlo algorithms, the state space partitions into many ergodicity classes, with the largest class being exponentially smaller than the total number of configurations as the chain length increases.
Oberdorf et al.\ \cite{oberdorf2006secondary} used a backbite-based Monte Carlo algorithm for sampling Hamiltonian walks on the two-dimensional square lattice. Their simulations suggested that the underlying Markov chain was ergodic, but they could not provide a proof. Theorem 2.1 proves that any Hamiltonian cycle can be reconfigured into any other using double-switch moves in rectangular grid graphs, thereby establishing irreducibility of the associated Markov chain. Theorem 5.9 extends this result to Hamiltonian paths. Together, these results make progress toward proving the ergodicity of the Markov chain studied in \cite{oberdorf2006secondary} and provide an alternative, rigorous algorithm for generating the Hamiltonian cycles considered there. Analogous Monte Carlo algorithms for three-dimensional cubic lattices \cite{mansfield1982monte, jacobsen2008unbiased, deutsch1997long} similarly assume ergodicity without proof. Proving analogous reconfiguration results in three dimensions remains open.

\subsection{Resistant Hamiltonian cycles}

\noindent In this section we consider Hamiltonian cycles on an $m \times n$ grid graph $G$ with $ n \geq m$ for which RtCF requires $\Theta(mn^2)$ moves to complete. We call such cycles \textit{resistant} and denote their class by $\mathcal{H}_{\text{res}}(m,n)$. We give a heuristic characterization of resistant cycles (Conjecture 6.3) and argue that they comprise a vanishingly small fraction of all Hamiltonian cycles.

\null 

\noindent Recall how the upper bound of $mn^2$ moves arises for Theorem 2.1: essentially, the factor $m$ comes from the iterations of RtCF over the nested rectangles; one factor $n$ comes from the number of small cookies we may need to collect within a rectangle; and the other factor $n$ comes from the cascade length required for the 1LC algorithm (Proposition 3.10) to collect a small cookie $C$.  Part of this cascade's move cost is due to the sequence of flips needed to collect $C$, if $C$ is followed by a $j$-stack of $A_0$s, and the other part due to a sequence of turn weakenings that may be required. (we refer to $j$ as the  \index{height of a stack of A0s@height (of a stack of $A_0$s)|textbf}\textit{height} of the stack). We make the following assumption. 

\null 

\noindent \textbf{Assumption 6.1.} Let $H$ be a Hamiltonian cycle on an $m \times n$ grid graph $G$. If the cascade required to collect a small cookie $C$ uses a sequence of $j$ moves for turn weakenings, then $C$ is followed by a stack of $A_0$s of height at least $j$, and so it also requires an additional sequence of at least $j$ flips. 

\null

\noindent We will make two additional assumptions in this section and discuss them briefly at the end of the section.

\noindent It follows from Assumption 6.1 that a $j$-stack of $A_0$s following a small cookie $C$ contributes at least as many moves to the cost of collecting $C$ as a sequence of turn weakenings that may also be needed. Thus, up to a constant multiple, we may view the height of a $j$-stack of $A_0$s as the main contributor to the move cost of collecting $C$. So, if two Hamiltonian cycles each contain many stacks of $A_0$-types of large height, and these stacks do not coincide (their locations in $G$ are pairwise disjoint), then RtCF will require $\Theta(m n^2)$ moves to terminate. The total move cost from large cookies alone is only $O(mn)$, so the only way to reach $\Theta(m n^2)$ is to have enough small cookies that each demand many moves. We make this more concrete below.

By 1LC, collecting a cookie followed by a $j$-stack of $A_0$s costs at least $j$ moves. Moreover, after the collection, in the next iteration of RtCF the stack of $A_0$s persists, with one fewer $A_0$ (see Figure 3.14(d) in Chapter 3). Consider a single small cookie followed by a stack of $A_0$s with height $h=h(m,n)$. We make the following assumption:

\null 

\noindent \textbf{Assumption 6.2.} Let $H$ be a Hamiltonian cycle on an $m \times n$ grid graph $G$. A stack of $A_0$s with height $h$ contributes at least $h$ small cookies to the total number of cookies that need to be collected by RtCF, one for each $A_0$ in the stack.

\null

\noindent That is, if the first $A_0$ of the stack contains the small cookie $C_i$ on rectangle $R_i$, after $C_i$ has been collected via the flip move $X \mapsto C_i$, then on the iteration $i+1$, there will be at least one small cookie $C_{i+1}$ on $R_{i+1}$ adjacent to $X$. See Figure 6.1 (the smallest grid we could construct containing a $3$-stack). For $j \in \{i, i+1, \ldots, i+h\}$, the cookie $C_j$ on $R_j$ requires $h-j$ moves to be collected, so the total move contribution to the RtCF move counter is $\frac{(h)(h+1)}{2}=\Theta(h^2)$. Let $c=c(m,n)$ be the number of such stacks of $A_0$s with height $h$. Then, to reach the bound of $\Theta(mn^2)$ moves, we need $h$ and $c$ to satisfy $ch^2 =\Theta(mn^2)$. Thus, there are many possibilities for $h$ and $c$. 

We record the above as Conjecture 6.3.

\begin{center}

\end{center}

\noindent\textbf{Conjecture 6.3.} Let $G$ be an $m \times n$ grid graph and let $H$ be a Hamiltonian cycle of $G$. Then $H$ is resistant if and only if there exist functions $h(m,n)$ and $c(m,n)$ such that $H$ contains $c$ stacks of $A_0$-types, each with height $h$, where $ch^2 = \Theta(mn^2)$.

%pagemarker

\noindent A proof of Conjecture 6.3 would require proving Assumptions 6.1 and 6.2, and a detailed consideration of the converse. Note that the conjecture assumes the existence of Hamiltonian cycles with the required stacks, but we do not prove such cycles exist for all $(m,n)$. As the heights of the $A_0$ stacks rise, verifying that a candidate graph is indeed a Hamiltonian cycle becomes increasingly difficult. %Figure \textcolor{red}{6.1} shows a concrete example on a $15 \times 24$ grid containing a $4$-stack of $A_0$-types. 
While we have not produced examples with arbitrarily long stacks for sufficiently large grid dimensions, the existence of such resistant cycles seems plausible since the presence of $A_0$-stacks does not create obvious obstructions to Hamiltonicity.
Moreover, by constructing Hamiltonian cycles with $j$-stacks of $A_0$s for small values of $j$, we noticed that accommodating the stacks seemed to require both $m \gtrsim j^2$ and $n \gtrsim j^2$. While we have only observed this for $j \le 6$ (with the $6$-stack requiring a grid of over $1000$ vertices), the pattern suggests that both grid dimensions must grow substantially with stack height.

%       OLD VERSION of Prop 6.1.
{
}

\null 

\noindent \textbf{Conjecture 6.4.} The $n^2 m$ bound on the number of moves is tight. That is, there exists a constant $\alpha \in(0,1)$ such that for all sufficiently large $m$ and $n$, there are Hamiltonian cycles $H,K$ of an $m \times n$ grid graph $G$ for which any reconfiguration sequence from $H$ to $K$ requires at least $\alpha n^2 m$ moves.

\null

\noindent \textbf{Conjecture 6.5.} Resistant cycles are exponentially rare:
$\lim_{m,n \to \infty} \frac{|\mathcal{H}_{\text{res}}(m,n)|}{|\mathcal{H}(m,n)|} = 0.$

%\textcolor{red}{[I should make it more clear here that this is for a $fixed$ constant $a$, as $a$ is chosen below... maybe i should say ``they are exponentially rare" instead...]}.

\null

\noindent We give a plausibility argument in support of Conjecture 6.5, using a heuristic but unproven assumption (Assumption 6.6). Since Hamiltonian cycles are self-avoiding walks, $|\mathcal{H}(m,n)| \leq 3^{mn}$. In Note 6.7 at the end of the section, we show that $(2^{1/6})^{mn} \leq |\mathcal{H}(m,n)|$, so the per-vertex growth rate $\tau$ of $|\mathcal{H}(m,n)|$ satisfies $2^{1/6} \leq \tau \leq3$. The authors in \cite{bousquet2005self} give better bounds for the related case of Hamiltonian paths on square grid graphs. 

What is the smallest possible space occupied by stacks of $A_0$s in a resistant Hamiltonian cycle? Each stack of height $h$ occupies $O(h)$ vertices, so $c$ stacks occupy $ch$ vertices total. From the characterization $ch^2 = \Theta(mn^2)$, we have $ch = \Theta(mn^2/h)$. Since $m$ and $n$ are fixed, to minimize the space occupied we must maximize $h$. The maximum possible stack height is $h = \Theta(n)$, in which case $c = \Theta(m)$. Thus, stacks must occupy at least $\Theta(mn^2/n) = \Theta(mn)$. That is, there exists a constant $a \in (0,1)$ such that the stacks occupy at least $amn$ vertices -- a fraction $a$ of the total grid. We now make the following assumption about how this constraint affects the count of Hamiltonian cycles.

\null

\noindent \textbf{Assumption 6.6.} Constraining $amn$ vertices with $A_0$-stacks reduces the possible number of Hamiltonian cycles as if these vertices were removed from the grid. That is, for each fixed placement of the stacks there are $\tau^{mn(1-a)}$ possible resistant Hamiltonian cycles.

\null

\noindent In how many ways can we position the stacks? We need to choose $c = \Theta(m)$ positions from the $mn$ vertices in the grid. For an upper bound, we use $(mn)^c = (mn)^{m}$. Thus:

$$
|\mathcal{H}_{\text{res}}(m,n)| \lesssim (mn)^{m} \cdot \tau^{mn(1-a)}.
$$

\noindent The proportion of resistant cycles to all cycles is then at most

\begin{align*}
\frac{|\mathcal{H}_{\text{res}}(m,n)|}{|\mathcal{H}(m,n)|} 
&\lesssim \frac{(mn)^{m} \cdot \tau^{mn(1-a)}}{\tau^{mn}} \\
&= (mn)^{m} \cdot \tau^{-amn} \\
&= (mn)^{m} \cdot \left(\frac{1}{\tau}\right)^{amn}.
\end{align*}

\noindent Since $(mn)^m = e^{m \ln(mn)}$ grows sub-exponentially compared to $\tau^{amn} = e^{amn \ln \tau}$, and since $\tau \geq 2^{1/6} > 1$ so $1/\tau < 1$, we have 

$$
(mn)^m \cdot \left(\frac{1}{\tau}\right)^{amn} \to 0 \quad \text{as } m,n \to \infty.
$$

\noindent This heuristic suggests that the class of resistant Hamiltonian cycles is vanishingly small relative to all Hamiltonian cycles. 

Before showing that $(2^{1/6})^{mn} \leq |\mathcal{H}(m,n)|$, we need a definition. Let $G$ be a graph. A \index{disjoint cycle cover|textbf}\textit{disjoint cycle-cover} of $G$ is a spanning subgraph that is a disjoint union of cycles. This is also known as a \index{2-factor|textbf}\textit{2-factor} of $G$.

\null
\noindent \textbf{Note 6.7.} To see that the number of Hamiltonian cycles on an $m \times n$ rectangular grid grows exponentially in the area, assume $m$ is even and that $n$ is divisible by~3 (recall that $m$ is the horizontal dimension and $n$ is the vertical dimension, with $m \le n$). Partition $G$ into $n/3$ horizontal strips, each of size $3 \times m$ (that is, $3$ rows and $m$ columns). By Lemma~3.4(d) in~\cite{nishat2020reconfiguration}, each $3\times m$ strip contains exactly $2^{(m-2)/2}$ Hamiltonian cycles. Fix a disjoint cycle cover $\{H_1,\ldots,H_{n/3}\}$ of~$G$ consisting of one Hamiltonian cycle per strip. Observe that there are at least two switchable boxes between every adjacent pair $(H_i,H_{i+1})$. Switching one of those boxes merges the two cycles into one without affecting the others, and repeating this process merges all strips into a single Hamiltonian cycle of~$G$. See Figure 6.2. There are $(2^{(m-2)/2})^{n/3}=2^{n(m-2)/6}$ such Hamiltonian cycles, so
$$
|\mathcal{H}(m,n)|^{1/(mn)} \ge
\left(2^{n(m-2)/6}\right)^{1/(mn)}
= 2^{(m-2)/(6m)} \to 2^{1/6}
\quad\text{as } m\to\infty.
$$
\noindent Hence $\tau \ge 2^{1/6} \approx 1.122 > 1$.
When $n \not\equiv 0 \pmod{3}$, the construction can be extended by including the remaining $1$ or $2$ rows without affecting the exponential rate.

\noindent We conclude with remarks on the two key assumptions in this section. Proofs for Assumptions 6.1 and 6.2 seem within reach of our current tools, without requiring new techniques. We believe Assumption 6.1 would require substantial case analysis related to the structure of sturdy looping fat paths, and that Assumption 6.2 would be more straightforward. Assumption 6.6 on the other hand, seems more difficult and we have no idea how to prove it. We believe the insights about resistant cycles are valuable even with these assumptions remaining unproven.

\begin{center}

\end{center}

\subsection{Role of the boundary}

\noindent Throughout this dissertation, many results are stated for rectangular grid graphs. In the proofs, we do not explicitly invoke this assumption again, but the rectangular structure is used implicitly at nearly every step: in the design of canonical forms and the RtCF algorithm (Chapter 2); in the MLC algorithm and Lemma 3.5 (Chapter 3); in the turn-weakening lemmas near boundaries (Lemmas 4.23, 4.26, 4.29, and 4.30 in Chapter 4); in the description of some impossible configurations (Lemmas 5.11 and 5.12 in Chapter 5); and in the EtB and EtC algorithms (Chapter 5). In this section we consider what happens when we relax the rectangular boundary condition or embed rectangular grids on a cylinder or torus. We will see that reconfiguration can fail: Hamiltonian cycles may become frozen, or the space of cycles may split into disconnected components.

\null 

\noindent\textbf{Example 6.8.} Hamiltonian cycles of a simply connected polyomino that are frozen under the double-switch move. See Figure 6.3 (a). 

%pagemarker

\noindent\textbf{Example 6.9.} Hamiltonian paths of a simply connected polyomino that belong to distinct ergodicity classes under the backbite move. See Figure 6.3 (b). 

\null 

\begin{adjustbox}{trim=0.15cm 0cm 0cm 0cm} 

\end{adjustbox}

\null 

\noindent\textbf{Example 6.10.} A Hamiltonian cycle of a $6 \times 7$ grid graph embedded on a cylinder that is frozen under the double-switch move. This construction appears to generalize $k \times (k+1)$ grids, for each even $k\geq 6$. See Figure 6.3 (c).

\null 

\noindent\textbf{Example 6.11.} A Hamiltonian cycle of a $6 \times 12$ grid graph embedded on a torus that is frozen under the double-switch move. This construction appears to generalize to $k \times 2k$ grids for each even $k\geq 6$. See Figure 6.3 (d).

\null 

\noindent The counterexamples illustrated above raise the following question: can we reconfigure any Hamiltonian cycle on a polyomino to any other by expanding our move set? Can we do the same for rectangular grid graphs embedded on a cylinder or torus? We conclude this section with two conjectures. We first need a definition. Let $H=H_0$ be a Hamiltonian cycle of a polyomino $G$. A \index{quadruple-switch move|textbf}\textit{quadruple-switch move} on $H$ is defined as follows: for $i \in \{1, \ldots 4\}$, perform a switch $\sigma_i$ on a switchable box of $H_{i-1}$ to obtain the subgraph $H_i$. The move is \textit{valid} if $H_4$ is a Hamiltonian cycle of $G$.

\null

\noindent\textbf{Conjecture 6.12.} Let $G$ be a polyomino and let $H$ and $K$ be two distinct Hamiltonian cycles of $G$. Then there exists a sequence of valid quadruple-switch moves that reconfigures $H$ into $K$.

\null

\noindent\textbf{Conjecture 6.13.} Let $G$ be an $m \times n$ grid graph embedded on a cylinder or torus, and let $H$ and $K$ be two distinct Hamiltonian cycles of $G$. Then there exists a sequence of valid quadruple-switch moves that reconfigures $H$ into $K$.

\null 

\noindent Why might a quadruple-switch move be enough for these more general reconfiguration questions? It seems that these moves might be sufficient to break down the stubborn structures arising from combinations of $A$-types and stairs subgraphs (recall definition at the beginning of Chapter 4), and compensate for the loss of control incurred by relaxing the conditions that the boundary be rectangular and that the graph be embedded in the plane.

\null 

\subsection{Extensions and future work}

\noindent In this section we give an overview of possible extensions of our results and state open problems for future research. These extensions include restricting the move set, reconfiguration of cycle covers (of which Hamiltonian cycles are a special case), reconfiguration on higher dimensional square lattices, and reconfiguration on the triangular lattice. We extend the usual $m \times n$ rectangular grid graph to higher dimensions as follows. For integers $n_1,\dots,n_d \ge 2$, a \index{d-orthotope@$d$-orthotope|textbf}\textit{$d$-dimensional orthotope grid graph (or $d$-orthotope)} is the graph $G$ with vertex set

$$
V(G) = \{ (x_1,\dots,x_d) \in \mathbb{Z}^d : 0 \le x_i < n_i \text{ for all } i \},
$$
and with an edge between two vertices if they differ by $1$ in exactly one coordinate and are equal in all others.

Nishat and Whitesides in \cite{nishat2017bend} ask whether it is possible to reconfigure Hamiltonian cycles in rectangular grid graphs using just flip and transpose moves (recall the definitions from the Introduction chapter). Nishat in \cite{nishat2020reconfiguration} also asks if given two disjoint cycle covers in a grid graph, we can reconfigure them into one another. We restate these questions here as conjectures.

\null

\noindent\textbf{Conjecture 6.14.} Let $G$ be an $m \times n$ grid graph and let $H$ and $K$ be two distinct Hamiltonian cycles of $G$. Then there exists a sequence of flips and transposes that reconfigures $H$ into $K$.

\null 

\noindent We have frequently encountered the scenario consisting of  a leaf $L$ followed by an $A_1$-type with switchable middle-box $W$ with looping $H$-path $P(X,Y)$. We have shown that we can always find a cascade ending with the move $Z \mapsto W$, where $Z$ is a switchable box of $P(X,Y)$. We believe that Conjecture 6.14 hinges on showing that this last ``long-range" double-switch move $Z \mapsto W$ can be replaced by a cascade of flips and transposes. The tools in Chapter 4 seem useful to address this conjecture, however we believe additional results would be needed regarding the structure of fat paths and turns, including a proof of Assumption 6.1.

\null

\noindent\textbf{Conjecture 6.15.} Let $G$ be an $m \times n$ grid graph and let $H$ and $K$ be two distinct disjoint cycle covers of $G$. Then there exists a sequence of double-switch moves that reconfigures $H$ into $K$. Moreover, if $G$ is a general polyomino, then quadruple-switch moves are enough for reconfiguration.

\null

\noindent The rationale for the polyomino case is similar to that given in support of Conjectures 6.12 and 6.13: more flexible moves may compensate for the loss of the rectangular structure.

\null 

% THREE PARAGRAPHS ON HEX AND TRIANGULAR LATTICE
{
}

\noindent Our results are specific to square lattices. The reconfiguration problem remains open for triangular and hexagonal lattices. The backbite move may be particularly effective for reconfiguring Hamiltonian paths in the triangular lattice, since we have more edges to choose from: an interior end-vertex in a triangular lattice has five possible edges (vs three in the square lattice), available for initiating a backbite move. Moreover, for subgraphs with triangular boundaries, there are two possible edges toward each boundary side, suggesting that relocating an end-vertex to the boundary may be more straightforward. That is, an analogue of the EtB algorithm may be easier to prove for triangular-lattice graphs. This additional degree of freedom suggests that backbite-only reconfiguration may be more tractable in the triangular lattice.

\null

\noindent Extending our results to higher dimensions presents new challenges. The tools developed in this dissertation relied on the 2-dimensional planar structure, in particular, the consequences of Jordan's curve theorem for polygons (Theorem 1.1.4). In three dimensions and beyond, there is no good analogue of this theorem, and much of our toolbox becomes unusable. 

We have been researching Hamiltonian cycle reconfiguration in 3-dimensional grids with the help of GridLab, a tool with interactive visualization features that can be used to generate Hamiltonian paths and cycles in 3 dimensions, available online at \cite{graphsofwrath}. This research suggests that reconfiguration may still be possible with appropriate move sets, leading us to the following conjecture.

\null

\noindent\textbf{Conjecture 6.16.} Let $G$ be a $3$-dimensional cuboid grid graph (a $3$-orthotope) and let $H$ and $K$ be two distinct Hamiltonian cycles of $G$. Then there exists a sequence of double-switch and switch moves that reconfigures $H$ into $K$.

\null

\noindent For Hamiltonian paths in higher dimensions, the backbite move may suffice on its own:

\null

\noindent\textbf{Conjecture 6.17.} Let $G$ be a $d$-dimensional orthotope grid graph (for $d \geq 3$) and let $H$ and $K$ be two distinct Hamiltonian paths of $G$. Then there exists a sequence of backbite moves that reconfigures $H$ into $K$.

\null

\noindent Conjecture 6.16 may hinge on resolving the following subproblem: Let $H = P(v_1,v_s),$ $P(v_s,v_t),$ $P(v_t,v_r)$ be a Hamiltonian path on a $d$-orthotope, where $P(v_s,v_t)$ is already in a ``nice" form. Can we relocate an end-vertex, say $v_1$, to any vertex of $H \setminus P(v_s,v_t)$ that $v_1$ can reach using only backbite moves? We are currently researching this question for the $3$-dimensional case, with a move set that includes single-switches and double-switches.

\null

\newpage

\printindex

\newpage 

\cleardoublepage
\phantomsection
\bibliography{references}

\end{document}